\documentclass[a4paper,12pt]{article}
\usepackage{amscd,amsfonts,amssymb,amscd,amsmath,latexsym}
\usepackage[hypertex,backref]{hyperref}
\usepackage[all]{xy}

\makeatletter
\renewcommand{\subsubsection}{\@startsection
{subsubsection}
{1}
{0mm}
{0mm}
{0mm}
{\normalfont\normalsize}}
\makeatother

\usepackage{mathrsfs}
\newtheorem{theorem}{Theorem}[section] 
\newtheorem{prop}[theorem]{Proposition}
\newtheorem{lem}[theorem]{Lemma}
\newtheorem{ddd}[theorem]{Definition}
\newtheorem{kor}[theorem]{Corollary}

\newtheorem{fact}[theorem]{Fact}

\newcommand{\R}{{\mathbb{R}}}
\newcommand{\Z}{{\mathbb{Z}}}
\newcommand{\hB}{{$\hfill \square$}\newline}

\newcommand{\Aut}{{\tt Aut}}
\newcommand{\Hom}{{\tt Hom}}

\newcommand{\Map}{{\tt Map}}
\newcommand{\Pic}{{\mathcal{P}\tt ic}}

\newcommand{\cD}{{\mathcal{D}}}
\newcommand{\Mod}{{\tt mod}}

\newcommand{\cF}{{\mathcal{F}}}

\newcommand{\cT}{{\mathcal{T}}}
\newcommand{\id}{{\tt id}}

\newcommand{\pr}{{\tt pr}}
\newcommand{\cR}{{\mathcal{R}}}

\newcommand{\Ext}{{\tt Ext}}
  
  \newcommand{\Ob}{{\tt Ob}}

  \newcommand{\T}{{\mathbb{T}}}  
\newcommand{\cB}{{\mathcal{B}}}
\newcommand{\proof}{{\tt Proof:\mbox{$\:\:$}}}

\newcommand{\HOM}{{\tt HOM}}
\newcommand{\Sh}{{\tt Sh}}
\newcommand{\Ab}{{\tt Ab}}
\newcommand{\ch}{{\tt ch}}
\newcommand{\pch}{{\tt{}^p ch}}
\newcommand{\cPic}{{\tt PIC}}
\newcommand{\Q}{\mathbb{Q}}
\newcommand{\Top}{{\tt TOP}}

\title{Duality for topological abelian group stacks and $T$-duality}
\author{U. Bunke, Th. Schick, M. Spitzweck, and A. Thom}

\renewcommand{\Pic}{{\tt PIC}}
\renewcommand{\P}{{\mathbb{P}}}
\newcommand{\cov}{{\tt cov}}

\newcommand{\cC}{{\mathcal{C}}}
\newcommand{\uU}{{\underline{U}}}
\newcommand{\cI}{{\mathcal{I}}}
\newcommand{\cK}{{\mathcal{K}}}
\newcommand{\Cat}{{\tt Cat}}
\newcommand{\uG}{{\underline{G}}}
\newcommand{\uH}{{\underline{H}}}
\renewcommand{\top}{{\tt top}}
\newcommand{\uX}{{\underline{X}}}
\newcommand{\uExt}{{\underline{\tt Ext}}}
\newcommand{\uR}{{\underline{\R}}}
\newcommand{\Sets}{{\tt Sets}}
\newcommand{\N}{{\mathbb{N}}}
\newcommand{\coinvv}{{\tt coinv}}

\newcommand{\uE}{\underline{E}}

\renewcommand{\lim}{{\tt lim}}
\newcommand{\cW}{{\mathcal{W}}}
\newcommand{\res}{{\tt res}}
\newcommand{\ind}{{\tt ind}}
\newcommand{\Tor}{{\tt Tor}}
\newcommand{\colim}{{\tt colim}}
\newcommand{\rel}{{\tt rel}}
\newcommand{\uD}{{\underline{D}}}

\newcommand{\im}{{\tt im}}
\newcommand{\uC}{\underline{C}}
\newcommand{\uL}{{\underline{L}}}
\newcommand{\uT}{{\underline{\mathbb{T}}}}
\newcommand{\nat}{\mathbb{N}}
\newcommand{\uW}{{\underline{W}}}
\newcommand{\bS}{{\mathbf{S}}}
\newcommand{\Mor}{{\tt Mor}}
\newcommand{\uF}{{\underline{F}}}
\newcommand{\EXT}{{\tt EXT}}
\newcommand{\cone}{{\tt Cone}}
\newcommand{\coinv}{{\tt \coinv}}
\newcommand{\cE}{{\mathcal{E}}}
\newcommand{\cTorsor}{{\mathcal{T}\tt ors}}
\newcommand{\Torsor}{{\tt Tors}}
\newcommand{\groupoids}{{\tt groupoids}}
\newcommand{\cEXT}{{\mathcal{E}\tt XT}}
\newcommand{\cG}{\mathcal{G}}
\newcommand{\cH}{\mathcal{H}}
\newcommand{\cL}{\mathcal{L}}
\newcommand{\uB}{{\underline{B}}}
\newcommand{\uZ}{\underline{\mathbb{\Z}}}

\newcommand{\uHOM}{{\underline{\tt HOM}}}
\newcommand{\Triple}{{\tt Triple}}
\newcommand{\Prin}{{\tt Prin}}
\newcommand{\uK}{{\underline{K}}}
\newcommand{\Gerbe}{{\tt Gerbe}}
\newcommand{\can}{{\tt can}}
\newcommand{\Gr}{{\tt Gr}}
\newcommand{\coker}{{\tt coker}}
\newcommand{\ev}{{\tt ev}}
\newcommand{\diag}{{\tt diag}}
\newcommand{\antidiag}{{\tt antidiag}}
\newcommand{\uA}{{\underline{A}}}
\newcommand{\cU}{\mathcal{U}}
\newcommand{\coind}{{\tt coind}}
\newcommand{\cV}{\mathcal{V}}
\newcommand{\uP}{×\underline{P}}
\newcommand{\F}{\mathbb{F}}
\newcommand{\End}{{\tt End}}
\newcommand{\rk}{{\tt rk}}
\newcommand{\uHom}{{\underline{\Hom}}}

\begin{document}

\maketitle
 
\begin{abstract}
We extend Pontrjagin duality from topological abelian groups to certain locally compact group stacks.
To this end we develop a sheaf theory on the big site of topological spaces $\bS$ in order to prove that the sheaves
$\uExt^i_{\Sh_\Ab\bS}(\uG,\uT)$, $i=1,2$, vanish, where $\uG$ is the sheaf represented by a locally compact abelian
group and $\T$ is the circle. As an application of the theory we
interpret topological $T$-duality of principal $\T^n$-bundles in terms of Pontrjagin duality of abelian group stacks.
\end{abstract}

\tableofcontents
\section{Introduction}

\subsection{A sheaf theoretic version of Pontrjagin duality}

\subsubsection{}

A character of an abelian topological group $G$ is a continuous homomorphism $$\chi:G\to \T$$
from $G$ to the circle group $\T$.
The set of all characters of $G$ will be denoted by $\widehat G$. It is again a group under point-wise multiplication.

Assume that $G$ is locally compact. The compact-open topology on the space of continuous maps $\Map(G,\T)$
induces a compactly generated topology on this space of maps, and hence a topology
on its subset $\widehat G\subseteq \Map(G,\T)$.
The group $\widehat G$ equipped with this topology is called the dual group  of $G$. 

\subsubsection{}

An element $g\in G$ gives rise to  a character $\ev(g)\in \widehat{\widehat G}$ defined by $\ev(g)(\chi):=\chi(g)$ for $\chi\in \widehat G$. In this way we get a continuous homomorphism
$$\ev:G\to \widehat{\widehat G}\ .$$
The main assertion of Pontrjagin duality is
\begin{theorem}[Pontrjagin duality]
If $G$ is a locally compact abelian group, then
$\ev:G\to \widehat{\widehat G}$ is an isomorphism of topological groups.
\end{theorem}
Proofs of this theorem can be found e.g. in \cite{MR1397028} or \cite{MR1646190}.

\subsubsection{}

In the present paper we use the language of sheaves in order to encode the topology of spaces and groups.
Let $X,B$ be topological spaces.
The space $X$ gives rise to a sheaf of sets $\uX$ on $B$ which associates to an open subset
$U\subseteq B$ the set $\uX(U)=C(U,X)$ of continuous maps from $U$ to $X$.
If $G$ is a topological abelian group, then $\uG$ is a sheaf of abelian groups on $B$.

The space $X$ or the group $G$ is not completely determined by the sheaf it generates over $B$. As an extreme
example take $B=\{*\}$. Then $X$ and the underlying discrete space $X^\delta$ induce isomorphic sheaves
on $B$. 

For another example, assume that $B$ is totally disconnected. Then the sheaves generated by  the spaces $[0,1]$ and $\{*\}$
are isomorphic.

\subsubsection{}

But one can do better and consider sheaves which are defined on all topological spaces or at least on a sufficiently big subcategory $\bS\subset \Top$. We turn $\bS$ into a Grothendieck site determined by the pre-topology of open coverings (for details about the choice of $\bS$ see 
Section \ref{hwdqwwqdqww}).
We will see that the topology in $\bS$ is sub-canonical so that every object
$X\in \bS$ represents a sheaf $\uX\in \Sh\bS$. By the Yoneda Lemma the space $X$ can be recovered from the sheaf $\uX\in \Sh\bS$ represented by $X$. The evaluation of the sheaf $\uX$ on $A\in \bS$ is defined by $\uX(A):=\Hom_{\bS}(A,X)$.
Our site $\bS$ will in particular contain all locally compact spaces.

\subsubsection{}\label{dfdfa}

We can reformulate Pontrjagin duality in sheaf theoretic terms.
The circle group $\T$ belongs to $\bS$ and gives rise to a sheaf $\uT$.
Given a sheaf of abelian groups $F\in \Sh_\Ab\bS$ we define its dual by
$$D(F):=\uHom_{\Sh_\Ab\bS}(F,\uT)$$
and observe that
$$D(\uG)\cong \underline{\widehat G}$$
for an abelian group $G\in \bS$.

The image of the natural pairing
$$F\otimes_\Z D(F)\to \uT$$ under the isomorphism
$$\Hom_{\Sh_\Ab\bS}(F\otimes_\Z D(F),\uT)\cong \Hom_{\Sh_\Ab\bS}(F,\uHom_{\Sh_\Ab\bS}(D(F),\uT))=\Hom_{\Sh_\Ab\bS}(F,D(D(F)))$$
gives the evaluation map
$$\underline{\ev_F}:F\to D(D(F))\ .$$

\begin{ddd}\label{hqdqwdwqdiowqdoiqwdqwd}
We call a sheaf of abelian groups $F\in \Sh_\Ab\bS$ dualizable, if the evaluation map
$\underline{\ev_F}:F\to D(D(F))$ is an isomorphism of sheaves.
\end{ddd}
The sheaf theoretic reformulation of Pontrjagin duality is now:
\begin{theorem}[Sheaf theoretic version of Pontrjagin duality]\label{ttt57}
If $G$ is a locally compact abelian group, then $\uG$ is dualizable.
 \end{theorem}

\subsection{Picard stacks}

\subsubsection{}

A group gives rise to  a category $BG$ with one object so that the group appears as the group of automorphisms
of this object. A sheaf-theoretic analog is the notion of a gerbe.

\subsubsection{}

A set can be identified with a small category which has only identity morphisms.
In a similar way a sheaf of sets can be considered as a strict presheaf of categories. In the present paper a presheaf of categories on $\bS$ is a lax contravariant functor $\bS\to \Cat$. Thus a presheaf $F$ of categories
associates to each object $A\in \bS$ a category $F(A)$, and to each morphism $f:A\to B$ a functor
$f^*:F(B)\to F(A)$. The adjective lax means, that in addition for each pair of  composable morphisms $f,g\in \bS$ we have  specified an isomorphism of functors
$g^*\circ  f^*\stackrel{\phi_{f,g}}{\cong} (f\circ g)^*$ which satisfies higher associativity relations.
The sheaf of categories is called strict if these isomorphisms are identities.
\subsubsection{}

A category is called a groupoid if all its morphisms are isomorphisms.
A prestack on $\bS$ is a presheaf of categories on $\bS$ which takes values in
groupoids.
A prestack is a stack if it satisfies  in addition
descent conditions on the level of objects and morphisms. For details about stacks we refer to \cite{MR2223406}.

\subsubsection{}

A sheaf of groups $F\in \Sh_\Ab\bS$ gives rise to a prestack ${}^p\cB F$
which associates to $U\in \bS$ the groupoid ${}^p \cB F(U):=B F(U)$.
This prestack is not a stack  in general.
But it can be stackified (this is similar to sheafification) to the stack $\cB F$.

\subsubsection{}

For an abelian group $G$ the category $BG$ is actually a group object in $\Cat$. The group operation is implemented by the functor $BG\times BG\to BG$ which is obvious on the level of objects and given by the group structure of $G$
on the level of morphisms. In this case associativity and commutativity is strictly satisfied.
In a similar manner, for $F\in \Sh_\Ab\bS$ the prestack
${}^p\cB F$ on $\bS$ becomes a group object in prestacks on $\bS$.

\subsubsection{}\label{uqwddqwdwq}
A group $G$ can of course also be viewed as a category $G$ with only identity morphisms.
This category is again a strict group object in $\Cat$.
The functor  $G\times G\to G$ is given by the group operation on the level of objects, and in the obvious way on the level of morphisms.

\subsubsection{}

In general, in order to define group objects in two-categories like $\Cat$ or the category of stacks on $\bS$, 
one would relax the strictness of associativity and commutativity. For our purpose the appropriate relaxed notion is that of a Picard category which we will explain in detail in \ref{uwqdiwqdwdqdqdqwdwq4334}. The corresponding sheafified notion
is that of a Picard stack. We let $\Pic\bS$ denote the two-category of Picard stacks on $\bS$.

\subsubsection{}
\label{def_of_BF}

A sheaf of groups $F\in \bS_\Ab\bS$ gives rise to a Picard stack in two ways.

First of all it determines the Picard stack $F$ which associates to $U\in \bS$ the
Picard category $F(U)$ in the sense of \ref{uqwddqwdwq}.
The other possibility is the Picard stack
$\cB F$ obtained as the stackification of the Picard prestack ${}^p\cB F$.

 \subsubsection{}

Let $A\in \Pic\bS$ be a Picard stack. The presheaf of its isomorphism classes of objects generates a sheaf of abelian groups which will be denoted by
$H^{0}(A)\in \Sh_\Ab\bS$. The Picard stack $A$ gives furthermore rise to the sheaf of abelian groups
 $H^{-1}(A)\in \Sh_\Ab\bS$ of automorphisms of the unit object.

A Picard stack
$A\in \Pic\bS$ fits into an extension (see \ref{dedwhedjewdn66})
\begin{equation}\label{hdjqwhdqwjdqwdqwdd}
\cB H^{-1}(A)\to A \to H^0(A)\ .
\end{equation}

\subsection{Duality of Picard stacks}

\subsubsection{}

It is an essential observation by Deligne that the two-category of Picard stacks $\cPic\bS$ admits an interior
$\uHOM_{\Pic\bS}$. Thus for Picard stacks $A,B\in \Pic\bS$ we have a Picard stack
$$\uHOM_{\Pic\bS} (A,B)\in \Pic\bS$$
of additive morphisms from $A$ to $B$ (see \ref{zqdiqdqwqwodiowqdqwdwqd}).

\subsubsection{}

Since we can consider sheaves of groups as Picard stacks in two different ways one can 
now ask how Pontrjagin duality is properly reflected in the language of Picard stacks.
It turns out that the correct dualizing object is the stack $\cB \uT\in \Pic\bS$, and not $\uT\in \Pic\bS$ as one might 
guess first. For a Picard stack
$A$ we define its dual by
$$D(A):=\uHOM_{\Pic\bS}(A,\cB \uT)\ .$$
We hope that using the same symbol for the dual sheaf and dual Picard stack does not introduce to much confusion (see the two footnotes attached to the formulas (\ref{ttt54}) and (\ref{ttt56}) below).
One can ask what the duals of the Picard  stacks $F$ and $\cB F$ look like.
In general we have (see \ref{idosdc333})
\begin{equation}\label{ttt54}
D(\cB F)\cong D(F)\footnote{On the right-hand side $D(F)$ is the dual sheaf as in \ref{dfdfa}, considered as a Picard stack as in \ref{def_of_BF}.}\ .
\end{equation}
One could expect that
\begin{equation}\label{ttt56}
D(F)\cong \cB D(F)\footnote{The symbol $D(F)$ on the left-hand side of this equation denotes the dual of the Picard stack given by the sheaf $F$, while $D(F)$ on the right-hand side is the dual sheaf.}\ ,
\end{equation}
but this only holds under the condition that
$\uExt^1_{\Sh_\Ab\bS}(F,\uT)=0$
(see \ref{zuzede333}).
This condition is not always satisfied, e.g.
$$\uExt^1_{\Sh_\Ab\bS}(\underline{\oplus_{n\in \nat}\Z/2\Z},\uT)\not=0$$
(see \ref{zwdwqidqidwqpdopop}).
But the main effort of the present paper is made to show that 
$$\uExt^1_{\Sh_\Ab\bS}(\uG,\uT)=0$$ for a large class of locally compact abelian 
groups (this condition is part of admissibility \ref{hkdjkwjkqwkdjqwdkjqwdqwd}).

\subsubsection{}

Let $A\in \Pic \bS$ be a Picard stack.
There is a natural evaluation
$$\underline{\ev_A}:A\to D(D(A))\ .$$

In analogy to \ref{hqdqwdwqdiowqdoiqwdqwd} we make the following definition.
\begin{ddd}
We call a Picard stack dualizable if $\underline{\ev_A}:A\to D(D(A))$ is an equivalence of Picard stacks.
\end{ddd}

\subsubsection{}

If $G$ is a locally compact abelian group
and
\begin{equation}\label{hjwdhkwqdwqdqdq}
\Ext^1_{\Sh_\Ab\bS}(\uG,\uT)\cong \Ext^1_{\Sh_\Ab\bS}(\widehat \uG,\uT)\cong 0\ ,
\end{equation}
then the evaluation maps
$$\underline{\ev_\uG}:\uG\to D(D(\uG))\ ,\quad 
\underline{\ev_{\cB\uG}}:\cB \uG\to D(D(\cB\uG))$$
are isomorphisms by applying (\ref{ttt54}) and (\ref{ttt56}) twice, and using  the sheaf-theoretic version of
Pontrjagin duality \ref{ttt57}.
In other words, under the condition (\ref{hjwdhkwqdwqdqdq}) above $\uG$ and $\cB \uG$ are dualizable Picard stacks.

\subsubsection{}

Given the structure $$ \cB H^{-1}(A)\to A \to H^0(A)$$ of $A\in \Pic\bS$,
we ask for a similar description of the dual  $D(A)\in \Pic\bS$.

We will see under the crucial condition $$\uExt^1_{\Sh_\Ab\bS}(H^0(A),\uT)\cong \uExt^2_{\Sh_\Ab\bS}(H^{-1}(A),\uT)\cong 0\ ,$$  that $D(A)$ fits into 
$$ \cB D(H^{0}(A))\to D(A) \to D(H^{-1}(A))\ .$$
Without the condition the description of $H^{-1}(D(A))$ is more complicated, and we refer to (\ref{zsauxasx7}) for more details.

\subsubsection{}

This discussion now leads to one of the main results of the present paper.

\begin{ddd}\label{hkdjkwjkqwkdjqwdkjqwdqwd}
 We call the sheaf $F\in \Sh_\Ab\bS$  admissible, iff
$$\uExt^1_{\Sh_\Ab\bS}(F,\uT)\cong \uExt^2_{\Sh_\Ab\bS}(F,\uT)\cong 0\ .$$
\end{ddd}

\begin{theorem}[Pontrjagin duality for Picard stacks]\label{dwiduwidqdqwdqwd}
If $A\in \Pic\bS$ is a Picard stack such that $H^{i}(A)$ and $D(H^{i}(A))$ are dualizable and admissible for $i=-1,0$,
then $A$ is dualizable.
\end{theorem}

\subsection{Admissible groups}

\subsubsection{}

Pontrjagin duality for locally compact abelian  groups implies that the sheaf $\uG$
associated to a locally compact abelian  group $G$ is dualizable.
In order to apply Theorem \ref{dwiduwidqdqwdqwd} 
to Picard stacks $A\in \Pic\bS$ whose sheaves $H^{i}(A)$, $i=-1,0$ are represented by locally compact abelian groups we must know for a locally compact  group $G$ whether 
the sheaf $\uG$ is admissible.
\begin{ddd}
We call a locally compact group $G$ admissible, if the sheaf $\uG$ is admissible. 
\end{ddd}

\subsubsection{}

Admissibility of a locally compact group is a complicated property. Not every locally compact group
is admissible, e.g. a discrete group of the form
$\oplus_{n\in \nat}\Z/p\Z$ for some integer $p$ is not admissible (see \ref{zwdwqidqidwqpdopop}).

\subsubsection{}

Admissibility is a vanishing condition
$$\uExt^1_{\Sh_\Ab\bS}(\uG,\uT)\cong \uExt^2_{\Sh_\Ab\bS}(\uG,\uT)\cong 0\ .$$
This condition depends on the site $\bS$.  In the present paper we shall also consider the sub-sites
$\bS_{lc-acyc}\subset \bS_{lc}\subset \bS$ of locally compact locally acyclic spaces and locally compact spaces.

For these sites the $\uExt$-functor commutes with restriction (we verify this property in \ref{odiwqpwqdqwdqwdq}).
Admissibility thus becomes a weaker condition on a smaller site. We will refine our notion of admissibility by saying that a group $G$ is admissible on the site $\bS_{lc}$ (or similarly for $\bS_{lc-acyc}$), if the corresponding restrictions of the extension sheaves vanish, e.g.
$$\uExt^1_{\Sh_\Ab\bS}(\uG,\uT)_{|\bS_{lc}}\cong \uExt^2_{\Sh_\Ab\bS}(\uG,\uT)_{|\bS_{lc}}\cong 0$$
in the case $\bS_{lc}$.

\subsubsection{}

Some locally compact abelian  groups are admissible on the site $\bS$. This applies e.g. to
finitely generated groups like $\Z, \Z/n\Z$, but also to  $\T^n$ and $\R^n$.

In the case of profinite groups $G$ we need the technical assumption that it
does not have too much two-torsion and three-torsion.

\begin{ddd}[\ref{two-three-cond}]
We say that the topological abelian group  $G$ satisfies the two-three condition, if
\begin{enumerate}
\item it does not admit $\prod_{n\in \nat} \Z/2\Z$ as a sub-quotient,
\item  the multiplication by $3$ on the component $G_0$ of the identity has
  finite cokernel.
\end{enumerate}
\end{ddd}

We can show that a profinite abelian group which satisfies the two-three
condition is also admissible.
We conjecture that it is possible to remove the condition using other techniques.

 A compact connected abelian group is divisible. Hence, if $p\in \nat$ is a prime, then
the multiplication $p:G\to G$ is set-theoretically surjective. 

\begin{ddd} We say that $G$ is locally topologically $p$-divisible ($ltd_p$), if this map has a continuous local section. The group $G$ is locally topologically divisible, if it is $ltd_p$ for all primes $p$.
\end{ddd}

For a general  connected compact group which satisfies the two-three condition we can
only show that it is admissible 
on $\bS_{lc-acyc}$.  If it is locally topologically divisible, then
it is admissible on the larger site $S_{lc}$.

We think that the two-tree condition and the restriction to a locally compact site is of technical nature. 
The condition of local compactness enters the proof since at one place we want to calculate the cohomology of the sheaf $\uZ$ on the space $A\times G$ using a K{\"u}nneth formula. For this reason we want that $A$ is locally compact.

As our counterexample above shows, a general (infinitely generated) discrete group is not
admissible unless we restrict to the site $\bS_{lc-acyc}$. We do not know if the restriction to the  site  $\bS_{lc-acyc}$ is really necessary for a general compact connected groups.

Using that the class of admissible groups is closed under finite products and extensions we get the following general theorem. 

\begin{theorem}[\ref{wdwquidioqwdopqwdq}]\label{wdwquidioqwdopqwdq1}
\begin{enumerate}
\item  If $G$ is a locally compact abelian group which satisfies the two-three
  condition, then it is admissible over $\bS_{lc-acyc}$.
\item  Assume that 
\begin{enumerate}
\item $G$ satisfies the two-three condition,
\item $G$ admits an open subgroup of the form $C\times \R^n$ with $C$ compact such that $G/C\times \R^n$ is finitely generated
\item the connected component of the identity of $G$ is locally topologically divisible.
\end{enumerate} 
Then $G$ is admissible over 
$\bS_{lc}$.
\end{enumerate}
\end{theorem}
The whole Section \ref{squfceweqwdqwd} is devoted to the proof of this theorem.
This section is very long and technical. A reader who is interested in the extension of Pontrjagin duality to topological abelian group stacks and applications to $T$-duality is advised to skip this section in a first reading and to take Theorem \ref{wdwquidioqwdopqwdq} for granted. 

\subsection{$T$-duality}

\subsubsection{} 

The aim of topological $T$-duality is to model the underlying topology of mirror symmetry in algebraic geometry and $T$-duality in string theory. For a detailed motivation we refer to \cite{math.GT/0501487}. The objects of topological $T$-duality over a space $B$ are pairs $(E,G)$ of a $\T^n$-principal bundle $E\to B$ and a gerbe $G\to E$ with band
$\T$. A gerbe with band $\T$ over a space $E$ is a map of stacks $G\to E$ which is locally isomorphic to $\cB \T_{|E}\to \uE$. Topological $T$-duality associates to $(E,G)$ dual pairs $(\widehat E,\widehat G)$.
For precise definitions we refer to \cite{math.GT/0501487} and \ref{uqdidwqwqdwqdwqd}.

\subsubsection{}

The case $n=1$ ($n$ is the dimension of the fibre of $E\to B$) is quite easy to understand (see \cite{MR2130624}).
In this case every pair admits a unique $T$-dual (up to isomorphism, of course).
The higher dimensional case is more complicated since on the one hand not
every pair admits a $T$-dual, and on the other hand, in general a  $T$-dual is
not unique, compare \cite{MR2130624,math.GT/0501487}.

\subsubsection{}

The general idea is that the construction of  a $T$-dual pair of $(E,G)$ needs the choice of an additional structure.
This structure might not exist, and in this case there is no $T$-dual. On the other hand, if the additional structure exists, then it might not be unique, so that the $T$-dual is not unique, too.

\subsubsection{}

The additional structure in \cite{math.GT/0501487} was the extension of the pair to a $T$-duality triple.
We review this notion in \ref{uqdidwqwqdwqdwqd}.

One can also  interpret the approach to $T$-duality via non-commutative topology \cite{MR2116734}, \cite{MR2222224}
in this way. The gerbe $G\to E$ determines an isomorphism class of a bundle of compact operators (by equality of the Dixmier-Douady classes). The additional structure which determines a $T$-dual is an $\R^n$-action on this bundle of compact operators which lifts the $\T^n$-action on $E$. Let us call a bundle of algebras of compact operators on $E$ together with such a $\R^n$-action a dynamical triple. 

The precise relationship between $T$-duality triples and dynamical triples is studied in the thesis
of Ansgar Schneider \cite{SchneiderThesis}.

\subsubsection{}\label{hjhjdqdqdwdd}

The initial motivation of the present paper was a third choice for the additional structure of the pair which came to live in a discussion with T. Pantev in spring 2006.
It was motivated by the analogy with some constructions in algebraic geometry, see e.g.\cite{math.AG/0306213}.

The starting point is the observation that a $T$-dual exists if and only if
the restriction of the gerbe $G\to E$ to $E|_{B^{(1)}}$, the restriction of
$E$ to a one-skeleton of $B$, is trivial. In particular, the restriction of
$G$ to the fibres of $E$  has to be trivial. Of course, the $T^n$-bundle
$E\to B$ is locally trivial on $B$, too. Therefore, locally on the base $B$,
the sequence of maps
$G\to E\to B$ is equivalent to
$$(\cB \uT\times \uT^n)_{|B}\to \uT^n_{|B}\to \uB_{|B}\ ,$$
where we identify spaces over $B$ with the corresponding sheaves.
The stack $(\cB \uT\times \uT^n)_{|B}$ is a  Picard stack.

Our proposal for the additional structure on $G\to E\to B$ is that of a torsor over the
Picard stack $(\cB \uT\times \uT^n)_{|B}$.

\subsubsection{}

In order to avoid the definition of a torsor over a group object in the two-categorial world of stacks
we use the following trick (see \ref{hqdwjhqdwqjdwqdqw} for details). Note that 
a torsor $X$ over an abelian group $G$ can equivalently be described as an extension of abelian groups
$$0\to G\to U\stackrel{p}{\to} \Z\to 0$$ so that $X\cong p^{-1}(1)$.
We use the same trick to describe a torsor over the Picard stack
$(\cB \uT\times \uT^n)_{|B}$ as an extension
$$(\cB \uT\times \uT^n)_{|B}\to U\to \uZ_{|B}$$
of Picard stacks.

\subsubsection{}

The sheaf of sections of $E\to B$ is a torsor over $\uT^n$.
Let $$0\to \uT^n_{|B}\to \cE\to \uZ_{|B}\to 0$$
be the corresponding extension of sheaves of abelian groups.
The filtration (\ref{hdjqwhdqwjdqwdqwdd}) of the Picard stack $U$
has the form
$$\cB \uT_{|B}\to U\to \cE\ ,$$
and locally on $B$
\begin{equation}\label{hqedjqdhdhdwhjdqwd6674133}
U\cong (\cB\uT \times \uT^n\times \uZ)_{|B}\ .
\end{equation}

\subsubsection{}

The proposal \ref{hjhjdqdqdwdd} was motivated by the hope that the Pontrjagin dual 
$D(U)$ of $U$ determines the $T$-dual pair. In fact, this can not be true directly since
the structure of $D(U)$ (this uses admissibility of $\T^n$ and $\Z^n$)
is given locally on $B$ by
$$D(U)\cong (\cB \uT\times \cB \uZ^n\times \uZ)_{|B}\ .$$
Here
the factor $\uZ$ in (\ref{hqedjqdhdhdwhjdqwd6674133}) gives rise to $\cB \uT$, the factor $\cB \uT$ yields $\uZ$, and $\uT^n$ yields
$\cB \uZ^n$ according to the rules (\ref{ttt54}) and (\ref{ttt56}).
The problematic factor is $\cB \Z^n$ in a place where we expect a factor $\uT^n$.
The way out is to interpret the gerbe
$\cB \uZ^n$ as the gerbe of $\R^n$-reductions of the trivial  principal bundle
$\T^n\times B\to B$.

\subsubsection{}\label{dqhdjqdhwqdqwdqdw}

We have canonical isomorphisms  $H^{0}(D(U))\cong D(H^{-1}(U))\cong D(\uT_{|B})\cong \uZ_{|B}$
and let $D(U)_1\subset D(U)$ be the pre-image of $\{1\}_{|B}\subset \uZ_{|B}$ under the natural map
$D(U)\to H^0(D(U))$. The quotient
$D(U)_1/(\cB \uT)_{|B}$ 
is a $\Z^n$-gerbe over $B$ which we interpret as the gerbe of $\R^n$-reductions of a
$\T^n$-principal bundle $\widehat E\to B$ which is well-defined up to unique isomorphism.
The full structure of $D(U)$ supplies the dual gerbe $\widehat G\to \widehat E$.
The details of this construction are explained in \ref{hdjhdqjdwdqwdqwdwqdwqdwqd54433}.

\subsubsection{}

In this way we use Pontrjagin duality of Picard stacks in order to construct a $T$-dual pair to $(E,G)$.
Schematically the picture is
$$(E,G)\stackrel{1}{\leadsto} U \stackrel{2}{\leadsto} D(U) \stackrel{3}{\leadsto} (\widehat E,\widehat G)\ ,$$
where the steps are as follows:
\begin{enumerate}
\item choice of the structure on $G$ of a torsor over $(\cB \uT\times \uT^n)_{|B}$
\item Pontrjagin duality $D(U):=\uHom_{\Pic\bS}(U,\cB \uT)$
\item Extraction of the dual pair from $D(U)$ as explained in \ref{dqhdjqdhwqdqwdqdw}.
\end{enumerate}

\subsubsection{}

Consider a pair $(E,G)$.
In Subsection \ref{hdjhdqjdwdqwdqwdwqdwqdwqd54433} we provide two constructions:
\begin{enumerate}
\item The construction $\Phi$ starts with the choice of a $T$-duality triple extending $(E,G)$ and constructs
the structure of a  torsor over $(\cB \uT\times \uT^n)_{|B}$ on $G$.
\item The construction $\Psi$ starts with the structure on $G$ of a   torsor over $(\cB \uT\times \uT^n)_{|B}$ and constructs an extension of $(E,G)$ to a $T$-duality triple.
\end{enumerate}
Our main result Theorem \ref{main} asserts that the constructions $\Psi$ and $\Phi$ are inverses to each other. In other words, the theories of topological  $T$-duality via Pontrjagin duality of Picard stacks and via $T$-duality triples are equivalent.

\section{Sheaves of Picard categories}\label{uwqdiwqdwdqdqdqwdwq4334}
\subsection{Picard categories}
 
\subsubsection{}
In a cartesian closed category one can define the notion of a group object in the standard way.
It is given by an object, a multiplication, an identity, and  an inversion morphism. The group axioms can be written as a collection of commutative diagrams involving these morphisms.

Stacks on a site $\bS$ form a two-cartesian closed two-category.  A group object in a two-cartesian closed two-category is again given by an object and the multiplication, identity and inversion morphisms. In addition each of the commutative diagrams from the category case is now filled by a two-morphisms. These two-morphisms must satisfy higher associativity relations. 

Instead of writing out all  these relations we will follow the exposition of Deligne SGA4, expos\'e XVIII \cite{deli}, which gives a rather effective way of working with group objects in two-categories.
We will not be interested in the most general case. For our purpose it suffices to consider a notion which includes sheaves of abelian groups and gerbes with band given by a sheaf of abelian groups. 
We choose to work with sheaves of strictly commutative Picard categories. 

\subsubsection{}

 Let $C$ be a category, $$F:C\times C\to C$$ be a bi-functor, and $$\sigma:F(F(X,Y),Z)\stackrel{\sim}{\to} F(X,F(Y,Z))$$
 be a natural isomorphism of tri-functors.
 \begin{ddd}
 The pair $(F,\sigma)$ is called an associative functor if the following holds.
  For every family $(X_i)_{i\in I}$ of objects of $C$ let $e:I\to M(I)$ denote the canonical map into the free monoid (without identity) on $I$. We require the existence of a map
 $\underline{F}:M(I)\to \Ob(C)$ and isomorphisms $a_i:\underline{F}(e(i))\stackrel{\sim}{\to}X_i$, and isomorphisms $a_{g,h}:\underline{F}(gh)\stackrel{\sim}{\to} F(\underline{F}(g),\underline{F}(h))$ such that the following diagram commutes:
 \begin{equation}\label{er3}\xymatrix{\underline{F}(f(gh))\ar@{=}[d]\ar[r]^{\hspace*{-0.5cm}a_{f,gh}}&F(\underline{F}(f),\underline{F}(gh))\ar[r]^{\hspace*{-0.5cm}a_{g,h}}&F(\underline{F}(f),F(\underline{F}(g),\underline{F}(h))\\\underline{F}((fg)h)\ar[r]_{\hspace*{-0.5cm}a_{fg,h}}&F(\underline{F}(fg),\underline{F}(h))\ar[r]_{\hspace*{-0.5cm}a_{f,g}}&F(F(\underline{F}(f),\underline{F}(g)),\underline{F}(h))\ar[u]^{\sigma}}\ .\end{equation} 
  \end{ddd}

\subsubsection{}
Let $(F,\sigma)$ be as above. Let in addition be  given a natural transformation of bi-functors
 $$\tau:F(X,Y)\to F(Y,X)\ .$$
 \begin{ddd}
 $(F,\sigma,\tau)$ is called a commutative and associative functor if the following holds.
  For every  family $(X_i)_{i\in I}$ of objects in $C$ let $e:I\to N(I)$ denote the canonical map into the free abelian monoid (without identity) on $I$. We require that there exist $\underline{F}:N(I)\to C$, isomorphisms $a_i:\underline{F}(e(i))\to X_i$ and isomorphisms $a_{g,h}:\underline{F}(gh)\stackrel{\sim}{\to} F(\underline{F}(g),\underline{F}(h))$ such that (\ref{er3}) and 
  \begin{equation}\label{er4}
  \xymatrix{\underline{F}(gh)\ar@{=}[d]\ar[r]^{\hspace*{-0.5cm}a_{g,h}}& F(\underline{F}(g),\underline{F}(h))\ar[d]^\tau\\\underline{F}(hg)\ar[r]_{\hspace*{-0.5cm}a_{h,g}}&F(\underline{F}(h),\underline{F}(g))}
  \end{equation}
  commute.
 \end{ddd}
 
\subsubsection{}
 We can now define the notion of a strict Picard category. Instead of $F$ will use the symbol "$+$".
 \begin{ddd}
 A Picard category is a groupoid $P$ together with a bi-functor $+:P\times P\to P$ and natural isomorphisms $\sigma$ and $\tau$ as above such that $(+,\sigma,\tau)$ is a commutative and associative functor,
 and such that for all $X\in P$ the functor $Y\mapsto X+Y$ is an equivalence of categories.
 \end{ddd}

\subsection{Examples of Picard categories}

\subsubsection{}\label{tzezwd}

Let $G$ be an abelian group. We consider $G$ as a category with only identity morphisms.
The addition $+:G\times G\to G$ is a bi-functor. For $\sigma$ and $\tau$ we choose
the identity transformations. Then $(G,+,\sigma,\tau)$ is a strict Picard category which we will denote again by $G$.

 \subsubsection{}\label{tzezwd1}

Let us now consider the abelian group $G$ as a category $BG$ with one object $*$ and $\Mor_{BG}(*,*):=G$. We let $+:BG\times BG\to BG$ be the bi-functor which acts as addition
$\Mor_{BG}(*)\times \Mor_{BG}(*)\to \Mor_{BG}(*)$. For $\sigma$ and $\tau$ we again choose the identity transformations. We will denote this strict Picard category  $(BG,+,\sigma,\tau)$  shortly by  $BG$.

\subsubsection{}\label{tzezwd2}

Let $G,H$ be abelian groups. We define the category $\EXT(G,H)$ as the category of short exact sequences $$\cE:0\to H\to E\to G\to 0\ .$$ Morphisms in $\EXT(G,H)$ are isomorphisms of complexes which reduce to the identity on $G$ and $H$.  We define the bi-functor
$+:\EXT(G,H)\times \EXT(G,H)\to \EXT(G,H)$ as the Baer addition
$$\cE_1+\cE_2:0\to H\to F\to G\to 0\ ,$$
where $F:=\tilde F/H$, $\tilde F\subset E_1\times E_2$ is the pre-image of the diagonal in $G\times G$, and the action of  $H$ on $\tilde F$ is induced by the anti-diagonal action on $E_1\times E_2$. The associativity morphism
$\sigma:(\cE_1+\cE_2)+\cE_3\to \cE_1+(\cE_2+\cE_3)$ is induced by the canonical isomorphism
$(E_1\times E_2)\times E_3\to E_1\times(E_2\times E_3)$. Finally, the transformation
$\tau:\cE_1+\cE_2\to \cE_2+\cE_1$ is induced by the flip $E_1\times E_2\to E_2\times E_1$.
The triple $(+,\sigma,\tau)$ defines on $\EXT(G,H)$ the structure of a Picard category.

\subsubsection{}\label{weufiwe}

Let $G$ be a topological abelian group, e.g. $G:=\T$. On a space $B$ we consider the
category of $G$-principal bundles $\cB G(B)$. Given two
$G$-principal bundles $Q\to B$, $P\to B$ we can define a new $G$-principal bundle
$Q\otimes_G P:=Q\times_B P/G$. The quotient is taken by the anti-diagonal action, i.e. we identify
$(q,p)\sim (qg,pg^{-1})$. The $G$-principal structure on $Q\otimes_G P$ is induced by the action $[q,p]g:=[q,pg]$. 

We define the associativity morphism
$$\sigma:(Q\otimes_G P)\otimes_G R\to Q\otimes_G(P\otimes_G R)$$
as the map $[[q,p],r]\to [q,[p,r]]$. Finally we let $\tau:Q\otimes_GP\to P\otimes_GQ$ be given by
$[q,p]\to [p,q]$.

We claim that $(\otimes_G,\sigma,\tau)$ is a strict Picard category.
Let $I$ be a collection of objects of $\cB G(B)$.
We choose an ordering of $I$. Then we can write each element of $j\in N(I)$ in a unique way as ordered product $f=P_1P_2\dots P_r$ with $P_1\le P_2 \dots \le P_r$.
We define
$\uF(f):=P_1\otimes_G (P_2\otimes_G \dots (\dots\otimes_G P_r)\dots)$.
We let $a_P:\uF(P)\stackrel{\sim}{\to} P$, $P\in I$, be the identity.
Furthermore, for $f,g\in N(I)$ we define
$a_{f,g}:\uF(fg)\stackrel{\sim}{\to} \uF(f)\otimes_B \uF(g)$ as the map
induced the permutation which puts the factors of $f$ and $g$ in order (and so
that the order to repeated factors is not changed). Commutativity of
(\ref{er3}) is now clear. In order to check the commutativity of (\ref{er4})
observe that $\tau:Q\otimes_GQ\to Q\otimes_G Q$ is the identity, as $[qg,q]
\mapsto [q,qg]=[qg,q]$.

\subsection{Picard stacks and examples}
 
\subsubsection{}
We consider a site $\bS$. A prestack on $\bS$ is a lax associative functor $P:\bS^{op}\to \groupoids$,
where from now on by a groupoid we understand an essentially small category\footnote{Later, e.g. in \ref{jqwddqkwdnwd}, we need that the isomorphism classes in a groupoid form a set.} in which all morphisms are isomorphisms.
Lax associative means that as a part of the data  for each composeable pair of maps $f,g$ in $ \bS$ 
there is an isomorphism of functors $l_{f,g}:P(f)\circ
P(g)\stackrel{\sim}{\to} P(g\circ f)$, and these isomorphisms satisfy higher
associativity relation. A prestack is a stack if it satisfies  the usual
descent conditions one the level of objects and morphisms. For a reference of the language of stacks see e.g. \cite{MR2223406}. 

\subsubsection{} 

\begin{ddd}\label{zduqwgdwqdqdwd}
 A Picard stack $P$ on $\bS$ is a stack $P$ together with an operation $+:P\times P\to P$ and transformations
 $\sigma$, $\tau$ which induce for each $U\in \bS$ the structure of a Picard category on $P(U)$.
 \end{ddd}

\subsubsection{}\label{diewodewdw1}
Let $\Sh_\Ab \bS$ denote the category of sheaves of abelian groups on $\bS$.
Extending the example \ref{tzezwd} to sheaves we can view each object
$\cF\in \Sh_\Ab\bS$ as a Picard stack, which we will again denote by $\cF$.

\subsubsection{}\label{diewodewdw}
We can also extend the example \ref{tzezwd1} to sheaves. In this way every sheaf $\cF\in \Sh_\Ab \bS$ gives rise to a Picard stack $\cB\cF$.

\subsubsection{} \label{catExt}
 
We now extend the example \ref{tzezwd2} to sheaves. First of all note that one can define a Picard category $\EXT(\cG,\cH)$ of extensions of sheaves $0\to \cH\to \cE\to \cG\to 0$ as a direct generalization of \ref{tzezwd2}. 
Then we define a prestack $\cEXT(\cG,\cH)$ which associates to
$U\in \bS$ the Picard category $\cEXT(\cG,\cH)(U):=\EXT(\cG_{|U},\cH_{|U})$.
One checks that $\cEXT(\cG,\cH)$ is a stack.
 
\subsubsection{}   \label{d6327d2d}

Let $\cG\in \Sh_\Ab \bS$ be a sheaf of abelian groups. We consider $\cG$ as a group object in the category $\Sh\bS$ of sheaves of sets on $\bS$. A $\cG$-torsor is an object $\cT\in \Sh\bS$ with an action of $\cG$ such that $\cT$ is locally isomorphic to $\cG$.
We consider the category of $\cG$-torsors $\Torsor(\cG)$ and their isomorphisms. We define the functor
$+:\Torsor(\cG)\times \Torsor(\cG)\to \Torsor(\cG)$ by 
$\cT_1+\cT_2:=T_1\times \cT_2/\cG$, where we take the quotient by the anti-diagonal action. The structure of a $\cG$-torsor is induced by the action of $\cG$ on the second factor $\cT_2$. We let $\sigma$ and  $\tau$ be induced by the associativity transformation of the cartesian product and the flip, respectively. With these structures $\Torsor(\cG)$ becomes a Picard category.

We define a Picard stack $\cTorsor(\cG)$ by localization, i.e. we set
$\cTorsor(\cG)(U):=\Torsor(\cG_{|U})$.

We have a canonical map of Picard stacks
$$U:\cEXT(\uZ,\cG)\to \cTorsor(\cG)$$
which maps the extension
$\cE:0\to \cG\to \cE\stackrel{\pi}{\to} \uZ\to 0$ to the $G$-torsor $U(\cE)=\pi^{-1}(1)$.
Here $\uZ$ denotes the constant sheaf on $\bS$ with value $\Z$.
One can check that $U$ is an equivalence of Picard stacks (see \ref{iiiss} for precise definitions).

\subsubsection{}\label{jksakasdioioeqweqw}

The example \ref{weufiwe} has a sheaf theoretic interpretation.
We assume that $\bS$ is some small subcategory of the category of topological
spaces which is closed under taking open subspaces, and such that the topology
is induced by the coverings of spaces by families of open subsets.

Let $G$ be a topological abelian group.
We have a Picard prestack ${}^p\cB G$ on $\bS$ which associates to each space $B\in \bS$ the Picard category ${}^p\cB G(B)$ of $G$-principal bundles on $B$. 
As in \ref{weufiwe} the monoidal structure on ${}^p\cB G(B)$ is given by the tensor product of $G$-principal bundles (this uses that $G$ is abelian).
 We now define
$\cB G$ as the stackification of ${}^p\cB G$.

The topological group $G$ gives rise to a sheaf $\uG\in\Sh_\Ab\bS$ which associates to  $U\in \bS$ the abelian group  $\uG(U):=C(U,G)$. We have a canonical transformation
$$\Gamma: \cB G\to \cTorsor(\uG)\ .$$
It is induced by a transformation ${}^p\Gamma:{}^p\cB G\to   \cTorsor(\uG)$.
We describe the functor ${}^p\Gamma_B:{}^p\cB G(B)\to \cTorsor(\uG)(B)$ for all  $B\in \bS$.
Let $E\to B$ be a $G$-principal bundle.
Then we define ${}^p\Gamma_B(E)\in \Torsor(\uG_{|B})$ to be  the sheaf which associates to each
$(\phi:U\to B)\in \bS/B$ the set of sections of $\phi^*E\to U$.

One can check that $\Gamma$ is an equivalence of Picard stacks (see \ref{iiiss} for precise definitions).

\subsection{Additive functors}\label{zqdiqdqwqwodiowqdqwdwqd}

\subsubsection{}
We now discuss the notion of an additive functor between Picard stacks $F:P_1\to P_2$.
\begin{ddd}
An additive functor between Picard categories is a functor $F:P_1\to P_2$ and a natural transformation $F(x+y)\stackrel{\sim}{\to} F(x)+F(y)$ such that the following diagrams commute:
\begin{enumerate}
\item $$\xymatrix{F(x+y)\ar[r]\ar[d]^{F(\tau)}&F(x)+F(y)\ar[d]^\tau\\ F(y+x)\ar[r]&F(y)+F(x)}$$
\item $$\xymatrix{F((x+y)+z)\ar[r]\ar[d]^{F(\sigma)}&F(x+y)+F(z)\ar[r]&(F(x)+F(y))+F(z)\ar[d]^\sigma\\F(x+(y+z))\ar[r]&F(x)+F(y+z))\ar[r]&F(x)+(F(y)+F(z))}$$
\end{enumerate}
\end{ddd}

\begin{ddd}
An isomorphism between additive functors $u:F\to G$ is an isomorphism of functors such that
$$\xymatrix{F(x+y)\ar[r]^{u_{x+y}}\ar[d]&G(x+y)\ar[d]\\F(x)+F(y)\ar[r]^{u_x+u_y}&G(x)+G(y)}$$
commutes.
\end{ddd}

\begin{ddd}
We let $\Hom(P_1,P_2)$ denote the groupoid of additive functors from $P_1$ to $P_2$.
By $\cPic$ we denote the two-category of Picard categories.
\end{ddd}

\subsubsection{}

The groupoid $\Hom(P_1,P_2)$ has a natural structure of a Picard category.
We set $$(F_1+F_2)(x):=F_1(x)+F_2(x)$$ and define
the transformation $$(F_1+F_2)(x+y)\to (F_1+F_2)(x)+(F_1+F_2)(y)$$ such that
{\scriptsize $$\xymatrix{(F_1+F_2)(x+y)\ar@{=}[d]\ar[rr]&&(F_1+F_2)(x)+(F_1+F_2)(y)\ar@{=}[d]\\
F_1(x+y)+F_2(x+y)\ar[r]&(F_1(x)+F_1(y))+(F_2(x)+F_2(y))\ar[r]^{\alpha\circ\tau\circ\alpha}&F_1(x)+F_2(x)+F_1(y)+F_2(y)}$$}
commutes. The associativity and commutativity constraints are induced by those of $P_2$.

\subsubsection{}\label{iiiss} 

Let now $P_1,P_2$ be two Picard stacks on $\bS$.
\begin{ddd}
An additive functor $F:P_1\to P_2$ is a morphism of stacks together with a two-morphism
$$\xymatrix{P_1\times P_1\ar[r]^+\ar[d]^{F\times F}\ar@{=>}[dr]&P_1\ar[d]^F\\P_2\times P_2\ar[r]^+&P_2}$$
which induces for each $U\in \bS$
 an additive functor $F_U:P_1(U)\to P_2(U)$ of Picard categories. An isomorphism between additive
functors $u:F_1 \to F_2$ is a two-isomorphism of morphisms of stacks which induces for each $U\in \bS$ an isomorphism of additive functors $u_U:F_{1,U}\to F_{2,U}$. 
\end{ddd}
As in the case of Picard categories, the  additive functors between Picard stacks form again a Picard category.
\begin{ddd}
We let $\Hom(P_1,P_2)$ denote the groupoid of additive functors from $P_1,P_2$.We get a two-category $\cPic(\bS)$ of Picard stacks on $\bS$.
\end{ddd}

\begin{ddd}
We let $\HOM(P_1,P_2)$ denote the Picard category of additive functors between Picard stacks $P_1$ and $P_2$.
\end{ddd}

\subsubsection{}

Let $P,Q$ be Picard stacks on the site $\bS$.
By localization we define a Picard stack $\uHOM(P,Q)$.
\begin{ddd}\label{lwqdjkqwdjwqq89812}
The Picard stack $\uHOM(P,Q)$ is the sheaf of Picard categories given by
$$\uHOM(P,Q)(U):=\HOM(P_{|U},Q_{|U})\ .$$
\end{ddd}
A priori this describes a prestack, but one easily checks the stack conditions.

\subsection{Representation of Picard stacks by complexes of sheaves of groups}

\subsubsection{}

There is an obvious notion of a Picard prestack.
Furthermore, there is an associated Picard stack construction $a$ such that for a Picard prestack $P$ and a Picard stack $Q$ we have a natural equivalence 
\begin{equation}\label{ufweiucfwe}
\Hom(aP,Q)\cong \Hom(P,Q)
\end{equation}

\subsubsection{}

Let $\bS$ be a site.
We follow Deligne, SGA 4.3, Expose XVIII, \cite{deli}. Let $C(\bS)$ be the two-category of complexes of sheaves of abelian groups
$$K:0\to K^{-1} \stackrel{d}{\to} K^0\to 0$$
which live in degrees $-1,0$. Morphisms in $C(\bS)$ are morphisms of complexes, and two-isomorphisms are homotopies between morphisms.

To such a complex we associate a Picard stack
$\ch(K)$ on $\bS$ as follows. We first define a Picard prestack $\pch(K)$.
\begin{enumerate}
\item  For $U\in \bS$ we define the set of objects of $\pch(K)(U)$ as $K^0(U)$.
\item For $x,y\in K^0(U)$ we define the set of morphisms  by
$$\Hom_{\pch(K)(U)}(x,y) :=\{f\in K^{-1}(U)|df=x-y\}\ .$$
\item The composition of morphisms is addition in $K^{-1}(U)$.
\item The functor $+:\pch(K)(U)\times \pch(K)(U)\to \pch(K)(U)$ is given by addition of objects and morphisms.  The associativity and the commutativity constraints are the identities.
\end{enumerate}
\begin{ddd}\label{weiofwe}
We define $\ch(K)$ as the Picard stack associated to the prestack $\pch(K)$, i.e. $\ch(K):=a\pch(K)$.
\end{ddd}

\subsubsection{}\label{jqwddqkwdnwd}

If $P$ is a Picard stack, then we can define a presheaf ${}^p H^0(P)$ which associates to $U\in \bS$ the group of  isomorphism classes of $P(U)$. We let $H^0(P)$ be the associated sheaf.
 Furthermore we have a sheaf $H^{-1}(P)$ which associates to $U\in \bS$ the group of automorphisms $\Aut(e_U)$, where $e_U\in P(U)$ is a choice of a neutral object with respect to $+$. Since the neutral object $e_U$ is determined
up to unique isomorphism, the group $\Aut(e_U)$ is also determined up to unique isomorphism.
Note that $\Aut(e_U)$ is abelian, as shows the following diagram
\begin{equation*}
  \begin{CD}
    e+e @>{\id+f}>> e+e @>{g+\id}>> e+e\\
    @VV{\cong}V @VV{\cong}V @VV{\cong}V\\
    e @>{f}>> e @>{g}>> e
  \end{CD}
\end{equation*}
and the fact that $(\id+f)\circ (g+\id)=g+f=(g+\id)\circ(\id+f)$.

If $V\to U$ is a morphism in $\bS$, then we have a unique isomorphism $f:(e_U)_{|V}\stackrel{\sim}{\to }e_V$ which induces a group homomorphism
$f\circ \dots\circ f^{-1}:\Aut(e_U)_{|V}\to \Aut(e_V)$. 
This homomorphism is the structure map of the  
sheaf $H^{-1}(P)$ which is determined up to unique isomorphism in this way.

We now observe that 
\begin{equation}\label{uefiewfewf}
H^{-1}(\ch(K))\cong H^{-1}(K)\ ,\quad H^0(\ch(K))\cong H^0(K)\ .
\end{equation}

\subsubsection{}

A morphism of complexes $f:K\to L$ in $C(\bS)$, i.e. a diagram
$$\xymatrix{0\ar[r]&K^{-1}\ar[d]^{f^{-1}}\ar[r]^{d_K}&K^0\ar[d]^{f^0}\ar[r]&0\\0\ar[r]&L^{-1}\ar[r]^{d_L}&L^0\ar[r]&0
}\ ,$$ induces an additive functor of Picard prestacks $\pch(f):\pch(K)\to \pch(L)$ and, by the functoriality of the associated stack construction, an additive functor
$$\ch(f):\ch(K)\to \ch(L)\ .$$
In view of (\ref{uefiewfewf}), it is an equivalence if and only if $f$ is a
quasi isomorphism.

\subsubsection{}

We consider two morphisms $f,g:K\to L$ of complexes in $C(\bS)$.
A homotopy $H:f\to g$ is a morphism of sheaves $H:K^0\to L^{-1}$ such that
$f^0-g^0=d_L\circ H$ and $f^{-1}-g^{-1}= H\circ d_K$. It is easy to see that
$H$ induces an  isomorphism of additive functors
$\pch(H):\pch(f)\to \pch(g)$ and therefore
$$\ch(H):\ch(f)\to \ch(g)$$
by the associated stack construction and (\ref{ufweiucfwe}).

One can show that the isomorphisms $\pch(f)\to \pch(g)$ correspond precisely to homotopies
$H:f\to g$.  This implies that the  morphisms $\ch(f)\to \ch(g)$ also correspond bijectively to these homotopies.

\subsubsection{}

We have the following result.
\begin{lem}\cite[1.4.13]{deli}\label{etzwwe4}
\begin{enumerate}
\item For every Picard stack $P$ there exists a complex $K\in C(\bS)$ such that $P\cong \ch(K)$.
\item For every additive functor $F:\ch(K)\to \ch(L)$ there exists a quasi isomorphism $k:K^\prime\to K$ and a morphism $l:K^\prime\to L$ in $C(\bS)$ such that
$F\cong \ch(l)\circ \ch(k)^{-1}$.
\end{enumerate}
\end{lem}

\subsubsection{}

Let $\cPic^\flat(\bS)$ be the category of Picard stacks obtained from the
two-category $\Pic\bS$ by identifying isomorphic additive functors.
\begin{prop}\cite[1.4.15]{deli} \label{zudw77823d}
The construction $\ch$ gives an equivalence of categories
$$D^{[-1,0]}(\Sh_\Ab\bS)\to \cPic^\flat(\bS)\ ,$$ where $D^{[-1,0]}(\Sh_\Ab\bS)$ is the full subcategory of $D^+(\Sh_\Ab\bS)$ of objects whose cohomology is trivial in degrees $\not\in \{-1,0\}$.
\end{prop}

\subsubsection{}

\begin{lem}\cite[1.4.16]{deli}  Let $K,L\in C(\bS)$ and assume that $L^{-1}$ is injective.
\begin{enumerate}
\item $\pch(L)$ is already a stack.
\item For every $F\in \Hom_{\cPic(\bS)}(\ch(K),\ch(L))$ there exists a morphism
 $f\in \Hom_{C(\Sh_\Ab\bS)}(K,L)$ such that $\ch(f)\cong F$.
\end{enumerate}
\end{lem}

Let $C^\prime(\bS)\subseteq C(\bS)$ denote the full sub-two-category of complexes $0\to K^{-1}\to K^0\to 0$ with
$K^{-1}$ injective.

\begin{lem}\label{cprimeequi}\cite[1.4.17]{deli}   
The construction $\ch$ induces an equivalence of the two-categories
$\cPic(\bS)$ and $C^\prime(\bS)$. 
\end{lem}

\subsubsection{}
 
Finally we give a characterization of the Picard stack
$\uHom(\ch(K),\ch(L))$. 
\begin{lem}\cite[1.4.18.1]{deli} \label{hashjdjshadui} Assume that $L^{-1}$ is injective.
Then we have an equivalence
$$\ch(\tau_{\le 0} \uHom_{\Sh_\Ab\bS}(K,L))\stackrel{\sim}{\to } \uHOM_{\cPic(\bS)}(\ch(K),\ch(L))\ .$$
\end{lem}

\subsubsection{}
\begin{lem}\label{zuedzwe777}Let $K,L\in C(\bS)$.
Then we have  an isomorphism
$$H^{i}(\uHOM_{\cPic(\bS)}(\ch(K),\ch(L)))\cong R^i\uHom_{\Sh_\Ab\bS}(K,L)$$ for $i=-1,0$.
\end{lem}
\proof
First observe that by the discussion above the left hand side, and by the definition of $R\uHom$ the right side both  only depend on the quasi-isomorphism type of the complex $L$ of length $2$. Without loss of
generality we can therefore assume that $L^{-1}$ is injective.

We now choose an injective resolution $I:0\to L^{-1}\to I^0\to I^1 \dots$ of $L$ starting with the choice of an embedding $L^0\to I^0$. Then we have  $\uHom(K,I)\cong R\uHom(K,L)$.  We now observe that
$H^i\uHom(K,I)\cong H^i\uHom (K,L)$ for $i=0,-1$.
While the case $i=-1$ is obvious, for $i=0$ observe that a $0$-cycle in
$\uHom(K,I)$ is a morphism of complexes and necessarily factors over $L\to I$.
We thus have for $i\in\{-1,0\}$
$$R^i\uHom(K,L)\cong H^i\uHom (K,L)\cong H^{i}(\uHom(\ch(K),\ch(L)))\ .$$
\hB

\subsubsection{}\label{gdhasgdaui}

For $A,B\in \Sh_\Ab\bS$ we have canonical isomorphisms
$$\Ext^2_{\Sh_\Ab\bS}(B,A)\cong R^0\Hom_{\Sh_\Ab\bS}(B,A[2])\cong \Hom_{D(\Sh_\Ab\bS)}(B,A[2])\ .$$
In the following we recall two eventually equivalent ways how an exact complex 
$$\cK:0\to A\to X
\to Y\to B\to 0$$ represents an  element $$Y(\cK)\in \Hom_{D(\Sh_\Ab\bS)}(B,A[2])\cong \Ext_{\Sh_\Ab\bS}^2(B,A)$$
(the letter $Y$ stands for Yoneda who investigated this construction first).
Let $\cK_A$ be the complex $$\cK_A:0\to X\to Y\to B\to 0\ ,$$
where $B$ sits in degree $0$. The obvious inclusion $\alpha:A[2]\to \cK_A$ induced by $A\to X$ is a quasi-isomorphism. Furthermore, we have a canonical map $\beta:B\to \cK_A$.
The element
$Y(\cK)\in \Hom_{D(\Sh_\Ab\bS)}(B,A[2])$ is by definition the composition
\begin{equation}\label{ygdeduz83e}
Y(\cK):B\stackrel{\beta}{\to}\cK_A\stackrel{\alpha^{-1}}{\to} A[2]\ .
\end{equation}

We can also consider the complex $\cK_B$ given by
$$0\to A\to X\to Y\to 0$$ where $A$ is in degree $-2$.
The projection $Y\to B$ induces a quasi-isomorphism
$\gamma:\cK_B\to B$. We furthermore have a canonical map $\delta:\cK_B\to A[2]$.
We consider the composition $Y^\prime(\cK)\in \Hom_{D(\Sh_\Ab\bS)}(B,A[2])$ 
\begin{equation}\label{gf6f67whrjwr}
Y^\prime(\cK):B\stackrel{\gamma^{-1}}{\to} \cK_B\stackrel{\delta}{\to} A[2]\ .
\end{equation}
\begin{lem}\label{soundso}
In $\Hom_{D(\Sh_\Ab\bS)}(B,A[2])$ we have the equality $$Y(\cK)=Y^\prime(\cK)\ .$$
\end{lem}
\proof
We consider the morphism of complexes $\phi:\cK_B\to \cK_A$ given by obvious maps in the diagram
$$\xymatrix{X\ar[r]&Y\ar[r]& B\\
A\ar[u]\ar[r]&X\ar[r]\ar[u]&Y\ar[u]}
 \ .$$
It fits into
$$\xymatrix{B\ar[r]^\beta\ar@{=}[d]&\cK_A&A[2]\ar@{=}[d]\ar[l]^\alpha\\
B&\cK_B\ar[u]^\phi\ar[l]^\gamma\ar[r]^\delta&A[2]}\ .$$
It suffices to show that the two squares commute.
In fact we will show that
$\phi$ is homotopic to $\beta\circ \gamma$ and $\alpha\circ \delta$.
Note that $\phi-\beta\circ \gamma$ is given by  
$$\xymatrix{X\ar[r]&Y\ar[r]& B\\
A\ar[u]\ar[r]&X\ar@{.>}[ul]^{\id}\ar[r]\ar[u]&Y\ar[u]^0\ar@{.>}[ul]^{0}}
 \ .$$
The dotted arrows in this diagram gives the zero homotopy of
this difference.

Similarly $\phi-\alpha\circ \delta$ is given by 
$$\xymatrix{X\ar[r]&Y\ar[r]& B\\
A\ar[u]^0\ar[r]&X\ar@{.>}[ul]^{0}\ar[r]\ar[u]&Y\ar[u]\ar@{.>}[ul]^{\id}}
 \ ,$$
 and we have again indicated the required zero homotopy. 
 \hB

\subsubsection{}

The equivalence between the two-categories $\cPic(\bS)$ and $C^\prime(\bS)$ (see \ref{cprimeequi}) allows us to classify equivalence classes of Picard stacks with fixed $H^i(P)\cong A_i$, $i=-1,0$, 
$A_i\in \Sh_\Ab\bS$. An equivalence of such Picard stacks is an equivalence which induces the identity on the cohomology.

\begin{lem}\label{ewhdwejdew8}
The set $\Ext_{\cPic(\bS)}(A_0,A_{-1})$ of equivalence classes of Picard
stacks $P$ with given isomorphisms
$H^i(P)\cong A_i$, $i=0,-1$, 
$A_i\in \Sh_\Ab\bS$ is in bijection with $\Ext^2_{\Sh_\Ab\bS}(A_0,A_{-1})$.  
\end{lem}
\proof
We define a map
\begin{equation}\label{phidef632423}
\phi:\Ext_{\cPic(\bS)}(A_0,A_{-1})\to \Ext^2_{\Sh_\Ab\bS}(A_0,A_{-1})
\end{equation}
as follows.

Consider an exact complex 
 $$\cK:0\to A_{-1}\to K^{-1}\to K^0\to A_0\to 0$$
 such that $K^{-1}$ is injective, and let
  $K\in C(\bS)$ be the complex $$0\to K^{-1}\to K^0\to 0$$ 
 Then we define
$\phi(\ch(K)):=Y(\cK)$.

We must show that $\phi$ is well-defined. Indeed, if $\ch(K)\cong \ch(L)$, then (see \ref{cprimeequi})
$K\cong L$ by an equivalence which induces the identity on the level of cohomology.
It follows that 
$Y(\cK)=Y(\cL)$. 

Since $\ch:C^\prime(\bS)\to \cPic^\flat(\bS) $ is an equivalence, hence in particular surjective on the level of equivalence classes, the map $\phi$ is well-defined.
Since every element of $\Ext^2_{\Sh_\Ab\bS}(A_{-1},A_0)$ can be written as $Y(\cK)$ for a suitable complex $\cK$ as above
we conclude that $\phi$ is surjective. If $\phi(\ch(K))\cong \phi(\ch(L))$,
then there exists $M\in C(\bS)$ with given isomorphisms $H^i(M)\cong A_i$ together with quasi-isomorphisms
$K\leftarrow M\rightarrow L$ inducing the identity on cohomology. But this diagram induces an equivalence
$\ch(K)\cong \ch(M)\cong \ch(L)$.
\hB

\subsubsection{}\label{def_of_BF2}

We continue with a further description of $\cB F$ for $F\in\Sh_\Ab \bS$,
compare \ref{def_of_BF}. Note that the sheaf of groups of automorphisms of
$\cB F$ is 
$F$, whereas the group of objects is trivial. Therefore we can represent $\cB
F$ as $\ch(F[1])$, where $F[1]$ is the complex $0\to F\to 0\to 0$ concentrated in
degree $-1$.

\subsubsection{}\label{dedwhedjewdn66}

Let $P\in \cPic(\bS)$ be a Picard stack and $G\subseteq H^{-1}(P)$. Let us assume that $P = \ch(K)$ for a suitable complex $K:0\to K^{-1}\to K^0\to 0$. We have a natural injection $G\hookrightarrow \ker(K^{-1}\to K^0)$. 
We consider the quotient $\bar K^{-1}$ defined by the exact sequence
$0\to G\to K^{-1}\to \bar K^{-1}\to 0$ and obtain a new Picard stack
$\bar P\cong \ch(\bar K)$, where $\bar K:0\to \bar K^{-1}\to K^0\to 0$.
The following diagram
represents a sequence of morphisms of Picard stacks
$$\cB G\to P\to \bar P$$
$$\xymatrix{0\ar[r]&G\ar[r]\ar[d]&0\ar[d]\ar[r]&0\\0\ar[r]&K^{-1}\ar[r]\ar[d]&K^0\ar[d]\ar[r]&0\\0\ar[r]&\bar K^{-1}\ar[r]&K^0\ar[r]&0}\ .$$ 
The Picard stack $\bar P$ can be considered as a quotient of $P$ by $\cB G$.
We will employ this construction in \ref{dcgbsasjkscsklioiowqdqwdw}.

\section{Sheaf theory on big sites of topological spaces}\label{hwdqwwqdqww}

\subsection{Topological spaces and sites}

\subsubsection{}\label{euwfiwef}
In this paper a topological space will always be compactly generated and Hausdorff. We will define categorical limits and colimits in the category of compactly generated Hausdorff spaces. Furthermore, we will equip mapping spaces $\Map(X,Y)$ with the compactly generated topology obtained from the compact-open topology. 
In this category we have the exponential law
$$\Map(X\times Y,Z)\cong \Map(X,\Map(Y,Z))\ .$$
By $\Map(X,Y)^\delta$ we will denote the underlying set.
For details on this {\em convenient category of topological spaces} we refer to \cite{MR0210075}.

\subsubsection{} \label{firfore}

The sheaf theory of the present paper refers to the Grothendieck site  $\bS$. The underlying category of $\bS$  is the category of compactly generated topological Hausdorff spaces. The covering families of a space $X\in \bS$  are coverings by families of open subsets.

We will also need the sites $\bS_{lc}$ and $\bS_{lc-acyc}$ given by the full subcategories of locally compact and locally compact locally acyclic spaces (see \ref{locacyintro}).

\subsubsection{}

We let $\Pr\bS$ and $\Sh\bS$ denote the category of presheaves and sheaves of sets on $\bS$. Then we have an adjoint pair of functors
$$i^\sharp:\Pr\bS\Leftrightarrow \Sh\bS:i$$
where $i$ is the inclusion of sheaves into presheaves, and $i^\sharp$ is the sheafification functor.
If $F\in \Pr\bS$, then sometimes  we will write $F^\sharp:=i^\sharp F$. 
 
As before, by $\Pr_\Ab\bS$ and $\Sh_\Ab\bS$ we denote the categories of presheaves and sheaves of abelian groups.

\subsection{Sheaves of topological groups}

\subsubsection{}

In this subsection we collect some general facts about sheaves generated by spaces and topological groups.
We formulate the results for the site $\bS$. But they remain true if one replaces $\bS$ by $\bS_{lc}$ or $\bS_{lc-acyc}$.

\subsubsection{}\label{hdjqhwdwqd}

Every object $X\in \bS$ represents a presheaf $\uX\in \Pr\bS$ of sets defined by
$$\bS\ni U\mapsto \uX(U):=\Hom_{\bS}(U,X)\in \Sets$$ on objects, and by
$$\Hom_{\bS}(U,V)\ni f\mapsto f^*:\uX(V)\to \uX(U)$$ 
on morphisms. 
\begin{lem}\label{zzzuzuzr}
$\uX$ is a sheaf. 
\end{lem}
\proof Straightforward.\hB
Note that a topology is called sub-canonical if all representable presheaves are sheaves. 
Hence $\bS$ carries a sub-canonical topology.

\subsubsection{}\label{uiqqwq4432}

We will also need the relative version.
For $Y\in \bS$ we consider the relative site $\bS/Y$ of spaces over $Y$. Its objects are morphisms $A\to Y$, and its morphisms are commutative diagrams
$$\xymatrix{A\ar[dr]\ar[rr]&&B\ar[dl]\\&Y&}\ .$$
The covering families of $A\to Y$ are induced from the coverings of $A$ by open subsets.

An object  $X\to Y\in \bS/Y$ represents the presheaf $\underline{X\to Y}\in \Pr\bS/Y$ (by a similar definition
as  in the absolute case \ref{hdjqhwdwqd}). The induced topology on $\bS/Y$ is again sub-canonical. In fact, we have the following Lemma.
\begin{lem}
For all $(X\to Y)\in \bS$ the presheaf
$\underline{X\to Y}$ is a sheaf.
\end{lem}
\proof 
Straightforward. \hB

\subsubsection{}

Let $I$ be a small category and $X\in \bS^I$. Then we have
$$\lim_{i\in I} \underline{X(i)}\cong \underline{\lim_{i\in I} X(i)}\ .$$
One can not expect a similar property for arbitrary colimits. But we have the following result.
\begin{lem}\label{ehewfudwqdwqdwq}
Let $I$ be a directed partially ordered set and $X\in \bS^I$ be a direct system of discrete sets such that $X(i)\to X(j)$ is injective for all $i\le j$.
Then the canonical map
$$\colim_{i\in I}\underline{X(i)}\to \underline{\colim_{i\in I} X(i)}$$ is an isomorphism. 
\end{lem}
\proof
Let ${}^p\colim_{i\in I}\underline{X(i)}$ denote the colimit in the sense of presheaves.
Then we have an inclusion
$${}^p\colim_{i\in I}\underline{X(i)}\hookrightarrow \underline{\colim_{i\in I}X(i)}$$
(this uses injectivity of the structure maps of the system $X$).
Since the target is a sheaf and sheafification preserves injections it induces an inclusion
\begin{equation}\label{nachtwededs}
\colim_{i\in I}\underline{X(i)}\hookrightarrow \underline{\colim_{i\in I} X(i)}\ .
\end{equation}
Let now $A\in \bS$ and $f\in \underline{\colim_{i\in I}X(i)}(A)=\Hom_{\bS}(A,\colim_{i\in I}X(i))$.
Since $\colim_{i\in I}X(i)$ is discrete and $f$ is continuous the family of subsets $\{f^{-1}(x)\subseteq A|x\in \colim_{i\in I} X(i)\}$ is an open covering of $A$ by disjoint open subsets.
The family
$$\prod_{x} f_{|f^{-1}(x)} \in \prod_x {}^p\colim_{i\in I}\underline{X(i)}(f^{-1}(x))$$
represents a section over $A$ of the sheafification $\colim_{i\in I}\underline{X(i)}$ of ${}^p\colim_{i\in I}\underline{X(i)}$ which of course maps to $f$ under the inclusion (\ref{nachtwededs}).
Therefore  (\ref{nachtwededs}) is also surjective.
\hB

\subsubsection{}\label{multbla}

If $G\in \bS$ is a topological abelian group, then $\Hom_\bS(U,G):=\uG(U)$ has the structure of a group by point-wise multiplications. In this case $\uG$ is a sheaf of abelian groups $\uG\in \Sh_\Ab\bS$.

\subsubsection{}\label{azuzasudsad87}

We can pass back and forth between (pre)sheaves of abelian groups and (pre)sheaves of sets using the following adjoint pairs of functors
$${}^p\Z(\dots):\Pr\bS\Leftrightarrow{\Pr}_\Ab \bS:\cF\ ,\quad \Z(\dots):\Sh\bS\Leftrightarrow\Sh_\Ab \bS:\cF\ .$$
The functor  $\cF$ forgets the abelian group structure. The functors ${}^p\Z(\dots)$ and $\Z(\dots)$
 are called the linearization functors.
The presheaf linearization
associates to a presheaf  of sets $H\in \Pr\bS$ the presheaf ${}^p\Z(H)\in \Pr_\Ab\bS$ which sends  $U\in \bS$ to the free abelian group $\Z(H)(U)$ generated by the set $H(U)$. The linearization functor for sheaves is the composition
$\Z(\dots):=i^\sharp\circ {}^p\Z(\dots)$ of the presheaf linearization and the sheafification.

\subsubsection{}

\begin{lem}\label{dezuqwideqwd}
Consider an exact  sequence  of topological groups in $\bS$ $$1\to G\to H\to
L\to 1\ .$$ 
Then $1\to \uG\to \uH\to \uL$ is an exact sequence of sheaves of 
groups.

The map of topological spaces $H\to L$ has local sections if and only if
$$1\to \uG\to \uH \to \uL\to 1$$ is an exact sequence of sheaves of groups.
\end{lem} 
\proof
Exactness at $\uG$ and $\uH$ is clear. 

Assume the existence of local sections of $H\to L$.
This implies
exactness at $\uL$. 

On the other hand, evaluating the sequence at the object $L\in \bS$
we see that the exactness of the sequence of sheaves implies the existence of local sections to $H\to L$.
\hB

\subsubsection{}\label{kfff388355}

For sheaves $G,H\in \Sh\bS$ we define the sheaf
$$\uHom_{\Sh\bS}(G,H)\in \Sh\bS\ ,\quad \uHom_{\Sh\bS}(G,H)(U):=\Hom_{\Sh\bS}(G_{|U},H_{|U})\ .$$
Again, this a priori defines a presheaf, but one checks the sheaf conditions
in a straightforward manner.

\subsubsection{}

If $X,H\in \bS$ and $H$ is in addition a topological abelian group, then  
$\Map(X,H)$ is again a topological abelian group. For a topological abelian group $G\in \bS$ we let
$\Hom_{\top-\Ab}(G,H)\subseteq \Map(G,H)$ be the closed subgroup of homomorphisms.
Recall the construction of the internal $\uHom_{\Sh_\Ab\bS}(\dots,\dots)$, compare \ref{kfff388355}.
\begin{lem}\label{eudiwe33}
Assume that $G,H\in \bS$ are topological abelian groups.
Then we have a canonical isomorphism of sheaves of abelian groups
$$e: \underline{\Hom_{\top-\Ab}(G,H)}\stackrel{\sim}{\to}\uHom_{\Sh_\Ab\bS}(\uG,\uH)\ .$$
\end{lem}
\proof
We first describe the morphism $e$.
Let $\phi\in \underline{\Hom(G,H)}(U)$. We define the element $e(\phi)\in \uHom_{\Sh_\Ab\bS}(\uG,\uH)(U)$, i.e. a morphism of sheaves 
$e(\phi):\uG_{|U}\to \uH_{|U}$, such that it sends 
$f\in G(\sigma:V\to U)$ to $\{V\ni v\mapsto \phi(\sigma(v))f(v)\}\in H(\sigma:V\to U)$.  

Let us now describe the inverse $e^{-1}$. Let $\psi\in  \uHom_{\Sh_\Ab\bS}(\uG,\uH)(U)$ be given. 
Since $(\pr_U:U\times G\to U)\in \bS/U$
it gives rise to a map
$$\psi_G:\uG(\pr_U:U\times G\to U)\to \uH(\pr_U:U\times G\to U)\ .$$
We now consider $(\pr_G:U\times G\to G)\in \uG(\pr_U:U\times G\to U)$ and set
$\phi_G:=\psi_G(\pr_G)\in \uH(\pr_U:U\times G\to U)=\Hom_\bS(U\times G,H)$.
We now invoke the  exponential law  isomorphism $\exp:\Hom_\bS(U\times G,
H)\stackrel{\sim}{\to}\Hom_\bS(U,\Map(G,H))$. Since $\psi$ was a homomorphism
and using the sheaf property, we actually get an element
$e^{-1}(\psi):=\exp(\phi_G)\in
\Hom_\bS(U,\Hom_{\top-\Ab}(G,H))=\underline{\Hom_{\top-\Ab}(G,H)}(U)$.
(Details of the argument for this fact are left to the reader).\hB 

If $\bS$ is replaced by one of the sub-sites $\bS_{lc}$ or $\bS_{lc-acyc}$, then it may happen that
$\Hom_{\top-\Ab}(G,H)$ does not belong to the sub-site. In this case 
Lemma \ref{eudiwe33} remains true if one interprets
$\underline{\Hom_{\top-\Ab}(G,H)}$ as the restriction of the sheaf on $\bS$  represented by
$\Hom_{\top-\Ab}(G,H)$ to the corresponding sub-site.

\subsection{Restriction}

\subsubsection{}\label{hfjdsfc}

In this subsection we prove a general result in sheaf theory (Proposition \ref{hsdjfjsdf674}) which is probably
well-known, but which we could not locate in the literature. Let $\bS$ be a site. For $Y\in \bS$ we can consider the relative site $\bS/Y$. Its objects are maps
$(U\to Y)$ in $\bS$, and its morphisms are diagrams
$$\xymatrix{U\ar[dr]\ar[rr]&&V\ar[dl]\\&Y&}\ .$$ The covering families for $\bS/Y$ are induced by the covering families of $\bS$, i.e. for $(U\to Y)\in \bS/Y$ a covering family $\tau:=(U_i\to U)_{i\in I}$ induces a covering family
$$\tau_{|Y}:=\left(\xymatrix{U_i\ar[dr]\ar[rr]&&U\ar[dl]\\&Y&}\right)_{i\in I}$$
in $\bS/Y$.

\subsubsection{}\label{hjadhaduqiuwiue}

There is a canonical functor $f:\bS/Y\to \bS$ given on objects by
$f(U\to Y):=U$. It induces  adjoint pairs of functors
$$f_*:\Sh\bS/Y\Leftrightarrow \Sh\bS:f^* \ ,{}^p f_*:\Pr\bS/Y\Leftrightarrow \Pr\bS:{}^pf^*\ .$$
We will often write $f^*(F)=:F_{|Y}$ for the restriction functor. For
$F\in \Pr\bS$ it is given in explicit terms by 
$$F_{|Y}(U\to Y):=F(U)\ .$$ 
The functor $f^*$ is in fact the restriction of ${}^pf^*$ to the subcategory of sheaves $\Sh\bS\subseteq \Pr\bS$. 
If $F\in \Sh\bS$, then one must check that ${}^pf^* F$  is a sheaf.
To this end note that the descent conditions for ${}^pf^*F$ with respect to the induced coverings described in \ref{hfjdsfc} immediately follow from the descent conditions for $F$.

\subsubsection{}
Let us now study how the restriction acts on representable sheaves. We assume that $\bS$ has finite products. Let $X,Y\in \bS$.
\begin{lem}\label{udiqwdqww}
If $\bS$ has finite products, then  we have a natural isomorphism
$$\uX_{|Y}\cong \underline{X\times Y\to Y}\ .$$
\end{lem}
\proof
By the universal property of the product for $\sigma:U\to Y$ we have 
\begin{eqnarray*}
\underline{X\times Y\to Y}(U\to Y)&=& \{h\in \Hom_\bS(U,X\times Y)| \pr_Y\circ h= \sigma\}\\&\cong& \Hom_{\bS}(U,X)\\&=& \uX(U)\\&=& \uX_{|Y}(U\to Y)\ .
\end{eqnarray*}

\hB

\subsubsection{}

Let 
$i^\sharp:\Pr\bS\to \Sh\bS$ and $i^\sharp_Y:\Pr\bS/Y\to \Sh\bS/Y$ be the sheafification functors.
\begin{lem}\label{udiwqdiuiuqw}
We have a canonical isomorphism
$f^*\circ i^\sharp\cong i^\sharp_Y\circ {}^p f^*$. 
\end{lem}
\proof
The argument is similar to that in the proof of \cite[Lemma 2.31]{bss}.
The main point is that the sheafifications $i^\sharp F(U)$ and $i^\sharp_YF_{|Y}(U\to Y)$ can be expressed in terms of the category of covering families
$\cov_\bS(U)$ and $\cov_{\bS/Y}(U\to Y)$. Compared with  \cite[Lemma 2.31]{bss} the argument is simplified by the fact that the  induction of covering families
described in \ref{hfjdsfc} induces an isomorphism of categories
$\cov_\bS(U)\to \cov_{\bS/Y}(U\to Y)$.
\hB

\subsubsection{}

Recall that a functor is called exact if it commutes with limits and colimits.

\begin{lem}\label{teezwdc}
The restriction functor $f^*:\Sh\bS\to \Sh\bS/Y$ is exact.
\end{lem}
\proof
The functor $f^*$ is a right-adjoint and therefore commutes with limits.
Since colimits of presheaves are defined object-wise, i.e. for a diagram
$$\cC\to \Pr\bS\ ,\quad c\mapsto F_c$$ of presheaves we have 
$$({}^p\colim_{c\in \cC} F_c)(U)=\colim_{c\in \cC} F_c(U)\ ,$$
 it follows from the explicit
description of ${}^pf^*$, that it commutes with colimits of presheaves.
Indeed, for a diagram of sheaves
$F:\cC\to \Sh\bS$ we have $$\colim_{c\in \cC} F_c=i^\sharp {}^p\colim_{c\in\cC} F_c\ .$$
By Lemma \ref{udiwqdiuiuqw} we get
$f^*\colim_{c\in \cC} F_c\cong f^*i^\sharp {}^p\colim_{c\in\cC} F_c\cong i^\sharp_Y {}^pf^*{}^p\colim_{c\in\cC} F_c\cong i^\sharp {}^p\colim_{c\in\cC} {}^p f^*F_c\cong \colim_{c\in \cC} f^*F_c\ .$
\hB

\subsubsection{}

Let $H\in \Sh\bS$ and $X\in \bS$.
\begin{lem}\label{hgeuidwqdqwdwqd}
For $U\in \bS$ we have a natural bijection
$$\uHom_{\Sh\bS}(\uX,H)(U)\cong H(U\times X)\ .$$
\end{lem}
\proof

We have the following chain of isomorphisms
\begin{eqnarray*}
\uHom_{\Sh\bS}(\uX,H)(U)&\cong& \Hom_{\Sh\bS/U}(\uX_{|U},H_{|U})\\&\stackrel{Lemma \:\ref{udiqwdqww}}{\cong}&
\Hom_{\Sh\bS/U}(\underline{U\times X\to U},H_{|U})\\&\cong& H(U\times X)\ .\end{eqnarray*}
\hB

We will need the explicit description of this bijection.
We define a map
$$\Psi:\uHom_{\Sh\bS}(\uX,H)(U)\to H(U\times X)$$ as follows.
An element $h\in \uHom_{\Sh\bS/U}(\uX,H)(U)=\Hom_{\Sh\bS}(\uX_{|U},H_{|U})$
induces a map $\tilde h:\uX_{|U}(U\times X\to U)\to H_{|U}(U\times X\to U)$.
Now we have $\uX_{|U}(U\times X\to U)=\uX(U\times X)=\Map(U\times X,X)$ and
$H_{|U}(U\times X\to U)=H(U\times X)$. Let $\pr_X:U\times X\to X$ be the projection.
We define $\Psi(h):=\tilde h(\pr_X)$.

We now define
$\Phi:H(U\times X)\to \uHom_{\Sh\bS}(\uX,H)(U)$.
Let $g\in H(U\times X)$. 
For each $(e:V\to U)\in \bS/U$ we must define a map
$\widehat g_{(V\to U)}:\uX_{|U}(V\to U)\to H_{|U}(V\to U)$.
Note that $\uX_{|U}(V\to U)= \uX(V)=\Map(V,X)^\delta$.
Let $x\in \Map(V,X)$. Then we have an induced map
$(e,x):V\to U\times X$.
We define $\Psi(g):=H(e,x)(g)$.
One checks that this construction is functorial in $e$ and therefore defines
a morphism of sheaves.

A straight forward calculation shows that $\Psi$ and $\Phi$ are inverses to each other and induce the bijection above.

\subsubsection{}\label{ajdskasda}

For a map $U\to V$ we have an isomorphism of sites
$\bS/U\cong (\bS/V)/(U\to V)$. Several formulas below implicitly contain this identification. For example, we have a canonical isomorphism
$$F_{|U}\cong (F_{|V})_{|U\to V}$$ for $F\in \Sh\bS$. For $G,H\in \Sh\bS$ this induces isomorphisms
\begin{equation}\label{cdhdsc}\uHom_{\Sh\bS}(G,H)_{|V}\cong \uHom_{\Sh\bS/V}(G_{|V},H_{|V})\ .\end{equation}

\subsubsection{}

We now consider a sheaf of abelian groups $H\in \Sh_\Ab\bS$.
If $\tau\in \cov_\bS(U)$ is a covering family of $U$ and $H\in \Sh_\Ab \bS$, then we can define the \v{C}ech complex $\check{C}^*(\tau;H)$ (see \cite[2.3.5]{bss}).
\begin{ddd}
The sheaf $H$ is called flabby if $H^i(\check{C}^*(\tau;H))\cong 0$ for all $i\ge 1$, $U\in \bS$ and $\tau\in \cov_\bS(U)$.
\end{ddd}

For $A\in \bS$ we have the section functor 
$$\Gamma(A;\dots):\Sh_\Ab\bS\to \Ab\ ,\quad \Gamma(A;H):=H(A) .$$ This functor
is left exact and admits a right-derived functor $R\Gamma(A;\dots)$ which can be calculated using injective resolutions. Note that $\Gamma(A;\dots)$ is acyclic on flabby sheaves, i.e. $R^i\Gamma(A;H)\cong 0$ for $i\ge 1$ if $H$ is flabby. Hence, the derived functor $R\Gamma(A;\dots)$ can also be calculated using flabby resolutions. Note that an injective sheaf is flabby.

\subsubsection{}
 
\begin{lem}\label{hajasdfaf}
The restriction functor $\Sh_\Ab\bS\to \Sh_\Ab\bS/Y$ preserves flabby sheaves.
\end{lem}
\proof
Let $H\in \Sh_\Ab\bS$ be flabby.
For $(U\to Y)\in \bS/Y$ and $\tau\in \cov_\bS(U)$ we let
$\tau_{|Y}\in \cov_{\bS/Y}(U\to Y)$ be the induced covering family as in 
\ref{hfjdsfc}. Note that
$\check{C}^*(\tau_{|Y};H_{|Y})\cong \check{C}^*(\tau;H)$ is acyclic.
Since $\cov_{\bS/Y}(U\to Y)$ is exhausted by families of the form
$\tau_{|Y}$, $\tau\in \cov_{\bS}(U)$ it follows that
$H_{|Y}$ is flabby. \hB

\subsubsection{}

We now consider sheaves of abelian groups $F,G\in \Sh_\Ab\bS$. 
\begin{prop}\label{hsdjfjsdf674}
Assume that the site $\bS$ has the property that for all $U\in \bS$ the
restriction $\Sh\bS\to \Sh\bS/U$ preserves representable sheaves. 
If $Y\in \bS$, then 
in $D^+(\Sh_\Ab\bS/Y)$ there is a natural isomorphism 
$$R\uHom_{\Sh_\Ab\bS}(F,G)_{|Y}\cong   R\uHom_{\Sh_\Ab\bS/Y}(F_{|Y},G_{|Y}) .$$
\end{prop}
\proof
We choose an injective resolution $G\to I$ and furthermore an injective
resolution $I_{|Y}\to J$. Note that by Lemma \ref{hajasdfaf} restriction is
exact and therefore $G_{|Y}\to  J$ is a resolution of $G_{|Y}$.
Then we have
$$R\uHom_{\Sh_\Ab\bS}(F,G)_{|Y}\cong \uHom_{\Sh_\Ab\bS}(F,I)_{|Y}\stackrel{(\ref{cdhdsc})}{\cong} \uHom_{\Sh_\Ab\bS/Y}(F_{|Y},I_{|Y})\ .$$
Furthermore,
$$R\uHom_{\Sh_\Ab\bS/Y}(F_{|Y},G_{|Y})\cong \uHom_{\Sh_\Ab\bS/Y}(F_{|Y},J)$$
and the map
$I_{|Y}\to J$ induces 
\begin{eqnarray}\label{chdsjew77823}R\uHom_{\Sh_\Ab\bS}(F,G)_{|Y}&\cong& \uHom_{\Sh_\Ab\bS/Y}(F_{|Y},I_{|Y})\\&\stackrel{*}{\to}& \uHom_{\Sh_\Ab\bS/Y}(F_{|Y},J) \nonumber\\&\cong& R\uHom_{\Sh_\Ab\bS/Y}(F_{|Y},G_{|Y}) \nonumber\ .\end{eqnarray}
If the restriction $f^*:\Sh\bS\to \Sh\bS/Y$ would preserve injectives, then 
$I_{|Y}\to J$ would be a homotopy equivalence and the marked map would be a 
quasi-isomorphism.

In the generality of the present paper we do not know whether
$f^*$ preserves injectives. Nevertheless we show that our assumption on
$\bS$ implies that the marked map in (\ref{chdsjew77823}) is a quasi isomorphism
for all sheaves $F\in\Sh_\Ab\bS$.

\subsubsection{}

We first show a special case.
\begin{lem}\label{hasjsahdjashd}
If $F$ is representable, then the marked map in (\ref{chdsjew77823}) is a quasi isomorphism. 
\end{lem}
\proof
Let $(U\to Y)\in \bS/Y$ and $(A\to Y)\in \bS$.
Then on the one hand we have
\begin{eqnarray*}
\uHom_{\Sh_\Ab\bS/Y}(\Z(\uA)_{|Y},I_{|Y})(U\to Y)&\stackrel{\ref{kfff388355}}{=}&\Hom_{(\Sh\bS/Y)/(U\to Y)}((\uA_{|Y})_{|(U\to Y)},(I_{|Y})_{|U\to Y})\\&\stackrel{\ref{ajdskasda}}{\cong}& \Hom_{\Sh\bS/U}(\uA_{|U},I_{|U})\ .
\end{eqnarray*}
Since $\uA_{|U}$ is representable by our assumption on $\bS$ there exists $(B\to U)\in\bS/U$ such that $\uA_{|U}\cong \underline{B\to U}$.

Since the restriction functor is exact (Lemma \ref{teezwdc}) and the restriction of an injective sheaf from $\bS$ to $\bS/Y$ is still flabby 
(Lemma \ref{hajasdfaf}) we see that $G_{|U}\to I_{|U}$ is a flabby resolution.
It follows that
$$\Hom_{\Sh\bS/U}(\uA_{|U},I_{|U})\cong \Hom_{\Sh\bS/U}(\underline{B\to U},I_{|U}) \cong I_{|U}(B\to U) \cong R\Gamma(B\to U,G_{|U})\ .$$
On the other hand
\begin{eqnarray*}
\uHom_{\Sh_\Ab\bS/Y}(\Z(\uA)_{|Y},J)(U\to Y)&\cong& 
\Hom_{(\Sh\bS/Y)/(U\to Y)}((\uA_{|Y})_{|U\to Y},J_{|U\to Y})\\&\cong& \Hom_{\Sh\bS/U}(\uA_{|U} ,J_{|U})\\&\cong& 
\Hom_{\Sh\bS/U}(\underline{B\to U},J_{|U})\\&\cong& J(B\to U)\\&\cong& R\Gamma(B\to U,G_{|U}) .\end{eqnarray*}
\hB

\subsubsection{}

We consider the functor which associates to $F\in \Sh_\Ab\bS$ the cone
$$C(F):=\cone\left(\uHom_{\Sh_\Ab\bS/Y}(F_{|Y},I_{|Y})\to \uHom_{\Sh_\Ab\bS/Y}(F_{|Y},J)\right)\ .$$
In order to show that the marked map in (\ref{chdsjew77823}) is a quasi isomorphism we must show that $C(F)$ is acyclic for all $F\in \Sh_\Ab\bS$.

The restriction functor $(\dots)_{|Y}$ is exact by Lemma \ref{teezwdc}. 
For $H\in \Sh_\Ab\bS/Y$ the functor
$\uHom_{\Sh_\Ab\bS/Y}(\dots,H)$ is a right-adjoint. As a contravariant
functor it transforms colimits into limits and is left exact.
In particular, the functors
$$\uHom_{\Sh_\Ab\bS/Y}((\dots)_{|Y},I_{|Y}),\quad \uHom_{\Sh_\Ab\bS/Y}((\dots)_{|Y},J) $$ transform coproducts of sheaves into products of complexes. It follows that
$F\to C(F)$ also transforms coproducts of sheaves into products of complexes.
If $P\in \Sh_\Ab\bS$ is a coproduct of representable sheaves, then by Lemma 
\ref{hasjsahdjashd} $C(A)$ is a product of acyclic complexes and hence acyclic.

We claim that
$F\to C(F)$ transforms short exact sequences of sheaves to short exact sequences of complexes. Since $J$ is injective and
$(\dots)_{|Y}$ is exact, the functor
$F\mapsto \uHom_{\Sh_\Ab\bS/Y}(F_{|Y},J)$ is a composition of 
exact functors and thus has this property.
Furthermore we have
$\uHom_{\Sh_\Ab\bS/Y}(F_{|Y},I_{|Y})\cong \uHom_{\Sh_\Ab\bS}(F,I)_{|Y}$.
Since $I$ is injective, the functor 
$F\to \uHom_{\Sh_\Ab\bS/Y}(F_{|Y},I_{|Y})$ has this property, too.
This implies the claim.

We now argue by induction. Let $n\in \nat$ and assume that we have already shown that $H^i(C(F))\cong 0$ for all  $F\in \Sh_\Ab\bS$ and $i<n$. 
 
Consider   $F\in \Sh_\Ab\bS$. Then there exists a coproduct of representables
$P\in \Sh_\Ab\bS$ and an exact sequence $$0\to K\to P\to F\to 0\ .$$
In order to construct $P$ observe that
in general one has a canonical isomorphism
$$F\cong \colim_{\uA\to F} \uA\ .$$ We let $P:=\bigoplus_{\uA\to F} \uA$.
The collection of maps  $\uA\to F$ in the index category  induces a canonical surjection 
$$P=\bigoplus_{\uA\to F} \uA\to \colim_{\uA\to F} \uA \cong F\ .$$

The short exact sequence of complexes 
$$0\to C(F)\to C(P)\to C(K)\to 0$$ induces
a long exact sequence in cohomology. Since $H^i(P)\cong 0$ for all $i\in \Z$
we conclude that $H^{i}(C(F))\cong H^{i-1}(C(K))$ for all $i\in \nat$.
By our induction hypothesis we get
$H^n(C(F))\cong 0$.

By induction on $n$ we show that $C(F)$ is acyclic for all $F\in \Sh_\Ab\bS$, and this implies Proposition \ref{hsdjfjsdf674}.
\hB

\subsubsection{}

Note that a site $\bS$ which has finite products  satisfies the assumption of  Proposition \ref{hsdjfjsdf674}. In fact, for $A\in \bS$ we have by Lemma \ref{udiqwdqww} that 
$$\uA_{|Y}\cong \underline{A\times Y\to Y}\ .$$
Therefore the restriction functor $(\dots)_{|Y}$ preserves representables.

\begin{kor}\label{relativecaseall}
Let $\bS$ be as in  Proposition \ref{hsdjfjsdf674}.
For sheaves $F,G\in \Sh_\Ab\bS$ and $Y\in \bS$ we have
$$\uExt^i_{\Sh_\Ab\bS}(F,G)_{|Y}\cong \uExt^i_{\Sh\Ab\bS/Y}(F_{|Y},G_{|Y})\ .$$
\end{kor}
\proof
Since $(\dots)_{|Y}$ is exact (Lemma \ref{teezwdc})
we have
\begin{eqnarray*}
\uExt^i_{\Sh_\Ab\bS}(F,G)_{|Y}&\cong&(H^iR\uHom_{\Sh_\Ab\bS}(F,G))_{|Y}\\
&\cong&H^i(R\uHom_{\Sh_\Ab\bS}(F,G)_{|Y})\\
&\stackrel{Proposition \ref{hsdjfjsdf674}}{\cong}&H^i(R\uHom_{\Sh_\Ab\bS/Y}(F_{|Y},G_{|Y}))\\
&\cong&\uExt^i_{\Sh\Ab\bS/Y}(F_{|Y},G_{|Y})\ .
\end{eqnarray*}
\hB

\subsubsection{}\label{timeserkl}

Let us assume that $\bS$ has finite fibre products.
Fix $U\in \bS$. Then we can define a morphism of sites
${}^U\times:\bS\to \bS/U$ which maps
$A\in \bS$ to $(U\times A\to U)\in \bS/U$.
One easily checks the conditions given in \cite[1.2.2]{MR1317816}.
This morphism of sites induces an adjoint pair of functors
$${}^U\times_*:\Sh\bS\Leftrightarrow \Sh\bS/U:{}^U\times^*\ .$$
Let $f:\bS/U\to \bS$ be the restriction defined in \ref{hjadhaduqiuwiue} (with $Y$ replaced by $U$) and let
$$f_*:\Sh\bS/U\Leftrightarrow \Sh\bS:f^*$$
be the corresponding pair of adjoint functors.
\begin{lem}\label{hjhfwefweifew}
We assume that $\bS$ has finite products.
Then we have a canonical isomorphism $f^*\cong {}^U\times_*$.
\end{lem}
\proof
We first define this isomorphism on representable sheaves.
Since every sheaf can be written as a colimit of representable sheaves, 
${}^U\times_*$ commutes with colimits as a left-adjoint, and $f^*$ commutes with colimits
by Lemma \ref{teezwdc}, the isomorphism then extends to all sheaves.
Let $W\in \bS$ and $F\in \Sh\bS/U)$.
Then we have
\begin{eqnarray*}
\Hom_{\Sh\bS/U}({}^U\times_*\uW,F)&\cong &\Hom_{\Sh\bS}(\uW,{}^U\times^*(F))\\
&\cong& {}^U\times^*(F)(W)\\
&\cong&F(W\times U\to U)\\
&\stackrel{Lemma \ref{udiqwdqww}}{\cong}& \Hom_{\Sh\bS/U}(\underline{(W\times U\to U)},F)\\
&\cong&\Hom_{\Sh\bS/U}(f^*\uW,F)
\end{eqnarray*}
This isomorphism is represented by a canonical isomorphism
$${}^U\times_*(\uW)\cong f^*(\uW)$$ in
$\Sh\bS/U$.
\hB

 \subsubsection{}

Let $f:\bS\to \bS^\prime$ be a morphism of sites. It induces an adjoint pair
$$f_*:\Sh\bS\Leftrightarrow\Sh\bS^\prime:f^*\ .$$
\begin{lem}\label{wdiudwqdwq}
For $F\in \Sh\bS$ the functor $f_*$ has the explicit description
$$f_*(F):=\colim_{(\uU\to F)\in \Sh(\bS/F)} \underline{f(U)}\ .$$
\end{lem}
\proof
In fact, for  $G\in \Sh\bS^\prime$ we have
\begin{eqnarray*}
\Hom_{\Sh\bS^\prime}(\colim_{(\uU\to F)\in \Sh(\bS)/F} \underline{f(U)},G)&\cong&\lim_{(\uU\to F)\in \Sh(\bS)/F} \Hom_{\Sh\bS^\prime}( \underline{f(U)},G)\\&\cong &\lim_{(\uU\to F)\in \Sh(\bS)/F} G(f(U))\\&\cong&\lim_{(\uU\to F)\in \Sh(\bS)/F}  f^*G(U)\\&\cong& \lim_{(\uU\to F)\in \Sh(\bS)/F} \Hom_{\Sh\bS^\prime}(\uU,f^*G)\\&\cong&
\Hom_{\Sh\bS}(\colim_{(\uU\to F)\in \Sh(\bS)/F}\uU,f^*G)\\&\cong& 
\Hom_{\Sh\bS}(F,f^*G)\ .
\end{eqnarray*}

\begin{lem}\label{zuuzud}
 For $A\in \bS$
we have
$$f_*(\uA)=\underline{f(A)}\ .$$
\end{lem}
\proof
Since $(\id_A:A\to A)\in \bS/A$ is final we have
\begin{eqnarray*}
f_*(\uA)&\cong&\colim_{(\uU\to \uA)\in \Sh(\bS)/\uA} \underline{f(U)} \\
&\cong&\colim_{(U\to A)\in \bS/A} \underline{f(U)}\\
&\cong&\underline{f(A)}\ .
\end{eqnarray*}
\hB

\subsubsection{}

\begin{lem}
If $f:\bS\to \bS^\prime$ is fully faithful, then we have for $A\in \bS$ that 
$f^* f_*\uA\cong \uA$.
\end{lem}
\proof
We calculate for $U\in \bS$ that
\begin{eqnarray*}
f^*f_* \uA(U)&\stackrel{Lemma \:\ref{zuuzud} }{\cong}& f^*\underline{f(A)}(U)\\&\cong&
\underline{f(A)}(f(U))\\&\cong&\Hom_{\bS^\prime}(f(U),f(A))\\&\cong&
\Hom_{\bS}(U,A)\\&\cong&\uA(U)
\end{eqnarray*}
\hB

\subsubsection{}

For $U\in \bS$ we have an induced morphism of sites
${}_U f:\bS/U\to \bS^\prime/f(U)$ which maps
$(A\to U)$ to $(f(A)\to f(U))$.

If $f$ is fully faithful, then ${}_Uf$ is fully faithful for all $U\in \bS$.

\begin{lem}\label{wqdiwqduiqdwqdoqqd}
For $G\in \Sh\bS^\prime$ we have in  $\bS/U$ the identity
$${}_Uf^*(G_{|f(U)}) \cong (f^* G)_{|U}\ .$$
If we know in addition that $\bS,\bS^\prime$ have products which are preserved by  $f:\bS\to \bS^\prime$, then
for $F\in \bS$ we have
in $\bS^\prime/f(U)$ the identity
$${}_U f_* F_{|U}\cong (f_*F)_{|f(U)}\ .$$
\end{lem}
\proof
Indeed, for
$(V\to U)\in \bS/U$ we have
$$(f^*G)_{|U}(V\to U)=G(f(V))=G_{|f(U)}(f(V)\to f(U))=({}_U f^*G_{|f(U)})(V\to U)\ .$$
In order to see the second identity note that
we can write
$$F_{|U}\cong\colim_{(\uA\to F_{|U})\in \Sh(\bS/U)/F_{|U}} \uA\ .$$
Since ${}_Uf_*$ is a left-adjoint it commutes with colimits.
Furthermore the restriction functors $(\dots)_{|U}$ and $(\dots)_{f(U)}$ are exact by Lemma \ref{teezwdc} and therefore also commute with colimits.
Writing
$$F\cong \colim_{(\uA\to F)\in \Sh(\bS)/F}\uA$$ we get
\begin{eqnarray*}
F_{|U}\cong (\colim_{(\uA\to F)\in \Sh(\bS)/F}\uA)_{|U}&\cong& \colim_{(\uA\to F)\in \Sh(\bS)/F}\uA_{|U}\\&\stackrel{Lemma \:\ref{udiqwdqww}}{\cong}&
\colim_{(\uA\to F)\in \Sh(\bS)/F}\underline{A\times U\to A}\ .\end{eqnarray*}
Using that $f$ preserves products in the isomorphism marked by $(!)$ we calculate
\begin{eqnarray*}
{}_Uf_*F_{|U}&\cong&{}_Uf_*\colim_{(\uA\to F)\in \Sh(\bS)/F}\underline{A\times U\to U}\\
&\cong&\colim_{(\uA\to F)\in \Sh(\bS)/F}\:{}_Uf_*\underline{A\times U\to U}\\
&\stackrel{Lemma \:\ref{zuuzud}}{\cong}&\colim_{(\uA\to F)\in \Sh(\bS)/F}\underline{f(A\times U)\to f(U)
}\\&\stackrel{(!)}{\cong}&\colim_{(\uA\to F)\in \Sh(\bS)/F}\underline{f(A)\times f(U)\to f(U)}\\
&\cong&\colim_{(\uA\to F)\in \Sh(\bS)/F}\underline{f(A)}_{|f(U)}\\
&\cong&(\colim_{(\uA\to F)\in \Sh(\bS)/F}\underline{f(A)})_{|f(U)}\\
&\stackrel{Lemma \:\ref{wdiudwqdwq}}{\cong}& (f_*F)_{|f(U)}\ .
\end{eqnarray*}
 \hB

\subsubsection{}

\begin{lem}\label{k12uiwdwuqdqdwqdqwdd}
Assume that $\bS,\bS^\prime$ have finite products which are preserved by $f:\bS\to \bS^\prime$.
For $F\in \Sh\bS$ and $G\in \Sh\bS^\prime$
we have a natural isomorphism 
$$\uHom_{\Sh\bS}(F,f^*G)\cong f^*\uHom_{\Sh\bS^\prime}(f_*F,G)\ .$$
\end{lem}
\proof
For $U\in \bS$ we calculate
\begin{eqnarray*}
\uHom_{\Sh\bS}(F,f^*G)(U)&=&
\Hom_{\Sh\bS/U}(F_{|U},(f^*G)_{|U})\\&\stackrel{Lemma\ref{wqdiwqduiqdwqdoqqd}}{\cong}&
\Hom_{\Sh\bS/U}(F_{|U},{}^Uf^*(G_{|f(U)}))\\
&\cong&
\Hom_{\Sh\bS^\prime/f(U)}({}^Uf_*(F_{|U}),G_{|f(U)})\\
&\stackrel{Lemma \:\ref{wqdiwqduiqdwqdoqqd}}{\cong}&
\Hom_{\Sh\bS^\prime/f(U)}((f_*F)_{|f(U)},G_{|f(U)})\\
&=&\uHom_{\Sh\bS^\prime}(f_*F,G)(f(U))\\
&=&f^*\uHom_{\Sh\bS^\prime}(f_*F,G)(U)
\end{eqnarray*}
\hB

\subsubsection{}

We now consider the derived version of Lemma \ref{k12uiwdwuqdqdwqdqwdd}.
Let $F\in \Sh_\Ab\bS$ and $G\in \Sh_\Ab\bS^\prime$ be sheaves of abelian groups. 
\begin{prop}\label{neuhsdjfjsdf674133}
We make the following assumptions:
\begin{enumerate}
\item The sites $\bS,\bS^\prime$ have finite products.
\item The morphism of sites $f:\bS\to \bS^\prime$ preserves finite products.
 \item We assume that $f^*: \Sh\bS^\prime\to \Sh\bS$ is exact.
\end{enumerate}
 Then in  $D^+(\Sh_\Ab\bS)$ there is a natural isomorphism 
$$R\uHom_{\Sh_\Ab\bS}(F,f^*G)\cong  f^* R\uHom_{\Sh_\Ab\bS^\prime}(f_*F,G) .$$
\end{prop}
\proof
By \cite[Proposition 3.6.7]{MR1317816} the conditions on the sites and $f$ imply that the functor $f_*$ is exact. It follows that $f^*$ preserves injectives (see e.g. \cite[Proposition 3.6.2]{MR1317816}).
We choose an injective resolution $G\to I^\bullet$.  Then $f^*G\to f^*I^\bullet$
is an injective resolution of $f^*G$. It follows that
$$f^*R\uHom_{\Sh_\Ab\bS^\prime}(f_*F,G)\cong f^*\uHom_{\Sh_\Ab\bS^\prime}(f_*F,I^\bullet)\cong \uHom_{\Sh_\Ab\bS}(F,f^*I^\bullet)\cong R\uHom_{\Sh_\Ab\bS}(F,f^*G)\ .$$
\hB

Proposition \ref{neuhsdjfjsdf674133} is very similar in spirit to Proposition \ref{hsdjfjsdf674}.
On the other hand, we can not deduce \ref{hsdjfjsdf674} from \ref{neuhsdjfjsdf674133}
since the functor $f:\Sh\bS/U
\to \bS$ in general does not preserve products.
In fact, if $\bS$ has fibre products, then 
$(A\to U)\times (B\to U)=(A\times_UB\to U)$.
Then $f((A\to U)\times (B\to U))\cong A\times_UB$ while
$f(A\to U)\times f(B\to U)\cong A\times B$.

\subsubsection{}
Recall the notations ${}^W\times, {}_W\times$ from \ref{timeserkl}.

\begin{lem}Let $W\in \bS$. For every $F\in \Sh\bS$ and $A\in \bS$ we have
$$\uHom_{\Sh\bS}(\uW,F)(A)\cong {}^W\times^* {}^W\times_*F(A)\ .$$
\end{lem}
\proof
For $A\in \bS$ we consider $f:\bS/A\to \bS$ and get
\begin{eqnarray}\label{hjuziqwuxq}
\uHom_{\Sh\bS}(\uW,F)(A)&\cong \Hom_{\Sh\bS/A}(\uW_{|A},F_{|A})\\
&\cong&\Hom_{\Sh\bS/A}(f^*\uW,f^*F)\nonumber\\
&\stackrel{Lemma\:\ref{hjhfwefweifew}}{\cong}&\Hom_{\Sh\bS/A}({}^A\times_*(\uW),f^*F)\nonumber\\
&\cong&\Hom_{\Sh\bS}(\uW,{}^A\times^*f^*(F))\nonumber\\
&\cong&{}^A\times^*(f^*F)(W)\nonumber\\
&\cong& f^*F(A\times W\to A)\nonumber\\
&\cong&F(A\times W)\nonumber\\
&\cong&{}^W\times^*(f^*F)(A)\nonumber\\
&\stackrel{Lemma\:\ref{hjhfwefweifew}}{\cong}&{}^W\times^*({}^W\times_*(F))(A)\nonumber\ .
\end{eqnarray}
\hB

\subsubsection{}\label{hdfjcsdcuis}

To abbreviate, let us introduce the following notation.
\begin{ddd}\label{jhbqwddwqd}
If $\bS$ has finite products, then for $W\in \bS$ we introduce the functor
$$\cR_W:={}^W\times^*\circ  {}^W\times_*:\Sh\bS\to \Sh\bS\ .$$
\end{ddd}
We denote its restriction to the category of sheaves of abelian groups by the same symbol.
 Since $ {}^W\times_*$ is exact and
${}^W\times^*$ is left-exact, it admits a right-derived functor
$R\cR_W:D^+(\Sh_\Ab\bS)\to D^+(\Sh_\Ab\bS)$.

\subsubsection{}

\begin{lem}\label{udhiqwudqwdwqd}
Let $F\in \Sh_\Ab\bS$ and $W\in \bS$. Then 
we have a canonical isomorphism $R\uHom_{\Sh_\Ab\bS}(\Z(\uW),F)\cong R\cR_W(F)$.
\end{lem}
 \proof
This follows from the isomorphism of functors
$\uHom_{\Sh_\Ab\bS}(\Z(\uW),\dots)\cong  \cR_W(\dots)$ from
$\Sh_\Ab\bS$ to $\Sh_\Ab\bS$.
In fact, we have for $F\in \Sh_\Ab\bS$ that
\begin{eqnarray}\label{uidwdqdiwqdqdwq09}
\uHom_{\Sh_\Ab\bS}(\Z(\uW),F)&\cong&\uHom_{\Sh\bS}(\uW,\cF(F))\\
&\stackrel{(\ref{hjuziqwuxq})}{\cong}&{}^W\times^*({}^W \times_*(F))\nonumber\\
&\stackrel{\ref{jhbqwddwqd}}{\cong}&\cR_W(F)\ .\nonumber
\end{eqnarray} 
\hB

\subsection{Application to sites of topological spaces}\label{odiwqpwqdqwdqwdq}

\subsubsection{}

In this subsection we consider the site $\bS$ of compactly generated
topological spaces as in \ref{firfore} and
some of its sub-sites. We are interested in proving that restriction to sub-sites preserve
$\uExt^i$-sheaves.

We will further study properties of the functor $\cR_W$. In particular, we are interested in results asserting that the higher derived functors $R^i\cR_W(F)$, $i\ge 1$ vanish under certain conditions on $F$ and $W$.

\subsubsection{}

\begin{lem}\label{ddwqdiwduiwqdiud}
If $C\in \bS$ is compact and $H$ is a discrete space, 
then $\Map(C,H)$ is discrete, and
$$\cR_C(\uH)=\underline{\Map(C,H)}\ .$$
\end{lem}
\proof
We first show that
 $\Map(C,H)$ is a discrete space in the compact-open topology. Let $f\in \Map(C,H)$.
Since $C$ is compact, the image $f(C)$ is compact, hence finite. We must show that $\{f\}\subseteq \Map(C,H)$ is open.  Let $h_1,\dots,h_r$ be the finite set of values of $f$.
The sets $f^{-1}(h_i)\subseteq C$ are closed and therefore compact and their
union is $C$. The sets
$\{h_i\}\subseteq H$ are open.
Therefore
$U_i:=\{g\in \Map(C,H)|g(f^{-1}(h_i))\subseteq\{h_i\}\}$ are open subsets of $\Map(C,H)$.
We now see that
$\{f\}=\cap_{i=1}^r U_i$ is open. 

We have by the exponential law 
$$\cR_C(\uH)(A)\cong \uH(A\times C)\cong \Hom_{\bS}(A\times C,H)\cong \Hom_{\bS}(A,\Map(C,H))\cong \underline{\Map(C,H)}(A)\ .$$
\hB

\subsubsection{}

A sheaf $G\in \Sh_\Ab\bS$ is called $\cR_W$-acyclic if $R^i\cR_W(G)\cong 0$ for $i\ge 1$.

\begin{lem}\label{gfrgh3z244r82r2324}
If $G$ is a discrete group,
then $\uG$ is $\cR_{\R^n}$-acyclic.
\end{lem}
\proof
Let $\uG\to I^\bullet$ be an injective resolution.
Then we have for $A\in \bS$ that
$$\cR_{\R^n}(I^\bullet)(A)\stackrel{(\ref{hjuziqwuxq})}{\cong} I^\bullet(A\times \R^n)\ .$$
Therefore $R^i\cR_{\R^n}(\uG)$ is the sheafification of the presheaf
$$A\mapsto H^i(I^\bullet(A\times \R^n))\ .$$
Since $G$ is discrete, the sheaf cohomology of $\uG$ is homotopy invariant, and therefore
$$H^i(I^\bullet(A\times \R^n))\cong H^i(A\times \R^n;\uG)\stackrel{!}{\cong} H^i(A;\uG)\ .$$
To be precise this can be seen as follows.
Let $(A)$ denote the site of open subsets of $A$. It comes with a natural map
$\nu_A:(A)\to \bS$. 
The sheaf cohomology functor is the derived functor of the evaluation functor. In order to indicate on which category this evaluation functor is defined  we temporarily use subscripts. If $I\in \Sh_\Ab\bS$ is injective, then
$\nu_A^*I\in \Sh_\Ab(A)$ is still flabby (see \cite[Lemma 2.32]{bss}).
This implies that $H^*_{\Sh_\Ab\bS}(A;\uG)\cong H_{(A)}^*(A,\nu_A^*(\uG))$.

The diagram
$$\xymatrix{A\ar[rr]^{\id_A\times0}\ar[dr]^{\id_A}&&A\times \R^n\ar[dl]^{\pr_A}\\&A&}$$
is a homotopical isomorphism (in the sense of \cite[2.7.4]{MR1299726}) $A\to A\times \R^n$ over $A$.
We now apply \cite[2.7.7]{MR1299726} which says that the natural map
$$R(\pr_A)_*\pr_A^*(\nu_{A}^*(\uG))\to R(\id_A)_*\id_A^*(\nu^*_A(\uG))$$
is an isomorphism in $D^+(\Sh_\Ab(A))$. But since $G$ is discrete we get
$\pr_A^*(\nu_A^*\uG)\cong \nu_{A\times \R^n}^*(\uG)$.
If we apply the functor $R\Gamma_{(A)}(A,\dots)$ to this isomorphism and take cohomology we get the desired isomorphism marked by $!$\footnote{It is tempting to apply a K{\"u}nneth formula to calculate $H^*(A\times \R^n,\uG)$. But since $A$ is not necessarily compact it is not clear that the K{\"u}nneth formula holds.}.

Now, the sheafification of the presheaf $\bS\ni A\to H^i(A;\uG)\in \Ab$ is exactly the $i$th cohomology sheaf
of $I^\bullet$ which vanishes for $i\ge 1$. \hB 

\subsubsection{}

\begin{lem}
The sheaf $\uR^n$ is $\cR_W$-acyclic for every compact $W\in \bS$.
\end{lem}
\proof
Let $\uR^n\to I^\bullet$ be an injective resolution. Then $R^i\cR_W(\uR^n)$ is the sheafification of the presheaf
$$\bS\ni A \mapsto H^i(\uHom_{\Sh_\Ab\bS}(\Z(\uW_{|A}),I^\bullet_{|A}))\cong H^i(I^\bullet(A\times W))\ .$$
Let $[s]\in H^i(I^\bullet(A\times W))$ be represented by  $s\in I^i(A\times W)$, and $a\in A$. Then we must find a neighbourhood $U\subseteq A$ of $a$ such that $[s]_{|U\times W}=0$, i.e.
$s_{|U\times W}=dt$ for some $t\in I^{i-1}(U\times W)$, where $d:I^{i-1}\to I^i$ is the boundary operator of the resolution.

The ring structure of $\R$ induces on   $\uR^n$ the structure of a sheaf of rings.
In order to distinguish this sheaf of rings from the sheaf of groups $\uR^n$ we will use the notation $\cC$.
Note that $\uR^n$ is in fact a sheaf of $\cC$-modules.

The forgetful functor $\res:\Sh_{\cC-\Mod}\bS\to \Sh_\Ab\bS$ fits into an adjoint pair
$$\ind:\Sh_\Ab\bS\Leftrightarrow \Sh_{\cC-\Mod}\bS:\res\ ,$$
where $\ind$ is given by $\Sh_\Ab\bS\ni V\to V\otimes_\Z \cC\in \Sh_{\cC-\Mod}$.
Since $\cC$ is a torsion-free sheaf it is flat. It follows that $\ind$ is exact and $\res$ preserves injectives.

We can now choose an injective resolution $\uR^n\to J^\bullet$ in $\Sh_{\cC-\Mod}\bS$ and assume
that $I^\bullet=\res(J^\bullet)$.

Since the complex of sheaves $I^\bullet$ is exact we can find an open covering
$(V_r)_{r\in R}$ of $A\times W$ such that $s_{|V_r}=dt_r$ for some $t_r\in I^{i-1}(V_r)$.
Since $W$ is compact (locally compact suffices), by \cite[Thm. 4.3]{MR0210075})
the product topology on $A\times W$ is the compactly generated topology used for the product in $\bS$.
Hence after refining the covering $(V_r)$ we can assume that $V_r=A_r\times W_r$ for open subsets $A_r\subseteq A$ and $W_r\subseteq W$ for all $r\in R$.

We define 
$R_a:=\{r\in R|a\in A_r\}$.
The family $(W_r)_{r\in R_a}$ is an open covering of $W$. Since $W$ is compact we can choose a finite set
$r_1,\dots,r_k\in R_a$ such that $\cW:=(W_{r_1},\dots,W_{r_k})$ is still an open covering of $W$.
The subset $U:=\cap_{j=1}^k A_{r_j}$ is an open neighbourhood of $a\in A$.

Since $I^{i-1}$ is injective we can choose\footnote{Let $U\subseteq X$ be an open subset. Then we have an injection $\uU\to \uX$ and hence an injection $\Z(\uU)\to \Z(\uX)$. For an injective sheaf $I$ we get a surjection $\Hom_{\Sh_\Ab\bS}(\Z(\uX),I)\to \Hom_{\Sh_\Ab\bS}(\Z(\uU),I)$. In other symbols,
$I(X)\to I(U)$ is surjective.} extensions $\tilde t_r\in I^{i-1}(A\times W)$ such that $(\tilde t_r)_{|V_r}=t_r$.
 
We choose a  partition of unity $(\chi_1,\dots,\chi_v)$ subordinate to the finite covering $\cW$.
We take advantage of  the fact that $I^\bullet=\res(J^\bullet)$ which implies that we can multiply sections by continuous functions, and that $d$ commutes with this multiplication. We define
$$t:=\sum_{k=1}^v \chi_k (\tilde t_{r_k})_{|U\times W} \in I^{i-1}(U\times W)\ .$$
Note that $\chi_k (s-d\tilde t_{r_k})_{|U\times W}=0$. In fact we have
$\chi_k (s-d\tilde t_{r_k})_{|(U\times W)\cap V_{r_k}}=\chi_k(s-dt_{r_k})_{|(U\times W)\cap V_{r_k}}=0$.
Furthermore, there is a neighbourhood $Z$ of the complement of $(U\times W)\cap V_{r_k}$ in $U\times W$ where $\chi_k$ vanishes.
Therefore the restrictions $\chi_k (s-d\tilde t_{r_k})$ vanish on the open covering
$\{Z,(U\times W)\cap V_{r_k}\}$ of $U\times W$, and this implies the assertion.
We get
\begin{eqnarray*}
dt&=&\sum_{j=1}^v d(\chi_k (\tilde t_{r_k})_{|U\times W})\\
&=&\sum_{j=1}^v \chi_k d(\tilde t_{r_k})_{|U\times W}\\
&=&\sum_{j=1}^v \chi_k (d\tilde t_{r_k})_{|U\times W}\\
&=&\sum_{j=1}^v \chi_k s_{|U\times W}\\
&=&s_{|U\times W}\ .
\end{eqnarray*}
\hB

\begin{kor}\label{sofdmadd}
\begin{enumerate}
\item 
If $G$ is discrete, then we have 
$\uExt^i_{\Sh_\Ab\bS}(\Z(\uR^n),\uG)\cong 0$ for all $i\ge 1$.
\item
For every compact $W\in \bS$ and $n\ge 1$ we have
$\uExt^i_{\Sh_\Ab\bS}(\Z(\uW),\uR^n)\cong 0$ for all $i\ge 1$.
\end{enumerate}
\end{kor}

\subsubsection{}

\begin{lem}\label{hdjfcdcuisc}
If  $W$ is a profinite space, then every sheaf $F\in \Sh_\Ab\bS$ 
 is $\cR_W$-acyclic. Consequently,
$\uExt^i_{\Sh_\Ab\bS}(\Z(\uW),F)\cong 0$ for $i\ge 1$.
\end{lem}
\proof 
We first show the following intermediate result which is used in \ref{uiqwdgqwdwqdqdwqdwqd} in order
to finish the proof of Lemma \ref{hdjfcdcuisc}.

\subsubsection{}

\begin{lem}\label{wzdqdwqd}
If $W$ is a profinite space, then $\Gamma(W;\dots)$ is exact.
\end{lem}
\proof
A profinite topological space can be written as limit, $W\cong \lim_I W_n$, for an inverse  system of finite spaces
$(W_n)_{n\in I}$. Let $p_n:W\to W_n$ denote the projections.
First of all, $W$ is compact. Every covering of $W$ admits
a finite subcovering. Furthermore, a finite covering admits a refinement  to a covering
by pairwise disjoint open subsets of the form $\{p_n^{-1}(x)\}_{x\in W_n}$ for an appropriate $n\in I$. 
This implies the vanishing $\check{H}^p(W,F)\cong 0$ of the \v{C}ech cohomology groups  for $p\ge 1$ and every presheaf $F\in \Pr_\Ab\bS$.

Let $\cH^q=R^q i$ be the derived functor
of the embedding $i:\Sh_\Ab\bS\to \Pr_\Ab\bS$ of sheaves into presheaves.
We now consider the  
\v{C}ech cohomology spectral sequence \cite[3.4.4]{MR1317816}
$(E_r,d_r)\Rightarrow R^*\Gamma(W,F)$ with $E_2^{p,q}\cong \check{H}^{p}(W;\cH^q(F))$ and use
 \cite[3.4.3]{MR1317816} to the effect that
$\check{H}^{0}(W;\cH^q(F))\cong 0$ for all $q\ge 1$.
Combining these two vanishing results we see that
the only non-trivial term of the second page of the spectral sequence is
$E_2^{0,0}\cong  \check{H}^{0}(W;\cH^0(F))\cong F(W)$.
Vanishing of $R^i\Gamma(W;\dots)$ for $i\ge 1$ is equivalent to the exactness of $\Gamma(W;\dots)$.
\hB
 
\subsubsection{}\label{uiqwdgqwdwqdqdwqdwqd}

We now prove Lemma \ref{hdjfcdcuisc}. Let $i\ge 1$.
 We use that
$R^i\cR_W(F)$ is the sheafification of the presheaf $\bS\ni A\mapsto
R^i\Gamma(A\times W;F)\in \Ab$.
For every sheaf $F\in \Sh_\Ab\bS$ we have by some intermediate steps in (\ref{hjuziqwuxq}) 
$$\Gamma(A\times W;F)\cong \Gamma(W; \cR_A(F))\ .$$
Let us choose an injective resolution $F\to I^\bullet$. 
Using  Lemma \ref{wzdqdwqd} for the second isomorphism we get
\begin{equation}\label{hequcqx}
H^i\Gamma(A;\cR_W(I^\bullet))\stackrel{s\mapsto \tilde s}{\cong} H^i\Gamma(A\times W;I^\bullet)\stackrel{\tilde s\mapsto \bar s}{\cong}\Gamma(W; H^i\cR_A (I^\bullet))\ .
\end{equation}
Consider a point
$a\in A$ and $s\in H^i\Gamma(A;\cR_W(I^\bullet))$. 
We must show that there exists a neighbourhood $a\in U\subseteq A$ of
$a$ such that $s_{|U}=0$.
Since $I^\bullet$ is an exact sequence of sheaves in degree $\ge 1$ there exists an open covering
$\{Y_r\}_{r\in R}$ of $A\times W$ such that $\tilde s_{|Y_r}=0$. 
After refining this covering we can assume that $Y_r\cong U_r\times V_r$ for suitable open subsets $U_r\subseteq A$ and $V_r\subseteq W$.
Consider the set
$R_a:=\{r\in R|a\in U_r\}$.
Since $W$ is compact the covering 
$\{V_r\}_{r\in R_a}$ of $W$ admits a finite subcovering
indexed by $Z\subseteq R_a$. The set
$U:=\cap_{r\in Z } U_r\subseteq A$ is an open neighbourhood of $a$.
By further restriction we get $\tilde s_{|U\times V_r}=0$ for all $r\in Z$.
By (\ref{hequcqx}) this means that $0=\overline{s_{|U}}_{|V_r}\in \Gamma(W; H^i\cR_U (I^\bullet))$. Therefore $\overline{s_{|U}}$ vanishes locally on $W$ and therefore globally.
This implies $s_{|U}=0$.
\hB

\subsubsection{}

Let $f:\bS_{lc}\to \bS$ be the inclusion of the full subcategory of locally compact topological spaces.
Since an open subset of a locally compact space is again locally we can define the topology on $\bS_{lc}$ by
$$\cov_{\bS_{lc}}(A):=\cov_{\bS}(A)\ ,\quad A\in \bS_{lc}\ .$$
 
The compactly generated topology and the product topology on products of
locally compact spaces coincides. The same applies to fibre products.
Furthermore, a fibre product of locally compact spaces is locally compact.
The functor preserves fibre products. 
In view of the definition of the topology $\bS_{lc}$ the inclusion functor
$f: \bS_{lc}\to \bS$ is a morphism of sites.

\begin{lem}\label{ioieeqwe}
Restriction to the site $ \bS_{lc}$ commutes with sheafification, i.e.
$$i^\sharp\circ {}^p f^*\cong f^*\circ i^\sharp\ .$$
\end{lem}
\proof
(Compare with the proof of Lemma \ref{udiwqdiuiuqw}.)
This follows from
$\cov_{\bS_{lc}}(A)\cong\cov_{\bS}(A)$ for all $A\in \bS$ and
the explicit construction of $i^\sharp$ in terms of the
set $\cov_{\bS_{lc}}$ (see \cite[Sec. 3.1]{MR1317816}).
\hB

\begin{lem}\label{jkasjsaoxsacxs}
The restriction $f^*:\Sh\bS\to \Sh\bS_{lc}$ is exact.
\end{lem}
\proof (Compare with the proof of Lemma \ref{teezwdc}.)
The functor $f^*$ is a right-adjoint and thus commutes with limits.  Colimits of presheaves are defined object-wise, i.e. for a diagram
$$\cC\to \Pr\bS\ ,\quad c\mapsto F_c$$ of presheaves we have 
$$({}^p\colim_{c\in \cC} F_c)(U)=\colim_{c\in \cC} F_c(U)\ .$$
 It follows from the explicit
description of ${}^pf^*$ that this functor commutes with colimits of presheaves.
For a diagram of sheaves
$F:\cC\to \Sh\bS$ we have 
$$\colim_{c\in \cC} F_c=i^\sharp {}^p\colim_{c\in\cC} F_c\ .$$
By Lemma  \ref{ioieeqwe} we get $
f^*\colim_{c\in \cC} F_c\cong f^*i^\sharp {}^p\colim_{c\in\cC} F_c\cong i^\sharp {}^pf^*{}^p\colim_{c\in\cC} F_c\cong i^\sharp {}^p\colim_{c\in\cC} {}^p f^*F_c\cong \colim_{c\in \cC} f^*F_c\ .$
\hB

\subsubsection{}

We have now verified that $f:\bS_{lc}\to \bS$ satisfies  the assumptions of Proposition \ref{neuhsdjfjsdf674133}. 
\begin{kor}\label{wuiwqdiuwiduwqdiuqwd78}
Let $f:\bS_{lc}\to \bS$ be the inclusion of the site of locally compact spaces.
For $F\in \Sh_\Ab\bS_{lc}$ and $G\in \Sh_\Ab\bS$
we have
$$f^*R\uHom_{\Sh_\Ab\bS}(f_*F,G)\cong R\uHom_{\Sh_\Ab\bS_{lc}}(F,f^*G)\ .$$
In particular we have
$$f^*\uExt^k_{\Sh_\Ab\bS}(f_*F,G)\cong \uExt^k_{\Sh_\Ab\bS_{lc}}(F,f^*G)$$
for all $k\ge 0$.
\end{kor}
In fact, the first assertion implies the second since $f^*$ is exact.

We need this result in the following special case.
If $G\in \bS$ is a locally compact group, then by abuse of notation we
write $\uG$ for the sheaves of abelian group represented by
$G$ in both categories $\Sh_\Ab\bS$ and $\Sh_\Ab\bS_{lc}$.

We have
 $f^*\uG=\uG$. By Lemma \ref{zuuzud} we also have
$f_*\uG\cong \uG$.

\begin{kor}\label{uiwddiqwdwqdddqd}
Let $G,H\in \bS_{lc}$ be locally compact abelian groups. We have
$$f^*\uExt^k_{\Sh_\Ab\bS}(\uG,\uH)\cong \uExt^k_{\Sh_\Ab\bS_{lc}}(\uG,\uH)$$
for all $k\ge 0$.
\end{kor}

\subsubsection{}\label{locacyintro}

In some places we will need a second sub-site of $\bS$, the site
$\bS_{loc-acyc}$ of locally acyclic spaces. 
\begin{ddd}
A space $U\in \bS$ is called acyclic, if $H^i(U;\uH)\cong 0$
for all discrete abelian groups $H$ and $i\ge 1$.
\end{ddd}
By Lemma \ref{wzdqdwqd} all profinite spaces are acyclic.
The space $\R^n$ is another example of an acyclic space.
In fact, the homotopy invariance used in the proof of \ref{gfrgh3z244r82r2324}
shows that the inclusion $0\to \R^n$ induces an isomorphism 
$H^i(\R^n;\uH)\cong H^i(\{0\};\uH)$, and  a one-point space is clearly acyclic.

\begin{ddd}
A space $A\in \bS$ is called locally acyclic if it admits an open covering by acyclic
spaces.
\end{ddd}

In general we do not know if the product of two locally acyclic spaces
is again locally acyclic (the K{\"u}nneth formula needs a compactness assumption).
In order to ensure the existence of finite products we consider the combination of the conditions
locally acyclic and locally compact.

Note that all finite-dimensional manifolds   are locally acyclic and locally compact.
An open subset of a locally acyclic locally compact space is again locally acyclic.
We let $\bS_{lc-acyc}\subset \bS$ be the full subcategory of locally acyclic locally compact
spaces. The topology on 
$\bS_{lc-acyc}$ is given by
$$\cov_{\bS_{lc-acyc}}(A):=\cov_{\bS}(A)\ .$$
Let $g:\bS_{lc-acyc}\to \bS$ be the inclusion.
The proofs of Lemma \ref{ioieeqwe}, Lemma \ref{jkasjsaoxsacxs}
and Corollary \ref{wuiwqdiuwiduwqdiuqwd78}
apply verbatim.
\begin{kor}
\begin{enumerate}
\item The restriction $g^*:\Sh\bS\to \Sh\bS_{lc-acyc}$ is exact.
\item For $F\in \Sh_\Ab\bS_{lc-acyc}$ and $G\in \Sh_\Ab\bS$
we have
$$g^*R\uHom_{\Sh_\Ab\bS}(g_*F,G)\cong R\uHom_{\Sh_\Ab\bS_{lc-acyc}}(F,g^*G)\ .$$
and
$$g^*\uExt^k_{\Sh_\Ab\bS}(g_*F,G)\cong \uExt^k_{\Sh_\Ab\bS_{lc-acyc}}(F,g^*G)$$
for all $k\ge 0$.
\end{enumerate}
\end{kor}

\begin{kor}\label{uiwddiqwdwqdddqd2121}
Let $G,H\in \bS_{lc-acyc}$ be locally acyclic abelian groups. We have
$$g^*\uExt^k_{\Sh_\Ab\bS}(\uG,\uH)\cong \uExt^k_{\Sh_\Ab\bS_{lc-acyc}}(\uG,\uH)$$
for all $k\ge 0$.
\end{kor}

Note that the product
$\prod_{\nat} \T$ is compact but not locally acyclic.

\subsection{$\Z_{mult}$-modules}\label{udwwuqidqqw434}

\subsubsection{}\label{dkwqwjdkqkdjqwdw}

We consider the multiplicative semigroup $\Z_{mult}$ of non-zero integers. 
Every abelian group $G$ (written multiplicatively) has a tautological action of $\Z_{mult}$ by homomorphisms given by
$$\Z_{mult}\times G\to G\ ,\quad (n,g)\mapsto g^n\ .$$
\begin{ddd}
When we consider $G$ with this action we write $G(1)$.
\end{ddd}

The semigroup $\Z_{mult}$ will therefore act on all objects naturally constructed from an abelian group $G$.
We will in particular use the action of $\Z_{mult}$ on the group-homology and group-cohomology of $G$.

If $G$ is a topological group, then $\Z_{mult}$ acts by continuous maps. It therefore also acts on the cohomology of $G$ as a topological space. Also this action will play a role later.
 
The remainder of the present subsection sets up some language related to $\Z_{mult}$-actions.

\subsubsection{}

\begin{ddd}
A $\Z_{mult}$-module is an abelian group with an action of $\Z_{mult}$ by homomorphisms.
\end{ddd}
We will write this action as 
$\Z_{mult}\times G\ni (m,g)\mapsto \Psi^m(g)\in G$.
Thus in the case of $G(1)$ we have $\Psi^m(g)=g^m$.

\subsubsection{}

We let $\Z_{mult}-\Mod$ denote the category of $\Z_{mult}$-modules. An equivalent description of this category is as the category of modules under the 
commutative semigroup ring $\Z[\Z_{mult}]$. The category $\Z_{mult}-\Mod$ is an abelian tensor category.

We have an exact inclusion of categories
$$\Ab\hookrightarrow (\Z_{mult}-\Mod)\ , \quad G\mapsto G(1)$$ 
By $$\cF:(\Z_{mult}-\Mod)\to \Ab$$ we denote the forgetful functor.

\begin{ddd}Let $G$ be an abelian group. 
For $k\in \Z$ we let $G(k)\in \Z_{mult}-\Mod$ denote the $\Z_{mult}$-module given by the action
$\Z_{mult}\times G\ni (p,g)\mapsto \Psi^p(g):=g^{p^k}\in G$.  
\end{ddd}
Observe that for abelian groups $V,W$ we have a natural isomorphism
\begin{equation}\label{hwqdwqdwqd}V(k)\otimes_\Z W(l)\cong (V\otimes_\Z W)(k+l)\ .\end{equation}

\subsubsection{}

Let $V$ be a $\Z_{mult}$-module.
\begin{ddd}
We say that $V$ has weight $k$ if there exists an isomorphism 
$V\cong \cF(V)(k)$ of $\Z_{mult}$-modules.
\end{ddd} 
If $V$ has weight $k$, then every sub-quotient of $V$ has weight $k$.
Note that a $\Z_{mult}$-module can have many weights. We have e.g. isomorphisms of $\Z_{mult}$-modules
$(\Z/2\Z)(1)\cong (\Z/2\Z)(k)$ for all $k\not=0$. 

\subsubsection{}

Let $V\in \Z_{mult}-\Mod$. We say that $v\in V$ has weight $k$ if it generates a submodule
$\Z<v>\subset V$ of weight $k$.
For $k\in \nat$ we let $W_k:(\Z_{mult}-\Mod)\to \Ab$ be the functor which associates to $V\in  \Z_{mult}-\Mod$ its subgroup of vectors of weight $k$. Then we have an adjoint pair
of functors
 $$(k):\Ab\Leftrightarrow (\Z_{mult}-\Mod):W_k\ .$$
 Observe that the functor $W_k$ is not exact.
Consider for example a prime number $p\in \nat$ and 
the sequence
$$0\to \Z(1)\stackrel{p}{\to} \Z(1)\to (\Z/p\Z)(p)\to 0\ .$$
The projection map is indeed $\Z_{mult}$-equivariant since
$m^p\equiv m\:\mbox{mod $p$}$ for all $m\in \Z$.
Then $$0\cong W_p(\Z(1))\to W_p((\Z/p\Z)(p))\cong \Z/p\Z$$ is not surjective.
 
\subsubsection{}\label{fbefewfwef}

 Let $V\in \Ab$ and $V(1)\in \Z_{mult}-\Mod$. Then we can form the graded tensor algebra 
 $$T^*_\Z(V(1))=\Z\oplus V(1)\oplus V(1)\otimes_\Z V(1)\oplus\dots\ .$$
 We see that $T^k(V(1))$ has weight $k$.
 The elements $x\otimes x\in V(1)\otimes V(1)$
 generate a homogeneous ideal $I$. Hence the graded algebra
 $\Lambda_\Z^* (V(1)):=T^*_\Z(V(1))/I$ has the property that
 $\Lambda^k_\Z(V(1))$  has weight $k$.

\subsubsection{}

It makes sense to speak of a sheaf or presheaf of $\Z_{mult}$-modules on the site $\bS$. We let
$\Sh_{\Z_{mult}-\Mod}\bS$ and $\Pr_{\Z_{mult}-\Mod}\bS$  denote the corresponding abelian categories of sheaves and presheaves. 
\begin{ddd}
Let $V\in \Sh_{\Z_{mult}-\Mod}\bS$ (or $V\in \Pr_{\Z_{mult}-\Mod}\bS$) and $k\in \Z$.
We say that $V$ is of weight $k$ if the map
$V\stackrel{m^k-m}{\to} V$ vanishes for all $m\in \Z_{mult}$.
\end{ddd}
We define the functors
$(k):\Sh_\Ab\bS\to \Sh_{\Z_{mult}-\Mod}\bS$, $\cF:\Sh_{\Z_{mult}-\Mod}\bS\to \Sh_\Ab\bS$,  and
$W_k:\Sh_{\Z_{mult}-\Mod}\bS\to \Sh_\Ab\bS$ (and their presheaf versions) object-wise.
We also have  a pair of adjoint functors
$$(k):\Sh_\Ab\bS\Leftrightarrow \Sh_{\Z_{mult}-\Mod}\bS:W_k$$
(and the corresponding  presheaf version).
A sheaf $V\in \Sh_{\Z_{mult}-\Mod}\bS$ of $\Z_{mult}$-modules has weight $k\in \Z$ if
$V\cong \cF(V)(k)\cong W_k(V)(k)$.

\section[Admissible topological groups]{Admissibility of sheaves represented by topological abelian groups}\label{squfceweqwdqwd}

\subsection{Admissible sheaves and groups}

\subsubsection{}

The main topic of the present paper is a duality theory for abelian group stacks (Picard stacks, see \ref{zduqwgdwqdqdwd}) on the site $\bS$. A Picard stack $P\in \cPic(\bS)$ gives rise to the sheaf
of objects $H^0(P)$ and the sheaf of automorphisms of the neutral object $H^{-1}(P)$.
These are sheaves of abelian groups on $\bS$.

We will define the notion of a dual Picard stack $D(P)$ (see \ref{wdkwqdidqiddqw}). With the intention to generalize
the Pontrjagin duality for locally compact abelian groups to group stacks we study the question under which conditions the natural map $P\to D(D(P))$ is an isomorphism. In Theorem \ref{hjjhjad78234} we see that this is the case if
the sheaves are dualizable (see \ref{fuhwefiwfewfewf}) and admissible (see
\ref{iuefefewffwf}). Dualizability is a sheaf-theoretic   generalization of
the classical Pontrjagin duality and is satisfied e.g.~for the sheaves $\uG$ for locally compact groups $G\in \bS$ (see \ref{hj444cq}). Admissibility is more exotic and will be defined below (\ref{iuefefewffwf}). One of the main results
of the present paper asserts that the sheaves $\uG$ are admissible for a large (but not exhaustive) class  of locally compact groups $G\in \bS$, and for an even larger class if $\bS$ is replaced by $\bS_{lc-acyc}$.

\subsubsection{}
\begin{ddd}\label{iuefefewffwf}
We call a sheaf of groups $F$ admissible if $\uExt^i_{\Sh_\Ab \bS}(F,\uT)\cong 0$
for $i=1,2$. A topological abelian group $G\in \bS$ is called admissible if $\uG$ is an admissible sheaf.
\end{ddd}

We will also consider the sub-sites $\bS_{lc-acyc}\subset \bS_{lc}\subset \bS$ of locally compact spaces.
\begin{ddd}\label{iuefefewffwf1}
 A locally compact abelian topological group $G\in \bS_{lc}$ is called admissible on $\bS_{lc}$ if 
$\uExt^i_{\Sh_\Ab \bS_{lc}}(\uG,\uT)\cong 0$
for $i=1,2$ (and correspondingly for $\bS_{lc-acyc}$). 
\end{ddd}

\begin{lem}Let $f:\bS_{lc}\to \bS$ be the inclusion.
If  $F\in \Sh_\Ab\bS_{lc}$ and  $f_*F$ is admissible, then $F$ is admissible on $\bS_{lc}$. 
\end{lem}
\proof
This is an application of Corollary \ref{wuiwqdiuwiduwqdiuqwd78}.\hB

\begin{kor}\label{ioodqwdqwd}
If $G\in \bS_{lc}$ is admissible, then it is admissible on $\bS_{lc}$.
\end{kor} 
This is an application of Corollary \ref{uiwddiqwdwqdddqd}.
\subsubsection{}

\begin{lem}\label{wuihfewfefwfw}
 The class of admissible sheaves is closed under finite products and extensions.
\end{lem}
\proof
For finite products the assertions follows from the fact that
$\uExt_{\Sh_\Ab\bS}(\dots,\uT)$ commutes with finite products.  Given an extension of sheaves
$$0\to F\to G\to H\to 0$$
such that $F$ and $H$ are admissible, then also $G$ is admissible.
This follows immediately from the long exact sequence obtained by applying
$\uExt_{\Sh_\Ab\bS}^*(\dots,\uT)$. 
\hB

\subsubsection{}

In this paragraph we formulate one of the main theorems of the present paper.
Let $G$ be a locally compact topological abelian group. We first recall some technical conditions.
\begin{ddd}\label{two-three-cond}
We say that $G$ satisfies the two-three condition, if
\begin{enumerate}
\item it does not admit $\prod_{n\in \nat} \Z/2\Z$ as a subquotient,
\item  the multiplication by $3$ on the component $G_0$ of the identity has
  finite cokernel.
\end{enumerate}
\end{ddd}

\begin{ddd}\label{ltdp} We say that $G$ is locally topologically divisible if for all primes $p\in \nat$ the multiplication map
$p:G\to G$ has a continuous local section.
\end{ddd}

Using that the class of admissible groups is closed under finite products and extensions we get the following general theorem. 

\begin{theorem}\label{wdwquidioqwdopqwdq}
\begin{enumerate}
\item  If $G$ is a locally compact abelian group which satisfies the two-three
  condition, then it is admissible over $\bS_{lc-acyc}$.
\item  Assume that 
\begin{enumerate}
\item $G$ satisfies the two-three condition,
\item $G$ has an open subgroup of the form $C\times \R^n$ with $C$ compact such that $G/C\times \R^n$ is finitely generated
\item the connected component of the identity of $G$ is locally topologically divisible.
\end{enumerate} 
Then $G$ is admissible over 
$\bS_{lc}$.
\end{enumerate}
\end{theorem}
 \proof
By \cite[Thm. 7.57(i)]{MR1646190} the group $G$ has a splitting
$G\cong H\times \R^n$ for some $n\in \nat_0$, where $H$ has a compact open subgroup $U$.
The quotient $G/(U\times \R^n)\cong H/U$ is therefore discrete.
Using Lemma \ref{wuihfewfefwfw} conclude that
$G$ is admissible over $\bS_{lc}$ if $\R^n$ and $H$ are so.

Admissibility of $\R^n$ follows from Theorem \ref{udiqwdwqdwd} (which we prove later) in conjunction
with \ref{wuihfewfefwfw} and \ref{ioodqwdqwd}.

Since $D:=H/U$ is discrete, the exact sequence
$$0\to U\to H\to H/U\to 0$$
has local sections and therefore induces an exact sequence of associated sheaves 
$$0\to \uU\to \uH\to \uD \to 0$$
by Lemma 
\ref{dezuqwideqwd}.
If $D$ is finitely generated, then it is admissible by Theorem \ref{wjebdqwdwdqdwqd}, and hence admissible on $\bS_{lc}$  by Lemma \ref{ioodqwdqwd}.
Otherwise it is admissible on $\bS_{lc-acyc}$ by Theorem \ref{uwdiqwdqwdwqdwqd}.

Therefore $H$ is admissible on $\bS_{lc}$  (or $\bS_{lc-acyc}$,
respectively) by Lemma \ref{wuihfewfefwfw} if $U$ is admissible over 
$\bS_{lc}$ (or $\bS_{lc-acyc}$,
respectively). 
The compact group $U$ fits into an exact  sequence
$$0\to U_0\to U\to P\to 0$$ where $P$ is profinite and $U_0$ is closed and connected.
By assumption $U$, $U_0$ and $P$ satisfy the two-three condition.
The connected compact group $U_0$ is admissible on $\bS_{lc}$ (or on  $\bS_{lc}$ without the assumption that it is locally topologically divisible), by Theorem \ref{wqidduiwdwddwduiqwdiu}, and
the profinite $P$ is admissible by Theorem \ref{wqdiqwdiudwqdd}, and hence admissible
on $\bS_{lc}$ by  Lemma \ref{ioodqwdqwd}.

Note that 
$$0\to \uU_0\to \uU\to \uP\to 0$$ is exact by Lemma \ref{jhbuisaicascsa1}.
Now it follows from  Lemma \ref{wuihfewfefwfw}
that $U$ is admissible on $\bS_{lc}$ (or $\bS_{lc-acyc}$, respectively).
\hB 

\subsubsection{}

We conjecture that the assumption that $G$ satisfies the two-three condition is only technical and forced by our technique
to prove admissibility  of compact groups.

In the remainder of this section, we will mainly be concerned with the proof
of the statements used in the proof of Theorem \ref{wdwquidioqwdopqwdq} above.

\subsubsection{}

Our proofs of admissibility for a sheaf $F\in \Sh_\Ab\bS$ will usually be based on the following argument.
\begin{lem}\label{udiuddwqwqd}
Assume that $F\in \Sh_\Ab\bS$ satisfies
\begin{enumerate}
\item $\uExt^i_{\Sh_\Ab\bS}(F,\uZ)\cong 0$ for $i=2,3$
\item $\uExt^i_{\Sh_\Ab\bS}(F,\uR)\cong 0$ for  $i= 1,2$.
\end{enumerate}
Then $F$ is admissible.
\end{lem}
 \proof
We apply the functor $\uExt_{\Sh_\Ab\bS}^*(F,\dots)$ to the sequence 
$$0\to \uZ\to \uR\to \uT\to 0$$
and get the following segments of the long exact sequence
$$\dots\to  \uExt^i_{\Sh_\Ab\bS}(F,\uR)\to
\uExt^i_{\Sh_\Ab\bS}(F,\uT)\to \uExt^{i+1}_{\Sh_\Ab\bS}(F,\uZ)\to\dots\ .$$
We see that the assumptions on $F$ imply that $\uExt^i_{\Sh_\Ab\bS}(F,\uT)\cong 0$ for $i=1,2$.
\hB 

\subsubsection{}

A space $W\in \bS$ is called profinite if it can be written as the limit of
an inverse system of finite spaces. Lemma \ref{hdjfcdcuisc} has as a special
case the following theorem.
\begin{theorem}\label{uihdqewddwqdwqd}
If $W\in \bS$ is profinite, then $\Z(\uW)$ is admissible.
\end{theorem}

\subsection{A double complex}

\subsubsection{}\label{ufiwefuew}

Let $H\in \Sh_\Ab\bS$ be a sheaf of groups with underlying sheaf of sets $\cF(H)\in \Sh\bS$.
Applying the linearization functor $\Z$ (see \ref{azuzasudsad87}) we get again a sheaf of groups
$\Z(\cF(H))\in \Sh_\Ab\bS$. The group structure of $H$ induces on $\Z(\cF(H))$ a ring structure.
We denote this sheaf of rings by $\Z[H]$. It is the sheafification of the presheaf ${}^p\Z[H]$ which 
associates to  $A\in \bS$ the integral group ring ${}^p\Z[H](A)$ of the group $H(A)$.

We consider the category $\Sh_{\Z[H]-\Mod}\bS$ of sheaves of $\Z[H]$-modules. The trivial action of the sheaf of groups $H$ on the sheaf $\uZ$ induces the structure of a sheaf of $\Z[H]$-modules on $\uZ$.
With the notation introduced below we could (but refrain from doing this) write this sheaf of $\Z[H]$-modules
as $\coind(\uZ)$.

 \subsubsection{}

The forgetful functor $\res:\Sh_{\Z[H]-\Mod}\bS\to \Sh_\Ab\bS$ fits into an adjoint pair of functors
$$\ind:\Sh_\Ab\bS\Leftrightarrow\Sh_{\Z[H]-\Mod}\bS:\res\ .$$
Explicitly, the functor
$\ind$ is given by
$$\ind(V):=\Z[H]\otimes_\Z V\ .$$
Since $\Z[H]$ is  a torsion-free sheaf and therefore a sheaf of flat $\uZ$-modules the functor $\ind$ is exact. Consequently the functor
$\res$ preserves injectives.

\subsubsection{} 

We let $\coinvv:\Sh_{\Z[H]-\Mod}\bS\to \Sh_\Ab\bS$ denote the coinvariants functor given by
$$\Sh_{\Z[H]-\Mod}\bS\ni V \mapsto  \coinvv(V):=V\otimes_{\Z[H]}\uZ\in \Sh_\Ab\bS\ .$$
This functor fits into the adjoint pair
$$\coinvv:\Sh_{\Z[H]-\Mod}\bS\Leftrightarrow\Sh_\Ab\bS:\coind$$
with the coinduction functor which maps $W\in \Sh_\Ab\bS$ to the
sheaf of $\Z[H]$-modules induced by the trivial action of $H$ in $W$, formally this can be written as
 $$\coind(W):=\uHom_{\Sh_{\Ab}\bS}(\uZ,W)\ .$$
 
\subsubsection{}\label{wjkeffqfqwfwqfwqf}

\begin{lem}
Every sheaf $F\in \Sh_{\Z[H]-\Mod}\bS$ is a quotient of a flat sheaf.
\end{lem}
\proof
Indeed, the counits of the adjoint pairs $(\Z(\dots),\cF)$ and $(\ind,\res)$ induce a surjection
$\ind(\Z(\cF(\res (F))))\to F$.  Explicitly it is given by the composition of the sum and action map
(omitting to write some forgetful functors) $$\Z[H]\otimes_\Z\Z [F] \to  \Z[H]\otimes F\to F\ .$$
Since $\Z[H]$ is  a sheaf of unital rings this action is surjective.
Moreover, for $A\in \Sh_{\Z[H]-\Mod}\bS$ we have
$$A\otimes_{\Z[H]}(\Z[H]\otimes_\Z\Z [F])\cong A \otimes_\Z  \Z[F]\ .$$
Since $\Z[F]$ is a torsion-free sheaf of abelian groups the operation $A\to A\otimes_{\Z[H]}(\Z[H]\otimes_\Z\Z [F])$ preserves exact sequences in $\Sh_{\Z[H]-\Mod}\bS$.
Therefore $\Z[H]\otimes_\Z\Z [F]$ is a flat sheaf of $\Z[H]$-modules.
 \hB

\subsubsection{}
\begin{lem}
The class of flat $\Z[H]$-modules is $\coinvv$-acyclic.
\end{lem}
\proof
Let $F^\bullet$ be an exact complex of flat $\Z[H]$-modules.
We choose a flat resolution $P^\bullet\to \uZ$ in $\Sh_{\Z[H]-\Mod}\bS$
which exists by \ref{wjkeffqfqwfwqfwqf}\footnote{Note that $P^\bullet$ is a homological complex, i.e. the differentials have degree $-1$.}.
Since $F^\bullet$ consists of flat modules  the induced map
$$F^\bullet\otimes_{\Z [H]} P^\bullet\to F^\bullet\otimes_{\Z [H]} \uZ=\coinvv(F^\bullet)$$
is a quasi-isomorphism. Since $P^\bullet$ consists of flat modules, tensoring by $P^\bullet$ commutes with taking cohomology, so that we have
$$H^*( F^\bullet\otimes_{\Z [H]} P^\bullet)\cong H^*(H^*(F^\bullet)\otimes _{\Z [H]} P^\bullet)\cong 0\ .$$
Therefore $\coinvv(\dots)$ maps acyclic complexes of flat $\Z[H]$-modules to acyclic complexes of $\uZ$-modules.
\hB
\begin{kor}
We can calculate $L^*\coinvv(\uZ)$ using a flat resolution.
\end{kor}

\subsubsection{}\label{iuqwhdiqdwqdwqd}

We will actually work with a very special flat resolution of $\uZ$.
The bar construction on the sheaf of groups $H$ gives a sheaf $H^\bullet$ of simplicial sets with an action of $H$.
We let $C^\bullet(H):=C(\Z(H^\bullet))$ be the sheaf of homological chain complexes associated to the sheaf of simplicial groups $\Z(H^\bullet)$. The $H$-action on $H^\bullet$ induces a $H$-action on $C^\bullet(H)$ and therefore the structure of a sheaf of $\Z[H]$-modules.
In order to understand the structure of $C^\bullet(H)$ we first consider the presheaf version
${}^pC^\bullet(H):=C({}^p\Z(H^\bullet))$.
In fact, we can write $$H^i\cong H\times \underbrace{H\times\dots\times H}_{i\:factors}\ ,$$
as a sheaf of $H$-sets, 
and therefore
\begin{equation}\label{uibiquwdwdqwd1213}
{}^pC^i(H)\cong {}^p\Z[H]\otimes^p_\Z {}^p\Z(\underbrace{H\times\dots\times H}_{i\:factors})
\end{equation}
The cohomology of the complex ${}^pC^\bullet(H)$  is given by
$$H^i({}^pC^\bullet(H))\cong \left\{\begin{array}{cc}
{}^p\uZ&i=0\\0&i\ge 1
\end{array}\right\}\ ,$$
where ${}^p\uZ\in \Pr_\Ab\bS$ denotes the constant presheaf with value $\Z$.
Since sheafification $i^\sharp$ is an exact functor and by definition $C^\bullet(H)=i^\sharp {}^pC^\bullet(H)$ we get
$$H^i(C^\bullet(H))\cong \left\{\begin{array}{cc}
\uZ&i=0\\0&i\ge 1
\end{array}\right\}\ .$$
Furthermore,
\begin{equation}\label{uibiquwdwdqwd}
C^i(H)\cong \Z[H]\otimes_\Z \Z(\underbrace{H\times\dots\times H}_{i\:factors})=\ind(\Z(\underbrace{H\times\dots\times H}_{i\:factors}))
\end{equation}
 shows that 
$C^i(H)$ is flat.
Let us write
$C^\bullet:=C^\bullet(H)=\ind(D^\bullet)$
with $$D^i:=\Z(\underbrace{H\times\dots\times H}_{i\:factors})\ .$$
 
\begin{ddd}\label{iudhuiqwdqwdwqd}
For a sheaf $H\in \Sh_\Ab\bS$ we define
the complex $U^\bullet:=U^\bullet(H):=\coinvv(C^\bullet)$.
\end{ddd}

It follows from the construction that $U^\bullet$ depends functorially on $H$. In particular, for  $\alpha:H\to H^\prime$ we have a map of complexes
$U^\bullet(\alpha):U^\bullet(H)\to U^\bullet(H^\prime)$.

\subsubsection{}

The main tool in our proofs of admissibility of a sheaf $F$ is the study of  the sheaves
$$R^*\uHom_{\Sh_{\Z[H]-\Mod}\bS}(\uZ,\coind(W))$$
for $W=\uZ,\uR$.
Let us write this in a more complicated way using the special flat resolution $C^\bullet(H)\to \uZ$ constructed in \ref{iuqwhdiqdwqdwqd}. We choose an injective resolution
$\coind(W)\to I^\bullet$ in $\Sh_{\Z[H]-\Mod}\bS$. 
 Using that $\res\circ \coind=\id$, and that $\res(I^\bullet)$ is injective in $\Sh_\Ab\bS$, we get
\begin{eqnarray*}
R\uHom_{\Sh_{\Z[H]-\Mod}\bS}(\uZ,\coind(W))&\cong&\uHom_{\Sh_{\Z[H]-\Mod}\bS}(\uZ,I^\bullet)\\
&\cong&\uHom_{\Sh_{\Z[H]-\Mod}\bS}(C^\bullet(H),I^\bullet)\\
&\cong&\uHom_{\Sh_{\Z[H]-\Mod}\bS}(\ind(D^\bullet),I^\bullet)\\
&\cong&\uHom_{\Sh_\Ab\bS}(D^\bullet,\res(I^\bullet))\\
&=&\uHom_{\Sh_\Ab\bS}(D^\bullet,\res(\coind(\res(I^\bullet))))\\
&=&\uHom_{\Sh_{\Z[H]-\Mod}\bS}(\ind(D^\bullet),\coind(\res(I^\bullet)))\\
&\cong&\uHom_{\Sh_{\Z[H]-\Mod}\bS}(C^\bullet,\coind(\res(I^\bullet)))\\
&\cong&\uHom_{\Sh_\Ab\bS}(\coinvv(C^\bullet),\res(I^\bullet))\\
&\cong&R\uHom_{\Sh_\Ab\bS}(L\coinvv(\uZ),W)\ .
\end{eqnarray*}

\subsubsection{}

In general, the coinvariants functor $\coinvv(\dots)=\dots\otimes_{\Z[H]}\uZ$ can be written in terms of the tensor product  in the sense of presheaves composed with a sheafification. Furthermore, $C^\bullet(H)$ is the sheafification
of ${}^pC^\bullet(H)$.
Using the fact that the tensor product of presheaves commutes with sheafification
\begin{eqnarray}\label{nbmncm87324234asc}U^i&\stackrel{Definition \ref{iudhuiqwdqwdwqd}}{=}&
\coinvv(C^i)\\&\cong&(C^i(H)\otimes^p_{\Z[H]} \uZ)^\sharp\nonumber\\
&=&(({^p}C^i(H))^\sharp\otimes^p_{({}^p\Z[H])^\sharp} ({}^p\uZ)^\sharp)^\sharp\nonumber\\
&\cong&({}^pC^i(H)\otimes^p_{{}^p\Z[H]} {}^p\uZ)^\sharp\nonumber\\
&\stackrel{(\ref{uibiquwdwdqwd1213})}{\cong}&({}^p\Z[H]\otimes^p_\Z {}^p\Z(\underbrace{H\times\dots\times H}_{i\:factors})\otimes^p_{{}^p\Z[H]} {}^p\uZ)^\sharp\nonumber\\
&\cong&({}^p\Z(\underbrace{H\times\dots\times H}_{i\:factors}))^\sharp\nonumber\\
&=&({}^pD^i)^\sharp\nonumber
\end{eqnarray}
with 
\begin{equation}\label{uebiwqdwq}
{}^pD^i=:{}^p\Z(\underbrace{H\times\dots\times H}_{i\:factors})\ .
\end{equation}
In particular, we have 
\begin{equation}\label{uieiefefewfw}
U^i=\Z(H^i)\ .
\end{equation}

\subsubsection{}\label{iudiqwdwqdwqd}
Let $G$ be a group.
Then we can form the standard reduced bar complex for the group homology with integer coefficients
$$G^\bullet \dots \to \Z(G^n)\to \Z(G^{n-1})\to \dots\to \Z \to 0\ .$$
The abelian group $\Z(G^n)$ (sitting in degree $n$) is freely generated by the underlying set of
$G^n$, and we write the generators in the form $[g_1|\dots|g_n]$.
The differential is given by 
$$d=\sum_{i=0}^{n}(-1)^id_i:\Z(G^n)\to \Z(G^{n-1})\ ,$$ where
$$d_i[g_1|\dots|g_n]:=\left\{
\begin{array}{cc} 
{}[g_2|\dots|g_n]&i=0\\
{}[g_1|\dots|g_ig_{i+1}|\dots|g_n]&1\le i\le n-1\\
{}[g_1|\dots|g_{n-1}]&i=n\end{array}\right. \ .$$
The cohomology of this complex is the group homology $H_*(G;\Z)$. 

\subsubsection{}

For $A\in \bS$ the complex
${}^pD^\bullet(A)$ (see \ref{uebiwqdwq}) is exactly the standard complex (see \ref{iudiqwdwqdwqd}) for the group homology $H_*(H(A),\uZ)$ of the group $H(A)$.  
The cohomology sheaves $H^*(U^\bullet)$ are thus the sheafifications of the cohomology presheaves
$$\bS\ni A\mapsto H_*(H(A),\uZ)\in \Ab\ .$$

\subsubsection{}\label{lazdzu32d898}

In this paragraph we collect some facts about the homology of abelian groups.
An abelian group $V$ is the same thing as a $\Z$-module. We define the graded $\Z$-algebra $\Lambda^*_\Z V$ as the quotient of the tensor algebra $$T_\Z V :=\bigoplus_{n\ge 0}\underbrace{ V\otimes_\Z\dots\otimes_\Z V}_{n\:factors}$$
by the graded ideal $I\subseteq T_\Z V $ generated by the elements $x\otimes x$, $x\in V$.

\subsubsection{}

Let $G$ be an abelian group.  We refer to \cite[V.6.4]{MR672956} for the following fact.
\begin{fact}\label{eidoewd}
There exists a canonical map $$m:\Lambda^i_\Z G\to H_i(G;\Z)\ .$$
It is an isomorphism for $i=0,1,2$, and it becomes an isomorphism after tensoring with $\Q$ for all $i\ge 0$. If $G$ is torsion-free, then it is  an isomorphism $m:\Lambda^i_\Z G\stackrel{\sim}{\to} H_i(G;\Z)$ for all $i\ge 0$.
\end{fact}

\subsubsection{}\label{diqwodwqdwqdq}

Let $U^\bullet=U^\bullet(H)$ (see Definition \ref{iudhuiqwdqwdwqd}) for $H\in \Sh_\Ab\bS$.
The cohomology sheaves 
$L^*\coinvv(\uZ))\cong H^*(U^\bullet)$ are the sheafifications of the presheaves
$H^*({}^pD^\bullet)$. By the fact \ref{eidoewd} we have
$H^i({}^pD^\bullet)\cong \Lambda^i_\Z H$ for $i=0,1,2$ for  the presheaves $\bS\ni U\mapsto \Lambda^i_\Z H(U)$.
In particular, we have $H^0(U^\bullet)\cong \uZ$ and $H^1(U^\bullet)\cong H$.
If $H$ is a torsion-free sheaf, then $H^i(U^\bullet)\cong (\Lambda^i_\Z H)^\sharp$ for all $i\ge 0$.
Finally, for an arbitrary sheaf $H\in \Sh_\Ab\bS$ we have $$(\Lambda^*_\Z H)^\sharp\otimes_\Z\Q\cong H^*(U^\bullet)\otimes_\Z\Q\ .$$

\subsubsection{}\label{uiwdiqwdwqdq}

Our application of this relies on the study of the two spectral sequences converging to
$$\uExt^*_{\Sh_\Ab\bS}(L\coinvv(\uZ),\uZ)\cong \uExt^*_{\Sh_{\Z[H]-\Mod}\bS}(\uZ,\uZ)\\ .$$
We choose an injective resolution
$\uZ\to I^\bullet$ in $\Sh_\Ab\bS$.
Then 
$$\uExt^*_{\Sh_\Ab\bS}(L\coinvv(\uZ),\uZ)\cong H^*(\uHom_{\Sh_\Ab\bS}(U^\bullet,I^\bullet))\ .$$

The first spectral sequence denoted by
$(F_r,d_r)$ is obtained by taking the cohomology in the $U^\bullet$-direction first. Its second page is given by 
$$F_2^{p,q}\cong \uExt^p_{\Sh_\Ab\bS}(L^q\coinvv(\uZ),\uZ)\ .$$
This page contains the object of our interest, namely by \ref{diqwodwqdwqdq} the sheaves of groups
$$F_2^{p,1}\cong \uExt^p_{\Sh_\Ab\bS}(H,\uZ)\ .$$

The other spectral sequence
$(E_r,d_r)$ is obtained by taking the cohomology in the $I^\bullet$-direction first.
In view of (\ref{nbmncm87324234asc}) its  first page is
given by
$$E_1^{p,q}\cong \uExt^q_{\Sh_\Ab\bS}(\Z(H^p),\uZ)\ ,$$
which can be evaluated easily in many cases.

Let us note that
$$H^*R\uHom_{\Sh_{\Z[H]-\Mod}\bS}(\uZ,\uZ)\cong \uExt^*_{\Sh_{\Z[H]-\Mod}\bS}(\uZ,\uZ)$$
has the structure of a graded ring with multiplication given by the Yoneda product.

\subsubsection{}

We now verify Assumption $2$ of Lemma 
\ref{udiuddwqwqd} for all compact groups.

\begin{prop}\label{hdfjsdfiuiuwef}
Let $H\in\bS$ be a compact group.
Then we have
$\uExt^i_{\Sh_\Ab\bS}(\uH,\uR)\cong 0$ for $i= 1,2$.
\end{prop}
As in Definition \ref{iudhuiqwdqwdwqd} let
 $U^\bullet:=U^\bullet(\uH)$.
 Let $\uR\to I^\bullet$ be an injective resolution.
Then  we get a double complex
$\uHom_{\Sh_\Ab\bS}(U^\bullet,I^\bullet)$.

We first take the cohomology in the $I^\bullet$-, and then in the $U^\bullet$-direction.
We get a spectral sequence with first term 
$$E^{p,q}_1\cong \uExt^q_{\Sh_\Ab\bS}(\Z(\cF \uH^p),\uR)\ .$$
It follows from Corollary \ref{sofdmadd}, 2., that   $E_1^{p,q}\cong 0$ for $q\ge 1$.

We consider the complex
$$C^\bullet(H,\R):0\to \Map(H,\R)\to \dots\to \Map(H^{p-1},\R)\to \Map(H^p,\R)\to\dots$$
of topological groups
which calculates the continuous group cohomology $H_{cont}^*(H;\R)$ of $H$ with coefficients in $\R$. 
Now observe that by the exponential law for $A\in \bS$
$$\uExt^0_{\Sh_\Ab\bS}( \Z(\cF \uH^p),\uR)(A)\stackrel{Lemma \:\ref{hgeuidwqdqwdwqd}}{\cong} \Hom_\bS(H^p\times A,\R)\cong \Hom_\bS(A,\Map(H^p,\R))\ .$$ Hence
the complex $(E_1^{*,0},d_1)(A)$ is isomorphic to the complex 
$$\Hom_\bS(A,C^\bullet(H,\R))=\underline{C^\bullet(H,\R)}(A)\ .$$

Since $H$ is a compact group we have
$H^i_{cont}(H;\R)\cong 0$ for $i\ge 1$.
Of importance for us is a particular continuous chain contraction 
$h^p:\Map(H^p,\R)\to \Map(H^{p-1},\R)$, $p\ge 1$, 
which is given by the following explicit formula. 
If $c\in \Map(H^p, \R)$ is a cocycle, then we can define
$h^p(c):=b\in \Map(H^{p-1}, \R)$ by the formula
$$b(t_1,\dots,t_{p-1}):=(-1)^p\int_{H} c(t_1,\dots,t_{p-1},t)dt\ ,$$
where $dt$ is the normalized Haar measure.
Then we have $db=c$.
The maps $(h^p)_{p>0}$ induce a chain contraction $(h^p_*)_{p>0}$
of the complex $\underline{C^\bullet(H,\R)}$.
Therefore $H^i\Map(A,\cC^\bullet(H,\R))\cong H^i\underline{C^\bullet(H,\R)}(A)$ for $i\ge 1$, too.

The spectral sequence thus degenerates from the second page on.
We conclude that
$$H^i\uHom_{\Sh_\Ab\bS}(U^\bullet,I^\bullet)\cong 0\ , \quad i\ge 1\ .$$

We now take the cohomology of the double complex
$\uHom_{\Sh_\Ab\bS}^i(U^\bullet,I^\bullet)$ in the other order, first in the
$U^\bullet$-direction and then in the $I^\bullet$-direction.
In order to calculate the cohomology of $U^\bullet$ in degree $\le 2$ we use the fact \ref{eidoewd}.
We get again a spectral sequence
with second term (for $q\le 2$, using \ref{diqwodwqdwqdq})
$$F^{p,q}_2\cong\uExt^p_{\Sh_\Ab\bS}((\Lambda^q_\Z\uH)^\sharp,\uR)\ .$$
We know by Corollary \ref{sofdmadd}, 2., that $$F_2^{p,0}\cong \uExt^p_{\Sh_\Ab\bS}(\uZ,\uR)\cong 0\ ,\quad p\ge 1$$
(note that we can write $\uZ=\Z(\underline{\{*\}})$ for a one-point space). 
Furthermore note that
$$F_2^{0,1}\cong \uHom_{\Sh_\Ab\bS}(\uH,\uR)\cong \underline{\Hom_{\top-\Ab}(H,\R)}\cong 0$$
since there are no continuous homomorphisms $H\to \R$.
The second page of the spectral sequence thus has the structure
$$
\begin{array}{|c|c|c|c|c|c|}\hline
2&\uHom_{\Sh_\Ab\bS}((\Lambda_\Z^2\uH)^\sharp,\uR)&&&&\\\hline
1&0&\uExt^1_{\Sh_\Ab\bS}(\uH,\uR)&\uExt^2_{\Sh_\Ab\bS}(\uH,\uR)&&\\\hline
0&\uR&0&0&0&0\\\hline
&0&1&2&3&4\\\hline
\end{array}\ .$$
Since the spectral sequence must converge to zero in positive degrees we see that
$$\uExt^1_{\Sh_\Ab\bS}(\uH,\uR)\cong 0\ .$$

We claim that $\uHom_{\Sh_\Ab\bS}((\Lambda_\Z^2\uH)^\sharp,\uR)\cong 0$.
Note that for $A\in \bS$ we have
 $$\uHom_{\Sh_\Ab\bS}((\Lambda_\Z^2\uH)^\sharp,\uR)(A)\cong
\uHom_{\Pr_\Ab\bS}(\Lambda_\Z^2\uH,\uR)(A)\cong \Hom_{\Pr_\Ab\bS/A}(\Lambda_\Z^2\uH_{|A},\uR_{|A}) \ .$$
An element $\lambda\in \Hom_{\Pr_\Ab\bS/A}(\Lambda_\Z^2\uH_{|A},\uR_{|A})$ induces a family of 
 a  biadditive (antisymmetric) maps
$$\lambda^W:\uH(W)\times \uH(W)\to \uR(W)$$ for $(W\to A)\in \bS/A$ which is compatible with restriction. Restriction to points gives continuous biadditive maps
$H\times H\to \R$. Since $H$ is compact the only such map is the constant map to zero.
Therefore $\lambda^W$ vanishes for all $(W\to A)$.
 This proves the claim.
 
Again, since the spectral sequence $(F_r,d_r)$  must converge to zero in higher degrees  we see that
$$F_2^{2,1}\cong \uExt^2_{\Sh_\Ab\bS}(\uH,\uR)\cong 0\ .$$
This finishes the proof of the Lemma. \hB

\subsection{Discrete groups}

\subsubsection{}

In this subsection we study admissibility of discrete abelian groups.
First we show the easy fact that a finitely generated discrete abelian group is admissible.
In the second step we try to generalize this result using the representation of an arbitrary
discrete abelian group as a colimit of its finitely generated subgroups.
The functor $\uExt^*_{\Sh_\Ab\bS}(\dots,\uT)$ does not commute with colimits because of the presence
of higher $R\lim$-terms in the spectral sequence (\ref{uihdiwqdqwdd434}), below.

And in fact, not every discrete abelian group is admissible.

\subsubsection{}

\begin{lem}\label{ehjidqwwqd}
For $n\in \nat$ we have
$\uExt^i_{\Sh_\Ab\bS}(\uZ/n\uZ,\uZ)\cong 0$ for $i\ge 2$.
\end{lem}
 \proof
We apply the functor $\uExt^*_{\Sh_\Ab\bS}(\dots,\uZ)$ to the exact sequence
$$0\to \uZ\to \uZ\to \uZ/n\uZ\to 0$$
and get the long exact sequence
$$\uExt^{i-1}_{\Sh_\Ab\bS}(\uZ,\uZ)\to \uExt^{i}_{\Sh_\Ab\bS}(\uZ/n\uZ,\uZ)\to \uExt^{i}_{\Sh_\Ab\bS}(\uZ,\uZ)\to \dots\ .$$
Since by Theorem \ref{uihdqewddwqdwqd} $\uExt^{i}_{\Sh_\Ab\bS}(\uZ,\uZ)\cong 0$ for $i\ge 1$, the assertion follows.
\hB

\begin{theorem}\label{wjebdqwdwdqdwqd}
A finitely generated abelian group is admissible.
\end{theorem}
\proof
The group $\Z/n\Z$ is admissible, since Assumption $2$ of Lemma \ref{udiuddwqwqd} follows from
Proposition \ref{hdfjsdfiuiuwef}, while Assumption $1$ follows from Lemma \ref{ehjidqwwqd}.
The group $\Z$ is admissible by Theorem \ref{uihdqewddwqdwqd} since we can write $\uZ\cong \Z(\underline{\{*\}})$ for a point $\{*\}\in \bS$.
A finitely generated abelian group is a finite product of groups of the form $\Z$ and $\Z/n\Z$ for various $n\in \nat$.
A finite product of admissible groups is admissible. \hB

\subsubsection{}

We now try to extend this result to general discrete abelian groups using colimits.
Let $I$ be a filtered category (see \cite[0.3.2]{MR1317816} for definitions). The category $\Sh_\Ab\bS$ 
is an abelian Grothendieck category (see \cite[Thm. I.3.2.1]{MR1317816}). The category $\Hom_{\Cat}(I,\Sh_\Ab\bS)$ is again abelian and a Grothendieck category (see \cite[Prop. 0.1.4.3]{MR1317816}). We have the adjoint pair
$$\colim:\Hom_{\Cat}(I,\Sh_\Ab\bS)\Leftrightarrow \Sh_\Ab\bS:C_?\ ,$$
where the functor $C_?$ associates to $F\in \Sh_\Ab\bS$ the constant functor $C_F\in \Hom_{\Cat}(I,\Sh_\Ab\bS)$
with value $F$. 
The functor $C_?$ also has a right-adjoint $\lim$.
$$C_{?}:\Sh_\Ab\bS\Leftrightarrow\Hom_{\Cat}(I,\Sh_\Ab\bS):\lim$$
This functor is left-exact and admits a right-derived functor
$R\lim$.
Finally, for $F\in \Hom_{\Cat}(I,\Sh_\Ab\bS)$  we have the functor
$${}^I\uHom_{\Sh_\Ab\bS}(F,\dots): \Sh_\Ab\bS\to \Hom_{\Cat}(I^{op},\Sh_\Ab\bS)$$
which fits into the adjoint pair
\begin{equation}\label{udiwqdwqudwqdpojk}
\int^{I}\dots\otimes F:\Hom_{\Cat}(I^{op},\Sh_\Ab\bS)\Leftrightarrow: \Sh_\Ab\bS:{}^{I}\uHom_{\Sh_\Ab\bS}(F,\dots)\ ,
\end{equation}
where for $G\in \Hom_{\Cat}(I^{op},\Sh_\Ab\bS)$ the symbol
$\int^{I} G\otimes F\in \Sh_\Ab\bS$ denotes the coend of the functor
$I^{op}\times I\to \Sh_\Ab$, $(i,j)\mapsto G(i)\otimes F(j)$; correspondingly
$\int_I$ denotes the end of the appropriate functor.
Indeed we have for all $A\in \Hom_{\Cat}(I^{op},\Sh_\Ab\bS)$ and $B\in \Sh_\Ab\bS$ a natural isomorphism
\begin{eqnarray*}
\Hom_{\Sh_\Ab\bS}(\int^{I}A\otimes F,B)&\cong &\int_I{}^{I^{op}\times I}\uHom_{\Sh_\Ab\bS}(A\otimes F,B)\\&\cong&
\int_I{}^{I^{op}}\uHom_{\Sh_\Ab\bS}(A,{}^I\uHom_{\Sh_\Ab\bS}(F,B))\\
&\cong&\Hom_{\Hom_{\Cat}(I^{op},\Sh_\Ab\bS)}(A,{}^I\uHom_{\Sh_\Ab\bS}(F,B))
\end{eqnarray*}

The functor
${}^I\uHom_{\Sh_\Ab\bS}(F,\dots)$ is therefore also left-exact and admits a right-derived version.

\subsubsection{}

We say that $P\in \Hom_{\Cat}(I,\Sh_\Ab\bS)$ is $I$-free if there exists a collection of flat
sheaves $(U(l))_{l\in I}$, $U(l)\in \Sh_\Ab\bS$, such that $P(j)=\bigoplus_{l\to j} U(l)$, and
$P(i\to j):\bigoplus_{l\to i} U(l)\to \bigoplus_{l\to j} U(l)$ maps
the summand $U(l)$  at $l\to i$ identically to the summand $U(l)$ at $l\to i\to j$.

\begin{lem}\label{hjabddwodqwdq}
If $P\in \Hom_{\Cat}(I,\Sh_\Ab\bS)$  is $I$-free, and if $J\in \Sh_\Ab\bS$ is injective, then
${}^I\uHom_{\Sh_\Ab\bS}(P,J)\in \Hom_{\Cat}(I^{op},\Sh_\Ab\bS)$ is injective.
\end{lem}
\proof
We consider an exact sequence
$$(A^\bullet:0\to A^0\to A^1\to A^2\to 0)$$ in $\Hom_{\Cat}(I^{op},\Sh_\Ab\bS)$. Then we have

$$\Hom_{\Hom_{\Cat}(I^{op},\Sh_\Ab\bS)}(A^\bullet,{}^I\uHom_{\Sh_\Ab\bS}(P,J))\stackrel{(\ref{udiwqdwqudwqdpojk})}{\cong}
\Hom_{\Sh_\Ab\bS}(\int^I A^\bullet\otimes P,J)\ .$$

Exactness in $\Hom_{\Cat}(I^{op},\Sh_\Ab\bS)$ is defined object-wise \cite[Thm. 0.1.3.1]{MR1317816}.
Therefore
$0\to A^0(i)\to A^1(i)\to A^2(i)\to 0$ is an exact complex of sheaves for all $i\in I$.
Since $P(j)$ is flat for all $j\in I$
the complex
$A^\bullet(i)\otimes P(j)$ is exact for all pairs $(i,j)\in I\times I$.
The complex of sheaves
$\int^I A^\bullet\otimes P$ is the complex of push-outs along the exact sequence of diagrams
$$
\begin{CD}
  \bigsqcup_{(i\to j)\in \Mor(I)}A^\bullet(j)\otimes P(i) @>{\id\otimes
  P(i\to j)}>> \bigsqcup_{j\in I} A^\bullet(j)\otimes P(j)\\
@VV{A^\bullet(i\to j)\otimes \id}V\\
    \bigsqcup_{i\in I} A^\bullet(i)\otimes P(i)\ .
\end{CD}
$$
We claim that this $\int^I A^\bullet\otimes P\cong \bigoplus_{i\in I} A^\bullet(i)\otimes U(i)$.
Since $U(i)$ is flat for all $i\in I$ this is a sum of exact sequences and hence exact.
Since $J$ is injective we conclude that
$\Hom_{\Sh_\Ab\bS}(\int^I A^\bullet\otimes P,J)$ is exact.
Therefore
$\Hom_{\Hom_{\Cat}(I^{op},\Sh_\Ab\bS)}(\dots,{}^I\uHom_{\Sh_\Ab\bS}(P,J))$
preserves exactness and hence  ${}^I\uHom_{\Sh_\Ab\bS}(P,J)$ is an injective object in $\Hom_{\Cat}(I^{op},\Sh_\Ab\bS)$.

In order to finish the proof it remains to show the claim. We expand the push-out diagram by inserting the structure of $P$.
$$
\begin{CD}
  \bigoplus_{i\to j}A^\bullet(j)\otimes \bigoplus_{l\to i}U(l) @>{a:=\id\otimes
  (i\to j)}>> \bigoplus_{l\to j} A^\bullet(j)\otimes U(l)\\
@VV{b:=A^\bullet(i\to j)\otimes \id}V @VV{c:=A^\bullet(l\to i)\otimes \id}V\\
    \bigoplus_{l\to i} A^\bullet(i)\otimes U(l) @>{d:=A^\bullet (l\to i)\otimes \id}>> \bigoplus_{l\in I} A^\bullet(l)\otimes U(l)\ .
\end{CD}
$$
It suffices to show the the lower right corner has the universal property.
The lower horizontal map has a split
$s:\bigoplus_{l\in I} A^\bullet(l)\otimes U(l)\to  \bigoplus_{l\to i} A^\bullet(i)\otimes U(l)$
given by the inclusion of summands induced by $l\mapsto l\stackrel{\id}{\to}l$.
Two maps $\psi_l,\psi_r$ from the lower left and upper right corner to some sheaf $V$ which satisfy
$\psi_l\circ b=\psi_r\circ a$
must induce a unique map $\psi:\bigoplus_{l\in I} A^\bullet(l)\otimes U(l)\to V$
such that $\psi\circ d=\psi_l$ and $\psi\circ c=\psi_r$. We have now other choice than to define $\psi:=\psi_l\circ s$, and this map has the required properties as can be checked by an easy diagram chaise.
\hB

\begin{lem}\label{resoexists}
Every sheaf $F\in \Hom_{\Cat}(I,\Sh_\Ab\bS)$ has an $I$-free resolution.
\end{lem}
\proof
It suffices to show that there exists a surjection
$P\to F$ from an $I$-free sheaf.
We start with the surjection
$\Z(F)\to F$ and note that $\Z(F)(i)$ is flat for all $i\in I$. Then we define the $I$-free sheaf $P$ by  $P(i):=\oplus_{l\to i} \Z(F)(l)$, and the surjection $P\to Z(F)$ by
$(l\to i)_*:Z(F)(l)\to \Z(F)(i)$ on the summand of $P(i)$ with index $(l\to i)$.
\hB

\begin{lem}\label{zuiwqdwqidwqd}
We have a natural isomorphism in $D^+(\Sh_\Ab\bS)$
$$R\uHom_{\Sh_\Ab\bS}(\colim(F),H)\cong R\lim R{}^I\uHom_{\Sh_\Ab\bS}(F,H)\ .$$
\end{lem}
\proof
%
%
%
Let now $F\in \Hom_{\Cat}(I,\Sh_\Ab\bS)$. By Lemma \ref{resoexists} we can choose an $I$-free resolution $P^\bullet$ of $F$.
Let $H\in \Sh_\Ab\bS$ and $H\to I^\bullet$ be an injective resolution.
Then we have 
$$
R\uHom_{\Sh_\Ab\bS}(\colim\: F,H)\cong \uHom_{\Sh_\Ab\bS}(\colim \:F,I^\bullet)\ .$$
Since the category $\Sh_\Ab\bS$ is a Grothendieck abelian category \cite[Thm. I.3.2.1]{MR1317816}
the functor $\colim$ is exact \cite[Thm. 0.3.2.1]{MR1317816}.
Therefore
$\colim \:P^\bullet(F)\to \colim \:F$ is a quasi-isomorphism.
It follows that
$$\uHom_{\Sh_\Ab\bS}(\colim \:F,I^\bullet)\cong \uHom_{\Sh_\Ab\bS}(\colim \:P^\bullet(F),I^\bullet) \ .$$

By Lemma \ref{teezwdc} the restriction functor
is exact and therefore commutes with colimits.
We get
\begin{eqnarray*}
\uHom_{\Sh_\Ab\bS}(\colim\: P^\bullet(F),I^\bullet)(A)&\cong& \Hom_{\Sh_\Ab\bS/A}((\colim\: P^\bullet(F))_{|A},I_{|A}^\bullet)\\&\cong &\Hom_{\Sh_\Ab\bS/A}(\colim (P^\bullet(F)_{|A}),I_{|A}^\bullet) 
 \\&\cong& \lim \:{}^I\Hom_{\Sh_\Ab\bS/A}(P^\bullet(F)_{|A},I_{|A}^\bullet) \\&\cong&\lim \:({}^I\uHom_{\Sh_\Ab\bS}(P^\bullet(F),I^\bullet)(A)) \\&\cong&
(\lim \:{}^I\uHom_{\Sh_\Ab\bS}(P^\bullet(F),I^\bullet))(A)
\ ,
\end{eqnarray*}
where in the last step we use that a limit of sheaves is defined object-wise.
 In other words
$$\uHom_{\Sh_\Ab\bS}(\colim \:F,I^\bullet)\cong \uHom_{\Sh_\Ab\bS}(\colim \:P^\bullet(F),I^\bullet)\cong \lim \:{}^I\uHom_{\Sh_\Ab\bS}(P^\bullet(F),I^\bullet) 
 $$

Applying Lemma \ref{hjabddwodqwdq} we see that
${}^I\uHom_{\Sh_\Ab\bS}(P^\bullet(F),I^\bullet)\in\Hom_{\Cat}(I,\Sh_\Ab\bS) $ is injective, and
\begin{eqnarray*}
R\uHom_{\Sh_\Ab\bS}(\colim\:F,H)&\cong& \lim
\:{}^I\uHom_{\Sh_\Ab\bS}(P^\bullet(F),I^\bullet)\\& \cong&
R\lim\:{}^I\uHom_{\Sh_\Ab\bS}(P^\bullet(F),I^\bullet)\\& \cong& R\lim\:
R{}^I\uHom_{\Sh_\Ab\bS}(F,H)\ . 
\end{eqnarray*}

\hB 

\subsubsection{}
Lemma \ref{zuiwqdwqidwqd} implies the existence of s spectral sequence with second page
\begin{equation}\label{uihdiwqdqwdd434}
E_2^{p,q}\cong R^p\lim\:R^q{}^I\uHom_{\Sh_\Ab\bS}(F,H)
\end{equation}
 which converges to
$\Gr R^*\uHom_{\Sh_\Ab\bS}(\colim\:F,H)$.

\subsubsection{}

\begin{lem}\label{uiddqwwqdwd}
Let $I$ be a  category and $U\in \bS$.
Then we have a commutative diagram
$$\xymatrix{
D^+((\Sh_\Ab\bS)^{I^{op}})\ar[d]^{R\Gamma(U,\dots)}\ar[r]^{R\lim}&D^+(\Sh_\Ab\bS)\ar[d]^{R\Gamma(U,\dots)}\\
D^+(\Ab^{I^{op}})\ar[r]^{R\lim}&D^+(\Ab)}\ .$$
\end{lem}
\proof
Since the  limit of a diagram of sheaves is defined object-wise we have
$\Gamma(U,\dots)\circ \lim=\lim\circ \Gamma(U,\dots)$.
We show that the two compositions
$R\Gamma(U,\dots)\circ R\lim$ ,$R\lim\circ R\Gamma(U,\dots)$ are both isomorphic to
$R(\Gamma(U,\dots)\circ \lim)$ by showing that $\lim$ and  $\Gamma(U,\dots)$ preserve injectives.

The left-adjoint of the functor
$\lim:(\Sh_\Ab\bS)^{I^{op}}\to \Sh_\Ab\bS$ is the constant diagram functor
$C_{\dots}: \Sh_\Ab\bS\to (\Sh_\Ab\bS)^{I^{op}}$ which is exact. Therefore $\lim$ preserves injectives.

The left-adjoint of the functor $\Gamma(U,\dots)$ is the functor
$\dots\otimes_\Z\Z(U)$. Since $\Z(U)$ is a sheaf of flat $\Z$-modules (it is torsion-free)
this functor is exact (object-wise and therefore on diagrams).
It follows that $\Gamma(U,\dots)$ preserves injectives. \hB 
 
\subsubsection{}

The general remarks on colimits above are true for all sites with finite
products (because of the use of Lemma \ref{teezwdc}). In particular we can
replace $\bS$ by the site $\bS_{lc-acyc}$ of locally acyclic locally compact
spaces.

\begin{lem}\label{uiqwdwqdwdwqd}
Let $F\in (\Sh_\Ab\bS_{lc-acyc})^{I^{op}}$.
If for all acyclic $U\in \bS_{lc-acyc}$ the canonical map
$\lim \:F(U)\to R\lim (F(U))$ is a quasi-isomorphism 
of complexes of abelian groups, then
$\lim \:F\to R\lim F$ is a quasi-isomorphism.
 \end{lem}
\proof
We choose an injective resolution $F\to I^\bullet$ in $(\Sh_\Ab\bS_{lc-acyc})^{I^{op}}$.
As in the proof of Lemma \ref{uiddqwwqdwd} for each $U\in \bS_{lc-acyc}$ the complex 
$ I^\bullet(U)$ is injective. If we assume that $U$ is acyclic, the map
$F(U)\to I^\bullet(U)$ is a quasi-isomorphism in $(\Ab)^{I^{op}}$, and $R\lim\:(F(U))\cong\lim(I^\bullet(U))$.
By assumption
$$(\lim F)(U)\cong \lim (F(U))\to \lim (I^\bullet(U))\cong (\lim I^\bullet)(U)\cong (R\lim F)(U)$$ is a quasi-isomorphism for all acyclic $U\in \bS_{lc-acyc}$.
An arbitrary  object $A\in \bS_{lc-acyc}$ can be covered by acyclic open subsets.
Therefore
$\lim F\to \lim I^\bullet$ is an quasi-isomorphism locally on $A$ for each $A\in \bS_{lc-acyc}$.
Hence 
$\lim F\to R\lim F$ is an isomorphism in $D^+(\Sh_\Ab\bS_{lc-acyc})$ . \hB

\subsubsection{}

Let $D$ be a discrete group.
Then we let $I$ be the category of all finitely generated subgroups of $D$.
This category is filtered.
Let $F:I\to \Ab$ be the "identity" functor.
Then we have a natural isomorphism
$$D\cong \colim \:F\ .$$
By Theorem \ref{wjebdqwdwdqdwqd} we know that a finitely generated group $G$ is admissible,
i.e.
$$R^q\uHom_{\Sh_\Ab\bS}(\uG,\uT)\cong 0\ ,\quad q=1,2\ .$$
Note that by Lemma \ref{ehewfudwqdwqdwq} we have ${\colim \uF}=\underline{D}$. 
Using the spectral sequence 
(\ref{uihdiwqdqwdd434}) we get for $p=1,2$ that
\begin{equation}
R^p\uHom_{\Sh_\Ab\bS}(\uD,\uT)\cong R^p\lim \:\uHom_{\Sh_\Ab\bS}(\uF,\uT)\ .
\end{equation}
Let $g\colon \bS_{lc-acyc}\to \bS$ be the inclusion.
Since $\T$ and $D$ belong to $\bS_{lc-acyc}$ using Lemma \ref{uiwddiqwdwqdddqd2121} we also have
\begin{eqnarray}\label{uwhihdiwqduwqdwqodop}
g^*R^p\uHom_{\Sh_\Ab\bS}(\uD,\uT)&\cong &R^p\uHom_{\Sh_\Ab\bS_{lc-acyc}}(\uD,\uT)\\&\cong& 
R^p\lim \:\uHom_{\Sh_\Ab\bS_{lc-acyc}}(\uF,\uT)\ .\nonumber
\end{eqnarray}

\subsubsection{}

We now must study the $R^p\lim$-term.
First we make the index category $I$ slightly smaller.
Let $\bar D\subseteq D\otimes_\Z\Q=:D_\Q$ be the image of the natural map
$i:D\to D_\Q$, $d\mapsto d\otimes 1$.
We observe that $\bar D$ generates $D_\Q$.
We choose a  basis $B\subseteq \bar D$ of the $\Q$-vector space $D_\Q$ and let
$\Z B\subseteq \bar D$ be the $\Z$-lattice generated by $B$.
For a subgroup $A\subseteq D$ let $\bar A:=i(A)\subset D_\Q$.
We consider the partially ordered (by inclusion) set
$$J:=\{A\subseteq D|\mbox{$A$ finitely generated and $\Q\bar A\cap \Z B\subset
  \bar A$ and $[\bar A:\bar A\cap \Z B]<\infty $}\}\ .$$
Here $[G:H]$ denotes the index of a subgroup $H$ in a group $G$.
We still let $A$ denote the "identity" functor
$A:J\to \Ab$.
\begin{lem}
We have $D\cong \colim \: A$.
\end{lem}
\proof
 
It suffices to show that $J\subseteq I$ is cofinal.
Let $A\subseteq D$ be a finitely generated subgroup. Choose a finite subset
$\bar b\subset B$ such that the $\Q$-vectorspace $\Q\bar b\subset D_\Q$
generated by $\bar b$ contains $\bar A$ and choose representatives $b$ in $D$
of the elements of $\bar b$. They generate a finitely generated group
$U\subset D$ such that $\bar U=i(U)=\Z\bar b$. 
We now consider the group $G:=<U,A>$. This group is still finitely generated.
Similarly, since
$\Q\bar G=\Q\bar U+\Q\bar A=\Q\bar b$
we have
$$\Q\bar G\cap \Z B = \Q\bar b\cap \Z B=\bar b\subseteq \bar U\subseteq \bar
G\ .$$
Moreover, since $\Q\bar G=\Q\bar b=\Q(\bar G\cap \Z B)$, $[\bar G:\bar G\cap
\Z B]<\infty$,
i.e. $G\subset J$. On the other hand, by construction $A\subseteq G$.
\hB 

On $J$ we define the grading
$w:J\to \nat_0$ by
$$w(A):=|A_{tors}|+ \rk\bar A + [\bar A:\bar A\cap \Z B]\qquad\text{for $A\in J$}\ .$$
\begin{lem}\label{uiqdiqwdwdqwdqd434}
The category $J$ together with the grading $w:J\to \nat_0$
is a direct category in the sense of \cite[Def. 5.1.1]{MR1650134}
\end{lem}
\proof
We must show that
$A\subset G$ implies $w(A)\le w(G)$, and that
$A\subsetneqq G$ implies $w(A)  <w(G)$.
First of all we have
$A_{tors}\subset G_{tors}$ and therefore
$|A_{tors}|\le |G_{tors}|$.
Moreover we have
$\bar A\subseteq \bar G$, hence
$\rk\bar A \le \rk \bar G$.
Finally, we claim that the canonical map
$$\bar A/(\bar a\cap \Z B)\to \bar G/(\bar G\cap \Z B)$$
is injective. In fact, 
we have
$\bar A\cap (\bar G\cap \Z B)=(\bar A\cap \bar G)\cap \Z B=\bar A\cap \Z B\ .$
It follows that
$$|\bar A:\bar A\cap\Z B|\le |\bar G:\bar G\cap\Z B| \ .$$

Let now $A\subseteq G$ and $w(A)=w(G)$.
We want to see that this implies $A=G$. First note that
the inclusion of finite groups $A_{tors}\to G_{tors}$ is an isomorphism since both groups have the same
number of elements.
It remains to see that $\bar A\to \bar G$ is an isomorphism.
We have $\rk \bar A=\rk \bar G$. Therefore $\bar A\cap \Z B=\Q\bar A\cap \Z
B=\Q \bar G\cap B\Z =\bar G\cap B\Z$.  Now the equality 
$[\bar A:\bar A\cap \Z B]=[\bar G:\bar G\cap \Z B]$ implies that
$\bar A=\bar G$. \hB

\subsubsection{}

The category $C(\Ab)$ has the projective model structure whose
weak equivalences are the quasi-isomorphisms, and whose fibrations are level-wise surjections.
By Lemma \ref{uiqdiqwdwdqwdqd434} the category $J^{op}$ with the grading $w$ is an inverse category.
On $C(\Ab)^{J^{op}}$ we consider the inverse model structure whose weak equivalences are
the quasi-isomorphisms, and whose cofibrations are the object-wise ones.
The fibrations are characterized by a matching space condition which we will explain in the following.

For $j\in J$ let $J_j\subseteq J$ be the category with non-identity maps all
non-identity maps in $J$ with codomain $j$.
Furthermore consider the functor
$$M_j:C(\Ab)^{J^{op}}\xrightarrow{restriction} C(\Ab)^{(J_j)^{op}}\xrightarrow{\lim}C(\Ab)$$
There is a natural morphism $F(j)\to M_jF$ for all $F\in C(\Ab)^{J^{op}}$.
The matching space condition asserts that a map
$F\to G$ in $C(\Ab)^{J^{op}}$ is a fibration if and only if
$$F(j)\to G(j)\times_{M_jG}M_j F$$ is a fibration in $C(\Ab)$, i.e. a level-wise surjection, for all $j\in J$.
In particular, $F$ is fibrant if
$F(j)\to M_jF$ is a level-wise surjection for all $j\in J$.

For $X\in C(\Ab)^{J^{op}}$ the fibrant replacement $X\to RX$ induces the morphism 
$$\lim \:X\to \lim\:RX\cong R\lim\:X\ .$$ 

Let now again $A$ denote the identity functor $A\colon J\to \Ab$ for the
discrete abelian group $D$.

\begin{lem}\label{jkdjkwdqwdwdwqdqwdqd}
\begin{enumerate}
\item 
For all $U\in \bS$ which are acyclic the diagram of abelian groups
$\Hom_{\Sh_\Ab\bS}(\uA,\uT)(U)\in C(\Ab)^{J^{op}}$ is fibrant.
\item $$\lim \:\Hom_{\Sh_\Ab\bS_{lc-acyc}}(\uA,\uT)\to
  R\lim\:\Hom_{\Sh_\Ab\bS_{lc-acyc}}(\uA,\uT) \in C(\Ab)^{J^{op}}$$ is a quasi-isomorphism.
\end{enumerate}
\end{lem}
\proof
The Assertion $1.$ for locally compact acyclic $U$ verifies the assumption
of Lemma \ref{uiqwdwqdwdwqd}. Hence 2. follows from 1.
We now concentrate on 1.
We must show that
\begin{equation}\label{uidqwdqwdwqdw}
\uHom_{\Sh_\Ab\bS}(\uA(j),\uT)(U)\to M_j\uHom_{\Sh_\Ab\bS}(\uA,\uT)(U)
\end{equation}
is surjective.
We have
$$M_j\uHom_{\Sh_\Ab\bS}(\uA,\uT)\cong \lim \ \uHom_{\Sh_\Ab\bS}(\uA_{|J_j},\uT)\cong
\uHom_{\Sh_\Ab\bS}(\colim\  \uA_{|J_j},\uT)\ .$$
The map (\ref{uidqwdqwdwqdw}) is induced by
the map
\begin{equation}\label{edkjqdqwddlödwqd0922}
\colim \uA_{|J_j}\to \uF(j)\ .
\end{equation}

By Lemma \ref{ehewfudwqdwqdwq} we have
$\colim \ \uA_{|J_j}\cong \underline{\colim \ A_{|J_j}}$.
The map $\colim \ A_{|J_j}\to A(j)$  is of course an injection.
 
We finish the argument by the following observation. Let $H\to G$ be an injective map of finitely generated  groups (we apply this with $H:=\colim A_{|J_j}$ and $G:=F(j)$). Then for $U\in \bS$
$$\uHom_{\Sh_\Ab\bS}(\uG,\uT)(U) \to \uHom_{\Sh_\Ab\bS}(\uH,\uT)(U)$$
is a surjection.
In fact we have
\begin{eqnarray*}
\uHom_{\Sh_\Ab\bS}(\uG,\uT)(U)&\stackrel{Lemma \:\ref{eudiwe33}}{\cong}&
\underline{\Hom_{\top-\Ab}(G,\T)}(U)\\
&\cong \widehat \uG(U)
\end{eqnarray*}
(and a similar equation for $H$).
Since $H\to G$ is injective, its Pontrjagin dual
 $\widehat G\to \widehat H$ is surjective. Because of the classification of
 discrete finitely generated abelian groups, $\hat G$ and $\hat H$ both are
 homeomorphic to finite unions of finite dimensional tori. Because $U$ is
 acyclic, every map $U\to \hat H$ lifts to $\hat G$.
Therefore
$\widehat \uG(U)\to \widehat \uH(U)$ is surjective.
\hB

\begin{kor}\label{wuqdiwqdwqdqwdqd}
We have $R^p\lim\  \uHom_{\Sh_\Ab\bS_{lc-acyc}}(\uA,\uT)\cong 0$ for all $p\ge 1$.
\end{kor}
 
\begin{theorem}\label{uwdiqwdqwdwqdwqd}
Every discrete abelian group is admissible on $\bS_{lc-acyc}$.
\end{theorem}
\proof
This follows from (\ref{uwhihdiwqduwqdwqodop}) and Corollary 
\ref{wuqdiwqdwqdqwdqd}.
\hB

\subsubsection{}

We now present an example which shows that not every discrete group
$D$ is admissible on $\bS$ or $\bS_{lc}$, using Corollary 
\ref{aposteuhefkjwfe} to be established later.

\begin{lem}\label{zwdwqidqidwqpdopop}
Let $I$ be an infinite set, $1\not=n\in \nat$ and $D:=\oplus_{I} \Z/n\Z$.
Then $D$ is not admissible on $\bS$ or $\bS_{lc}$.
\end{lem}
\proof
We consider the sequence
\begin{equation}\label{uquwdiuqwwqdmm8}
0\to \Z/n\Z\to \T\stackrel{n}{\to} \T\to 0
\end{equation}
which has no global section,
and the product of $I$ copies of it
\begin{equation}\label{hjqjdzuqdzuwquzdw}
0\to \prod_I \Z/n\Z\to \prod_{I}\T\to \prod_I \T\to 0\ .
\end{equation}
This sequence of compact abelian groups does not have local sections. 
In fact, an open subset of $\prod_I \T$ always contains a subset of the form
$U=\prod_{I^\prime}\T \times V$, where $I^\prime \subset I$ is cofinite and
$V\subset \prod_{I\setminus I^\prime}\T$ is open.
A section $s:U\to \prod_I \T$ would consist of sections of the sequence
(\ref{uquwdiuqwwqdmm8}) at the entries labeled with  $I^\prime$.

By Lemma \ref{dezuqwideqwd} the sequence of sheaves associated to (\ref{hjqjdzuqdzuwquzdw}) is not exact. 
In view of Corollary \ref{aposteuhefkjwfe} below, the group
$\widehat{\prod_I \Z/n\Z} \cong \oplus_I \Z/n\Z$
is not admissible on $\bS$ or $\bS_{lc}$.
\hB

\subsection{Admissibility of the groups $\R^n$ and $\T^n$}\label{ewdwedewdii7}

\subsubsection{}\label{gqwefbkeje}

\begin{theorem}\label{tqwzdqwmmq6672}
The group $\T$ is admissible. 
\end{theorem}
\proof
Since $\T$ is compact, Assumption $2.$ of Lemma \ref{udiuddwqwqd} follows from Proposition \ref{hdfjsdfiuiuwef}.
It remains to show the first assumption  of Lemma \ref{udiuddwqwqd}.
As we will see this follows from the following result.
 \begin{lem}\label{iewofw}
We have $\uExt^i_{\Sh_\Ab\bS}(\uR,\uZ)=0$ for $i=1,2,3$.
\end{lem}
Let us assume this lemma for the moment.
We apply $\uExt^*_{\Sh_\Ab\bS}(\dots,\uZ)$ to the exact sequence
$$0\to \uZ\to \uR\to \uT\to 0$$ 
and consider the following part of the resulting long exact sequence for $i=2,3$:
$$\uExt^{i-1}_{\Sh_\Ab\bS}(\uZ,\uZ)\to \uExt^i_{\Sh_\Ab\bS}(\uT,\uZ)\to \uExt^i_{\Sh_\Ab\bS}(\uR,\uZ)\to \uExt^i_{\Sh_\Ab\bS}(\uZ,\uZ)\ .$$
The outer terms vanish  since $R\uHom_{\Sh_\Ab\bS}(\uZ,F)\cong F$ for all $F\in \Sh_\Ab\bS$. Therefore by Lemma \ref{iewofw}
$$\uExt^i_{\Sh_\Ab\bS}(\uT,\uZ)\cong  \uExt^i_{\Sh_\Ab\bS}(\uR,\uZ)\cong 0\ ,$$
and this is Assumption 1. of \ref{udiuddwqwqd}.

It remains to prove Lemma \ref{iewofw}.

\subsubsection{Proof of Lemma \ref{iewofw}.}\hfill
 
We choose an injective resolution $\uZ\to I^\bullet$. The sheaf $\uR$ gives rise to the homological complex 
$U^\bullet$ introduced in \ref{iudhuiqwdqwdwqd}. We get
 the double complex
$\uHom_{\Sh_\Ab\bS}(U^\bullet,I^\bullet)$ as in
 \ref{uiwdiqwdwqdq}. 
As before we discuss the associated two spectral sequences which compute $\uExt^*_{\Sh_{\Z[\uR]-\Mod}\bS}(\uZ,\uZ)$.
 
\subsubsection{}

 We first take the cohomology in the $I^\bullet$-, and then in the $U^\bullet$-direction.
 In view of (\ref{uieiefefewfw}), the first page of the resulting spectral sequence  is given  by 
$$E^{p,q}_1\cong \uExt_{\Sh_\Ab\bS}^q(\Z(\cF\uR^p),\uZ)\ ,$$
where $\cF\uR$ denotes the underlying sheaf of sets of $\uR$. 
By Corollary  \ref{sofdmadd}, 1. we have $E_1^{p,q}\cong 0$ for $q\ge 1$.

We now consider the case $q=0$. For $A\in \bS$ we have
 $E_1^{p,0}(A)\cong \Gamma(A\times \R^p;\uZ)\cong \Gamma(A;\uZ)$ for all $p$, i.e.
$E_1^{p,0}\cong \uZ$.
 We can easily calculate the cohomology of the complex
 $(E_1^{*,0},d_1)$ which is isomorphic to   
 $$0\to \uZ\stackrel{0}{\to} \uZ\stackrel{\id}{\to}\uZ\stackrel{0}{\to} \uZ\stackrel{\id}{\to}\uZ\to \dots \ .$$
 We get 
 $$H^i(E^{*,0}_1,d_1)\cong \left\{\begin{array}{cc} \uZ&i=0\\ 0&i\ge 1\end{array}\right.\ .$$
 The spectral sequence  $(E_r,d_r)$ thus degenerates at the second term, and
\begin{equation}\label{fueieded77}
H^i\uHom_{\Sh_\Ab\bS}(U^\bullet,I^\bullet)\cong  \left\{\begin{array}{cc} \uZ&i=0\\ 0&i\ge 1\end{array}\right.\ .\end{equation}

\subsubsection{}
We now consider the second spectral sequence $(F^{p,q}_r,d_r)$ associated to the double complex $\uHom_{\Sh_\Ab\bS}(U^\bullet,I^\bullet)$ which 
takes cohomology first in the $U^\bullet$ and then in the $I^\bullet$-direction.  Its
second page is given by
$$F_2^{p,q}\cong \uExt^p_{\Sh_\Ab\bS}(H^q U^\bullet,\uZ)\ .
$$ 
Since $\uR$ is torsion-free we can apply \ref{diqwodwqdwqdq} and get  
$$F^{p,q}_2\cong \uExt^p_{\Sh_\Ab\bS}((\Lambda^q_\Z\uR)^\sharp,\uZ)\ .$$
In particular we have $\Lambda^0_\Z\uR\cong \uZ$ and thus
$$F_2^{p,0}\cong \uExt^p_{\Sh_\Ab\bS}(\uZ,\uZ)\cong 0$$ for $p\ge 1$
(recall that $\uExt^p_{\Sh_\Ab\bS}(\uZ,H)\cong 0$ for every $H\in \Sh_\Ab\bS$ and $p\ge 1$).
Furthermore, since $\Lambda^1_\Z\uR\cong \uR$ we have
$$F_2^{0,1}\cong \uHom_{\Sh_\Ab\bS}(\uR,\uZ)\stackrel{Lemma \: \ref{eudiwe33}}{\cong} \underline{\Hom_{\top-\Ab}(\R,\Z)}\cong0$$
since $\Hom_{\top-\Ab}(\R,\Z)\cong 0$.

Here is a picture of the relevant part of the second page. 
$$
\begin{array}{|c|c|c|c|c|c|c|}\hline3&\uHom_{\Sh_\Ab\bS}((\Lambda_\Z^3\uR)^\sharp,\uZ)&&&&&\\\hline
2&\uHom_{\Sh_\Ab\bS}((\Lambda_\Z^2\uR)^\sharp,\uZ)&\uExt^1_{\Sh_\Ab\bS}((\Lambda_\Z^2\uR)^\sharp,\uZ)&&&&\\\hline
1&0&\uExt^1_{\Sh_\Ab\bS}(\uR,\uZ)&\uExt^2_{\Sh_\Ab\bS}(\uR,\uZ)&\uExt^3_{\Sh_\Ab\bS}(\uR,\uZ)&&\\\hline
0&\uZ&0&0&0&0&0\\\hline
&0&1&2&3&4&5\\\hline
\end{array}\ .$$

\subsubsection{}

Let $V$ be an abelian group.
Recall the definition of 
  $\Lambda^*V$ from  \ref{lazdzu32d898}.
If $V$ has the structure of a $\Q$-vector space, then $T^{\ge 1}_\Z V$ has the structure of a graded $\Q$-vector space, and $I\subseteq T^{\ge 1}_\Z V$ is a graded $\Q$-vector subspace. Therefore $\Lambda_\Z^{\ge 1}V$ has the structure of a graded $\Q$-vector space, too.

\subsubsection{}

We claim that $$F_2^{0,i}\cong \uHom_{\Sh_\Ab\bS}((\Lambda_\Z^i\uR)^\sharp,\uZ)\cong 0$$ for $i\ge 1$.  Note that $$\uHom_{\Sh_\Ab\bS}((\Lambda_\Z^i\uR)^\sharp,\uZ)\cong
\uHom_{\Pr_\Ab\bS}(\Lambda_\Z^i\uR,\uZ)\ .$$
Let $A\in \bS$. An element $$\lambda\in \uHom_{\Pr_\Ab\bS}(\Lambda_\Z^i\uR,\uZ)(A)\cong \Hom_{\Pr_\Ab\bS/A}(\Lambda_\Z^i\uR_{|A},\uZ_{|A})$$ induces a homomorphism of groups
$\lambda^W: \Lambda^i_\Z\uR(W) \to \uZ(W)$ for every $(W\to A)\in \bS/A$.
Since $\Lambda^i_\Z\uR(W) $ is a $\Q$-vector space and $\uZ(W)$ does not contain  divisible elements we see that $\lambda^W=0$.
This proves the claim.

The claim implies that $F_2^{2,1}\cong \uExt^2_{\Sh_\Ab\bS}(\uR,\uZ)$ survives to the limit of the spectral sequence. Because of  (\ref{fueieded77}) it must vanish. This proves the Lemma \ref{iewofw}  in the case $i=2$.

The term $F_2^{1,1}\cong \uExt^1_{\Sh_\Ab\bS}(\uR,\uZ)$ also survives to the limit and therefore also vanishes  because of (\ref{fueieded77}). This proves Lemma \ref{iewofw}  in the case $i=1$.

\subsubsection{}

Finally, since $F_3^{0,3}\cong 0 $, we see that $d_2 :F^{1,2}_2\to F_2^{3,1}$ must be an isomorphism, i.e.
$$d_2:\uExt^1_{\Sh_\Ab\bS}((\Lambda_\Z^2\uR)^\sharp,\uZ)\stackrel{\sim}{\to} \uExt^3_{\Sh_\Ab\bS}(\uR,\uZ)\ .$$
We will finish the proof of Lemma 
 \ref{iewofw} for $i=3$ by showing that $d_2=0$.

\subsubsection{}

We consider the natural action of $\Z_{mult}$ on $\R$ and hence on $\uR$.
This turns $\uR$ into a sheaf of $\Z_{mult}$-modules of weight $1$ (see \ref{dkwqwjdkqkdjqwdw}).
It follows that $\uHom_{\Sh_\Ab\bS}(\uR,I^\bullet)$ is a complex of sheaves of $\Z_{mult}$-modules of weight $1$.
Finally we see that $\uExt^i_{\Sh_\Ab\bS}(\uR,\uZ)$ are sheaves of $\Z_{mult}$-modules of weight $1$ for $i\ge 0$.

Now observe (see \ref{fbefewfwef}) that $\Lambda_\Z^2\uR$ is a presheaf of $\Z_{mult}$-modules of weight $2$.
Hence $(\Lambda_\Z^2\uR)^\sharp$ and thus 
$\uExt^1_{\Sh_\Ab\bS}((\Lambda_\Z^2\uR)^\sharp,\uZ)$ are sheaves of $\Z_{mult}$-modules of weight $2$.

Since $\uR\to U^\bullet(\uR)=:U^\bullet$ is a functor we get an action of $\Z_{mult}$ on $U^\bullet$ and hence on the double complex $\uHom_{\Sh_\Ab\bS}(U^\bullet,I^\bullet)$. This implies that the differentials of the associated spectral  sequences commute with the $\Z_{mult}$-actions (see \ref{jhwefiewfewfe} for more details). This in particular applies to $d_2$.
The equality $d_2=0$ now directly follows from the following Lemma.
\begin{lem}
Let $V,W\in \Sh_\Ab\bS$ be sheaves of $\Z_{mult}$-modules of weights $k\not= l$.
Assume that $W$ has the structure of a sheaf of $\Q$-vector spaces.
If  $d\in \Hom_{\Sh_\Ab\bS}(V,W)$ is $\Z_{mult}$-equivariant, then $d=0$.
\end{lem}
\proof
Let $\alpha:\Z_{mult}\to \End_{\Sh_\Ab\bS}(V)$ and $\beta :\Z_{mult}\to\End_{\Sh_\Ab\bS}(W)$ denote the actions.
Then we have
$d\circ \alpha(q) - \beta(q)\circ d=0$
for all $q\in \Z_{mult}$. We consider $q:=2$.
Then $(2^k-2^l) \circ d=0$.
Since $W$ is a sheaf of $\Q$-vector spaces and $(2^k-2^l)\not=0$   this implies that $d=0$.
\hB

\subsubsection{}

\begin{theorem}\label{udiqwdwqdwd}
The group $\R$ is admissible.
\end{theorem}
\proof
The outer terms of the exact sequence
 $$0\to \uZ\to \uR\to \uT\to 0$$
are admissible by  Theorem \ref{tqwzdqwmmq6672} 
and Theorem \ref{wjebdqwdwdqdwqd}. We now apply  Lemma \ref{wuihfewfefwfw}, stability of admissibility under extensions.
\hB

\subsection{Profinite groups}

\begin{ddd}
A topological group $G$ is called profinite if there exists a small left-filtered\footnote{i.e. for every pair $i,k\in I$ there exists $j\in I$ with $j\le i$ and $j\le k$} poset\footnote{A poset is considered here as a category.}  $I$ and a system
$F\in \Ab^I$ such that
\begin{enumerate}
\item for all $i\in I$ the group $F(i)$ is finite,
\item for all $i\le  j$ the morphism $F(i)\to F(j)$ is a surjection,
\item there exists an isomorphism $G\cong \lim_{i\in I}F(i)$ as topological groups.
\end{enumerate}
\end{ddd}
For the last statement we consider finite abelian groups as topological groups with the discrete topology.
Note that the homomorphisms $G\to F(i)$ are surjective for all $i\in I$. We call the system $F\in \Ab^I$ an inverse system.
\begin{lem}[\cite{MR0156915}]
The following assertions on a topological abelian  group $G$ are equivalent:
\begin{enumerate}
\item  $G$ is compact and totally disconnected.
\item Every neighbourhood $U\subseteq G$ of the identity contains a compact subgroup $K$ such that
$G/K$ is a finite abelian group.
\item $G$ is profinite.
\end{enumerate}
 \end{lem}

\subsubsection{}

\begin{lem}\label{hjsahduidudiqdwq}
Let $G$ be a profinite abelian group and $n\in \Z$.
 We define the groups
$K,Q$ as the kernel and cokernel of the multiplication map by $n$, i.e. by the exact sequence
\begin{equation}\label{jhsabcasciscaca}
0\to K\to G\stackrel{n}{\to} G\to Q\to 0\ .
\end{equation}
Then $K$ and $Q$ are again profinite.
\end{lem}
\proof
We write $G:=\lim_{j\in J} G_j$ for an inverse  system $(G_j)_{j\in J}$ of finite abelian groups.
We define the system of finite subgroups
$(K_j)_{j\in J}$ by the sequences
$$0\to K_j\to G_j\stackrel{n}{\to} G_j\ .$$
Since taking kernels commutes with limits the natural projections
$K\to K_j$, $j\in J$ represent $K$ as the limit $K\cong \lim_{j\in J}K_j$.

Since cokernels do not commute with limits we will use a different argument for $Q$.
Since $G$ is compact and the multiplication by $n$ is continuous, $nG\subseteq G$ is a closed subgroup. Therefore the group theoretic  quotient $Q$ is a topological group in the quotient topology.

A quotient of a profinite group is again profinite \cite{MR1646190}, Exercise E.1.13.
Here is a solution of this exercise, using the following general structural result about 
compact abelian groups.
\begin{lem}[\cite{MR0156915}]
If $H$ is a compact abelian group, then for every open neighbourhood $1\in U\subseteq H$
there exists a compact subgroup $C\subset U$ such that $H/C\cong \T^a\times F$ for some $a\in \nat_0$ and a finite group $F$.
\end{lem}
Since $G$ is compact, its quotient $Q$ is compact, too.
This Lemma in particular implies that $Q$ is the limit of the system of these quotients $Q/C$.
In our case, since $G$ is profinite, it can not have $T^a$ as a quotient, i.e  given a surjection
$$G\to Q\to \T^a\times F$$
we conclude that $a=0$. Hence we can write $Q$ as a limit of an inverse  system of finite quotients.
This implies that $Q$ is profinite.
\hB

\subsubsection{}

Let $$0\to K\to G\to H\to 0$$ be an exact sequence of profinite groups, where $K\to G$ is the inclusion of a closed subgroup. 
\begin{lem}\label{jhbuisaicascsa444}
The sequence of sheaves
$$0\to \uK\to \uG\to \uH\to 0$$ is exact.
\end{lem}
\proof
By \cite[Proposition 1]{MR1867431} every surjection between profinite groups has a section.
Hence we can apply Lemma \ref{dezuqwideqwd}. \hB
 
In this result one can in fact drop the assumptions that $K$ and $G$ are profinite.
In our basic example the group $K$ is the connected component $G_0\subseteq G$ of the identity of $G$.
 
 \begin{lem}\label{jhbuisaicascsa1}
Let $$0\to K\to G\to H\to 0$$
be an exact sequence of compact abelian groups such that
$H$ is profinite.
Then the sequence of sheaves
$$0\to \uK\to \uG\to \uH\to 0$$ is exact.
\end{lem}
\proof  
We can apply \cite[Thm. 10.35]{MR1646190} which says that the projection
$G\to H$ has a global section. Thus the sequence
of sheaves is exact by  \ref{dezuqwideqwd} (even as sequence of presheaves).\hB

\subsubsection{}

In order to show that certain discrete groups are not admissible we use the following Lemma.
\begin{lem}\label{jhbuisaicascsa} Let $$0\to K\to G\to H\to 0$$
be an exact sequence of compact abelian groups. If the discrete abelian group $\widehat K$ is admissible,
then sequence of sheaves
$$0\to \uK\to \uG\to \uH\to 0$$ is exact.
\end{lem}
\proof  
If $U$ is a compact group, then its Pontrjagin dual $\widehat U:=\Hom_{\top-\Ab}(U,\T)$ is a discrete group.
Pontrjagin duality preserves exact sequences. Therefore we have the exact sequence of discrete groups
$$0\to \widehat H\to \widehat G\to \widehat K\to 0\ .$$
The surjective map of discrete sets $\widehat G\to \widehat K$ has of course a section.
Therefore by Lemma \ref{dezuqwideqwd} the sequence of sheaves of abelian groups 
$$0\to \underline{\widehat H}\to \underline{\widehat G}\to \underline{\widehat K}\to 0$$
is exact. We apply $\uHom_{\Sh_\Ab\bS}(\dots,\uT)$ and get the exact sequence
$$0\to \uHom_{\Sh_\Ab\bS}( \underline{\widehat K},\uT)\to\uHom_{\Sh_\Ab\bS}( \underline{\widehat G},\uT)\to \uHom_{\Sh_\Ab\bS}( \underline{\widehat H},\uT)\to\uExt^1_{\Sh_\Ab\bS}( \underline{\widehat K},\uT)\to\dots\ .$$
By Lemma \ref{eudiwe33} we have $\uHom_{\Sh_\Ab\bS}(\underline{\widehat G},\uT)\cong \uG$ etc. Therefore this sequence translates to
\begin{equation}\label{hbewdwdiwqdoiwqdwq}
0\to \uK\to\uG\to \uH\to\uExt^1_{\Sh_\Ab\bS}( \underline{\widehat K},\uT)\to\dots\ .
\end{equation}
By our assumption $\uExt^1_{\Sh_\Ab\bS}( \underline{\widehat K},\uT)=0$ so that
$\uG\to \uH$ is surjective. \hB

\subsubsection{}

The same argument does apply for the site $\bS_{lc}$.
Evaluating the surjection $\uG\to \uH$ on $H$ we conclude the following fact.
\begin{kor}\label{aposteuhefkjwfe}
If $K\subset G$ is a closed subgroup of a compact abelian group such that $\widehat K$ is admissible or admissible over $\bS_{lc}$, then the projection
$G\to G/K$ has local sections.
\end{kor}

\subsubsection{}

If $G$ is an abelian group, then let $G_{tors}\subseteq G$ denote the subgroup of torsion elements.
We call $G$ a torsion group, if $G_{tors}=G$. If $G_{tors}\cong 0$, then we say that $G$ is torsion-free.
If $G$ is a torsion group and $H$ is torsion-free, then
$\Hom_\Ab(H,G)\cong 0$.

A presheaf $F\in \Pr_\bS\bS$ is called a presheaf of torsion groups if $F(A)$ is a torsion group for every $A\in \bS$.
The notion of a torsion sheaf is more complicated.
\begin{ddd}
A sheaf $F\in \Sh_\Ab\bS$ is called a torsion sheaf if
for each $A\in \bS$ and $f\in F(A)$ there exists an open covering
$(U_i)_{i\in I}$ of $A$ such that $f_{|U_i}\in F(U_i)_{tors}$. 
\end{ddd}
The following Lemma provides equivalent characterizations of torsion sheaves.
\begin{lem}[\cite{MR1317816}, (9.1)]
Consider a sheaf $F\in \Sh_\Ab\bS$.
The following assertions are equivalent.
\begin{enumerate}
\item $F$ is a torsion sheaf.
\item $F$ is the sheafification of a presheaf of torsion groups.
\item The canonical morphism $\colim_{n\in \nat}( {}_nF)\to F$ is an isomorphism, where $({}_nF)_{n\in \nat}$ is the direct system ${}_nF:=\ker(F\stackrel{n!}{\to} F)$.
\end{enumerate}
\end{lem}
Note that a subsheaf or a quotient of a torsion sheaf is again a torsion sheaf.

\begin{lem}\label{qwuiwdiudwd}
If $H$ is a discrete torsion group, then $\uH$ is a torsion sheaf.
\end{lem}
\proof
We write $H=\colim_{n\in \nat} ({}_n H)$, where ${}_nH:=\ker(H\stackrel{n!}{\to} H)$.
By Lemma  \ref{ehewfudwqdwqdwq} we have
$\uH=\colim_{n\in \nat}\underline{{}_n H}$.
Now $\colim_{n\in \nat}\underline{{}_n H}$ is the sheafification of the presheaf ${}^p\colim_{n\in \nat}\underline{{}_n H}$ of torsion groups and therefore a torsion sheaf.
\hB
In general $\uH$ is not a presheaf of torsion groups. Consider e.g. $H:=\oplus_{n\in \nat}\Z/n\Z$.
Then the element $\id\in \uH(H)$ is not torsion.

\subsubsection{}

\begin{ddd}
A sheaf $F\in \Sh_\Ab\bS$ is called torsion-free if the group
$F(A)$ is torsion-free for all $A\in \bS$. 
\end{ddd}
A sheaf  $F$ is torsion-free if and only if for all  $n\in \nat$ the map
$F\stackrel{n}{\to} F$ is injective. 
It suffices to test this for all primes.

\begin{lem}\label{uihfuieddqwd}
If $F\in \Sh_\Ab\bS$ is a torsion sheaf and $E\in \Sh_\Ab\bS$ is torsion-free, then
$\uHom_{\Sh_\Ab\bS}(F,E)\cong 0$.
\end{lem}
\proof
We have
$$\uHom_{\Sh_\Ab\bS}(F,E)\cong \uHom_{\Sh_\Ab\bS}(\colim_{n\in \nat} ({}_n F),E)\cong \lim_{n\in \nat} 
\uHom_{\Sh_\Ab\bS}({}_n F,E)\ .$$
On the one hand, via the first entry multiplication by $n!$ induces on $\uHom_{\Sh_\Ab\bS}({}_n F,E)$ the trivial map.
On the other hand it induces an injection since $E$ is torsion-free. Therefore $\uHom_{\Sh_\Ab\bS}({}_n F,E)\cong 0$ for all $n\in \nat$. This implies the assertion.\hB

\subsubsection{}

If $F\in \Sh\bS$ and $C\in \bS$, then we can form the sheaf
$\cR_C(F)\in \Sh\bS$ by the prescription
$\cR_C(F)(A):=F(A\times C)$ for all $A\in \bS$ (see \ref{hdfjcsdcuis}).
We will show that $\cR_C$ preserves torsion sheaves provided $C$ is compact.
 \begin{lem}\label{hdwdqqdqwdd}
If $C\in \bS$ is compact and $H\in \Sh_\Ab\bS$ is a torsion sheaf, then
$\cR_C(H)$ is a torsion sheaf.
\end{lem}
\proof 
Given $s\in \cR_C(H)(A)=H(A\times C)$ and $a\in A$ we must find $n\in \nat$ and a neighbourhood
$U$ of $a$ such that $(ns)_{|U\times C}=0$. Since $C$ is  compact and $A\in \bS$ is compactly generated by assumption on $\bS$,
the compactly generated topology (this is the topology we use here) on 
the product $A\times C$ coincides with the product topology (\cite[Thm. 4.3]{MR0210075}). Since $H$ is torsion
there exists an open covering $(W_i=A_i\times C_i)_{i\in I}$ of $A\times C$ and a family of non-zero  integers  $(n_i)_{i\in I}$ such that $(n_i s)_{|W_i}=0$.

The set of subsets of $C$ 
$$\{ C_i|i\in I \ ,\:a\in A_i\}$$
forms an open covering of $C$. Using the compactness of $C$ we choose a finite set
$i_1,\dots,i_r\in I$ with $a\in A_{i_k}$ such that
$\{C_{i_k}|k=1,\dots r\}$ is still an open covering of $C$.
Then we define the open neighbourhood $U$ of $a\in A$ by  
$U:=\cap_{k=1}^r A_{i_k}$.
 Set $n:=\prod_{k=1}^r n_{i_k}$. Then we have $(ns)_{U\times C}=0$.\hB

\subsubsection{}

\begin{lem}\label{hjkduidqedqwdwqd}
Let $H$ be a discrete torsion group and $G\in \bS$  be  a compactly generated\footnote{i.e. generated by a compact subset} group. 
Then the sheaf $\uHom_{\Sh_\Ab\bS}(\uG,\uH)$ is a torsion sheaf.
\end{lem}
\proof
Let $C\subseteq G$ be a compact generating set of $G$.
Precomposition with the inclusion $C\to G$ gives an inclusion
$\Hom_{\top-\Ab}(G,H)\to \Map(C,H)$.
We have for $A\in \bS$ 
\begin{eqnarray*}
\uHom_{\Sh_\Ab\bS}(\uG,\uH)(A)&\stackrel{Lemma \:\ref{eudiwe33}}{\cong}& \underline{\Hom_{\top-\Ab}(G,H)}(A)\\&\cong& \Hom_{\bS}(A,\Hom_{\top-\Ab}(G,H))\\&\subseteq&
\Hom_\bS(A,\Map(C,H))\\&=&
\Hom_\bS(A\times C,H)\\
&\cong&\uH(A\times C)\\
&\cong&\cR_C(\uH)(A)\ .
\end{eqnarray*}
By Lemma \ref{qwuiwdiudwd} the sheaf $\uH$ is a torsion sheaf. It follows from  Lemma \ref{hdwdqqdqwdd}
that $\cR_C(\uH)$ is a torsion sheaf. By the calculation above
$\uHom_{\Sh_\Ab\bS}(\uG,\uH)$ is a sub-sheaf of the torsion sheaf $\cR_C(\uH)$ and therefore itself a torsion sheaf. \hB

\subsubsection{}
 
\begin{lem}\label{hdjasduiwdiuqwd11}
If $G$  is a compact abelian group, then
$$\uExt^1_{\Sh_\Ab\bS}(\uG,\uZ)\cong \uHom_{\Sh_\Ab\bS}(\uG,\uT)\cong \underline{\widehat G}\ .$$
Moreover, if $G$ is profinite, then $\underline{\widehat G}$ is a torsion sheaf.
\end{lem}
\proof
We apply the functor $\Ext^*_{\Sh_\Ab\bS}(\uG,\dots)$ to $$0\to \uZ\to \uR\to \uT\to 0$$
and get the following segment of a long exact sequence
$$\dots \to \uHom_{\Sh_\Ab\bS}(\uG,\uR)\to \uHom_{\Sh_\Ab\bS}(\uG,\uT)
\to  \uExt^1_{\Sh_\Ab\bS}(\uG,\uZ)\to \uExt^1_{\Sh_\Ab\bS}(\uG,\uR)\to\dots\ .$$
Since $G$ is compact we have
$0=\underline{\Hom_{\top-\Ab}(G,\R)}\stackrel{Lemma \:\ref{eudiwe33}}{\cong}\uHom_{\Sh_\Ab\bS}(\uG,\uR)$.
Furthermore, by Proposition 
\ref{hdfjsdfiuiuwef}  
we have
$\uExt^1_{\Sh_\Ab\bS}(\uG,\uR)=0$.
Therefore we get
$$\uExt^1_{\Sh_\Ab\bS}(\uG,\uZ)\cong \uHom_{\Sh_\Ab\bS}(\uG,\uT)
\stackrel{Lemma \:\ref{eudiwe33}}{\cong} \underline{\Hom_{\top-\Ab}(G,\T)} \cong \underline{\widehat G}\ .$$
If $G$ is profinite, then (and only then) its dual $\widehat G$ is a discrete torsion group by
\cite[Cor. 8.5]{MR1646190}.
In this case, by Lemma \ref{qwuiwdiudwd} the sheaf $\widehat \uG$ is a torsion sheaf.
\hB

\subsubsection{}\label{hdbqwdwqdddqwd}

 Let $H$ be a discrete group. For $A\in \bS$ we consider the continuous group
cohomology $H^i_{cont}(G;\Map(A,H))$, which is defined as the cohomology of the group cohomology  complex
$$0\to \Map(A,H)\to \Hom_{\bS}(G,\Map(A,H))\to\dots\to\Hom_\bS(G^i ,\Map(A,H))\to $$
with the differentials dual to the ones given by \ref{iudiqwdwqdwqd}.
The map
$$\bS\ni A\mapsto H^i_{cont}(G;\Map(A,H))$$ defines a presheaf whose sheafification we will denote by  $\cH^i(G,H)$.
\begin{lem}\label{torzqwdguqwdqwd1ededee}
If $G$ is profinite and  $H$ is a discrete group,  then $\cH^i(G,H)$ is a torsion sheaf for $i\ge 1$. 
\end{lem}
\proof

We must show that for each section
$s\in H^i_{cont}(G;\Map(A,H))$ and $a\in A$ there exists a neighbourhood $U$ of $a$ and a number $l\in \Z$ 
such that $(ls)_{|U}=0$. This additional locality is important.
Note that by the exponential law we have
 $$C^i_{cont}(G,\Map(A,H)):=\Hom_\bS(G^i ,\Map(A,H))\cong  \Hom_\bS(G^i\times A,H)\ .$$
Let $s\in H^i_{cont}(G;\Map(A,H))$ and $a\in A$. Let $s$  be represented by a cycle
$\widehat s\in C^i_{cont}(G,\Map(A,H))$. Note that $\widehat s:G^i\times A\to H$ is locally constant.
The sets $\{\widehat s^{-1}(h)|h\in H\}$ form an open covering of $G^i\times A$.

Since $G$ and therefore $G^i$ are compact and $A\in \bS$ is compactly generated by assumption on $\bS$,
the compactly generated topology (this is the topology we use here) on 
the product $G^i\times A$ coincides with the product topology (\cite[Thm. 4.3]{MR0210075}).

Since $G^i\times \{a\}\subseteq G^i\times A$ is compact we can choose
a finite set $h_1,\dots,h_r\in H$ such that $\{\widehat s^{-1}(h_i)|i=1,\dots,r\}$ covers
$G^i\times \{a\}$. Now there exists an open neighbourhood
$U\subseteq A$ of $a$ such that
$G^i\times U\subseteq \bigcup_{s=1}^r\widehat s^{-1}(h_i)$.
On $G^i\times U$ the function $\widehat s$ has at most a finite number of values belonging to the set
$\{h_1,\dots,h_r\}$.

Since $G$ is profinite there exists a finite quotient group
$G\to \bar G$ such that
$\widehat s_{|G^i\times U}$ has a factorization
$\bar s:\bar G^i\times U\to H$.
Note that $\bar s$ is a cycle in $C^i_{cont}(\bar G,\Map(U,H))$.

Now we use the fact that the higher (i.e. in degree $\ge 1$) cohomology of a finite group with arbitrary coefficients is annihilated by the order of the group. Hence 
 $|\bar G| \bar s$ is the boundary of some $\bar t\in C^{i-1}_{cont}(\bar G,\Map(U,H))$.
Pre-composing $\bar t$ with the projection
$G^i\times U\to \bar G^i\times U$ we get
$\widehat t\in C^i_{cont}(G,\Map(U,H))$ whose boundary is $\widehat s$.
This shows that $(|\bar G|s)_{|U}=0$.
\hB

\subsubsection{}\label{uidqwiudwqdqwd}

Let $G$ be a profinite group. We consider the complex $U^\bullet:=U^\bullet(\uG)$ as defined in \ref{iudhuiqwdqwdwqd}.
Let $\uZ\to I^\bullet$ be an injective resolution.
Then we get a double complex
$\uHom_{\Sh_\Ab\bS}(U^\bullet,I^\bullet)$ as in \ref{uiwdiqwdwqdq}. 

\begin{lem}\label{fweifreiwrfowef}
For $i\ge 1$
\begin{equation}\label{hdfjdsf7823r21}
H^i\uHom_{\Sh_\Ab\bS}(U^\bullet,I^\bullet) \:\mbox{is a torsion sheaf}\ .
\end{equation}
\end{lem}
\proof
We first take the cohomology in the $I^\bullet$-, and then in the $U^\bullet$-direction.
 The first page of the resulting spectral sequence  is given  by
$$E^{r,q}_1\cong \uExt_{\Sh_\Ab\bS}^q(\Z(\uG^r),\uZ)\ .$$ 
Since the $G^r$ are profinite topological spaces, by Lemma \ref{hdjfcdcuisc} we get
$E_1^{r,q}\cong 0$ for $q\ge 1$. 
 We now consider the case $q=0$. For $A\in \bS$ we have
$$E_1^{r,0}(A)\cong \uZ(A\times G^r)\cong \Hom_{\bS}(G^r\times A,\Z)\cong \Hom_\bS(G^r,\Map(A,\Z))$$ for all $r\ge 0$.  The differential of the complex  
 $(E_1^{*,0}(A),d_1)$ is exactly the differential of the complex
$C^i_{cont}(G,\Map(A,\Z))$ considered in \ref{hdbqwdwqdddqwd}. By  Lemma
 \ref{torzqwdguqwdqwd1ededee} we conclude that the cohomology sheaves $E_2^{i,0}$ are torsion sheaves.
The spectral sequence degenerates at the second term.  We thus have shown that 
$H^i\uHom_{\Sh_\Ab\bS}(U^\bullet,I^\bullet)$ is a torsion sheaf. \hB

\subsubsection{}

For abelian groups $V,W$ let 
$V*_\Z W:=\Tor^\Z_1(V,W)$ denote the $\Tor$-product.
If $V$ and $W$ in addition are $\Z_{mult}$-modules then 
$V*W$ is a $\Z_{mult}$-module by the functoriality of the $\Tor$-product.
\begin{lem}\label{jhdbwqidwqqd}
If $V, W$ are $\Z_{mult}$-modules of weight  $k,l$, then
$V\otimes_\Z W$ and $V*_\Z W$ are of weight $k+l$.
\end{lem}
\proof
The assertion for the tensor product follows from (\ref{hwqdwqdwqd}).
Let $P^\bullet \to V^\prime$ and $Q^\bullet \to W^\prime$ be projective resolutions of the underlying $\Z$-modules $V^\prime,W^\prime$ of $V,W$. Then
$P^\bullet(k)\to V$ and $Q^\bullet(l)\to W$ are $\Z_{mult}$-equivariant resolutions of
$V$ and $W$. We have by (\ref{hwqdwqdwqd})
$$V*_\Z W\cong H^1(P^\bullet(k)\otimes_\Z Q^\bullet(l))\cong
H^1(P^\bullet\otimes_\Z Q^\bullet)(k+l)\ .$$\hB

\subsubsection{}

The tautological action of $\Z_{mult}$ on an abelian group $G$ (we write $G(1)$ for this $\Z_{mult}$-module) (see \ref{dkwqwjdkqkdjqwdw}) induces an action
of $\Z_{mult}$ on the  bar complex $\Z(G^\bullet)$ by diagonal action on the generators,
and therefore on the group cohomology $H^{*}(G(1);\Z)$ of $G$.

In the following Lemma we calculate the cohomology of the group $\Z/p\Z(1)$ as a $\Z_{mult}$-module.
\begin{lem}\label{uidiqwdwqdwqdwqd}
Let $p\in \N$ be a prime. Then
we have
$$H^{i}(\Z/p\Z(1);\Z)\cong \left\{\begin{array}{cc} \Z(0)&i=0\\0&\mbox{$i$ \em odd}\\
\Z/p\Z(k)&i=2k\ge 2\end{array}\right\}\ .$$
\end{lem}
\proof
The group cohomology of $\Z/p\Z$ can be identified as a ring with the ring
$\Z\oplus c\Z/p\Z [c]$, where $c$ has degree $2$. 
Furthermore, for every group there is a canonical map
$\hat G\to H^2(G;\Z)$ which in the case $G\cong \Z/p\Z$ happens to be an isomorphism.
This implies that $H^2(\Z/p\Z(1);\Z)$ has weight $1$ as a $\Z_{mult}$-module.
Since the cup product in group cohomology is natural it is $\Z_{mult}$-equivariant.
Therefore, the power $c^k$ generates a module of weight $k$. Hence $H^{2k}(\Z/p\Z(1);\Z)$ has weight $k$. 
\hB

\subsubsection{}\label{iwfehewjkfnwelfalekfeaf}

  For a $\Z_{mult}$-module $V$ (with $\Z_{mult}$-action $\Psi$) let $P^v_k$ be the
  operator $x\mapsto \Psi^vx-v^kx$. Note that this is a commuting family of operators.
We let $M_{23}^v$ be the monoid generated by $P^v_2,P^v_3$.
\begin{ddd}
  A $\Z_{mult}$-module $V$ is called a {weight $2$-$3$-extension} if for all
  $x\in V$ and all $v\in \Z_{mult}$ there is $P^v\in M_{23}^v$  such that 
  \begin{equation}\label{eq:2-3-ext}
  P^vx=0.
  \end{equation}

  A sheaf of $\Z_{mult}$-modules is called a {weight $2$-$3$-extension}, if
  every section locally satisfies Equation \eqref{eq:2-3-ext} (with $P^v$
  depending on the section and the neighborhood).
\end{ddd}

In the tables below we will mark weight $2$-$3$-extensions with attribute
weight $2-3$.

\begin{lem}\label{ex:weight-2-3-ex}
  Every $\Z_{mult}$-module of weight $2$ or of weight $3$ is a weight
  $2$-$3$-extension. The class of weight $2$-$3$-extensions is closed under extensions.
\end{lem}
\proof
  This is a simple diagram chase, using the fact that for a $\Z_{mult}$-module $W$ of weight
  $k$, $\Psi^vx=v^kx$ for all $k\in \Z_{mult}$ and $x\in W$.
\hB

\begin{lem}\label{heddwqdqwd}
Let $n\in \nat$ and $V$ be a torsion-free $\Z$-module. For $i\in \{0,2,3\}$ the cohomology 
$H^{i}((\Z/p\Z)^n(1);V)$ is a $\Z/p\Z$-module whose weight is given by the table
$$\begin{array}{|c|c|c|c|c|c|c|}\hline i&0 &2&3&4\\ 
\hline
weight&0&1&2&2\mbox{-}3\\\hline\end{array}\ .$$
Moreover, $H^{1}((\Z/p\Z)^n;V)\cong 0$.
\end{lem}
\proof
We first calculate 
$H^{i}((\Z/p\Z)^n;\Z)$ using the K{\"u}nneth formula and induction by 
$n\ge 1$. The start is Lemma \ref{uidiqwdwqdwqdwqd}.
Let us assume the assertion for products with less than $n$ factors.
The cases $i=0,1,2$ are straightforward.
We further get
$$H^3((\Z/p\Z)^n(1);\Z)\cong H^2(\Z/p\Z(1);\Z)*_\Z H^2((\Z/p\Z)^{n-1}(1);\Z)\ .$$
By Lemma \ref{jhdbwqidwqqd} the $*_\Z$-product of two modules of weight $1$ is of weight $2$.
Similarly, we have an exact sequence
\begin{multline*}
  0\to \bigoplus_{j=0}^4 H^j(\Z/p\Z^{n-1}(1);\Z)\otimes H^{4-j}(\Z/p\Z;\Z)\to
  H^4((\Z/p\Z)^n(1);\Z)\to\\
  \to \bigoplus_{j=2}^3 H^j(\Z/p\Z(1);\Z) *_\Z
  H^{5-j}((\Z/p\Z)^{n-1}(1);\Z)\to 0\ .
\end{multline*}
By Lemma \ref{ex:weight-2-3-ex} and induction we conclude that
$H^4((\Z/p\Z)^n(1);\Z)$ is a weight $2$-$3$-extension.

We can calculate the cohomology  $H^{i}(G;V)$ of a group $G$ in a trivial $G$-module by the standard complex
$C^\bullet(G;V):=C(\Map(G^\bullet,V))$. If $G$ is finite, then  we have
$C^\bullet(G;V)\cong C^\bullet(G;\Z)\otimes_\Z V$.
If $V$ is torsion-free, it is a flat $\Z$-module and therefore
$H^i(C^\bullet(G;\Z)\otimes_\Z V)\cong H^i(C^\bullet(G;\Z))\otimes_\Z V$.
Applying this to $G=(\Z/p\Z)^n$ we get the assertion from the special case $\Z=V$. \hB

\subsubsection{}

An abelian group $G$ is a $\Z$-module. Let $p\in \Z$ be a prime. If $pg=0$ for every $g\in G$, then
we say that $G$ is a $\Z/p\Z$-module.

\begin{ddd}
A sheaf $F\in \Sh_\Ab\bS$ is a sheaf of $\Z/p\Z$-modules if
$F(A)$ is a $\Z/p\Z$-module for all $A\in \bS$.

\end{ddd}

\subsubsection{}

Below we consider profinite abelian groups which are also $\Z/p\Z$ modules. The following lemma describes their  structure.
\begin{lem}\label{hbdwqdqdwqd}
Let $p$ be a prime number. If $G$ is compact and a $\Z/p\Z$-module, then there
exists a set $S$ such that $G\cong \prod_{S}\Z/p\Z$. 
In particular, $G$ is profinite.
\end{lem}
\proof
The dual group $\widehat G$ of the compact group $G$ is discrete and also a $\Z/p\Z$-module, hence an $\F_p$-vector space.
Let $S\subseteq \widehat G$ be an $\F_p$-basis. Then we can write
$\widehat G\cong \oplus_{S} \Z/p\Z$. Pontrjagin duality interchanges sums and products.
We get $$G\cong \widehat{\widehat G}\cong \widehat{\oplus_{S}\Z/p\Z}\cong \prod_{S}\widehat{\Z/p\Z}\cong \prod_{S}\Z/p\Z\ .$$
\hB

\subsubsection{}

Assume that $G$ is a compact group and a $\Z/p\Z$-module.  Note that the construction $G\to U^\bullet:=U^\bullet(G)$ is functorial in $G$ (see \ref{jhwefiewfewfe} for more details).
The tautological action of $\Z_{mult}$ on $G$ induces an action of $\Z_{mult}$ on 
$U^\bullet$. We can improve Lemma \ref{fweifreiwrfowef} as follows.

\begin{lem}\label{wehjdwqdqwdq}
Let $\uZ\to I^\bullet$ be an injective resolution. For $i\in\{2,3,4\}$ the sheaves
$$
H^i\uHom_{\Sh_\Ab\bS}(U^\bullet,I^\bullet)$$
are sheaves of $\Z/p\Z$-modules  and $\Z_{mult}$-modules whose weights are given by  the following table:
$$\begin{array}{|c|c|c|c|c|}\hline i&0&2&3&4\\
\hline
\mbox{weight}&0&1&2&2\mbox{-}3\\\hline\end{array}\ .$$
Moreover,   $H^1\uHom_{\Sh_\Ab\bS}(U^\bullet,I^\bullet)\cong 0$.
\end{lem}
\proof
We argue with the spectral sequence as in the proof of \ref{fweifreiwrfowef}.
Since $E_1^{p,q}=0$ for $q\ge 1$, it suffices to show that the sheaves
$E_2^{i,0}$ are sheaves of $\Z/p\Z$-modules and $\Z_{mult}$-modules of weight $k$,
where $k$ corresponds to $i$  as in the table (and with the appropriate
modification for $i=4$).

We have for $A\in \bS$
$$E_1^{i,0}(A)\cong \uZ(A\times G^i)\cong \Hom_{\bS}(G^i,\Map(A,\Z))=C^i_{cont}(G;\Map(A,\Z))\ ,$$
where for a topological $G$-module $V$ the complex $C_{cont}^\bullet(G;V)$
denotes the continuous group cohomology complex.  
We now consider the presheaf $\tilde X^i$ defined by
$$A\mapsto \tilde X^i(A):=H^i_{cont}(G,\Map(A,\Z)).$$
Then by definition $E_2^{i,0}:=X^i:=(\tilde X^i)^\sharp$ is the sheafification of $\tilde X^i$. 

The action of $\Z_{mult}$ on $G$ induces an action $q\mapsto [q]$ on the presheaf $ \tilde X^i$,
which descends to an action of $\Z_{mult}$ on the associated sheaf $X^i$.

We fix $i\in \{0,1,2,3\}$ and let $k$ be the associated weight as in the table.
For $i\ge 2$ we must show that for each section
$s\in H^i_{cont}(G;\Map(A,\Z))$ and $a\in A$ there exists a neighbourhood $U$ of $a$  
such that $((q^k-[q]) s)_{|U}=0$ and $(ps)_{|U}=0$ for all $q\in \Z_{mult}$. For
$i=4$, instead of $((q^k-[q])s)_{|U}=0$ we must find $P^v\in M^v_{23}$ (depending on $U$
and $s$) (see \ref{iwfehewjkfnwelfalekfeaf} for notation)   such that
$P^vs_{|U}=0$. 
 For $i=1$ we must show that
$s_{|U}=0$.

We perform the argument in detail for $i=2,3$. It is a refinement of the proof
of Lemma \ref{torzqwdguqwdqwd1ededee}. The cases $i=1$ and $i=4$ are very similar. 

Let $s$  be represented by a cycle
$\widehat s\in C_{cont}^i(G,\Map(A,\Z))$. Note that $\widehat s:G^i\times A\to \Z$ is locally constant.
The sets $\{\widehat s^{-1}(h)|h\in \Z\}$ form an open covering of $G^i\times A$.
Since $G$ and therefore $G^i$ is compact and $A\in \bS$ is compactly generated by assumption on $\bS$,
the compactly generated topology (this is the topology we use here) on 
the product $G^i\times A$ coincides with the product topology (\cite[Thm. 4.3]{MR0210075}).
This allows the following construction.

The set of subsets
$$\{\widehat s^{-1}(h)|h\in \Z\}$$
forms an open covering of $G^i\times A$. Using the compactness of the subset  $G^i\times\{a\}\subseteq G^i\times A$ 
we choose a finite set $h_1,\dots,h_r\in \Z$ such that
$G^i\times\{a\}\subseteq \bigcup_{i=1}^r\widehat s^{-1}(h_i)$.
There exists a neighbourhood $U\subseteq A$ of $a$ such that
$G^i\times U\subseteq \bigcup_{i=1}^r \widehat s^{-1}(h_i)$.

We now use that $G$ is profinite by Lemma \ref{hbdwqdqdwqd}.
Since $\widehat s(G^i\times U)$ is a finite set (a subset of $\{h_1,\dots,h_r\}$)
there exists a finite quotient group
 $G\to \bar G$ such that
$\widehat s_{|G^i\times U}$ has a factorization
$\bar s:\bar G^i\times U\to \Z$.
In our case $\bar G\cong (\Z/p\Z)^r$ for a suitable $r\in \nat$.
Note that $\bar s$ is a cycle in $C^i_{cont}(\bar G,\Map(U,\Z))$.
Now by Lemma \ref{heddwqdqwd} we know that $(q^k-[q]) \bar s$ and $p\bar s$ are the boundaries of some $\bar t,\bar t_1\in C^{i-1}_{cont}(\bar G,\Map(U,\Z))$. 
Pre-composing $\bar t,\bar t_1$ with the projection
$G^i\times U\to \bar G^i\times U$ we get
$\widehat t,\widehat t_1\in C^i_{cont}(G,\Map(U,\Z))$ whose boundaries are  $(q^k-[q])\widehat s$ and $p\widehat s$, respectively.
This shows that $(q^k-[q]s)_{|U}=0$ and $(ps)_{|U}=0$.
\hB

 \subsubsection{}\label{dgbsaddiuidiud}

We consider a compact group $G$ and form the double complex  $\Hom_{\Sh_\Ab\bS}(U^\bullet,I^\bullet)$ defined in \ref{uidqwiudwqdqwd}.
Taking the cohomology first in the $U^\bullet$ and then in the $I^\bullet$-direction we get a second spectral sequence $(F^{p,q}_r,d_r)$
with second page
$$F_2^{p,q}\cong \uExt^p_{\Sh_\Ab\bS}(H^q U^\bullet,\uZ)\ .
$$ 
In the following we calculate the term $F_2^{1,2}$ and show the vanishing of several other entries.
\begin{lem}\label{udwqidqwdcwqcc}
The left lower corner of $F_2^{p,q}$ has the form
$$
\begin{array}{|c|c|c|c|c|c|c|}\hline3&0&&&&&\\\hline
2&0&\uExt^1_{\Sh_\Ab\bS}((\Lambda_\Z^2\uG)^\sharp,\uZ)&&&&\\\hline
1&0&\uExt^1_{\Sh_\Ab\bS}(\uG,\uZ)&\uExt^2_{\Sh_\Ab\bS}(\uG,\uZ)&\uExt^3_{\Sh_\Ab\bS}(\uG,\uZ)&&\\\hline
0&\uZ&0&0&0&0&0\\\hline
&0&1&2&3&4&5\\\hline
\end{array}\ .$$
\end{lem}
\proof
The proof is similar to the corresponding argument in the proof of \ref{hdfjsdfiuiuwef}.
Using \ref{diqwodwqdwqdq} we get  
$$F^{p,q}_2\cong \uExt^p_{\Sh_\Ab\bS}((\Lambda^q_\Z\uG)^\sharp,\uZ)$$
for $q=0,1,2$. In particular we have $\Lambda^0_\Z\uG\cong \uZ$ and thus
$$F_2^{p,0}\cong \uExt^p_{\Sh_\Ab\bS}(\uZ,\uZ)\cong 0\ .$$
Furthermore, since $\Lambda^1_\Z\uG\cong \uG$ we have
$$F_2^{0,1}\cong \uHom_{\Sh_\Ab\bS}(\uG,\uZ)\stackrel{\ref{eudiwe33}}{\cong} \underline{\Hom_{\top-\Ab}(G,\Z)}\cong0$$
since $\Hom_{\top-\Ab}(G,\Z)\cong 0$ by compactness of $G$.

We claim that $$F_2^{0,2}\cong \uHom_{\Sh_\Ab\bS}((\Lambda_\Z^2\uG)^\sharp,\uZ)\cong 0\ .$$
Note that for $A\in \bS$ we have
 $$\uHom_{\Sh_\Ab\bS}((\Lambda_\Z^2\uG)^\sharp,\uZ)(A)\cong
\uHom_{\Pr_\Ab\bS}(\Lambda_\Z^2\uG,\uZ)(A)\cong \Hom_{\Pr_\Ab\bS/A}(\Lambda_\Z^2\uG_{|A},\uZ_{|A}) \ .$$
An element $\lambda\in \Hom_{\Pr_\Ab\bS/A}(\Lambda_\Z^2\uG_{|A},\uZ_{|A})$ induces a family of 
biadditive (antisymmetric) maps
$$\lambda^W:\uG(W)\times \uG(W)\to \uZ(W)$$ for $(W\to A)\in \bS/A$ which is compatible with restriction. Restriction to points gives biadditive maps
$G\times G\to \Z$. Since $G$ is compact the only such map is the constant map to zero.
Therefore $\lambda^W$ vanishes for all $(W\to A)$.
This proves the claim. The same kind of argument shows that
$\uHom_{\Sh_\Ab\bS}((\Lambda^q_\Z\uG)^\sharp,\uZ)(A)=0$ for $q\ge 2$.
 
Let us finally show that $F_2^{0,3}\cong \uHom_{\Sh_\Ab\bS}((H^3U^\bullet)^\sharp,\uZ)\cong 0$.
By \ref{eidoewd} we have an exact sequence (see (\ref{uebiwqdwq}) for the notation ${}^pD^\bullet$)
$$0\to K\to \Lambda^3_\Z \uG\to H^3({}^pD^\bullet)\to C\to 0\ ,$$
where $K$ and $C$ are defined as the kernel and cokernel presheaf, and the middle map 
becomes an isomorphism after tensoring with $\Q$.
This means that $0\cong K\otimes \Q$ and $0\cong C\otimes \Q$.
Hence, these $K$ and $C$ are presheaves of torsion groups.
Let $A\in \bS$. Then an element
$s\in \uHom_{\Sh_\Ab\bS}((H^3U^\bullet)^\sharp,\uZ)(A)$ induces
by precomposition an element
$\tilde s\in \uHom_{\Sh_\Ab\bS}((\Lambda^3_\Z\uG)^\sharp,\uZ)(A)\cong 0$.
Therefore $s$ factors over
$\bar s\in \uHom_{\Sh_\Ab\bS}(C^\sharp,\uZ)(A)$.
Since $C$ is torsion we conclude that $\bar s=0$ by Lemma
\ref{uihfuieddqwd}. The same argument shows that $\uHom_{\Sh_\Ab\bS}((H^q
U^\bullet)^\sharp,\uZ)\cong 0$ for $q\ge 3$. \hB

\subsubsection{}

\begin{lem}\label{gashdsasda}
Assume that $p$ is an odd prime number, $p\ne 3$.
If $G$ is a compact group and a $\Z/p\Z$-module, then
$\uExt^i_{\Sh_\Ab\bS}(\uG;\uZ)\cong 0$
for $i=2,3$.

If $p=3$, then at least $\uExt^2_{\Sh_\Ab\bS}(\uG,\uZ)\cong 0$.
\end{lem}
\proof
We consider the spectral sequence $(F_r,d_r)$ introduced in \ref{dgbsaddiuidiud}.
It converges to the  graded sheaves associated to a certain filtrations of the cohomology sheaves of the total complex of the double complex
 $\uHom_{\Sh_\Ab\bS}(U^\bullet,I^\bullet)$ defined in \ref{uidqwiudwqdqwd}.
In degree $2,3,4$ these cohomology sheaves are sheaves of $\Z/p\Z$-modules carrying actions of $\Z_{mult}$ with weights
determined in Lemma \ref{wehjdwqdqwdq}.

The left lower corner of the second page of the spectral sequence was evaluated in Lemma \ref{udwqidqwdcwqcc}.
Note that $\uExt^i_{\Sh_\Ab\bS}(\uG,\uZ)$ is a sheaf of $\Z/p\Z$-modules with an action of $\Z_{mult}$ of weight $1$ for all $i\ge 0$.
The term $F_2^{2,1}\cong \uExt^2_{\Sh_\Ab\bS}(\uG,\uZ)$ survives to the limit of the spectral sequence and is a submodule of  a sheaf of $\Z/p\Z$-modules of weight $2$. 
On the other hand it has weight $1$. 

A $\Z/p\Z$-module $V$ with an action $\Z_{mult}$ which has weights $1$ and $2$ at the same time must be trivial\footnote{Note that in general a $\Z/p\Z$-module can very well have different weights. E.g a module of weight $k$
has also weight $pk$.}.
In fact, for every  $q\in \Z_{mult}$ we get the identity $q^2-q=(q-1)q=0$ in
$\End_{\Z}(V)$, and this implies that $q\equiv 1 (\mod p)$ for all $q\not\in
p\Z$. From this follows $p=2$, and this case was excluded.

Similarly, the sheaf $\uExt^1_{\Sh_\Ab\bS}((\Lambda_\Z^2\uG)^\sharp,\uZ)$ is a sheaf of $\Z/p\Z$-modules of weight $2$ (see \ref{fbefewfwef}), and $\uExt^3_{\Sh_\Ab\bS}(\uG,\uZ)$ is a sheaf of $\Z/p\Z$-modules of weight $1$.
Since $p$ is odd, this implies that the differential
$d_2^{1,3}:\uExt^1_{\Sh_\Ab\bS}((\Lambda_\Z^2\uG)^\sharp,\uZ)\to \uExt^3_{\Sh_\Ab\bS}(\uG,\uZ)$ is trivial.
Hence, $\uExt^3_{\Sh_\Ab\bS}(\uG,\uZ)$ survives to the limit. It is a subsheaf
of a sheaf of $\Z/p\Z$-modules which is a  weight $2$-$3$-extension. On the
other hand it has weight $1$. Substituting the weight $1$-condition into
equation \eqref{eq:2-3-ext}, because then $P^v_2=v-v^2=(1-v)v$ and
$P^v_3=v-v^3=(1-v)(1+v)v$ this implies that locally every section satisfies 
$(1-v)^n(1+v)^jv^ns=0$ for suitable $n,j\in\nat$ (depending on $s$ and the
neighborhood) and for all $v\in \Z_{mult}$. If $p>3$ we can choose $v$ such that
$(1-v)$, $(1+v)$, $v$ are simultaneously units, and in this case the equation implies
that locally every section is zero.

We conclude that $\uExt^3_{\Sh_\Ab\bS}(\uG,\uZ)\cong 0$.
\hB 
%

\subsubsection{}



\begin{lem}\label{wdqwidwqdqdqwwwwq}
Let $G$ be a compact group which satisfies the two-three condition (see \ref{two-three-cond}).
Assume further that \begin{enumerate}
\item $G$ is profinite, or
\item $G$ is connected and locally topologically divisible \ref{ltdp}.
\end{enumerate}
 Then the sheaves
$\uExt^i_{\Sh_\Ab\bS}(\uG,\uZ)$ are torsion-free for $i=2,3$.
\end{lem}
\proof
We must show that for all primes $p$ and $i=2,3$ the maps of sheaves
$$\uExt^i_{\Sh_\Ab\bS}(\uG,\uZ)\stackrel{p}{\to} \uExt^i_{\Sh_\Ab\bS}(\uG,\uZ)$$ are injective.
The multiplication by $p$ can be induced by the multiplication 
$$p:\uG\to \uG\ .$$
We consider the exact sequence
$$0\to K\to G\stackrel{p}{\to} G\to C\to 0\ .$$
The groups $K$ and $C$ are groups of $\Z/p\Z$-modules and therefore profinite
by \ref{hbdwqdqdwqd}. 

We claim that the sequence of sheaves
\begin{equation*}
  0\to \uK\to \uG\xrightarrow{p}\uG \to \uC\to 0
\end{equation*}
is exact. We first show the  claim in the case that $G$ is profinite.
Since $C$ is profinite, by Lemma \ref{jhbuisaicascsa1} and Lemma \ref{dezuqwideqwd} 
we know that
$$0\to \underline{G/K}\to \uG\to C\to 0$$ is exact.
Furthermore $G/K$ is profinite so that  the projection $G\to G/K$ has local sections.
This implies that
$\underline{G/K}\cong \uG/\uK$, and hence the claim.

We now discuss the case of a connected locally topologically divisible $G$.
In this case $C=\{1\}$. Since by assumption $p:G\to G$ has local sections
we can again use Lemma  \ref{dezuqwideqwd}  in order to conclude.

We decompose this sequence into two short exact sequences $$0\to \uK\to \uG\to Q\to 0\ , \quad 0\to Q\to
\uG\to \uC\to 0\ .$$
Now we apply the functor $\uExt^*_{\Sh_\Ab\bS}(\dots,\uZ)$
to the first sequence and study the associated long exact sequence
$$\uExt^{i-1}_{\Sh_\Ab\bS}(\uG,\uZ)\to \Ext^{i-1}_{\Sh_\Ab\bS}(\uK,\uZ)\to \uExt^i_{\Sh_\Ab\bS}(Q,\uZ)\to \uExt^i_{\Sh_\Ab\bS}(\uG,\uZ) \to \dots$$
Note that $\Ext^{2}_{\Sh_\Ab\bS}(\uK,\uZ)\cong 0$ by Lemma \ref{gashdsasda} if $p$ is odd, and by Lemma \ref{ehjidqwwqd} if $p=2$ (the two-three condition ensures that $K$ in this case is an at most finite product of copies of $\Z/2\Z$.).
This implies that
$$\uExt^3_{\Sh_\Ab\bS}(Q,\uZ)\to \uExt^3_{\Sh_\Ab\bS}(\uG,\uZ) $$ is injective.

Next we show that
$$\uExt^2_{\Sh_\Ab\bS}(Q,\uZ)\to \uExt^2_{\Sh_\Ab\bS}(\uG,\uZ) $$ is injective. For this it suffices to see that
$$\uExt^{1}_{\Sh_\Ab\bS}(\uG,\uZ)\to \uExt^{1}_{\Sh_\Ab\bS}(\uK,\uZ)$$ is surjective. But this follows from the diagram
$$\xymatrix{\uExt^1_{\Sh_\Ab\bS}(\uG,\uZ)\ar[r] \ar[d]^\cong& \uExt^1_{\Sh_\Ab\bS}(\uK,\uZ)\ar[d]^\cong\\
\underline{\widehat G}\ar[r]^{surjective}&\underline{\widehat K}}\ ,$$
where the vertical maps are the isomorphisms given by Lemma \ref{hdjasduiwdiuqwd11}, and surjectivity
of $\widehat G\to \widehat K$ follows from the fact that Pontrjagin duality preserves exact sequences.

Next we apply $\uExt^*_{\Sh_\Ab\bS}(\dots,\uZ)$ to the sequence
$$0\to Q\to \uG\to \uC\to 0$$ and get the long exact sequence
$$ \Ext^{i}_{\Sh_\Ab\bS}(\uC,\uZ)\to \uExt^i_{\Sh_\Ab\bS}(\uG,\uZ)\to \uExt^i_{\Sh_\Ab\bS}(Q,\uZ) \to \dots\ .$$
Again we use that 
$\Ext^{i}_{\Sh_\Ab\bS}(\uC,\uZ)\cong 0$ for $i=2,3$ by Lemma  \ref{gashdsasda} for $p>3$, and by 
 Lemma \ref{ehjidqwwqd} if $p\in\{2,3\}$ (in this case $C$ is an at most  finite product of copies of $\Z/p\Z$ by the two-three condition) in order to conclude that
$$\uExt^i_{\Sh_\Ab\bS}(\uG,\uZ)\to \uExt^i_{\Sh_\Ab\bS}(Q,\uZ)$$ is injective for $i=2,3$.
Therefore the composition
$$\uExt^i_{\Sh_\Ab\bS}(\uG,\uZ)\to \uExt^i_{\Sh_\Ab\bS}(Q,\uZ)\to \uExt^i_{\Sh_\Ab\bS}(\uG,\uZ) $$ is injective for $i=2,3$.
This is what we wanted to show.
\hB

\subsubsection{}

\begin{lem}\label{wdqwidwqdqdqwwwwq1}
Let $G$ be a compact connected group which satisfies the two-three condition.
Then the sheaves
$\uExt^i_{\Sh_\Ab\bS_{lc-acyc}}(\uG,\uZ)$ are torsion-free for $i=2,3$.
\end{lem}
\proof
The only point in the proof of  Lemma \ref{wdqwidwqdqdqwwwwq} where we have used the condition of local topological divisibility was the exactness of the sequence
$$0\to \uK\to \uG\stackrel{p}{\to} \uG\to 0$$
of sheaves on $\bS$.
We show that this sequence is  exact with this condition if we consider the sheaves on $\bS_{lc-acyc}$ (by restriction).
Let us start with the dual sequence
$$0\to \hat G\stackrel{p}{\to} \hat G\to \hat K\to 0$$
of discrete groups. It gives an exact sequence of sheaves
$$0\to \hat \uG\stackrel{p}{\to} \hat \uG\to \hat \uK\to 0\ .$$
We apply
$\Hom_{\Sh_\Ab\bS_{lc-acyc}}(\dots,\uT)$ and get, using Lemma  \ref{eudiwe33}, the long exact sequence
$$0\to \uK\to \uG\stackrel{p}{\to} \uG \to \uExt^1_{\Sh_\Ab\bS_{lc-acyc}}(\hat \uK,\uT)\to \dots\ .$$
By Theorem \ref{uwdiqwdqwdwqdwqd} we have  $\uExt^1_{\Sh_\Ab\bS_{lc-acyc}}(\hat \uK,\uT)=0$.
Now we can argue as in the proof of \ref{wdqwidwqdqdqwwwwq}.
\hB

\subsubsection{}

\begin{lem}\label{djekwdkjwqdwq}
Let  $G$ be a profinite abelian group.
Then the following assertions are equivalent.
\begin{enumerate}
\item $\uExt_{\Sh_\Ab\bS}^i(\uG,\uZ)$ is torsion-free for $i=2,3$
\item $\uExt_{\Sh_\Ab\bS}^i(\uG,\uZ)\cong 0$ for $i=2,3$.
\end{enumerate}
\end{lem}
 \proof
It is clear that 2. implies 1. Therefore let us show that
$\uExt_{\Sh_\Ab\bS}^i(\uG,\uZ)\cong 0$ for $i=2,3$ under the assumption that we already know that it is torsion-free.

We consider the double complex $\uHom_{\Sh_\Ab\bS}(U^\bullet,I^\bullet)$ introduced in \ref{uidqwiudwqdqwd}.
By Lemma \ref{fweifreiwrfowef} we know that the cohomology sheaves
$H^i(\uHom_{\Sh_\Ab\bS}(U^\bullet,I^\bullet))$ of the associated total complex are torsion sheaves.

The spectral sequence $(F_r,d_r)$ considered in \ref{dgbsaddiuidiud} calculates the associated graded sheaves
of a certain filtration of $H^i(\uHom_{\Sh_\Ab\bS}(U^\bullet,I^\bullet))$. The left lower corner of its second page was already evaluated in Lemma \ref{udwqidqwdcwqcc}.

The term $F_2^{2,1}\cong \uExt^2_{\Sh_\Ab\bS}(\uG,\uZ)$ survives to the limit of the spectral sequence. On the one hand by our  assumption it is torsion-free.
On the other hand by Lemma \ref{fweifreiwrfowef}, case $i=2$, it is a subsheaf of a torsion sheaf. It follows that
$$\uExt^2_{\Sh_\Ab\bS}(\uG,\uZ)\cong 0\ .$$
This settles the implication $1.\Rightarrow  2.$ in case $i=2$ of Lemma \ref{djekwdkjwqdwq}.

We now claim that $\uExt^1_{\Sh_\Ab\bS}((\Lambda_\Z^2\uG)^\sharp,\uZ)$ is a torsion sheaf.
Let us assume the claim and finish the proof of Lemma \ref{djekwdkjwqdwq}.
Since $F_2^{3,1}\cong \uExt^3_{\Sh_\Ab\bS}(\uG,\uZ)$ is torsion-free by assumption the differential $d_2$ must be trivial by Lemma \ref{uihfuieddqwd}. Since also $F_3^{0,3}\cong 0 $
the sheaf $\uExt^3_{\Sh_\Ab\bS}(\uG,\uZ)$ survives to the limit of the spectral sequence. By Lemma 
\ref{fweifreiwrfowef}, case $i=3$, it is a subsheaf of a torsion sheaf and therefore itself a torsion sheaf. It follows that 
$\uExt^3_{\Sh_\Ab\bS}(\uG,\uZ)$ is a torsion sheaf and a torsion-free sheaf at the same time, hence trivial. This is the  assertion $1.\Rightarrow  2.$ of Lemma \ref{djekwdkjwqdwq} for $i=3$.

 
We now show the claim.
We start with some general preparations.
If $F,H$ are two sheaves on some site, then one forms the presheaf $\bS\ni A\mapsto (F\otimes^p_\Z H)(A):=F(A)\otimes_\Z H(A)\in \Ab$. The sheaf $F\otimes_\Z H$
is by definition the sheafification of $F\otimes^p_\Z H$.
We can write the definition of the presheaf $\Lambda_\Z^2\uG$
in terms of the following exact sequence of presheaves
$$ 0\to K\to {}^p\Z(\cF(\uG))\stackrel{\alpha}{\to} \uG\otimes_\Z^p \uG \to \Lambda_\Z^2\uG\to 0\ ,$$
where $K$ is by definition the kernel. 
For $W\in \bS$ the map
$\alpha_W:{}^p\Z(\uG(W))\to \uG(W)\otimes_\Z\uG(W)$ is defined on generators by
$\alpha_W(x)=x\otimes x$, $x\in \uG(W)$.
Sheafification is an exact functor and thus gives  
$$0\to ({}^p\Z(\cF(\uG))/K)^\sharp\to \uG\otimes_\Z \uG \to (\Lambda_\Z^2\uG)^\sharp\to 0\ .$$
 We now apply the functor
$\uHom_{\Sh_\Ab\bS}(\dots,\uZ)$ and consider
the following segment of the associated long exact sequence
\begin{multline}\label{t23dui2dzzzq}
  \to\uHom_{\Sh_\Ab\bS}(\uG\otimes_\Z \uG,\uZ)\to
  \uHom_{\Sh_\Ab\bS}(({}^p\Z(\cF(\uG))/K)^\sharp,\uZ)\\
  \to \uExt^1_{\Sh_\Ab\bS}((\Lambda_\Z^2\uG)^\sharp,\uZ) \to
\uExt^1_{\Sh_\Ab\bS}(\uG\otimes_\Z \uG,\uZ)\to 
\end{multline}

The following two facts imply the claim:
\begin{enumerate}
\item $\uExt^1_{\Sh_\Ab\bS}(\uG\otimes_\Z \uG,\uZ) 
 $ is torsion. 
\item  $\uHom_{\Sh_\Ab\bS}(({}^p\Z(\cF(\uG))/K)^\sharp,\uZ)$ vanishes.
\end{enumerate}
Let us start with 1. Let $\uZ\to I^\bullet$ be an injective resolution.
We study
$$H^1\uHom_{\Sh_\Ab\bS}(\uG\otimes_\Z \uG,I^\bullet)\ . $$
We have
$$\uHom_{\Sh_\Ab\bS}(\uG\otimes_\Z \uG,I^\bullet)\cong \uHom_{\Sh_\Ab\bS}(\uG,\uHom_{\Sh_\Ab\bS}(\uG,I^\bullet))\ .$$
Let $K^\bullet=\uHom_{\Sh_\Ab\bS}(\uG,I^\bullet)$.
Since $\uHom_{\Sh_\Ab\bS}(\uG,\uZ)\cong \underline{\Hom_{\top-\Ab}(G,\Z)}\cong 0$ by the compactness of $G$ the map
$d^0:K^0\to K^1$ is injective. Let $d^1:K^1\to K^2$ be the second differential.
Now 
$$H^1\uHom_{\Sh_\Ab\bS}(\uG\otimes_\Z \uG,I^\bullet)\cong \uHom_{\Sh_\Ab\bS}(\uG,\ker(d^1))/\im(d^0_*)\ ,$$
where $d^0_*:\uHom_{\Sh_\Ab\bS}(\uG,K^0)\to \uHom_{\Sh_\Ab\bS}(\uG,\ker(d^1))$ is induced by $d^0$.
Since $d^0$ is injective we have
$$\im(d^0_*)=\uHom_{\Sh_\Ab\bS}(\uG,\im(d^0))\ .$$
Applying $\uHom_{\Sh_\Ab\bS}(\uG,\dots)$ to the exact sequence
$$0\to K^0\to \ker(d^1)\to \ker(d^1)/\im(d^0)\to 0$$
and using $\ker(d^1)/\im(d^0)\cong \uExt^1_{\Sh_\Ab\bS}(\uG,\uZ)$ we get
\begin{eqnarray*}
0\to \uHom_{\Sh_\Ab\bS}(\uG,K^0)\stackrel{d^0_*}{\to} \uHom_{\Sh_\Ab\bS}(\uG,\ker(d^1))\to &&\\\to \uHom_{\Sh_\Ab\bS}(\uG,\uExt^1_{\Sh_\Ab\bS}(\uG,\uZ))\to \uExt^1_{\Sh_\Ab\bS}(\uG,K^0)\to \dots&&\ .\end{eqnarray*}
In particular,\begin{eqnarray*}
\uExt^1_{\Sh_\Ab\bS}(\uG\otimes_\Z \uG,\uZ) &\cong& \uHom_{\Sh_\Ab\bS}(\uG,\ker(d^1))/\uHom_{\Sh_\Ab\bS}(\uG,K^0)
\\&\subseteq& \uHom_{\Sh_\Ab\bS}(\uG,\uExt^1_{\Sh_\Ab\bS}(\uG,\uZ))\ .\end{eqnarray*}
By Lemma \ref{hdjasduiwdiuqwd11} we know that $\uExt^1_{\Sh_\Ab\bS}(\uG,\uZ)\cong \underline{\widehat G}$. Since
$G$ is compact and profinite, the group $\widehat G$ is discrete and torsion.
It follows  by Lemma \ref{hjkduidqedqwdwqd} that $$\uHom_{\Sh_\Ab\bS}(\uG, \uExt^1_{\Sh_\Ab\bS}(\uG,\uZ))$$ is a torsion sheaf.
A subsheaf of a torsion sheaf is again a torsion sheaf. This finishes the argument for the first fact.

We now show the fact 2. 
 
Note that
\begin{equation}\label{qwhdwquidqdqd}
\uHom_{\Sh_\Ab\bS}(({}^p\Z(\cF(\uG))/K)^\sharp,\uZ)\cong \uHom_{\Sh_\Ab\bS}(\Z(\cF(\uG))/K^\sharp,\uZ)\ .
\end{equation}
Consider $A\in \bS$ and 
$\phi\in \uHom_{\Sh_\Ab\bS}(\Z(\cF(\uG))/K^\sharp,\uZ)(A)$.
We must show that each $a\in A$ has a neighbourhood $U\subseteq A$ such that $\phi_{|U}=0$.
 Pre-composing $\phi$ with the projection
$\Z(\cF(\uG))\to \Z(\cF(\uG))/K^\sharp$ we get
 an element
$$\bar \phi\in\uHom_{\Sh_\Ab\bS}(\Z(\cF(\uG)),\uZ)(A)\cong \uHom_{\Sh\bS}(\cF(\uG),\uZ)(A)\stackrel{Lemma\: \ref{hgeuidwqdqwdwqd}}{\cong} 
\uZ(A\times G)\ .$$

We are going to show that $\bar \phi=0$ after restriction to a suitable neighbourhood of  $a\in A$ using the fact that it annihilates $K^\sharp$.
Let us start again on the left-hand side of (\ref{qwhdwquidqdqd}).
We have
\begin{equation*}
  \begin{split}
    \uHom_{\Sh_\Ab\bS}(({}^p\Z(\cF(\uG))/K)^\sharp,\uZ)(A) &\cong
    \uHom_{\Pr_\Ab\bS}(({}^p\Z(\cF(\uG))/K),\uZ)(A)\\
    &\cong
    \Hom_{\Pr_\Ab\bS/A}({}^p\Z(\cF(\uG_{|A}))/K_{|A},\uZ_{|A})\ .
\end{split}
\end{equation*}
For $(A\times G\to A)\in \bS/A$ the morphism $\phi$ gives rise to a group homomorphism
$$\widehat \phi:\Z(\cF(\uG(A\times G)))/K(A\times G)\to \uZ(A\times G)\ .$$
The symbol $\cF(\uG(A\times G))$ denotes the underlying set of the group $\uG(A\times G)$.
We have $(\pr_G:A\times G\to G)\in \uG(A\times G)$ and by the explicit description of the element $\bar \phi$
given after the proof of Lemma \ref{hgeuidwqdqwdwqd} we see
$\bar \phi=\widehat \phi\left(\{[\pr_G]\}\right)$, where $[\pr_G]\in \Z(\cF(\uG(A\times G)))$ denotes the generator
corresponding to $\pr_G\in \uG(A\times G)$, and $\{\dots\}$ indicates that we take the class modulo
$K(A\times G)$.

The homomorphism
$\widehat \phi\in \Hom_\Ab(\Z\cF(\uG(A\times G))/ K(A\times G),\uZ(A\times G))$ is represented by a homomorphism
$$\tilde \phi:\Z\cF(\uG(A\times G))\to \uZ(A\times G)\ ,$$
and we have $\bar \phi=\tilde \phi([\pr_G])$.
By definition of $K$ we have the exact sequence
$$0\to K(A\times G)\to \Z(\cF(\uG(A\times G)))\to \uG(A\times G)\otimes_\Z \uG(A\times G)\to \Lambda^2_\Z \uG(A\times G)\to 0\ .$$
For $x\in \uG(A\times G)$ let $[x]\in \Z( \cF(\uG(A\times G)))$ denote the corresponding generator.
Since $(nx)\otimes (nx)= n^2(x\otimes x)$ in $\uG(A\times G)\otimes_\Z \uG(A\times G)$
 we have $n^2[x]-[nx]\in K(A\times G)$ for $n\in \Z$.
It follows that $\tilde \phi$ must satisfy the relation
$\tilde \phi(n^2[x])=\tilde \phi([nx])$ for all $n\in \Z$ and $x\in \uG(A\times G)$.
Let us now apply this reasoning to $x:=[\pr_G]$.
For every $n\in \Z$ we have a map $A\times G\stackrel{\id_A\times n=:\psi_n}{\to}A\times G$, and by naturality the map $\tilde \phi$ respects this action.
Moreover, the element $[n x]\in \Z(\cF(\uG(A\times G)))$ is obtained from $[x]$ via this action, since
$$\xymatrix{A\times G\ar[d]^{\pr_G}\ar[r]^{\psi_n}&A\times G\ar[d]^{\pr_G}\\G\ar[r]^n&G}$$ commutes.
It follows that
\begin{eqnarray*}
n^2\bar \phi&=&\tilde \phi(n^2[\pr_G])\\&=&\tilde \phi([n \pr_G])\\&=&\tilde \phi(\Z\cF(\bar G(\psi_n))([\pr_G]))\\&=&\uZ(\psi_n)(\bar \phi)\\&=&n^*\bar \phi\ ,
\end{eqnarray*}
where we write $n^*:=\uZ(\psi_n)$.
For $k\in \Z$ we define the open subset $V_k:=\bar \phi^{-1}(k)\subseteq A\times G$.
The family $(V_k)_{k\in K}$ forms an open pairwise disjunct  covering of $A\times G$.
Since $G$ is compact the compactly generated topology of $A\times G$ is the product topology
(\cite[Thm. 4.3]{MR0210075}).
Since $G$ is compact, we can choose a finite sequence $k_1,\dots,k_r\in \Z$ such that
$G\times \{a\}\subseteq \bigcup_{i=1}^r V_k$.
Furthermore, we can find a neighbourhood $U\subseteq A$ of $a$ such that
$G\times U\subseteq \bigcup_{i=1}^r V_k$.

Note that $\bar \phi_{|G\times U}$ has at most finitely many values.
Therefore there exists a finite quotient
$G\to F$ such that $\bar \phi_{|U\times G}$ factors through $f:U\times F\to \Z$.

The action of $\Z$ on $G$ is compatible with the corresponding  action of $\Z$ on $F$, and we have an action
$n^*:\Hom_{\bS}(U\times F,\Z)\to \Hom_\bS(U\times F,\Z)$.
The element $f\in \Hom_{\bS}(U\times F,\Z)$ still satisfies
  $n^2 f=n^* f$.
We now take $n:=|F|+1$. Then
$n^*=\id$ so that $(n^2-1)f=0$, hence $f=0$.
This implies $\bar \phi_{|U\times G}=0$.
 
 This finishes the proof of the second fact. \hB

\subsubsection{}

The Lemma \ref{wdqwidwqdqdqwwwwq} verifies Assumption  1. in \ref{djekwdkjwqdwq} for a large class of profinite groups.
\begin{kor}
If $G$ is a profinite group which satisfies the two-three condition \ref{two-three-cond}, then we have
$\Ext_{\Sh_\Ab\bS}^i(\uG,\uZ)\cong 0$ for $i=2,3$.
\end{kor}

\begin{theorem}\label{wqdiqwdiudwqdd}
A profinite abelian group which satisfies the two-three condition is admissible. 
\end{theorem}

\subsection{Compact connected abelian groups}

\subsubsection{}

Let $G$ be a compact abelian group.
We shall use the following fact shown in
\cite[Corollary 8.5]{MR1646190}.
\begin{fact}\label{thisfact}
$G$ is connected if and only if $\widehat G$ is torsion-free.
\end{fact}

\subsubsection{}\label{hjbdqwdqwuidwqdq}

For a space $X\in \bS$ and $F\in \Pr_\Ab\bS$ let $\check{H}^*(X;F)$ denote the \v{C}ech cohomology of $X$
with coefficients in $F$. It is defined as follows. To each open covering $\cU$ one associates the \v{C}ech complex
$\check{C}^\bullet(\cU;F)$. The open coverings form a left-filtered category whose morphisms are refinements.
The \v{C}ech complex depends functorially on the covering, i.e.~if
 $\cV\to \cU$ is a refinement, then we have a functorial chain map
$\check{C}^\bullet(\cU;F)\to \check{C}^\bullet(\cV;F)$.
We  define 
$$\check{C}^\bullet(X;F):=\colim_{\cU}\check{C}^\bullet(\cU;F)$$
and 
$$\check{H}^*(X;F):=\colim_{\cU}H^*(\check{C}^\bullet(\cU;F))\cong H^*(\colim_{\cU}\check{C}^\bullet(\cU;F))\ \cong H^*(\check{C}^\bullet(X;F))\ .$$

\subsubsection{}
Fix a discrete group $H$. Let ${}^p\uH$ be the constant presheaf with values
$H$. Note that then the sheafification ${}^p\uH^\sharp$ is isomorphic to $\uH$
as defined in \ref{multbla}. Moreover
\begin{equation}\label{checkcomp}
  \check{C}^\bullet(X;{}^p\uH)\cong \check{C}^\bullet(X,\uH)\ .
\end{equation}

In general \v{C}ech cohomology differs from sheaf cohomology. The relation between these two is given by the
\v{C}ech cohomology spectral sequence $(E_r,d_r)$ (see \cite[3.4.4]{MR1317816}) converging to $H^*(X;F)$.
Let $\cH^*(F):=R^*i(F)$ denote the right derived functors of the inclusion $i:\Sh_\Ab\bS\to \Pr_\Ab\bS$. 
Then the second page of the spectral sequence is given by
$$E_2^{p,q}\cong \check{H}^p(X;\cH^q(F))\ .$$
By  \cite[3.4.7]{MR1317816} the edge homomorphism
$$\check{H}^p(X;F)\to H^p(X;F)$$ is an isomorphism for $p=0,1$ and injective for $p=2$.

\subsubsection{}

We now observe that \v{C}ech cohomology transforms strict inverse limits of compact spaces into colimits of cohomology groups.
\begin{lem}\label{euihuqowdqwdwdwqd}
Let $H$ be a discrete abelian group.
If $(X_i)_{i\in I}$ is an inverse system of compact spaces in $\bS$ such that $X=\lim_{i\in I}X_i$, then
$$\check{H}^p(X;\uH)\cong \colim_{i\in I}\check{H}^p(X_i;\uH)\ .$$
\end{lem}
\proof 
We first show that the system of open coverings of $X$ contains a cofinal system of coverings  which are pulled back
from the quotients  $p_i:X\to X_i$.
 Let $\cU=(U_r)_{r\in R}$ be a covering by open subsets. For each $r$ there exists a family
$I_r\subset I$ and subsets $U_{r,i}\subseteq X_i$, $i\in I_r$ such that
$U_r=\cup_{i\in I_r}p_i^{-1}(U_{r,i})$.
The set of open subsets
$\{p_i^{-1}(U_{r,i})|r\in R,\: i\in I_r\}$
is an open covering of $X$.
Since $X$ is compact we can choose a finite subcovering
$\cV:=\{U_{r_1,i_1},\dots U_{r_k,i_k}\}$ which can naturally be viewed
as a refinement of $\cU$.
Since $I$ is left filtered we can  choose $j\in I$ such that
$j<i_d$ for $d=1,\dots,k$. For $j\le i$ let 
$p_{ji}:X_j\to X_i$ be the structure map of the system which are all surjective by the strictness assumption on the inverse system. Therefore
$\cV^\prime:=\{p_{j,i_d}^{-1}(U_{r_d,i_d})\:|\:d=1,\dots,k\}$
is an open covering of $X_j$, and $\cV=p^{-1}_j(\cV^\prime)$.

Now we
observe that $\check{C}^\bullet(p^*\cV';{}^p\uH)\cong \check{C}^\bullet(\cV';{}^p\uH)$.
Therefore
\begin{eqnarray*}\check{C}^\bullet(X;\uH)&\stackrel{\eqref{checkcomp}}{\cong}&
  \check{C}^\bullet(X;{}^p\uH)\\
  &\cong & \colim_{\cU}\check{C}^\bullet(\cU;{}^p\uH)\\
&=&\colim_{i\in I}\colim_{\small \mbox{coverings $\cU$ of $X_i$}}\check{C}^\bullet(\cU;{}^p\uH)\\
&\cong&\colim_{i\in I}\check{C}^\bullet(X_i;{}^p\uH)\\
&\stackrel{\eqref{checkcomp}}{\cong}&\colim_{i\in I}\check{C}^\bullet(X_i;\uH)
\end{eqnarray*}
since one can interchange colimits.
Since filtered colimits are exact and therefore commute with taking cohomology this implies
the Lemma. \hB

\subsubsection{}

Next we show that \v{C}ech cohomology has a K{\"u}nneth formula.
\begin{lem}\label{uihdqwdwdwdwqdq}
 Let $H$ be a discrete ring of finite cohomological dimension. 
Assume that $X,Y$ are compact.
Then there exists a K{\"u}nneth spectral sequence with second term
$$E^2_{p,q}:=\bigoplus_{i+j=q}\Tor^H_p(\check{H}^i(X;\uH),\check{H}^j(Y;\uH))$$
which converges to $\check{H}^{p+q}(X\times Y;\uH)$.
\end{lem}
\proof
Since $X$ is  compact the topology on $X\times Y$ is the product topology
\cite[Thm. 4.3]{MR0210075}.
Using again  the compactness of $X$ and $Y$ 
 we can find a cofinal system of coverings of $X\times Y$ of the form
$p^*\cU\cap q^*\cV$ for coverings $\cU$ of $X$ and $\cV$ of $Y$, where
$p:X\times Y\to X$ and $q:X\times Y\to Y$ denote the projections, and the intersection of coverings
is the covering by the collection of all cross intersections.
We have
$$\check{C}^*(p^*\cU\cap q^*\cV;{}^p\uH)\cong\check{C}^*(p^*\cU;{}^p\uH)\otimes_H
\check{C}^*(q^*\cV;{}^p\uH)\ .$$
Since the tensor product over $H$ commutes with colimits we get
$$\check{C}^*(X\times Y;{}^p\uH)\cong\check{C}^*(X;{}^p\uH)\otimes_H
\check{C}^*(Y;{}^p\uH)\ ,$$
and hence by \eqref{checkcomp}
$$\check{C}^*(X\times Y;\uH)\cong\check{C}^*(X;\uH)\otimes_H
\check{C}^*(Y;\uH)\ .$$

The K\"unneth spectral sequence is the spectral sequence associated to this
double complex of flat (as colimits of free) $H$-modules. 
\hB

\subsubsection{}

We now recall that  sheaf cohomology also transforms strict inverse limits of spaces 
into colimits, and that it has a K{\"u}nneth spectral sequence. The following is a specialization of \cite[II.14.6]{MR1481706} to compact spaces.
\begin{lem}
Let $(X_i)_{i\in I}$ be an inverse system of compact spaces and $X=\lim_{i\in
  I}X_i$, and $H$ let be a discrete abelian group.
Then 
$$H^*(X;\uH)\cong \colim_{i\in I} H^*(X_i;\uH)\ .$$
\end{lem}
 
For simplicity we formulate the K{\"u}nneth formula for the sheaf $\uZ$ only.
The following is a specialization of
 \cite[II.18.2]{MR1481706} to compact spaces and the sheaf $\uZ$.
\begin{lem}\label{qdhwqdqwdqwd}
For compact spaces $X,Y$ we
have a K{\"u}nneth spectral sequence with second term
$$E^2_{p,q}=\bigoplus_{i+j=q}\Tor^\Z_p(H^i(X;\uZ),H^j(Y;\uZ))$$
which converges to $H^{p+q}(X\times Y;\uZ)$.
\end{lem}
Of course, this formulation is much too complicated since $\Z$ has cohomological dimension one. In fact, the
 spectral sequence decomposes into a collection of short exact sequences.

\subsubsection{}

\begin{lem}\label{jkbeqwdwqdqwdw}
If $G$ is a compact connected abelian group, then
$H^*(G;\uZ)\cong \check{H}^*(G;\uZ)$.
\end{lem}
\proof
The \v{C}ech cohomology spectral sequence provides the map
$$\check{H}^*(G;\uZ)\to H^*(G;\uZ)\ .$$
We now use that fact that $G$ is the projective limit of groups isomorphic to $T^a$, $a\in \nat$.
Since the compact abelian group $G$  is connected, by Fact \ref{thisfact} we know that
$\widehat G$ is torsion-free.
Since  $\widehat G$ is torsion-free it is the filtered colimit of its finitely generated subgroups $\widehat F$.
If $\widehat F\subset \widehat G$ is finitely generated, then $\widehat F\cong \Z^a$ for some $a\in \nat$,
therefore $F:=\widehat{\widehat F}\cong T^a$. Pontrjagin duality transforms
the filtered colimit $\widehat G\cong \colim_{\widehat F}\widehat F$ into a strict limit $G\cong \lim_{\widehat F} F$.

Since $T^a$ is a manifold it admits a cofinal system of good open coverings where all multiple intersections
are contractible. Therefore we get the isomorphism
$$\check{H}^*(T^a;\uZ)\stackrel{\sim}{\to} H^*(T^a;\uZ)\ .$$
Since both cohomology theories commute with limits of compact spaces 
we get
$$\check{H}^*(G;\uZ)\stackrel{\sim}{\to} H^*(G;\uZ)\ .$$
\hB

\subsubsection{}

We now use the calculation of cohomology
of the underlying space of a connected compact abelian group $G$. 
By \cite[Theorem 8.83]{MR1646190} we have
\begin{equation}\label{eiu21hi21h211}
\check{H}^*(G,\uZ)\cong \Lambda_\Z^*\widehat G
\end{equation}
as a graded Hopf algebra.
This result uses the two properties \ref{euihuqowdqwdwdwqd} and
\ref{uihdqwdwdwdwqdq} of \v{C}ech cohomology in an essential way.

By Lemma \ref{jkbeqwdwqdqwdw} we also have
$$ H^*(G;\uZ)\cong \Lambda_\Z^*\widehat G\ .$$

Note that $\Lambda_\Z^*\widehat G$ is torsion free. In fact, $\Lambda_\Z^*\widehat G\cong\colim_{\widehat F} \Lambda^*_\Z\widehat F$,
where the colimit is taken over all finitely generated subgroups $\hat{F}$ of
$\hat G$. We therefore get a colimit of injections of
  torsion-free abelian groups, which is itself torsion-free.

\subsubsection{}

Recall the notation related to $\Z_{mult}$-actions introduced in \ref{udwwuqidqqw434}. 

 \begin{lem}\label{uifhweifwefwe}
If $V,W\in \Ab$ with $W$ torsion-free and $k,l\in \Z$, $k\not=l$, then $$\Hom_{\Z_{mult}-\Mod}(V(k),W(l))=0\ .$$
\end{lem}
\proof
Let $\phi\in \Hom_{\Ab}(V,W)$ and $v\in V$.
Then for all $m\in \Z_{mult}$ we have
$$m^k\phi(v)=\Psi^m(\phi(v))=\phi(\Psi^m (v))=\phi(m^l v)=m^l \phi(v)\ ,$$ i.e. each $w\in \im(\phi)\subseteq W(l)$ satisfies
$(m^k-m^l)w=0$ for all $m\in \Z_{mult}$.
This set of equations implies that $w=0$, since $W$ is torsion-free.
\hB

\subsubsection{}\label{jkdqwdwqdwqd}

\begin{lem}\label{uifhweifwefwe1}
If $V,W\in \Sh_\Ab\bS$ with $W$ torsion-free and $k,l\in \Z$, $k\not=l$, then $$\Hom_{\Sh_{\Z_{mult}-\Mod}\bS}(V(k),W(l))=0\ .$$
\end{lem}
We leave the easy proof of this sheaf version of Lemma \ref{uifhweifwefwe}
to the interested reader. Note that a subquotient of a sheaf of $\Z_{mult}$-modules
of weight $k$ also has weight $k$.

\subsubsection{}\label{jekjfefwefwef}

By $\bS_{lc}\subset\bS$ we denote the sub-site of locally compact objects.
The restriction to locally compact spaces becomes necessary because of the use of the K{\"u}nneth (or base change) formula below.

The first page of the spectral sequence
$(E_r,d_r)$ introduced in \ref{uiwdiqwdwqdq} is given by
$$E_1^{q,p}=\uExt^p_{\Sh_\Ab\bS}(\Z(\uG^q),\Z)\ .$$
This sheaf is the sheafification of the presheaf
$$\bS\ni A\mapsto H^p(A\times G^q,\uZ)\in \Ab\ .$$
In order to calculate this  cohomology we  use the K{\"u}nneth formula \cite[II.18.2]{MR1481706}
   and that $H^*(G,\uZ)$ is torsion-free (\ref{thisfact}). We get for $A\in \bS_{lc}$ that
$$H^*(A\times G^q,\uZ)\cong H^*(A,\uZ)\otimes_\Z  (\Lambda^*\widehat G)^{\otimes_\Z q}$$
for locally compact $A$.
The sheafification of $\bS_{lc}\ni A\to H^i(A,\uZ)$ vanishes for $i\ge 1$, and gives $\uZ$ for $i=0$.
Since sheafification commutes with the tensor product with a fixed group we get
\begin{equation}\label{uiqwdiqwdwqdw}
(E_1^{q,*})_{|\bS_{lc}} \cong \underline{(\Lambda^*\widehat G)^{\otimes_\Z q}}\ .
\end{equation}

\subsubsection{}\label{jdqwdwqdwqd}

We now consider the tautological action of the multiplicative semigroup $\Z_{mult}$ on $G$ of weight $1$. As before
we write $G(1)$ for the group $G$ with this action.
Observe that this action is continuous.
Applying the duality functor we get an action of $\Z_{mult}$ on the dual group $\widehat G$ which is
also of weight $1$. Therefore $$\widehat{G(1)}=\widehat G(1)\ .$$
 
The calculation of the cohomology
(\ref{eiu21hi21h211}) of the topological space $G$ with $\Z$-coefficients is functorial in $G$.
We conclude that
$H^*(G(1),\Z)\cong \Lambda_\Z^*(\widehat G(1))$ is a  decomposable $\Z_{mult}$-module.
By \ref{fbefewfwef} the group
$$H^p(G^q(1),\Z)\cong [(\Lambda^*\widehat G(1))^{\otimes_\Z q}]^p$$
is of weight $p$.

\subsubsection{}\label{jhwefiewfewfe}

We now observe that the spectral sequence $(E_r,d_r)$ is functorial in $G$.
To this end we use the fact that
$G\mapsto U^\bullet(\uG)$ is a covariant functor from groups $G\in \bS$ to complexes of sheaves
on $\bS$.  If $G_0\to G_1$ is a homomorphism of topological groups in $\bS$,
then we get an induced map $U^q(G_0)\to U^q(G_1)$, namely the map
$\Z(\uG_0^q)\to \Z(\uG^q_1)$.
Under the identification made in \ref{jekjfefwefwef} the induced map
 $$\uExt^p_{\Sh_\Ab\bS_{lc}}(\Z(\uG^q_1),\uZ)\to \uExt^p_{\Sh_\Ab\bS_{lc}}(\Z(\uG_0^q),\uZ)$$
goes to the map
$$\underline{H^p(G_1^q;\uZ)}\to \underline{H^p(G_0^q;\uZ)}$$
induced by the pull-back
 $$H^p(G_1^q;\uZ)\to H^p(G_0^q;\uZ)$$
associated to the map of spaces $G_0^q\to G_1^q$.

\subsubsection{}

The discussion in \ref{jhwefiewfewfe} and \ref{jdqwdwqdwqd} shows that
the $\Z_{mult}$-module structure $G(1)$ induces one on the spectral sequence
$(E_r,d_r)$, and we see that
$E_1^{q,p}$ has weight  $p$ (which is the number of factors $\widehat G(1)$ contributing to this term).

We introduce the notation $H^k:=H^k(\uHom_{\Sh_\Ab\bS}(U^\bullet,I^\bullet))_{|\bS_{lc}}$.
Note that  $H^k$ has a filtration 
$$0=F^{-1}H^k\subseteq F^{0}H^k\subseteq\dots\subseteq F^kH^k=H^k$$
which is preserved by the action of $\Z_{mult}$.
The spectral sequence $(E_r,d_r)$ converges to the associated 
graded sheaf $\Gr(H^k)$.

Note that $\Gr^p(H^k)$ is a subquotient of $E_2^{k-p,p}$.
Since a subquotient of a sheaf of weight $p$ also has weight $p$ (see \ref{jkdqwdwqdwqd}) we get the following conclusion.
\begin{kor}\label{gh12}
$\Gr^p(H^k)$ has weight $p$.
\end{kor}

\subsubsection{}

We now study some aspects of the second page
$E_2^{q,p}$ of the spectral sequence in order to show the following  Lemma.
\begin{lem}\label{gh122}
\begin{enumerate}
\item For $k\ge 1$ we have $F^0H^k\cong 0$.
\item For $k\ge 2$ we have $F^1H^k\cong 0$.
\end{enumerate}
\end{lem}
\proof
We start with $1$.
We see that
$E_2^{*,0}$ is the cohomology of the complex
$E_1^{*,0}\cong \uHom_{\Sh_\Ab\bS}(U^\bullet,\uZ)$.
Explicitly, using in the last step the fact that $G$ and therefore $G^q$ are connected,
$$E_1^{q,0}\cong \uHom_{\Sh_\Ab\bS}(U^q,\uZ)\stackrel{(\ref{uieiefefewfw})}{\cong} \uHom_{\Sh_\Ab\bS}(\Z(\uG^q),\uZ) \stackrel{(\ref{uidwdqdiwqdqdwq09})}{\cong} \cR_{G^q}(\uZ)  \stackrel{Lemma \: \ref{ddwqdiwduiwqdiud}}{\cong} \underline{\Map(G^q,\Z)}\cong  \uZ\ .$$
An inspection of the formulas \ref{iudiqwdwqdwqd} for the differential of the complex
$U^\bullet$ shows that
$d_{1}:E_1^{q,0}\to E_1^{q+1,0}$
vanishes for even $q$, and is the identity for odd $q$. In other words this complex is isomorphic to
$$0\to \uZ\stackrel{0}{\to} \uZ\stackrel{\id}{\to} \uZ \stackrel{0}{\to} \uZ \stackrel{\id}{\to}\uZ \stackrel{0}{\to}\dots\ .$$
This implies that $E_2^{q,0}=0$ for $q\ge 1$ and therefore the assertion of the Lemma.

We now turn to $2$. We calculate 
$$H^1(G^q;\uZ)\cong [(\Lambda^* \widehat G)^{\otimes_\Z q}]^1=\underbrace{\widehat G\oplus \dots\oplus \widehat G}_{\mbox{$q$ summands}}\ .$$
By (\ref{uiqwdiqwdwqdw}) we get
$$E_1^{q,1}\cong \underline{\widehat G^q}\ .$$
We see that
$(E_1^{*,1},d_1)$ is the sheafification of the complex of discrete groups
\begin{equation}\label{iqdiodqdwd}
\widehat G^\bullet:0\to \widehat G\to \widehat G^2 \to \widehat G^3 \to \dots \ ,
\end{equation}
where the differential is the dual of the differential of the complex
$$0\leftarrow G \leftarrow G^2\leftarrow G^3\leftarrow \dots$$ induced by the
maps given in 
\ref{iudiqwdwqdwqd}. 
Let us describe the differential more explicitly.
If $\chi\in \widehat G$, and $\mu:G\times G\to G$ is the multiplication map, then we have
$\mu^*\chi=(\chi,\chi)\in \widehat {G\times G}\cong \widehat G\times \widehat G$.
We can write
$\partial :\widehat G^q\to \widehat G^{q+1}$ as
$$\partial=\sum_{i=0}^{q+1} (-1)^i \partial_i\ .$$
Using the formulas of \ref{iudiqwdwqdwqd} we get
$$\partial_i(\chi_1,\dots,\chi_q)=(\chi_1,\dots,\chi_i,\chi_i,\dots,\chi_q)$$
for $i=1,\dots,q$.
Furthermore
$$\partial_0(\chi_1,\dots,\chi_q)=(0,\chi_1,\dots,\chi_q)$$ and
$$\partial_{q+1}(\chi_1,\dots,\chi_q)=(\chi_1,\dots,\chi_q,0)\ .$$
Note that $0=\partial:\widehat G^0\to \widehat G^1$.
We see that we can write the higher ($\ge 1$) degree part of the complex 
(\ref{iqdiodqdwd}) in the form
$$K^\bullet\otimes_\Z \widehat G\ ,$$
where
$K^q=\Z^q$ for $q\ge 1$ and
$\partial:\Z^q\to \Z^{q+1}$ is given by the same formulas as above.
One now shows\footnote{We leave this as an exercise in combinatorics to the interested reader.} that
$$H^q(K^\bullet)=\left\{\begin{array}{cc} \Z&q=1\\0&q\ge 2
                       \end{array}\right\}\ .$$
Since $\widehat G$ is torsion-free and hence a flat $\Z$-module we have
$$H^q(\widehat G^\bullet)\cong H^q(K^\bullet\otimes_\Z \widehat G)\cong H^q(K^\bullet)\otimes_\Z\widehat G\ .$$
Therefore $H^q(\widehat G^\bullet)=0$ for $q\ge 2$. This implies the assertion of the Lemma in the case $2$. 
\hB

\subsubsection{}

\begin{lem}\label{gh14}
Let $F^0\subseteq F^1\subseteq \dots\subseteq F^{k-1}\subseteq F^k=F$ be a filtered sheaf of $\Z_{mult}$-modules such that
$\Gr^l(F)$ has weight $l$, and such that $F^0=F^1=0$.
If $V\subseteq F$ is a torsion-free sheaf of weight $1$, then $V=0$.
\end{lem}
\proof
Assume that $V\not=0$.
We show by induction (downwards) that $F^l\cap V\not=0$ for all $l\ge 1$.
The case $l=1$ gives the contradiction.
Assume that $l>1$.
We consider the exact sequence
$$0\to F^{l-1}\cap V\to F^l\cap V\to \Gr_l(F)\ .$$
First of all, by induction assumption, the sheaf  $V\cap F^l$ is  non-trivial, and a 
as a subsheaf of $V$ it is torsion-free of weight $1$.
Since $\Gr_l(F)$ has weight $l\not= 1$, the map
$V\cap F_l\to \Gr_l(F)$ can not be injective.
Otherwise its image would be a torsion-free sheaf of two different weights $l$ and $1$, and this is impossible by Lemma \ref{jkdqwdwqdwqd}.
Hence
$F^{l-1}\cap V\not=0$. \hB

\subsubsection{}

We now show that Lemma \ref{djekwdkjwqdwq} extends to connected compact groups.
\begin{lem}\label{djekwdkjwqdwq1}
Let  $G$ be a compact connected abelian group.
Then the following assertions are equivalent.
\begin{enumerate}
\item $\uExt_{\Sh_\Ab\bS_{lc}}^i(\uG,\uZ)$ is torsion-free for $i=2,3$
\item $\uExt_{\Sh_\Ab\bS_{lc}}^i(\uG,\uZ)\cong 0$ for $i=2,3$.
\end{enumerate}
\end{lem}
\proof
The non-trivial direction is that 1. implies 2.. Thus let us assume that $\uExt_{\Sh_\Ab\bS_{lc}}^i(\uG,\uZ)$ is torsion-free for $i=2,3$.  
We now look at the spectral sequence $(F_r,d_r)$. The left lower corner of its  second page was calculated in
\ref{udwqidqwdcwqcc}. We see that $F_2^{1,2}=\uExt^1_{\Sh_\Ab\bS_{lc}}((\Lambda^2_\Z\uG)^\sharp,\uZ)$ has weight $2$, while $\uExt^3_{\Sh_\Ab\bS_{lc}}(\uG,\uZ)$ has weight $1$.
Since $\uExt^3_{\Sh_\Ab\bS_{lc}}(\uG,\uZ)$ is torsion-free by assumption, the
differential
$d_2^{1,2}:\uExt^1_{\Sh_\Ab\bS_{lc}}((\Lambda^2_\Z\uG)^\sharp,\uZ)\to 
\uExt^3_{\Sh_\Ab\bS_{lc}}(\uG,\uZ)$ is trivial by Lemma \ref{uifhweifwefwe1}.

We conclude that
$\uExt^2_{\Sh_\Ab\bS_{lc}}(\uG,\uZ)$ and $\uExt^3_{\Sh_\Ab\bS_{lc}}(\uG,\uZ)$ survive to the limit of the spectral sequence.
We see that $\uExt^i_{\Sh_\Ab\bS_{lc}}(\uG,\uZ)$ are torsion-free sub-sheaves  of weight $1$
of $H^{i+1}$ for $i=2,3$. Using the structure of $H^{i+1}$ given in \ref{gh12} in conjunction with the Lemmas \ref{gh122} and \ref{gh14} we get $\uExt^i_{\Sh_\Ab\bS_{lc}}(\uG,\uZ)=0$ for $i=2,3$.
\hB

\subsubsection{}

The combination of Lemma \ref{wdqwidwqdqdqwwwwq} (resp. Lemma \ref{wdqwidwqdqdqwwwwq})
and Lemma \ref{djekwdkjwqdwq1} gives the following result.
\begin{theorem}\label{wqidduiwdwddwduiqwdiu}
\begin{enumerate}
\item 
A compact connected abelian  group $G$  which satisfies the two-three condition is admissible
on the site $\bS_{lc-acyc}$.
\item If $G$ is in addition locally topologically divisible, then it is admissible on $\bS_{lc}$.\end{enumerate}
\end{theorem}

\section{Duality of locally compact group stacks}

\subsection{Pontrjagin Duality}

\subsubsection{}

In this subsection we extend Pontrjagin duality for locally compact groups to abelian
group stacks whose sheaves of objects and automorphisms are represented by locally compact groups.
In algebraic geometry a parallel theory has been considered in \cite{math.AG/0306213}.

The site $\bS$ denotes the site of compactly generated spaces as in \ref{firfore} or one of its sub-sites
$\bS_{lc}$, $\bS_{lc-acyc}$. The reason for considering these sub-sites lies in the fact that certain topological groups
are only admissible on these sub-sites (see the Definitions \ref{iuefefewffwf}, \ref{iuefefewffwf1} and Theorem  \ref{wdwquidioqwdopqwdq}).

\subsubsection{}

Let $F\in \Sh_\Ab \bS$. 
\begin{ddd}\label{dweodwe}
We define the dual sheaf of $F$ by
$$D(F):=\uHom_{\Sh_\Ab \bS}(F,\uT)\ .$$
\end{ddd}

\subsubsection{}

\begin{ddd}\label{fuhwefiwfewfewf}
We call $F$ dualizable, if the canonical evaluation morphism
\begin{equation}\label{zzuzuwedwdw78}c:F\to  D(D(F))\end{equation}
is an isomorphism of sheaves.
\end{ddd}

\begin{lem}\label{hj444cq}
If $G$ is a locally compact group which together with its Pontrjagin dual $\Hom_{\top-\Ab}(G,\T)$ is contained in $\bS$,  then $\uG\in \Sh_\Ab\bS$ is dualizable.
\end{lem}
\proof
By Lemma \ref{eudiwe33} we have  isomorphisms
$$\uHom_{\Sh_\Ab \bS}(\uG,\uT)\cong \underline{\Hom_{\top-\Ab}(G,\T)}$$
and \begin{eqnarray*}\uHom_{\Sh_\Ab \bS}(\uHom_{\Sh_\Ab \bS}(\uG,\uT),\uT)&\cong& \uHom_{\Sh_\Ab \bS}(\underline{\Hom_{\top-\Ab}(G,\T)},\T)\\&\cong& \underline{\Hom_{\top-\Ab}(\Hom_{\top-\Ab}(G,\T),\T)}\ .\end{eqnarray*}
The morphism $c$ in (\ref{zzuzuwedwdw78}) is induced by the
evaluation map
$$G\to  \Hom_{\top-\Ab}(\Hom_{\top-\Ab}(G,\T),\T)$$ which is an isomorphism by the classical Pontrjagin
duality of locally compact abelian groups \cite{MR1397028}, \cite{MR1646190}. \hB

\subsubsection{}

The sheaf of abelian groups $\uT\in \Sh_\Ab\bS$ gives rise to a Picard stack
$\cB\uT\in \cPic(\bS)$ as explained in \ref{diewodewdw}. Recall from \ref{def_of_BF2}
the following  alternative description of $\cB\uT$. Let $\uT[1]$ be the complex
with the only non-trivial entry $\uT[1]^{-1}:=\uT$. Then we have $\cB\uT\cong \ch(\T[1])$ in the notation of \ref{weiofwe}.

\subsubsection{}

Let now $P\in \cPic(\bS)$ be a Picard stack. Recall Definition \ref{lwqdjkqwdjwqq89812}, where we define the internal
$\uHOM$ between two Picard stacks.
\begin{ddd}\label{wdkwqdidqiddqw}
We define the dual stack by
$$D(P):=\uHOM_{\cPic(\bS)}(P,\cB\uT)\ .$$
\end{ddd}
We hope that using the same symbol $D$ for the dual in the case of Picard stacks
and the case of a sheaf of abelian groups will not cause confusion.

\subsubsection{}

This definition is compatible with Definition \ref{dweodwe} of the dual of a sheaf of abelian groups in the following sense.
\begin{lem}\label{idosdc333}
If $F\in \Sh_\Ab\bS$, then we have a natural isomorphism
$$\ch(D(F))\cong D(\cB F)\ .$$
\end{lem}
\proof
First observe that by definition $D(\cB
F)=\uHOM_{\cPic(\bS)}(\ch(F[1]),\ch(\uT[1]))$. 
We use Lemma \ref{zuedzwe777} in order to calculate $H^i(D(\cB F ))$.
We have $$H^{-1}(D(\cB F ))\cong R^{-1}\uHom_{\Sh_\Ab\bS}(F[1],\uT[1])\cong 0$$ and
$$H^0(D(\cB F))\cong R^0\uHom_{\Sh_\Ab\bS}(F[1],\uT[1])\cong \uHom_{\Sh_\Ab\bS}(F,\uT)\cong D(F)\ .$$
The composition of this isomorphism with the projection
$D(\cB F )\to \ch(H^0(D(\cB F)))$ from the stack $D(\cB F)$ onto its sheaf of isomorphism classes (considered as Picard stack) provides
the asserted natural isomorphism. 
\hB

\subsubsection{}

A sheaf of groups $F\in \Sh_\Ab\bS$ can also be considered as a complex $F\in C(\Sh_\Ab\bS)$  with non-trivial entry $F^0:=F$. It thus gives rise to a Picard stack $\ch(F)$.  Recall the definition of an admissible sheaf \ref{iuefefewffwf}.
\begin{lem}\label{zuzede333}
We have natural isomorphisms
$$H^{-1}(D(\ch(F)))\cong D(F)\ ,\quad H^0(D(\ch(F)))\cong \uExt^1_{\Sh_\Ab\bS}(F,\uT)\ .$$ In particular, if $F$ is admissible, then
$D(\ch(F))\cong \cB(D(F))$.
\end{lem}
\proof
We use again Lemma \ref{zuedzwe777}.
We have
$$H^{-1}(D(\ch(F)))\cong R^{-1}\uHom_{\Sh_\Ab\bS}(F,\uT[1])\cong \uHom_{\Sh_\Ab\bS}(F,\uT)\cong D(F)\ .$$
Furthermore,
$$H^0(D(\ch(F)))\cong   R^0\uHom_{\Sh_\Ab\bS}(F,\uT[1])\cong \Ext^1_{\Sh_\Ab\bS}(F,\uT)\ .$$
\hB

\subsubsection{}\label{zudwed4}

Assume now that $F$ is dualizable and admissible (at this point we only need $\uExt^1_{\Sh_\Ab\bS}(F,\uT)\cong 0$).
Then we have 
$$D(D(\ch(F)))\stackrel{Lemma \:\ref{zuzede333}}{\cong}  D(\cB(D(F))) \stackrel{Lemma \:\ref{idosdc333}}{\cong} \ch(D(D(F))\cong \ch(F)\ .$$  
Similarly, if $F$ is dualizable and $D(F)$ admissible (again we only need the weaker condition that $\uExt^1_{\Sh_\Ab\bS}(D(F),\uT)\cong 0$),  then we have 
$$D(D(\cB F)) \stackrel{Lemma \:\ref{idosdc333}}{\cong} D(\ch(D(F)))\stackrel{Lemma \:\ref{zuzede333}}{\cong}  \cB D(D(F))\cong \cB F\ .$$

\subsubsection{}
Let us now formalize this observation.
Let $P\in \cPic(\bS)$ be a Picard stack.
\begin{ddd}
We call $P$ dualizable if the natural evaluation morphism
$P\to D(D(P))$ is an isomorphism.
\end{ddd}
The discussion of \ref{zudwed4} can now be formulated as follows.
\begin{kor}
\begin{enumerate}
\item
If $F$ is dualizable and admissible, then $\ch(F)$ is dualizable.
\item
If $F$ is dualizable and $D(F)$ is admissible, then $\cB F $ is dualizable.
\end{enumerate}
\end{kor}
The goal of the present subsection is to extend this kind of result to more general Picard stacks.

\subsubsection{}

Let $P\in \cPic(\bS)$.
\begin{lem}
We have
$$H^{-1}(D(P))\cong D(H^0(P))$$
and an exact sequence
\begin{equation}\label{zsauxasx7}0\to \uExt^1_{\Sh_\Ab\bS}(H^{0}(P),\uT)\to H^{0}(D(P))\to D(H^{-1}(P))\to
\uExt^2_{\Sh_\Ab\bS}(H^{0}(P),\uT)\ .\end{equation}
\end{lem}
\proof 
By Lemma \ref{etzwwe4} we can choose $K\in C(\Sh_\Ab\bS)$ such that
$P\cong \ch(K)$.  We now get
\begin{eqnarray*}
H^{-1}(D(P))&\stackrel{Lemma \:\ref{zuedzwe777}}{\cong}& R^{-1}\uHom_{\Sh_\Ab\bS}(K,\uT[1])\\&\cong& R^0\uHom_{\Sh_\Ab\bS}(K,\uT)\\&\cong& \uHom_{\Sh_\Ab\bS}(H^0(K),\uT)\\&\cong& D(H^0(P))\ .
\end{eqnarray*}
Again by Lemma \ref{zuedzwe777} we have
$$H^{0}(D(P))\cong R^0\uHom_{\Sh_\Ab\bS}(K,\uT[1])\cong R^1\uHom_{\Sh_\Ab\bS}(K,\uT)\ .$$
In order to calculate this sheaf in terms of the cohomology sheaves $H^i(K)$ of $K$ we choose an injective resolution $\uT\to I$.
Then we have 
$$R\uHom_{\Sh_\Ab\bS}(K,\uT)\cong \uHom_{\Sh_\Ab\bS}(K,I)\ .$$
We must calculate the first total cohomology of this double complex. We first take the cohomology in the $K$-direction, and then in the $I$-direction. 
The second page of the associated spectral sequence has the form
$$\begin{array}{|c|c|c|c|c|}\hline
1&\uExt^0_{\Sh_\Ab\bS}(H^{-1}(P),\uT)&\uExt^1_{\Sh_\Ab\bS}(H^{-1}(P),\uT)&\uExt^2_{\Sh_\Ab\bS}(H^{-1}(P),\uT)&\uExt^3_{\Sh_\Ab\bS}(H^{-1}(P),\uT)\\\hline
0&\uExt^0_{\Sh_\Ab\bS}(H^{0}(P),\uT)&\uExt^1_{\Sh_\Ab\bS}(H^{0}(P),\uT)&\uExt^2_{\Sh_\Ab\bS}(H^{0}(P),\uT)&\uExt^3_{\Sh_\Ab\bS}(H^{0}(P),\uT)\\\hline
&0&1&2&3\\\hline
\end{array}\ .$$
The sequence (\ref{zsauxasx7}) is exactly the edge sequence for the total degree-$1$ term.
\hB

\subsubsection{}

The appearance in (\ref{zsauxasx7}) of the groups
$\Ext^i_{\Sh_\Ab\bS}(H^0(P),\uT)$ for $i=1,2$ was the motivation for the
introduction of the notion of an admissible sheaf in \ref{iuefefewffwf}.
\begin{kor}\label{dweido}
Let $P\in \cPic(\bS)$ be such that $H^0(P)$ is admissible. Then we have $$H^0(D(P))\cong D(H^{-1}(P))\ ,\quad H^{-1}(D(P))\cong D(H^0(P))\ .$$ 
\end{kor}

\subsubsection{}

\begin{theorem}\label{hjjhjad78234}
Let $P\in \cPic(\bS)$ and assume that
\begin{enumerate}
\item $H^0(P)$ and $H^{-1}(P)$ are dualizable
\item $H^0(P)$ and $D(H^{-1}(P))$ are admissible.
\end{enumerate}
Then $P$ is dualizable.
\end{theorem}
\proof
In view of \ref{zudw77823d} it suffices to show
that the evaluation map $c:P\to D(D(P))$ induces isomorphisms
$$H^i(c):H^i(P)\to H^i(D(D(P)))$$
for $i=-1,0$.
Consider first the case $i=0$.
Then by Corollary \ref{dweido} we have an isomorphism
$$h^0:H^0(D(D(P)))\stackrel{\sim}{\to} D(H^{-1}(D(F))\stackrel{\sim}{\to} D(D(H^0(P)))\ .$$
One now checks that the map
$$h^0\circ H^0(c):H^0(P)\to D(D(H^0(P)))$$ is the evaluation map
(\ref{zzuzuwedwdw78}). Since we assume that $H^0(P)$ is dualizable this map
is an isomorphism. Hence $H^0(c)$ is an isomorphism, too.

For $i=-1$ we use the isomorphism
$$h^{-1} :H^{-1}(D(D(P)))\stackrel{\sim}{\to} D(H^{0}(D(P)))\stackrel{\sim}{\to} D(D(H^{-1}(P)))$$
and the fact that $h^{-1}\circ H^{-1}(c):H^{-1}(P)\to D(D(H^{-1}(P)))$ is an isomorphism. \hB

\subsubsection{}

If $G$ is a locally compact group which together with its Pontrjagin dual belongs to $\bS$, then by \ref{hj444cq} the sheaf $\uG$ is dualizable. By Theorem \ref{wdwquidioqwdopqwdq} we know a large class of locally compact groups which are admissible on $\bS$ or at least after restriction to $\bS_{lc}$ or $\bS_{lc-acyc}$.
  
We get the following result.

\begin{theorem}\label{maoisosqiqwdq}
Let $G_0,G_{-1}\in \bS$ be two locally compact abelian groups. We assume that
 their Pontrjagin duals belong to $\bS$, and that $G_0
$ and $DG_{-1}$ are admissible on $\bS$.
 If $P\in \cPic(\bS)$ has $H^i(P)\cong \uG_i$ for $i=-1,0$, then
$P$ is dualizable.
\end{theorem}

Let us specialize to the case which is important for the application to $T$-duality.
Note that a group of the form $\T^n\times \R^m\times F$ for a finitely
generated abelian group $F$ is admissible by 
Theorem \ref{wdwquidioqwdopqwdq}. This class of groups is closed under forming Pontrjagin duals.

\begin{kor}
 If $P\in \cPic(\bS)$ has $H^i(P)\cong \T^{n_i}\times \R^{m_i}\times F_i$ for some finitely generated abelian groups $F_i$ for $i=-1,0$, then
$P$ is dualizable.
\end{kor}

For more general groups one may have to restrict to the  sub-site $\bS_{lc}$ or even to $\bS_{lc-acyc}$.

\subsection{Duality and classification}

\subsubsection{}

Let $A,B\in \Sh_\Ab\bS$. By Lemma \ref{ewhdwejdew8} we know that the
isomorphism classes $\Ext_{\cPic(\bS)}(A,B)$ of Picard stacks with $H^0(P)\cong
B$ and $H^{-1}(P)\cong A$ are classified by a characteristic class  
\begin{equation}\label{hhhs}\phi:\Ext_{\cPic(\bS)}(B,A)\stackrel{\sim}{\to} \Ext^2_{\Sh_\Ab\bS}(B,A)\ .\end{equation}

\subsubsection{}

\begin{lem}\label{dhjdqwhdw672} If $A$ is dualizable and $D(A),B$ are admissible,  then
there is a natural isomorphism $$\cD:\Ext^2_{\Sh_\Ab\bS}(B,A)\stackrel{\sim}{\to} \Ext^2_{\Sh_\Ab\bS}(D(A),D(B))$$ such that
\begin{equation}\label{hdhqwjdqwhjd}
\phi(D(P))=\cD(\phi(P)) \qquad\forall P\in \Ext_{\cPic(\bS)}(A,B)\ .
\end{equation}
\end{lem}
\proof
In order to define $\cD$ we use the identifications
\begin{eqnarray*}\Ext^2_{\Sh_\Ab\bS}(B,A)&\cong& \Hom_{D^+(\Sh_\Ab\bS)}(B,A[2])\ ,\\\ \Ext^2_{\Sh_\Ab\bS}(D(A),D(B))&\cong& \Hom_{D^+(\Sh_\Ab\bS)}(D(A),D(B)[2])\ .\end{eqnarray*}
We choose an injective resolution $\uT\to \cI$. For a complex of sheaves
$F$ we define $RD(F):=\uHom_{\Sh_\Ab\bS}(F,\cI)$. The map $\uT\to \cI$ induces a map
$D(F)\to RD(F)$. Note that $RD$ descends to a functor between derived categories
$RD:D^b(\Sh_\Ab\bS)^{op}\to D^+(\Sh_\Ab\bS)$. 
We now consider the following web of maps
$$\xymatrix{\Hom_{D^+(\Sh_\Ab\bS)}(B,A[2])\ar@{.>}[ddr]^{\hspace*{-0.2cm}\cD}\ar@{.>}[dr]^{\hspace{-0.2cm}\cD^\prime}\ar[r]^{\hspace*{-1cm}RD}&\Hom_{D^+(\Sh_\Ab\bS)}(RD(A[-2]),RD(B))\ar[d]^{D(A)\to RD(A)}\\
&\Hom_{D^+(\Sh_\Ab\bS)}(D(A),RD(B)[2])\\
&\Hom_{D^+(\Sh_\Ab\bS)}(D(A),D(B)[2])\ar[u]_{u}^{D(B)\to RD(B)}}\ .$$
Since $B$ is admissible the map $D(B)\to RD(B)$ is an isomorphism in
cohomology in degree $0,1,2$ (because $D(B)$ is concentrated in degree $0$ and
$H^k(RD(B))=R^k\Hom_{\Sh_\Ab\bS}(B,\uT)=\Ext_{\Sh_\Ab\bS}^k(B,\uT)$). Since $D(A)$ is
acyclic in non-zero degree, the map $u$ is an isomorphism. For this, observe
that $D(B)\to RD(B)$ can be replaced up to quasi-isomorphism by an embedding
$0\to \widetilde{RD(B)}\to RD(B)\to Q\to 0$ such that the quotient is zero in
degrees $0,1,2$. The statement then follows from the long exact Ext-sequence for
$\Hom_{\Sh_\Ab\bS}(D(A),-)$, because $\Ext_{\Sh_\Ab\bS}^i(D(A),B)=0$ for $i=0,1,2$. Therefore the
diagram defines the map
$$\cD:=u^{-1}\circ \cD^\prime:\Hom_{D^+(\Sh_\Ab\bS)}(B,A[2])\to \Hom_{D^+(\Sh_\Ab\bS)}(D(A),D(B)[2])\ .$$

\subsubsection{}

We now show the relation (\ref{hdhqwjdqwhjd}).
By Lemma \ref{etzwwe4} it suffices to show (\ref{hdhqwjdqwhjd}) for $P\in \cPic(\bS)$ of the form $P=\ch(K)$ for complexes
$$\cK:0\to A\to X\to Y\to B\to 0\ , \quad K:0\to X\to Y\to 0\ .$$
As in \ref{gdhasgdaui} we consider
$$\cK_A:0\to X\to Y\to B\to 0$$
with $B$ in degree $0$. Then by Definition (\ref{phidef632423})
of the map $\phi$ and the Yoneda map $Y$,  the element
$\phi(\ch(K))\in \Hom_{D^+(\Sh_\Ab\bS)}(B,A[2])$ is represented by 
the composition
(see (\ref{ygdeduz83e}))
$$Y(\cK)\colon B\stackrel{\beta}{\to}\cK_A\xleftarrow{\alpha} A[2]\ .$$
Since $RD$ preserves quasi-isomorphisms we get $RD(\alpha^{-1})=RD(\alpha)^{-1}$.
It follows that 
$$RD(Y(\cK)):RD(A[2])\stackrel{RD(\alpha)^{-1}}{\to} RD(\cK_A)\stackrel{RD(\beta)}{\to} RD(B)\ .$$ 
We read off that
$$\cD^\prime(Y(\cK)):D(A)\to RD(A)\stackrel{RD(\alpha)^{-1}}{\to} RD(\cK_A)[2]\stackrel{RD(\beta)}{\to} RD(B)[2]\ .$$ 
Let $\uT\to I^0\to I^1\to \dots$ be the injective resolution $\cI$ of $\uT$.
We define 
$J:=\ker(I^1\to I^2)$ and $I:=I^0$. Then we  have $\cB\uT=\ch(L)$ with
$L:0\to I\to J\to 0$ with $J$ in degree zero. Note that $I$ is injective.
Then by Lemma \ref{hashjdjshadui} we have
\begin{equation}\label{gdhfzzuuzwer}D(P)\cong \ch(H)
\ ,\end{equation} where $H:=\tau_{\le 0}\uHom_{\Sh_\Ab\bS}(K,L)$.
Let $$Q:=\ker(\uHom_{\Sh_\Ab\bS}(X,I)\oplus \uHom_{\Sh_\Ab\bS}(Y,J)\to \uHom_{\Sh_\Ab\bS}(X,J))\ .$$
Then $H$ is the complex
$$H:0\to \uHom_{\Sh_\Ab\bS}(Y,I)\to Q\to 0\ .$$
There is a natural map $D(B)\to \uHom_{\Sh_\Ab\bS}(Y,I)$ (induced by
$\uT\to I$ and $Y\to B$), and a projection
$Q\to D(A)$ induced by $A\to X$ and passage to cohomology.
Since $B$ is admissible the complex
$$\cH:0\to D(B)\to \uHom_{\Sh_\Ab\bS}(Y,I)\xrightarrow{d} Q\to D(A)\to 0$$ is
exact. Note that $\ker(d)=H^{-1}\uHom_{\Sh_\Ab\bS}(K,I)$ and
$\coker(d)=H^0\uHom_{\Sh_\Ab\bS}(K,I)$. 
We get 
$$\phi(D(P))\stackrel{(\ref{gdhfzzuuzwer})}{=}\phi(\ch(H))\stackrel{(\ref{phidef632423})}{=}Y(\cH)\stackrel{Lemma\:
  \ref{soundso}}{=}Y^\prime(\cH)\ .$$  
Explicitly, in view of (\ref{gf6f67whrjwr}) the map $Y^\prime(\cH)\in \Hom_{D^+(\Sh_\Ab\bS)}(D(A),D(B)[2])$ is given by the composition
$$Y^\prime(\cH):D(A)\stackrel{\gamma^{-1}}{\to} \cH_{D(A)}\stackrel{\delta}{\to} D(B)[2]\ ,$$
where
$$\cH_{D(A)}: 0\to D(B)\to \uHom_{\Sh_\Ab\bS}(Y,I)\to Q\to 0$$
with $Q$ in degree $0$, the map $\gamma:\cH_{D(A)}\to D(A)$ is the quasi-isomorphism induced by
the projection $Q\to D(A)$, and $\delta:\cH_{D(A)}\to D(B)[2]$ is the
canonical projection. Since $B$ is admissible we have
$R^1\uHom_{\Sh_\Ab\bS}(B,\uT)\cong 0$ and hence a quasi-isomorphism
$\cH_{D(A)}\to \cH_{D(A)}^\prime$ fitting into the following larger diagram 
{$${ \scriptsize\xymatrix{\cH_{D(A)}:\ar[d]^\sim&0\ar[r]&D(B)\ar[d]\ar[r]&\uHom_{\Sh_\Ab\bS}(Y,I)\ar[d]\ar[r]&Q\ar@{=}[d]\ar[r] &0\\
\cH_{D(A)}^\prime:\ar[d]&
0\ar[r]&\uHom_{\Sh_\Ab\bS}(B,I)\ar[d]\ar[r]&\uHom_{\Sh_\Ab\bS}(B,J)\oplus \uHom_{\Sh_\Ab\bS}(Y,I)\ar[d]\ar[r]&Q\ar[r]\ar[d] &0\\
RD(\cK_A):&0\ar[r]&\uHom_{\Sh_\Ab\bS}(B,I)\ar[r]& \uHom_{\Sh_\Ab\bS}(B,I^1)\oplus \uHom_{\Sh_\Ab\bS}(Y,I)\ar[r]&
 \uHom^2_{\Sh_\Ab\bS}(\cK_A,\cI)\ar[r]&\dots}\ .}$$
}
\subsubsection{}

The final step in the verification of (\ref{hdhqwjdqwhjd})
follows from a consideration of the diagram
$$\xymatrix{D(B)[2]\ar[rr]&&RD(B)[2]\\\cH_{D(A)}\ar[d]_\gamma^\sim\ar[u]^\delta\ar[r]^\sim&\cH^\prime_{D(A)}\ar[r]\ar[d]^\sim&RD(\cK_A)\ar[d]_{RD(\alpha)}^\sim\ar[u]^{RD(\beta)}\\ D(A)\ar@/_-1cm/[uu]^{Y^\prime(\cH)}\ar@/^1cm/[uurr]^{\hspace*{3cm}u\circ Y^\prime(\cH)}\ar@{=}[r]&D(A)\ar@/_3cm/[uur]_{\hspace*{-0.9cm}\cD^\prime(Y(\cK))}\ar[r]&RD(A)\ar@/_1cm/[uu]_{RD(Y(\cK))}}
$$
showing the marked equality in 
$$\phi(D(P))=Y^\prime(\cH)\stackrel{!}{=}\cD(Y(\cK))=\cD(\phi(P))\ .$$

It remains to show that 
$$\cD:\Ext^2_{\Sh_\Ab\bS}(B,A)\to \Ext^2_{\Sh_\Ab\bS}(D(A),D(B))$$ 
 is an isomorphism.

To this end we look at the following commutative web of maps
$${\scriptsize\xymatrix{R^2\Hom_{\Sh_\Ab\bS}(D(A),D(B))\ar[d]^u_\cong&R^2\Hom_{\Sh_\Ab\bS}(B,A)\ar@{.>}[l]^{\cD}\ar@{.>}[dl]^{\cD^\prime}\ar@{.>}[ldd]^{RD}\ar@{.>}[rdd]^{z}\ar[r]^\cong&R^2\Hom_{\Sh_\Ab\bS}(B,D(D(A)))\ar[d]_\cong^{v}\\R^2\Hom_{\Sh_\Ab\bS}(D(A),RD(B))\ar[r]^\cong&R^2\Hom_{\Sh_\Ab\bS}(D(A)\otimes^L B,\uT)&R^2\Hom_{\Sh_\Ab\bS}(B,RD(D(A))\ar[l]^\cong\\R^2\Hom_{\Sh_\Ab\bS}(RD(A),RD(B))\ar[u]\ar[r]^\cong&R^2\Hom_{\Sh_\Ab\bS}(RD(A)\otimes^L B,\uT)\ar[u]&R^2\Hom_{\Sh_\Ab\bS}(B,RD(RD(A)))\ar[l]^\cong\ar[u]}\ .}$$
The horizontal isomorphisms in the two lower rows are given by the derived
adjointness of the tensor product and the internal homomorphisms.
The horizontal isomorphism in the first row is induced by the isomorphism
$A\to D(D(A))$. The maps $u$ and $v$ are  isomorphisms since we assume that $D(A)$ is
admissible, compare the corresponding argument in the proof of Lemma
\ref{dhjdqwhdw672}. The map $z$ is induced by the canonical map $A\to
RD(RD(A))$. 
\hB

\section{$T$-duality of twisted torus bundles}

\subsection{Pairs and topological $T$-duality}\label{uqdidwqwqdwqdwqd}

\subsubsection{}
The goal of this subsection is to introduce the main objects of topological $T$-duality and review the structure of the theory.
Let $B$ be a topological space.
\begin{ddd}\label{wedewudwedwedewd}
A pair $(E,H)$ over $B$ consists of a locally trivial principal $\T^n$-bundle
$E\to B$ and a gerbe $H\to E$ with band $\T_{|E}$. An isomorphism of pairs is a diagram
$$\xymatrix{H\ar[d]\ar[r]^\psi&H^\prime\ar[d]\\E\ar[d]\ar[r]^\phi&E^\prime\ar[d]\\B\ar@{=}[r]&B}$$
consisting of an isomorphism of $\T^n$-principal bundles $\phi$ and an isomorphism of $\T$-banded gerbes $\psi$.
By $P(B)$ we denote the set of isomorphism classes of pairs over $B$. If $f:B\to B^\prime$ is a continuous map, then pull-back induces a functorial map
$P(f):P(B^\prime)\to P(B)$.
\end{ddd}

\subsubsection{}

The functor
$$P:\Top^{op}\to \Sets$$ has been introduced in \cite{MR2130624} for $n=1$ and in \cite{math.GT/0501487} for  $n\ge 1$ 
 in connection with the study of topological $T$-duality.

The main result of \cite{MR2130624} in the case $n=1$ is the construction of an involutive natural isomorphism $T:P\to P$, the $T$-duality isomorphism, which associates to each isomorphism class of pairs $t\in P(B)$ the class of the $T$-dual pair $\widehat t:=T(t)\in P(B)$.
 
In the higher-dimensional case $n\ge 2$ the situation is more complicated.
First of all, not every pair $t\in P(B)$ admits a $T$-dual. Furthermore, a $T$-dual, if it exists, may not be unique. The results of \cite{math.GT/0501487} are based on the notion of a $T$-duality triple.
In the following we recall the definition of a $T$-duality triple.

\subsubsection{}\label{wzewqeuwqzdewqd}

As a preparation we recall the following two facts.
Let $\pi:E\to B$ be a $\T^n$-principal bundle.
Associated to the decomposition of the global section functor $\Gamma(E,\dots)$ as  composition $\Gamma(E,\dots)=\Gamma(B,\dots)\circ \pi_*:\Sh_\Ab\bS/E\to \Ab$ there is a decreasing filtration
$$\dots \subseteq F^n H^*(E;\uZ)\subseteq F^{n-1}H^*(E;\uZ)\subseteq \dots\subseteq F^0 H^*(E;\uZ)$$ of the cohomology groups
$H^*(E;\uZ)\cong R^*\Gamma(E;\uZ)$ and a Serre spectral sequence with second term
$E_2^{i,j}\cong H^i(B; R^j\pi_*\uZ)$ which converges to
$\Gr H^*(E;\uZ)$.

\subsubsection{}

The isomorphism classes of gerbes over $X$ with band $\uT_{|X}$ are classified
by $H^2(X;\uT)\cong H^3(X;\uZ)$. The class associated to such a gerbe $H\to X$
is called  Dixmier-Douady class $d(X)\in H^3(X;\uZ)$.

If $H\to X$ is a gerbe with band $\uT_{|X}$ over some space $X$,
then the automorphisms of $H$ (as a gerbe with band $\uT_{|X}$) are classified by
$H^1(X;\uT)\cong H^2(X;\uZ)$.

\subsubsection{}

In the definition of a $T$-duality triple we furthermore need the following notation.
We let $y\in H^1(\T;\uZ)$ be the canonical generator. If $\pr_i:\T^n\to \T$ is the projection onto the $i$th factor, then we set $y_i:=\pr_i^*y\in H^1(\T^n;\uZ)$. Let $E\to B$ be a $\T^n$-principal bundle and $b\in B$. We consider its fibre $E_b$ over $b$. Choosing a base point $e\in E_b$ we use the $\T^n$-action in order to fix a homeomorphism $a_e:E\stackrel{\sim}{\to} \T^n$ such that
$a_{e}(et)=t$ for all $t\in \T^n$. The classes $x_i:=a_e^*(y_i)\in H^1(E_b;\uZ)$ are independent of the choice of the base point. Applying this definition to the bundle $\widehat E\to B$ below gives the classes $\widehat x_i\in H^1(\widehat E_b;\uZ)$ used in Definition \ref{zp00187265e3z12e}.

\begin{ddd}\label{zp00187265e3z12e}
A $T$-duality triple $t:=((E,H),(\widehat E,\widehat H),u)$ over $B$ consists of two
pairs $(E,H), (\widehat E,\widehat H)$  over $B$ and an isomorphism $u:\widehat p^*\widehat H\to p^*H$ of gerbes  with band $\T_{|E\times_B\widehat E}$ defined by the diagram
\begin{equation}\label{ediwedwe7823}
\xymatrix{&p^*H\ar[dl]\ar[dr]&&\widehat p^*\widehat H\ar[ll]^u\ar[dl]\ar[dr]&\\H\ar[dr]&&E\times_B\widehat E\ar[dl]^p\ar[dr]^{\widehat p}&&\widehat H\ar[dl]\\&E\ar[dr]^\pi&&\widehat E\ar[dl]^{\widehat \pi}&\\&&B&&}\ ,
\end{equation}
where all squares are two-cartesian. The following conditions are required:
\begin{enumerate}
\item The Dixmier-Douady classes of the gerbes satisfy $d(H)\in F^2H^3(E;\uZ)$ and
$d(\widehat H)\in F^2H^3(\widehat E;\uZ)$.
\item The isomorphism of gerbes $u$ satisfies the condition
\cite[(2.7)]{math.GT/0501487} which says the following.
If we restrict the diagram (\ref{ediwedwe7823}) to a point $b\in B$, then we can trivialize the restrictions of gerbes $H_{|E_b}$, $\widehat H_{|\widehat E_b}$ to the fibres $E_b,\widehat E_b$ of the $\T^n$-bundles over $b$ such that the induced isomorphism of trivial gerbes
$u_b$ over $E_b\times \widehat E_b$ is classified by $\sum_{i=1}^n \pr_{E_b}^*x_i\cup \pr_{\widehat E}^*\widehat x_i \in H^2(E_b\times \widehat E_b;\uZ)$.
\end{enumerate}
\end{ddd}

\subsubsection{}

There is a natural notion of an isomorphism of $T$-duality triples.
For a map $f:B^\prime\to B$ and a $T$-duality triple over $B$ there is a natural construction of a pull-back triple over $B^\prime$.
\begin{ddd}
We let $\Triple(B)$ denote the set of isomorphism classes of $T$-duality triples over $B$.
For $f:B^\prime\to B$ we let $\Triple(f):\Triple(B)\to \Triple(B^\prime)$ be the map induced by the pull-back of $T$-duality triples.
\end{ddd}

\subsubsection{}\label{ewjdhdwed78283324}
In this way we define a functor
 $$\Triple:\Top^{op}\to \Sets\ .$$  
This functor comes with specializations
$$s,\widehat s:\Triple\to P$$
given by $s((E,H),(\widehat E,\widehat H),u):=(E,H)$ and $\widehat s((E,H),(\widehat E,\widehat H),u)=(\widehat E,\widehat H)$.
\begin{ddd}\label{gwqhsqwsz}
A pair $(E,H)$ is called dualizable if there exists a triple $t\in \Triple(B)$ such that
$s(t)\cong (E,H)$. The pair $\widehat s(t)=(\widehat E,\widehat H)$ is called a $T$-dual of $(E,H)$.
\end{ddd}
Thus the choice of a triple $t\in s^{-1}(E,H)$ encodes the necessary choices in order to fix
a $T$-dual. One of the main results of \cite{math.GT/0501487} is the following
characterization of dualizable pairs.
\begin{theorem}\label{gdhqwdgqwuw}
A pair $(E,H)$ is dualizable in the sense of Definition
\ref{gwqhsqwsz} if and only if $d(H)\in F^2H^3(E;\uZ)$.
\end{theorem} 
Further results of  \cite[Thm. 2.24]{math.GT/0501487} concern the classification of the set of duals of a given pair $(E,H)$. 

\subsubsection{}

For the purpose of the present paper it is more natural to interpret the isomorphism of gerbes
$u:\widehat p^*\widehat H\to p^*H$ in a $T$-duality triple $((E,H),(\widehat E,\widehat H),u)$ in a different, but   equivalent manner. To this end we introduce the notion of a dual gerbe.

\subsubsection{}
First we recall the definition of the tensor product $w:H\otimes_X H^\prime\to X$ of gerbes $u:H\to X$ and $u^\prime:H^\prime\to X$ with band $\uT_{|X}$ over $X\in \bS$.
Consider first the fibre product of stacks $(u,u^\prime):H\times_XH^\prime\to X$. Let $T\in \bS/X$.
An object $s\in H\times_XH^\prime(T)$ is given by a triple $(t,t^\prime,\phi)$ 
of objects $t\in H(T)$, $t^\prime\in H^\prime(T)$ and an isomorphism
$\phi:u(t)\to u^\prime(t^\prime)$. By the definition of a $\T_{|X}$-banded gerbe
the group of automorphisms of $t$ relative to $u$
is the group  $\Aut_{H(T)/\rel(u)}(t)\cong \uT(T)$.  We thus have an  isomorphism $\Aut_{H\times_XH^\prime/\rel((u,u^\prime))}(s)\cong \uT(T)\times \uT(T)$.
Similarly, for $s_0,s_1\in H\times_XH^\prime(T)$ the set
$\Hom_{H\times_XH^\prime/\rel((u,u^\prime))}(s_0,s_1)$ is a torsor over $\uT(T)\times \uT(T)$.

We now define a prestack $H\otimes^p_XH^\prime$.
By definition the groupoid $H\otimes^p_XH^\prime(T)$ has the same objects as $H\times_XH^\prime(T)$, but the morphism sets are factored by the anti-diagonal
$\uT(T)\stackrel{\antidiag}{\subseteq} \uT(T)\times \uT(T)$, i.e. 
$$\Hom_{H\otimes^p_XH^\prime/\rel(w)}(s_0,s_1)=\Hom_{H\times_XH^\prime/\rel((u,u^\prime))}(s_0,s_1)/\antidiag(\uT(T))\ .$$
The stack $H\otimes_XH^\prime$ is defined as the stackification of the prestack
$H\otimes^p_XH^\prime$. It is again a gerbe with band $\uT_{|X}$. We furthermore have the following
relation of Dixmier-Douady classes.
$$d(H)+d(H^\prime)=d(H\otimes_X H^\prime)\ .$$

 \subsubsection{}

The sheaf $\uT_{|X}\in \Sh_\Ab\bS/X$ gives rise to the stack $\cB \uT_{|X}$ (see \ref{jksakasdioioeqweqw}).
Then $X\times  \cB \uT_{|X}\to X$ is the trivial $\T$-banded gerbe over $X$.
Let $H\to X$ be a gerbe with band $\uT_{|X}$. 
\begin{ddd}\label{hdbqwhjdwqduiwqdoiwqdopwqdwqd7}
A dual of the gerbe $H\to X$ is a pair $(H^\prime\to X,\psi)$ of a gerbe $H^\prime\to X$ and an isomorphism of gerbes $\psi:H\otimes_X H^\prime\to X\times \cB\uT_{|X}$.
\end{ddd}
Every gerbe $H\to X$ with band $\uT_{|X}$ admits a preferred dual $H^{op}\to X$ which we call the opposite gerbe.
The underlying stack of $H^{op}$ is $H$, but we change its structure of  a $\uT_{|X}$-banded gerbe using the inversion automorphism ${}^{-1}:\uT_{|X}\stackrel{\sim}{\to}  \uT_{|X}$.

If $(H^\prime_0\to X,\psi_0) $ and $(H^\prime_1\to X,\psi)$ are two duals, then there exists a unique isomorphism class of isomorphisms $H^\prime_0\to H_1^\prime$ of gerbes
such that the induced diagram
$$\xymatrix{H\otimes_X H_0^\prime\ar[dr]^{\psi_0}\ar[rr]&&H\otimes_X H_1^\prime\ar[dl]^{\psi_1}\\&X\times\cB \uT_{|X}&}$$
can be filled with a two-isomorphism.
Note that if $H^\prime\to X$ is underlying gerbe of the dual of $H\to X$, then we have the
relation of Dixmier-Douady classes
$$d(H)=-d(H^\prime)\ .$$

\begin{lem}\label{hedwqdqwdqw66712}
Let $\widehat H\to X$ and $H\to X$  be gerbes with band $\uT_{|X}$. 
There is a natural bijection between the sets of isomorphism classes of isomorphisms $\widehat H\to H$ of gerbes over $X$  and isomorphism classes of isomorphisms of $\uT_{|X}$-banded gerbes $\widehat H\otimes H^{op}\to X\times \cB\uT_{|X}$.
\end{lem}
\proof
An isomorphism
$\widehat H\to H$ induces an isomorphism $\widehat H\otimes_X H^{op}\to H\otimes_XH^{op}\stackrel{\psi}{\to} X\times \cB \uT_{|X}$.  
On the other hand, an isomorphism
$\widehat H\otimes_X H^{op}\to X\times \cB\uT_{|X}$ presents $\widehat H\to X$ as a dual of
$H^{op}\to X$. Since $\psi:H\otimes_XH^{op}\to X\times \cB\uT_{|X}$ presents
$H$ as a dual of $H^{op}$ we have a preferred isomorphism class of
isomorphisms $\widehat H\to H$, too.
\hB

\subsubsection{}\label{hdjdjhwqdddddd6676111}

In view of Lemma \ref{hedwqdqwdqw66712},
in Definition \ref{zp00187265e3z12e} of a $T$-duality triple
 $((E,H),(\widehat E,\widehat H),u)$ we can consider $u$ as an isomorphism
$$\widehat p^*\widehat H^{op}\otimes_X p^* H\to E\times_B\widehat E\times \cB\uT_{|E\times_B\widehat E}\ .$$
The condition  \ref{zp00187265e3z12e} 2. can be rephrased as follows. After restriction to a point
$b\in B$  we can find isomorphisms
\begin{equation}\label{hdjqwdhdwqdwqdw771}
v:\uT^n\times \cB\uT\stackrel{\sim}{\to } H_b\ ,\quad \widehat v: \uT^n\times \cB\uT\stackrel{\sim}{\to} \widehat H_b\ .
\end{equation}
After a choice of $*\to \cB \uT$ in order to define the map $s$ below
we get a map $r:\uT^n\times \uT^n\to \cB \uT$ by the diagram
$$\xymatrix{
\widehat H_b^{op}  \times H_b\ar[r]& \widehat H_b^{op}\otimes_{E_b\times
  \widehat E_b}  H_b\ar[r]^{u_b}&  E_b\times \widehat E_b\times \cB\uT
\ar[r]^{\hspace*{0.9cm}\pr_{\cB \uT}}& \cB \uT\\ 
(\uT^n\times \cB\uT^{op})\times (\uT^n\times \cB\uT)\ar[u]^{(\hat{v},v)}&&&&\\
\uT^n\times \uT^n\ar@/_1cm/[rrruu]^r\ar[u]^s&&&
}\ .$$ 
The condition \ref{zp00187265e3z12e}, 2. is now equivalent to the condition that we can 
choose the isomorphisms $v,\hat v$ in (\ref{hdjqwdhdwqdwqdw771}) such that
$$r^*(z)=\sum_{i=1}^n \pr_{1,i}^*x\cup \pr_{2,i}^* x\ ,$$
where
$z\in H^2(\cB \uT;\uZ)$ and $x\in H^1(\uT;\uZ)$ are the canonical generators, and 
$\pr_{k,i},:\uT^n\times \uT^n\stackrel{\pr_k}{\to} \T^n\stackrel{\pr_i}{\to} \uT$, $k=1,2$, $i=1,\dots,n$ are projections onto the factors.

\subsubsection{}\label{jqswijfoiwefjewwf}

Important topological invariants of $T$-duality triples are the Chern classes of the underlying $\T^n$-principal bundles. For a triple $t=((E,H),(\widehat E,\widehat H),u)$ we define
$$c(t):=c(E)\ ,\quad \widehat c(t):=c(\widehat E)\ .$$
These classes belong to $H^2(B;\uZ^n)$.

\subsection{Torus bundles, torsors and gerbes}\label{hqdwjhqdwqjdwqdqw}

\subsubsection{}

In this subsection we review various interpretations of the notion of a $\T^n$-principal bundle. 

\subsubsection{}
Let $G$ be a topological group.
Let us start with giving a precise definition of a $G$-principal bundle.
\begin{ddd}
A $G$-principal bundle $\cE$ over a space $B$ consists of a map of spaces $\pi:E\to B$ which admits local sections together with a fibrewise right action $E\times G\to E$ such that the natural map
$$E\times G\to E\times_B E\ ,\quad  (e,t)\mapsto (e,et)$$ is an homeomorphism.
An isomorphism of $G$-principal bundles $\cE\to \cE^\prime$ is a $G$-equivariant map $E\to E^\prime$ of spaces over $B$. 
\end{ddd}
By $\Prin_B(G)$ we denote the category of $G$-principal bundles over $B$.
The set $H^0(\Prin_B(G))$ of isomorphism classes in $\Prin_B(G)$ 
is in one-to-one correspondence with homotopy classes $[B,BG]$ of maps from $B$ to the classifying space $BG$ of $G$. If we fix a universal bundle $EG\to BG$,
then the bijection
$$[B,BG]\stackrel{\sim}{\to} H^0(\Prin_B(G))$$
is given by $$[f:B\to BG]\mapsto [B\times_{f,BG} EG\to B]\ .$$

\subsubsection{}

We now specialize to $G:=\T^n$.
The classifying space $B\T^n$ of $\T^n$ has the homotopy type of the Eilenberg-MacLane space $K(\Z^n,2)$. 
We thus have a natural isomorphism
$$H^2(B;\Z^n)\stackrel{Def.}{\cong} [B,K(\Z^n,2)]\cong [B,B\T^n]\cong H^0(\Prin_B(\T^n))\ .$$
Therefore, $\T^n$-principal bundles are classified by the characteristic class
$c(\cE)\in H^2(B;\Z^n)$.
Using the decomposition $$H^2(B;\Z^n)\cong \underbrace{H^2(B;\Z)\oplus\dots H^2(B;\Z)}_{n\:summands}$$
we can write $$c(\cE)=(c_1(\cE),\dots,c_n(\cE))\ .$$
\begin{ddd}\label{dqzwu55}
The class $c(\cE)$ is called the Chern class of $\cE$. The 
$c_i(\cE)$ are called the components of  $c(\cE)$.
\end{ddd}
In fact, if $n=1$, then $c(\cE)$ is the classical first Chern class of the 
$\T$-principal bundle $\cE$.

\subsubsection{} 

Let $\bS$ be a site and $F\in \Sh_\Ab \bS$ be a sheaf of abelian groups.
\begin{ddd}
An $F$-torsor $\cT$ is a sheaf of sets $\cT\in \Sh\bS$ together with an action
$\cT\times F\to \cT$ such that the natural map $\cT\times F\to \cT\times \cT$ is 
an isomorphism of sheaves. An isomorphism of $F$-torsors $\cT\to \cT^\prime$
is an isomorphism of sheaves which commutes with the action of $F$.
\end{ddd}
By $\Torsor(F)$ we denote the category of $F$-torsors.

\subsubsection{}\label{uegdzwdwe}

In \ref{catExt} we have introduced the category
$\EXT(\uZ,F)$ whose objects
are extensions of sheaves of groups
\begin{equation}\label{sgxhas98}\cW:0\to F\to W\stackrel{w}{\to} \uZ\to 0\ ,\end{equation}
and whose morphisms are isomorphisms of extensions.
We have furthermore defined the equivalence of categories
$$U:\EXT(\uZ,F)  \stackrel{\sim}{\to} \Torsor(F)$$
given by $U(\cW):=w^{-1}(1)$.

\subsubsection{}

Consider an extension $\cW$ as in (\ref{sgxhas98}) and apply $\Ext_{\Sh_\Ab\bS}(\uZ,\dots)$.
We get the following piece of the long exact sequence
$$\dots \to \Hom_{\Sh_\Ab\bS}(\uZ,W)\to \Hom_{\Sh_\Ab\bS}(\uZ,\uZ)\stackrel{\delta_\cW}{\to} \Ext^1_{\Sh_\Ab\bS}(\uZ,F)\to \Ext^1_{\Sh_\Ab\bS}(\uZ,W)\to \dots\ .$$
Let $1\in \Hom_{\Sh_\Ab\bS}(\uZ,\uZ)$ be the identity and set
$$e(\cW):=\delta_\cW(1)\in \Ext^1_{\Sh_\Ab\bS}(\uZ,F)\ .$$
Recall that $H^0(\EXT(\uZ,F))$ denotes the set of isomorphism classes of the category $\EXT(\uZ,F)$.
The following Lemma is well-known (see e.g. \cite{MR0225854}).
\begin{lem}
The map
$$\EXT(\uZ,F)\ni \cW\mapsto e(\cW)\in \Ext^1_{\Sh_\Ab\bS}(\uZ,F)$$
induces a bijection
$$e:H^0(\EXT(\uZ,F))\stackrel{\sim}{\to} \Ext^1_{\Sh_\Ab\bS}(\uZ,F)\ .$$
\end{lem}

In view of \ref{uegdzwdwe} we have natural bijections
\begin{equation}\label{eduwid89}H^0(\Torsor(F))\stackrel{U}{\cong} H^0(\EXT(\uZ,F))\stackrel{e}{\cong} \Ext^1_{\Sh_\Ab\bS}(\uZ,F)\ .\end{equation}

\subsubsection{} 

Let $G$ be an abelian topological group and consider  a principal $G$-bundle $\cE$ over $B$ with underlying map $\pi:E\to B$. 
By $\cT(\cE):=\underline{E\to B}\in \Sh\bS/B$ (see \ref{uiqqwq4432}) we denote its sheaf of sections. The right action of $G$ on $E$
induces an action $\cT(\cE)\times \uG_{|B}\to \cT(\cE)$.

\begin{lem}
$\cT(\cE)$ is a $\uG_{|B}$-torsor.
\end{lem}
\proof
This follows from the following fact. If $X\to B$ and $Y\to B$ are two maps, then
$$\underline{X\times_BY\to B}\cong \underline{X\to B}\times \underline{Y\to B}$$ in $\Sh\bS/B$. We apply this to the isomorphism 
$$E\times_B(B\times G)\cong E\times G\cong E\times_B E$$
of spaces over $B$ in order to get the isomorphism
$$\cT(\cE)\times \uG_{|B}\stackrel{\sim}{\to} \cT(\cE)\times\cT(\cE)\ .$$
\hB

\subsubsection{}

The construction $\cE\mapsto \cT(\cE)$ refines to a functor 
$$\cT:\Prin_B(G)\to \Torsor(\uG_{|B})$$
from the category of
$G$-principal bundles $\Prin_B(G)$ over $B$ to the category of $\uG_{|B}$-torsors $\Torsor(\uG_{|B})$ over $B$.

\begin{lem}\label{uefewieuefewnf72}
The functor
$$\cT:\Prin_B(G)\to \Torsor(\uG_{|B})$$
is an equivalence of categories.
\end{lem}
\proof
It is a consequence of the Yoneda Lemma that $\cT$ is an isomorphism on the level of morphism sets. It remains to show that the underlying sheaf $T$ of a
$\uG_{|B}$-torsor is representable by a $G$-principal bundle.  Since $T$ is locally isomorphic to
$\uG_{|B}$ this is true locally. The local representing objects can be glued to a global representing object.
\hB 

We can now prolong the chain of bijections (\ref{eduwid89}) to
\begin{eqnarray}\label{dehzuwdw82}\lefteqn{[B,BG]\cong H^0(\Prin_B(G))\cong H^0(\Torsor(\uG_{|B}))}&&\nonumber\\&&\hspace{1cm}\cong H^0(\EXT(\uZ_{|B},\uG_{|B}))\cong \Ext^1_{\Sh_\Ab\bS/B}(\uZ_{|B},\uG_{|B})\cong H^1(B;\uG)\ .
\end{eqnarray}

\subsubsection{}

Let $\bS$ be some site and $F\in \bS$ be a sheaf of abelian groups. By  $\Gerbe(F)$ we denote the two-category of gerbes with band $F$ over $\bS$.
It is well-known that isomorphism classes of gerbes with band $F$ are classified by $\Ext^2_{\Sh_\Ab\bS}(\uZ;F)$, i.e. there is a natural bijection
$$d:H^0(\Gerbe(F))\stackrel{\sim}{\to} \Ext^2_{\Sh_\Ab\bS}(\uZ;F)\ .$$
 
\subsubsection{}

Let $H\to G$ be a homomorphism of topological abelian groups with kernel
$K:=\ker(H\to G)$. Our main example is $\R^n\to \T^n$ with kernel $\Z^n$. 
\begin{ddd}
Let $T$ be a space.
An $H$-reduction of a $G$-principal bundle $\cE=(E\to B)$ on $T$ is a diagram
$$(F,\phi)\::\: \xymatrix{F\ar[r]^\phi\ar[d]&E\ar[d]\\T\ar[r]& B}\ ,$$
where $F\to T$ is a $H$-principal bundle, and $\phi$ is $H$-equivariant, where
$H$ acts on $E$ via $H\to G$. An isomorphism $(F,\phi)\to
(F^\prime,\phi^\prime)$ of $H$-reductions over $T$ is an isomorphism $f:F\to F^\prime$ of $H$-principal bundles such that $\phi^\prime\circ f=\phi$.
Let $R_H^\cE(T)$ denote the category of $H$-reductions of $\cE$.
\end{ddd}
Observe that the group of automorphisms of every object in
$R_H^\cE(T)$ is isomorphic to $\Map(T,K)^\delta$ (the superscript $\delta$ indicates that we take the underlying set).
For a map $u:T\to T^\prime$ over $B$ there is a natural pull-back functor
$R_H^\cE(u):R^\cE_H(T^\prime)\to R_H^\cE(T)$.

Let $\bS$ be a site as in \ref{firfore} and assume that $K,G,H,B$ belong to $\bS$. In this case it is easy to check that
$$T\mapsto R_H^\cE(T)$$ is a gerbe with band $\uK_{|B}$ on $\bS/B$. 

\subsubsection{}

Let $\pi:E\to B$ be a $G$-principal bundle. Note that the pull-back $\pi^*E\to E$ has a canonical section and is therefore trivialized. A trivialized $G$-bundle has a canonical $H$-reduction.
In other words, there is a canonical map of stacks over $B$
\begin{equation}\label{candef52}
\can:E\to R^\cE_H\ .
\end{equation}
Note that an object of $E(T)$ is a map $T\to E$; to this we assign the
pull-back of the canonical $H$-reduction of $\pi^*E$.

\subsubsection{}

The construction $\Prin_B(G)\ni \cE\mapsto R_H^\cE\in \Gerbe(\uK)$
is functorial in $\cE$ and thus induces a natural map of sets of isomorphism classes
$$r:H^0(\Prin_B(G)) \to H^0(\Gerbe(\uK))\ .$$

\begin{lem}\label{xhasxasxa}
If $H\to G$ is surjective and has local sections, and $H^1(B,\uH)\cong H^2(B,\uH)\cong 0$, then
$$r:H^0(\Prin_B(G)) \to H^0(\Gerbe(\uK))$$
is a bijection.
\end{lem}
\proof
The exact sequence $0\to K\to H\to G\to 0$
induces by \ref{dezuqwideqwd} an exact sequence
$$0\to \uK\to \uH\to \uG\to 0$$ of sheaves.
We consider the following segment of the associated long exact sequence in cohomology:
$$\dots\to H^1(B;\uH)\to H^1(B;\uG)\stackrel{\delta}{\to} H^2(B;\uK)\to H^2(B;\uH)\to\dots\ .$$
By our assumptions $\delta:H^1(B;\uG)\to H^2(B;\uK)$ is an isomorphism.
One can check that the following diagram commutes:
\begin{equation}\label{hjedweduuiiwe}
\xymatrix{H^0(\Prin_B(G))\ar[r]^r\ar[d]&H^0(\Gerbe(\uK_{|B}))\ar[d]^d\\
H^1(B;\uG)\ar[r]^\delta&H^2(B;\uK)}\ .
\end{equation}
This implies the result since the vertical maps are isomorphisms.
\hB

Note that we can apply this Lemma in our main example where $H=\R^n$ and $G=\T^n$.
In this case the diagram (\ref{hjedweduuiiwe}) is  the equality
\begin{equation}\label{hjedweduuiiwe1}
c(E)=d(R^E_{\uZ^n_{|B}})
\end{equation}
in $H^2(B;\uZ^n)$.
\subsection{Pairs and group stacks}

\subsubsection{}\label{uieudiqwd}

Let $\cE$ be a principal $\T^n$-bundle over $B$, or equivalently by (\ref{dehzuwdw82}), an extension $\cE\in \Sh_\Ab\bS/B$
\begin{equation}\label{dhwqjd712}
0\to \uT^n_{|B}\to \cE\to \uZ_{|B}\to 0
\end{equation}
 of sheaves of abelian groups. Let $\tilde c(\cE)\in  \Ext^1_{\Sh_\Ab\bS/B}(\uZ_{|B};\uT^n_{|B})$ be the class of this extension. 
Under the isomorphism 
$$\Ext^1_{\Sh_\Ab\bS/B}(\uZ_{|B};\uT^n_{|B})\cong H^1(B;\uT^n)\cong H^2(B;\uZ^n)$$ it corresponds to the Chern class $c(\cE)$ of the principal $\T^n$-bundle introduced in \ref{dqzwu55}. 

\subsubsection{}

We let $Q_\cE=\Ext_{\cPic(\bS)}(\cE,\uT_{|B})$ (see Lemma \ref{ewhdwejdew8}
for the notation) denote the set of equivalence classes of Picard stacks $P\in
\cPic(B)$ 
with isomorphisms 
$$H^0(P)\xrightarrow{\cong} \cE\ ,\quad H^{-1}(P)\xrightarrow{\cong} \uT_{B}\ .$$ By Lemma \ref{ewhdwejdew8} we have a bijection
$$\Ext^2_{\Sh_\Ab\bS/B}(\cE,\uT_{|B})\cong Q_\cE\ .$$
This bijection induces a group structure on $Q_\cE$ which we will use in the
discussion of long exact sequences below. We will not need a description of
this group structure in terms of the Picard stacks themselves. 

We apply 
$\Ext^*_{\Sh_\Ab\bS/B}(\dots,\uT_{|B})$ to the sequence (\ref{dhwqjd712}) and get the following segment of a long exact sequence
\begin{equation}\label{ewdwie66238p11}\scriptsize
\Ext^1_{\Sh_\Ab\bS/B}(\uT^n_{|B},\uT_{|B})\stackrel{\alpha}{\to}\Ext^2_{\Sh_\Ab\bS/B}(\uZ_{|B},\uT_{|B})\to Q_\cE\to\Ext^2_{\Sh_\Ab\bS/B}(\uT^n_{|B},\uT_{|B})\stackrel{\beta}{\to}\Ext^3_{\Sh_\Ab\bS/B}(\uZ_{|B},\uT_{|B})\ .
\end{equation}
The maps $\alpha,\beta$ are given by the left Yoneda product with the class
$$\tilde c(\cE)\in \Ext^1_{\Sh_\Ab\bS/B}(\uZ_{|B};\uT^n_{|B})\ .$$

\subsubsection{}

In this paragraph we identify the extension groups in the sequence (\ref{ewdwie66238p11}) with sheaf cohomology.
For a site $\bS$ with final object $*$ and $X,Y\in \Sh_\Ab\bS$  we have a local-global spectral sequence with second term
$$E_2^{p,q}\cong R^p\Gamma(*;\uExt^q_{\Sh_\Ab\bS}(X,Y))$$
which converges to $\Ext^{p+q}_{\Sh_\Ab\bS}(X,Y)$.

We apply this first to the sheaves $\uZ_{|B},\uT_{|B}\in \Sh_{\Ab}\bS/B$. The final object of the site $\bS/B$ is $\id:B\to B$.
We have
$\uExt^q_{\Sh_\Ab\bS/B}(\uZ_{|B},\uT_{|B})\cong 0$ for $q\ge 1$ so that this
spectral sequence degenerates at the second page  and gives 
\begin{eqnarray*}\Ext^p_{\Sh_\Ab\bS/B}(\uZ_{|B},\uT_{|B})\cong H^p(B;\uHom_{\Sh_\Ab\bS/B}(\uZ_{|B},\uT_{|B}))\hspace{1cm}&&\\\cong H^p(B;\uT_{|B})\cong H^p(B;\uT)\cong H^{p+1}(B;\uZ)&&\ .\end{eqnarray*}

Since the group $\T^n$ is admissible by Theorem \ref{tqwzdqwmmq6672}, and 
by Corollary \ref{relativecaseall}  we have
$\uExt^q_{\Sh_\Ab\bS/B}(\uT^n_{|B},\uT_{|B})\cong 0$ for $q=1,2$, we get
\begin{eqnarray*}\lefteqn{
\Ext^p_{\Sh_\Ab\bS/B}(\uT^n_{|B},\uT_{|B})\cong H^p(B;\uHom_{\Sh_\Ab\bS/B}(\uT^n_{|B},\uT_{|B}))}&&\\&&\hspace{2cm}\cong H^p(B;\uHom_{\Sh_\Ab\bS}(\uT^n,\uT)_{|B})\cong H^p(B;\uZ^n_{|B})\cong H^p(B;\uZ^n)\end{eqnarray*}
for $p=1,2$.
Therefore the sequence (\ref{ewdwie66238p11}) has the form
\begin{equation}\label{hcjwecwe778823eed}
H^1(B;\uZ^n)\stackrel{\alpha}{\to} H^3(B;\uZ)\to Q_\cE\stackrel{\widehat c}{\to} H^2(B;\uZ^n)\stackrel{\beta}{\to} H^4(B;\uZ)\ .
\end{equation}
In this picture the  maps $\alpha,\beta$ are both given by the cup-product
with the Chern class $c(\cE)\in H^2(B;\uZ^n)$, i.e. $\alpha( (x_i)) = \sum x_i\cup c_i(\cE)$.

\subsubsection{} 
Recall from \ref{wzewqeuwqzdewqd} the decreasing filtration $(F^kH^*(E;\uZ))_{k\ge 0}$ and the spectral sequence 
associated to the decomposition of functors $\Gamma(E;\dots)=\Gamma(B,\dots)\circ p_*:\Sh_{\Ab}\bS/E\to \Ab$.
This spectral sequence converges to $\Gr H^*(E;\uZ)$. Its second page is given by
$E_2^{i,j}\cong H^i(B; R^jp_*\uZ)$.
The edge sequence for
$F^2H^3(E;\uZ)$ has the form
\begin{eqnarray*}\lefteqn{\hspace{-1cm}
\ker(d_2^{0,2}:E_2^{0,2}\to E_2^{2,1})\stackrel{d_3^{1,1}}{\to} \coker(d_2^{1,1}:E_2^{1,1}\to E_2^{3,0})}&&\\&&\to F^2H^3(E;\uZ) \to (E_2^{2,1}/\im(d_2^{0,2}:E_2^{0,2}\to E_2^{2,1}))\stackrel{d_2^{2,1}}{\to} E_2^{4,0}\ .
\end{eqnarray*}

We now make this explicit. Since the fibre of $p$ is an $n$-torus,
$R^0p_*\uZ\cong\uZ$, $R^1p_*\uZ\cong \uZ^n$, $R^2p_*\uZ\cong \Lambda^2\uZ^n$.
The differential $d_2$ can be expressed in terms of the Chern class $c(\cE)$.
We get the following web of exact sequences, where $K$, $A$, $B$ are defined
as the appropriate kernels and cokernels.

\begin{equation}\label{dewidweudewdwed77283}
\xymatrix{&&&K\ar@/_3cm/[dddlll]^=\ar[d]&\\&H^1(B;\uZ^n)\ar[d]^\alpha&&H^0(B;\Lambda^2\uZ^n)\ar[d]^{i_{c(\cE)}}&\\&H^3(B;\uZ)\ar[d]&&H^2(B;\uZ^n)\ar[dr]^\beta\ar[d]&\\K\ar@{.>}[ur]^s\ar[r]&A\ar[r]\ar[d]&F^2H^3(E;\uZ)\ar[r]&B\ar[r]^{\bar \beta}\ar[d]&H^4(B;\uZ)\\&0&&0&}
\end{equation}
Since $K$ as a subgroup of the free abelian group $H^0(B,\Lambda^2\uZ^n)$ is free we can choose a lift $s$ as indicated.

\subsubsection{}\label{updef76}

We now define a map (the underlying pair map)
$$up:Q_\cE\to P(B)$$ as follows.
Any Picard stack $P\in \Pic(\bS/B)$ comes with a natural map of stacks $P\to H^0(P)$ (where on the right-hand side we consider a sheaf of sets as a stack). If $P\in Q_\cE$, then
$H^0(P)\cong \cE$ has a natural map to $\uZ_{|B}$ (see (\ref{dhwqjd712})).
We define the stack $G$ by the pull-back in stacks on $\bS/B$
$$\xymatrix{G\ar@/_0.8cm/[dd]\ar[d]\ar[r]&P\ar@/_-0.8cm/[dd]\ar[d]&\\ E\ar[r]\ar[d]&\cE\ar[d]\\
\underline{\{1\}}_{|B}\ar[r]&\uZ_{|B}&}\ .$$
All squares are two-cartesian, and the  outer square is the composition of the two inner
squares. We have omitted to write the canonical two-isomorphisms.
By construction $E$ is a sheaf of $\uT^n$-torsors (see \ref{uegdzwdwe}), i.e. by Lemma \ref{uefewieuefewnf72} a $\T^n$-principal bundle,
and $G\to E$ is a gerbe with band $\uT$. We set $$up(P):=(E,G)\ .$$

\subsubsection{}
Recall Definition \ref{gwqhsqwsz} of a dualizable pair.
\begin{lem}\label{dgqwhdwqdzwqud5}
If $P\in Q_\cE$, then $up(P)\in P(B)$ is dualizable.
\end{lem}
\proof
Let $(E,H):=up(E)$.
In view of Theorem \ref{gdhqwdgqwuw} we must show that $d(H)\in F^2H^3(E;\uZ)$.
For $k\in \Z$ we define $G(k)\to E(k)$ by the two-cartesian diagrams
\begin{equation}\label{hgt56767344}
\xymatrix{G(k)\ar@/_0.8cm/[dd]\ar[d]\ar[r]&P\ar@/_-0.8cm/[dd]\ar[d]&\\
 E(k)\ar[r]\ar[d]&\cE\ar[d]\\
\underline{\{k\}}_{|B}\ar[r]&\uZ_{|B}&}\ .
\end{equation}
The group structure of $\cE$ induces maps
\begin{equation}\label{wnnwqjdd6723}
\mu:E(k)\times_BE(m)\to E(k+m)
\end{equation}
On fibres appropriately identified with $\T^n$,  this map is the usual
group structure on $\T^n$. 
Since $P$ is a Picard stack these multiplications are covered by 
$\widehat \mu:G(k)\times G(m)\to G(k+m)$.
The isomorphism class
of $G(k)$ therefore must satisfy
\begin{equation}\label{e3e8923e}
\pr_{E(k)}^*G(k)\otimes \pr_{E(m)}^*G(m)\cong \mu^* G(k+m)\in \Gerbe(E(k)\times_BE(m))\ .
\end{equation}
We now write out this isomorphism in terms of Dixmier-Douady classes $d_k:=d(G(k))\in H^3(E(k);\uZ)$. 

We fix a generator of $H^1(\T;\Z)$. This fixes a choice
of generators of $x_i\in H^1(\T^n;\Z)$, $i=1,\dots,n$ via pullback along the
coordinate projections.
Let $a:\T^n\times \T^n\to \T^n$ be the group structure. Then
we have 
\begin{equation}\label{dedhewjde772834e324}
a^*(x_i)=\pr_1^*x_i+  \pr_2^* x_i\ .
\end{equation}
where $\pr_i:\T^n\times \T^n\to \T^n$, $i=1,2$, are the projections onto the factors.

Let us for simplicity assume that $B$ is connected.
Let $(E_r(k),d_r(k))$ be the Serre spectral sequence of the
composition $E(k)\to B \to *$ (see \ref{ewjdhdwed78283324}).
Then we can identify $E_2^{0,3}(k)\cong \Lambda^3_\Z H^1(\T^n;\Z)$.
The class $d_k$ has a symbol in $E_2^{0,3}(k)$ which can be written as 
 $\sum_{i<j<l}a_{i,j,l}(k)x_i\wedge x_j\wedge x_l$. Since the map $\mu$ (see \ref{wnnwqjdd6723}) on fibres can be written in terms of the
group structure we can use (\ref{dedhewjde772834e324}) in order to write out the
symbol of $\mu^*(d_k)$. Equation (\ref{e3e8923e}) implies the identity
\begin{eqnarray*}\lefteqn{
\sum_{i,j,l}a_{i,j,l}(k) \pr_1^*(x_i\wedge x_j\wedge x_l)+\sum_{s,t,u}a_{s,t,u}(m) \pr_2^*(x_s\wedge x_t\wedge x_u)}&&\\&&= \sum_{a,b,c}a_{a,b,c}(k+m)(\pr_1^*x_a+\pr_2^*x_a)\wedge (\pr_1^* x_b+\pr_2^* x_b)\wedge (\pr_1^* x_c+\pr_2^*x_c)\ .
\end{eqnarray*}
Because of the presence of mixed terms this is only possible if everything vanishes.
This implies that $d_k\in F^1 H^3(E(k);\uZ)$. 
We write $\sum_{i,j}x_i\wedge x_j\otimes u_{i,j}(k)$ for its symbol in $E_2^{1,2}$,
where $u_{i,j}\in H^1(B;\uZ)$.
As above equation (\ref{e3e8923e}) now implies
\begin{eqnarray*}\lefteqn{\sum_{i,j}\pr_1^*(x_i\wedge x_j)\otimes u_{i,j}(k)+\sum_{l,k}\pr_2^*(x_l\wedge x_r)\otimes u_{l,r}(m)}&&\\&&=\sum_{a,b}(\pr_1^*x_a+\pr_2^*x_a)\wedge(\pr_1^* x_b+\pr_2^*x_b)\otimes u_{a,b}(k+m)\ .
\end{eqnarray*}
Again the presence of mixed terms implies that everything vanishes.
This shows that $d_k\in F^2 H^3(E(k);\uZ)$.
The assertion of the Lemma is the case $k=1$. \hB

\subsubsection{}

Let us now combine (\ref{dewidweudewdwed77283}) and (\ref{hcjwecwe778823eed}) into a single diagram. We get the following web of horizontal and vertical exact sequences:
\begin{equation}\label{gdhgqwdqw691990290}
\xymatrix{&H^1(B;\Z^n)\ar@{=}[dl]\ar[d]^\alpha&&&\\H^1(B;\uZ^n)\ar[r]&H^3(B;\uZ)\ar[r]^h\ar[d]&Q_\cE\ar[d]^f\ar[r]^{\widehat c}&H^2(B;\uZ^n)\ar[d]\ar[r]^\beta&H^4(B;\uZ)\ar@{=}[d]\\
K\ar[r]\ar@{.>}[ur]^s&A\ar[r]\ar[d]&F^2H^3(E;\uZ)\ar[d]\ar[r]&B\ar[d]\ar[r]^{\bar \beta}&H^4(B;\uZ)\\&0&0&0&}\ .
\end{equation}
The map $f:Q_\cE\to F^2H^3(E;\uZ)$ associates
to $P\in Q_\cE$ the Dixmier-Douady class of the gerbe of the pair $up(P)\in P(B)$ with underlying $\T^n$-bundle $E$. Here we use Lemma \ref{dgqwhdwqdzwqud5}.
Surjectivity of $f$ follows by a diagram chase once we have shown the following Lemma.
\begin{lem}
The diagram (\ref{gdhgqwdqw691990290}) commutes.
\end{lem}
\proof
We have to check that the left and the right squares
\begin{equation}\label{hdjwhqdqwwqd77}
\xymatrix{H^3(B;\uZ)\ar[d]\ar[r]^h&Q_\cE\ar[d]^f\\A\ar[r]&F^2H^3(E;\uZ)}\ ,\quad 
\xymatrix{Q_\cE\ar[r]^{\widehat c}\ar[d]^f&H^2(B;\uZ^n)\ar[d]\\F^2H^3(E;\uZ)\ar[r]^e&B}
\end{equation} commute.
Let us start with the left square. Let $d\in H^3(B;\uZ)$. Under the identification
$$H^3(B;\uZ)\cong H^2(B;\uT)\cong \Ext^2_{\Sh_\Ab\bS/B}(\uZ;\uT)$$ it corresponds to a group stack
$P\in \Pic(\bS/B)$ with $H^0(P)\cong \uZ$ and $H^{-1}(P)\cong \uT$.
The group stack $h(d)\in Q_\cE$ is given by the pull-back (two-cartesian diagram)
$$\xymatrix{h(d)\ar[d]\ar[r]&\cE\ar[d]\\P\ar[r]&\uZ}\ .$$
In particular we see that the  gerbe  $H\to E$ of the pair $u(h(d))=:(E,H)$ is given by a pull-back
$$\xymatrix{H\ar[d]\ar[r]&E\ar[d]\\G\ar[r]&B}\ ,$$
where $G\in \Gerbe(B)$ is a gerbe with Dixmier-Douady class $d(G)=d$.
The composition $f\circ h:H^3(B;\uZ)\to H^3(E;\uZ)$
is thus given by the pull-back along the map $p:E\to B$, i.e. $p^*=f\circ h$.
By construction the composition $H^3(B;\uZ)\to A \to F^2H^3(E;\uZ)$ is a certain factorization of $p^*$. This shows that the left square commutes.

Now we show that the right square in (\ref{hdjwhqdqwwqd77}) commutes.
We start with an explicit description of $\widehat c$. Let $P\in Q_\cE$. The principal $\T^n$-bundle in
$G(0)\to E(0)\to B$ is trivial (see (\ref{hgt56767344}) for notation).  Therefore the
Serre spectral sequence $(E_r(0),d_r(0))$ degenerates at the second term.
We already know by Lemma \ref{dgqwhdwqdzwqud5} that the Dixmier-Douady class of $G(0)\to E(0)$ satisfies  $d_0\in F^2H^3(E(0);\uZ)$.
Its symbol can be written as
$\sum_{i=1}^n x_i\otimes \widehat c_i(P)$ for a uniquely determined sequence
$\widehat c_i(P)\in H^2(B;\uZ)$. These classes constitute the components of the class
$\widehat c(P)\in H^2(B;\uZ^n)$.  

We write the symbol of $f(P)=d_1$  as $\sum_{i} x_i\otimes a_i$ for a sequence $a_i\in H^2(B;\uZ)$. As in the proof of Lemma \ref{dgqwhdwqdzwqud5} the equation   (\ref{e3e8923e}) gives the identity
$$\sum_{i=1}^n \pr_2^* x_i\otimes a_i + \sum_{i=1}^n \pr_1^* x_i\otimes \widehat c_i(P)\equiv \sum_{i=1}^n (\pr_1^* x_i+\pr_2^* x_i)\otimes a_i$$
modulo the image of a second differential $d_2^{0,2}(1)$.
This relation is solved by $a_i:=\widehat c_i(P)$ and determines the image of the vector
$a:=(a_1,\dots,a_n)$  under $H^2(B;\uZ^n)\to H^2(B;\uZ^n)/\im(d_2^{0,2}(1))=:B$ uniquely.
Note that $e\circ f(P)$ is also represented by the image of the vector
$a$ in $B$. This shows that the right square in (\ref{hdjwhqdqwwqd77}) commutes.
\hB 

\subsection{$T$-duality triples and group stacks}\label{hdjhdqjdwdqwdqwdwqdwqdwqd54433}

\subsubsection{}
Let $R_n$ be the classifying space of $T$-duality triples introduced in \cite{math.GT/0501487}. It carries a universal $T$-duality triple
$t_{univ}:=((E_{univ},H_{univ}),(\widehat E_{univ},\widehat H_{univ}),u_{univ})$. 
Let $c_{univ},\widehat c_{univ}\in H^2(R_n;\Z^n)$ be the Chern classes of the bundles $E_{univ}\to R_n$, $\widehat E_{univ}\to R_n$. They satisfy the relation $c_{univ}\cup \widehat c_{univ}=0$.
Let $\cE_{univ}$ be the extension of sheaves corresponding to $E_{univ}\to R_n$ as  in \ref{uieudiqwd}.  In \cite{math.GT/0501487} we have shown that $H^3(R_n;\uZ)\cong 0$.
The diagram (\ref{gdhgqwdqw691990290}) now implies that there is a unique
Picard stack $P_{univ}\in Q_{\cE_{univ}}$ with
$\widehat c(P_{univ})=\widehat c_{univ}$ and underlying pair $up(P_{univ})\cong (E_{univ},H_{univ})$ (see \ref{updef76}).

\subsubsection{}

Let us fix a $\T^n$-principal bundle $E\to B$, or the corresponding extension of sheaves $\cE$. Let us furthermore fix a class $h\in F^2H^3(E;\uZ)$. In \cite{math.GT/0501487} we have introduced the set
$\Ext(E,h)$ of extensions of $(E,h)$ to a $T$-duality triple.
The main theorem about this extension set is \cite[Theorem 2.24]{math.GT/0501487}.

Analogously,  in the present paper we can consider the set of extensions of $(E,h)$ to a Picard stack $P$ with underlying $\T^n$-bundle $E\to B$ and $f(P)=h$, where $f:Q_\cE\to F^2H^3(E;\uZ)$ is as in (\ref{gdhgqwdqw691990290}). In symbols we can write  $f^{-1}(h)$ for this set.

The main goal of the present paper is to compare the sets $\Ext(E,h)$ and $f^{-1}(h)\subseteq Q_\cE$. In the following paragraphs we construct maps between these sets. 

\subsubsection{}\label{dwqzduwqdwqd663}

We fix a $\T^n$-principal bundle $E\to B$ and let $\cE\in \Sh_\Ab\bS/B$ be the corresponding extensions of sheaves. Let $\Triple_E(B)$ denote the set of isomorphism classes of triples $t$ such that $c(t)=c(E)$.
We first define a map
$$\Phi:\Triple_E(B)\to Q_\cE\ .$$
Let $t\in \Triple_E(B)$ be a triple which is classified by a map
$f_t:B\to R_n$. Pulling back the group stack $P_{univ}\in \Pic(\bS/R_n)$ we get an element $\Phi(t):=f_t^*(P_{univ})\in \Pic(\bS/B)$. If $t\in  \Triple_E(B)$, then  we have $\Phi(t)\in Q_\cE$ .
We further have
\begin{equation}\label{gghze673}
\widehat c (\Phi(t))=f_t^*\widehat c(P_{univ})=f^*_t\widehat c_{univ}= \widehat c(t)\ .
\end{equation}

\subsubsection{}\label{qudiuwqdiuwduwqdiwqd667}

In the next few paragraphs we  describe a map $$\Psi:Q_{\cE}\to \Triple_E(B)\ ,$$ i.e. a  construction of a $T$-duality triple $\Psi(P)\in \Triple_E(B)$ starting from a Picard stack
 $P\in Q_\cE$. 

\subsubsection{} \label{dcgbsasjkscsklioiowqdqwdw}
Consider $P\in \cPic(\bS/B)$. We have already constructed one pair $(E,H)$.
The dual $D(P):=\uHOM_{\cPic(\bS/B)}(P,\cB\uT_{|B})$
is a Picard stack with (see Corollary \ref{dweido})
\begin{eqnarray*}
H^0(D(P))\cong D(H^{-1}(P))\cong D(\uT_{|B})\cong \uZ_{|B}\\
H^{-1}(D(P))\cong D(H^0(P))\cong D(\cE)\ .
\end{eqnarray*}
In view of the structure (\ref{dhwqjd712}) of $\cE$,
the equalities $D(\uT^n_{|B})\cong \uZ^n_{|B}$, $D(\uZ_{|B})\cong \uT_{|B}$, and $\uExt^1_{\Sh_\Ab \bS/B}(\uT_{|B},\uT_{|B})\cong 0$ we have an exact sequence
\begin{equation}\label{hdwejdhqzu66623}
0\to \uT_{|B}\to D(\cE)\to \uZ^n_{|B}\to 0\ .
\end{equation}
Using the construction \ref{dedwhedjewdn66}
we can form a quotient $\overline{D(P)}$ which fits into the sequence of maps of Picard stacks
$$ \cB\uT_{|B}\to D(P)\to \overline{D(P)}\ ,$$
where $H^0(\overline{D(P)})\cong \uZ_{|B}$ and $H^{-1}(\overline{D(P)})\cong \uZ^n_{|B}$.
The fibre product
$$\xymatrix{R\ar[r]\ar[d]&\overline{D(P)}\ar[d]\\
\{1\}\ar[r]&\uZ_{|B}}$$
defines a gerbe $R\to B$ with band $\uZ^n_{|B}$.

\subsubsection{}

By Lemma \ref{xhasxasxa}
there exists a unique isomorphism class $\widehat E\to B$ of a $\T^n$-bundle
whose $\uZ^n_{|B}$-gerbe of $\R^n$-reductions 
$R_{\Z^n}^{\widehat E}$ is isomorphic to $R$. We fix such an isomorphism and obtain a canonical map of stacks $\can$ (see (\ref{candef52})) fitting into the diagram
\begin{equation}\label{fggdsad566677612e21}
\xymatrix{\widehat H^{op}\ar[rr]\ar[d]&&\tilde{\widehat H}\ar[d]\ar[r]&D(P)\ar[d]\\\widehat E\ar[r]^{\can}\ar[dr]& R_{\Z^n}^{\widehat E}\ar[r]^\cong\ar[d]& \ar[dl] R\ar[r]\ar[d]&\overline{D(P)}\ar[d]\\&B&\{1\}\ar[r]&\uZ_{|B}}\ .
\end{equation}
The gerbes $\widehat H^{op}\to \widehat E$ and $ \tilde{\widehat H}\to R$ are defined such that the squares
become two-cartesian (we omit to write the two-isomorphisms). 
The squares are cartesian, and the upper square is the definition of the gerbe $\widehat H\to \widehat E$.
In this way the Picard stack $P\in Q_\cE$ defines the second pair $(\widehat E,\widehat H)$
of the triple $\Psi(P)=((E,H),(\widehat E,\widehat H),u)$ whose construction has to be completed by providing $u$.

\subsubsection{}

\begin{lem}\label{hewdjewhdf89}
We have the equality $\widehat c(P)=\widehat c(\Psi(P))$ in $H^2(B;\uZ_{|B}^n)$.
\end{lem}
\proof
By the definition in \ref{jqswijfoiwefjewwf} we have $\widehat  c(\Psi(P))= c(\widehat E)$.  Furthermore,
by  (\ref{hjedweduuiiwe1}) we have
$c(\widehat E)=d(R^{\widehat E}_{\uZ_{|B}^n})=d(R)$.
By Lemma \ref{dhjdqwhdw672} we have
$\phi(D(P))=\cD(\phi(P))$, where 
$\phi$ is the characteristic class (\ref{hhhs}), and 
$$\cD:Q_\cE\cong \Ext^2_{\Sh_\Ab\bS/B}(\cE,\uT_{|B})\stackrel{\sim}{\to} \Ext^2_{\Sh_\Ab\bS/B}(D(\uT_{|B}),D(\cE))$$ is as in Lemma \ref{dhjdqwhdw672}.
The map $\widehat c:Q_\cE\to H^2(B;\uZ)$ by its definition fits into the
diagram
 $$\xymatrix{Q_\cE\ar@{.>}[rdd]^a\ar[dddr]^{\widehat c}\ar[r]^{\uT^n_{|B}\to
     \cE}\ar[d]^D&\Ext^2_{\Sh_\Ab\bS/B}(\uT^n_{|B};\uT_{|B})\ar[d]^{\cD}\\\Ext^2_{\Sh_\Ab\bS/B}(D(\uT_{|B}),D(\cE))\ar[r]^{D(\cE)\to D(\uT^n_{|B})}&\Ext^2_{\Sh_\Ab\bS/B}(D(\uT_{|B});D(\uT^n_{|B}))\ar[d]^\cong \\&\Ext^2_{\Sh_\Ab\bS/B}(\uZ_{|B},\uZ^n_{|B})\ar[d]^\cong\\&H^2(B;\uZ^n_{|B})}\ .$$ 
By construction the class $a(P)\in \Ext^2_{\Sh_\Ab\bS/B}(\uZ_{|B},\uZ^n_{|B})$
classifies the Picard stack $\overline{D(\cE)}$. This implies that
$\widehat c(P)$ classifies the gerbe $R\to B$\ , 
i.e. $\widehat c(P)=c(\widehat E)$. \hB 

\subsubsection{}\label{zdquwdzwq}

It remains to construct the last entry
\begin{equation}\label{ttztsw672}
u:H\otimes_{E\times_B\widehat E}\widehat H^{op}\to \cB\uT_{|E\times_B\widehat E}
\end{equation} of the triple
$\Psi(t)=((E,H),(\widehat E,\widehat H),u)$, where we use Picture \ref{hdjdjhwqdddddd6676111}.
Note that $H^{-1}(P)\cong \uT_{|B}$. Construction \ref{dedwhedjewdn66} gives rise to 
a sequence of morphisms of Picard stacks
$$\cB \uT_{|B}\to  P\to \bar P$$
(we could write $\bar P=\cE$).

We have an evaluation $\ev:P\times_B D(P)\to B\T_{|B}$.  For a pair $(r,s)\in \Z\times \Z$  we consider the following diagram with two-cartesian squares.
$$\xymatrix{H_r\times_B \tilde{\widehat H}_s\ar[r]\ar[d]&P\times D(P)\ar[d]\ar[r]^{\ev}&\cB\uT_{|B}\\
H_r\otimes_{E_r\times_B \widehat R_s}\tilde{\widehat H}_s\ar@{.>}[rru]^{u_{r,s}}\ar[r]\ar[d]^{w}&P\otimes_{\bar P\times \overline{D(P)}} D(P)\ar[d]&\\
E_r\times_B R_s\ar[r]\ar[d]&\bar P\times \overline{D(P)}\ar[d]&\\
\underline{\{r,s\}}_{|B}\ar[r]&\uZ_{|B}\times \uZ_{|B}&}\ .$$
By an inspection of the definitions one checks that the natural factorization $u_{r,s}$ exists if $r=s$.
Furthermore one checks that
$$(w,u_{r,s}):H_r\otimes_{E_r\times_B R_s}\tilde{\widehat H}_s\to (E_r\times_B R_s)\times \cB\uT_{|B}$$
is an isomorphism of gerbes with band $\uT_{|E_r\times_B R_s}$ if and only if
$r=s=1$. Both statements can be checked already when restricting to a point,
and therefore become clear when considering the argument in the proof of Lemma \ref{tdultrip}.

We define the map (\ref{ttztsw672}) by
$$\xymatrix{H\otimes_{E\times_B\widehat E}\widehat H^{op}\ar@/_-1cm/[rr]^u\ar[r]\ar[d]& H_1\otimes_{E_1\times_B \widehat R_1}\tilde{\widehat H}_1 \ar[d]  \ar[r]^{\hspace*{0.9cm}u_{1,1}}&\cB\uT_{|B}\\E\times_B\widehat E\ar[r]^{(\id,\can)}&E_1\times_B\widehat R_1&}$$
(note that $R_1=R$, $E_1=E$ and $H_1=H$).

\subsubsection{}

\begin{lem}\label{tdultrip}
The triple $\Psi(P)=((E,H),(\widehat E,\widehat H),u)$ constructed above is a $T$-duality triple.
\end{lem}
\proof
It remains to show that the isomorphism of gerbes $u$ satisfies the condition \ref{zp00187265e3z12e}, 2 in the version of \ref{hdjdjhwqdddddd6676111}. 

By naturality of the construction $\Psi$ in the base $B$ and the fact that
condition \ref{zp00187265e3z12e}, 2 can be checked at a single point $b\in B$, we can assume without loss of generality that $B$ is a point. 
We can further assume that $\cE=\uZ\times \uT^n$ and
$P=\cB\uT\times \uZ\times \uT^n$.
In this case
$$D(P)=D(\cB\uT)\times D(\uZ)\times D(\uT^n)\cong \uZ\times \cB\uT\times \cB\uZ^n\ .$$
We have
$$H\times \tilde{\widehat H}\cong (\cB\uT\times \uT^n)\times (\cB\uT\times \cB\uZ^n) \ .$$
The restriction of the evaluation map $H\times \tilde{\widehat H}\to H\otimes \tilde{\widehat H}\to \cB\uT$
 is the composition
$$ (\cB\uT\times \uT^n)\times (\cB\uT\times \cB\uZ^n) \cong  \cB\uT\times (\cB\uT)\times \uT^n\times \cB\uZ^n\stackrel{\id\times\id\times \ev}{\to}\cB\uT\times \cB\uT \times\cB\uT \stackrel{\Sigma}{\to} \cB\uT\ .$$
We are interested in the contribution
$\ev:\uT^n\times \cB\uZ^n\to \cB\T^n$.

It suffices to see that  in the case $n=1$ we have   $$\ev^*(z)=x\otimes y\ ,$$ where $x\in H^1(\T;\uZ)$,
$y\in H^1(\cB\Z;\uZ)$, and 
$z\in H^2(\cB T)\cong \Z$ are the canonical generators.
In fact this implies via the K{\"u}nneth formula that
$\ev^*(z)=\sum_{i=1}^n x_i\otimes y_i$,
where $x_i:=p_i^*(x)$, $y_i:=q_i^*(y)$ for the projections onto the components
$p_i:\uT^n\to \uT$ and $q_i:\cB \uZ^n\cong (\cB\uZ)^n\to \cB\uZ$.
Finally we use that
$\can^*(y_i)= x_i$, where $\can: \uT^n\to \cB\uZ^n$ is the canonical map  (\ref{candef52}) from a second copy of the torus to its gerbe of  $\R^n$-reductions (after identification of this gerbe with $\cB\uZ^n$). 
\hB

This finishes the construction of $\Psi$ which started in \ref{qudiuwqdiuwduwqdiwqd667}.

\subsubsection{}
 In \cite[2.11]{math.GT/0501487} we have seen that
the group $H^3(B;\uZ)$ acts on $\Triple(B)$ preserving the subsets $\Triple_E(B)\subseteq \Triple(B)$ for every $\T^n$-bundle $E\to B$. We will recall the description of the action in the proof of Lemma \ref{ghhwezfuwf78223} below. By (\ref{gdhgqwdqw691990290}) it also acts on $Q_\cE$.

\begin{lem}\label{ghhwezfuwf78223}
The map $\Psi$ is $H^3(B;\uZ)$-equivariant.
\end{lem}

The proof requires some preparations.

\subsubsection{}

Note that we have a canonical isomorphism $\uT\cong D(\uZ)$.
In order to work with canonical identifications we are going to use $D(\uZ)$ instead of $\uT$.

The isomorphisms classes of gerbes with band $D(\uZ)$ over a space $B\in \bS$ are
classified by 
\begin{equation}\label{ffffwsqwjwqswqsq}
H^2(B;D(\uZ)_{|B})\cong \Ext^2_{\Sh_\Ab\bS/B}(\uZ_{|B},D(\uZ)_{|B})\ .
\end{equation}
The latter group also classifies Picard stacks $P$ with fixed isomorphisms
$H^0(P)\cong \uZ_{|B}$ and $H^{-1}(P)\cong D(\uZ_{|B})$.

Given $P$ the gerbe $G\to B$ can be reconstructed
as a pull-back
$$\xymatrix{G\ar[d]\ar[r]&P\ar[d]\\\underline{\{1\}}_{|B}\ar[r]&\uZ_{|B}}\ .$$
If we want to stress the dependence of $G$ on $P$ we will write $G(P)$.

Recall that an object in $P(T)$ as a stack over $\bS$ consists of a map $T\to B$ and an object of
$P(T\to B)$. In a similar manner we interpret morphisms.

For example, the stack $\underline{{\{1\}}}_{|B}$ is the space $B$.

\subsubsection{}

Let $P\in \Ext_{\cPic(\bS/B)}(\uZ_{|B},D(\uZ)_{|B})$ be a Picard stack $P$ with a fixed isomorphisms
$H^0(P)\cong \uZ_{|B}$ and $H^{-1}(P)\cong D(\uZ_{|B})$.
Let $\phi: \Ext_{\cPic(\bS/B)}(\uZ_{|B},D(\uZ)_{|B})\to  \Ext^2_{\Sh_\Ab\bS/B}(\uZ_{|B},D(\uZ)_{|B})$
be the characteristic class.
\begin{lem}\label{iudwqodhqwdqwdqwd}
$\phi(D(P))\cong -\phi(P)$.
\end{lem}
\proof
First of all note that we have canonical isomorphisms
$$H^0(D(P))\cong D(H^{-1}(P)\cong D(D(\uZ_{|B}))\cong \uZ_{|B}$$ and
$$H^{-1}(D(P))\cong D(H^{0}(P))\cong D(\uZ_{|B})\cong D(\uZ)_{|B}\ .$$
Therefore we can consider
$D(P)\in  \Ext_{\cPic(\bS/B)}(\uZ_{|B},D(\uZ)_{|B})$ in a canonical way,
and
$\phi(P)$ and $\phi(D(P))$ belong to the same group.

Note that
$$d(G(P))=\phi(P)\ ,d(G(D(P)))=\phi(D(P))$$
under the isomorphism (\ref{ffffwsqwjwqswqsq}).
It suffices to show that
$d(G(P))=-d(G(D(P)))$.

In fact, as in \ref{zdquwdzwq} we have the following factorization of the evaluation map

\bigskip
$$\xymatrix{G(P)
\otimes_BG(D(P))\ar[d]\ar[r]&P\otimes_{\uZ_{|B}}D(P)\ar@/^1cm/[rrr]\ar[d]&P\times_{\uZ_{|B}}D(P)\ar[d]\ar[r]\ar[l]&P\times_BD(P)\ar[r]^{ev}\ar[d]&\uT_{|B}\\\underline{\{1\}}_{|B}\ar[r]
&\uZ_{|B}\ar@{=}[r]&\uZ_{|B}\ar[r]^{\diag}& \uZ_{|B}\times \uZ_{|B}& \\}
$$
This represents $G(D(P))$ as the dual gerbe $G(P)^{op}$ of $G(P)$ in the sense of 
Definition \ref{hdbqwhjdwqduiwqdoiwqdopwqdwqd7}.
The relation $d(G(P))=-d(G(D(P)))$ follows. \hB 

\subsubsection{} We now start the actual proof of Lemma \ref{iudwqodhqwdqwdqwd}.

In order to see that $\Psi$ is $H^3(B;\uZ)$-equivariant we will first describe the action of $H^3(B;\uZ)$ on the sets of isomorphism classes of $T$-duality triples with fixed underlying $\T^n$-bundles $E$ and $\widehat E$ on the one hand, and on the set of isomorphism classes of Picard stacks $Q_\cE$, on the other.

Consider $g\in H^3(B;\uZ)$. It classifies the isomorphism class of a gerbe $G\to B$ with band $\uT_{|B}$. If $t$ is represented by $((E,H),(\widehat E,\widehat H),u)$, then
$g+t$ is represented by $$((E,H\otimes\pi^*G),(\widehat E,\widehat H\otimes \widehat \pi^*G), u\otimes \id_{r^*G})\ , $$ where the maps $\pi,\widehat \pi,r$ are as in (\ref{ediwedwe7823}).

 Note the isomorphism
 $H^3(B;\uZ)\cong H^2(B;\uT)\cong \Ext^2_{\Sh_\Ab\bS/B}(\uZ_{|B},\uT_{|B})$. Therefore the class $g$ also classifies an isomorphism class of Picard stacks
 with $H^0(\tilde G)\cong \uZ_{|B}$ and $H^{-1}(\tilde G)\cong \uT_{|B}$. From $\tilde G$ we can derive the gerbe $G$ by the pull-back
$$\xymatrix{G\ar[r]\ar[d]&\tilde G\ar[d]\\
  \underline{\{1\}}_{|B}\ar[r]&\uZ_{|B}}\ .$$

\subsubsection{} 

Recall that we consider an extension $\cE\in \Sh_\Ab\bS/B$ of the form
$$0\to \uT_{|B}^n\to \cE\to \uZ_{|B}\to 0\ .$$
We consider a Picard stack
$P\in \Ext_{\cPic(\bS/B)}(\cE,D(\uZ)_{|B})$.

Let furthermore
$\tilde G\in  \Ext_{\cPic(\bS/B)}(\uZ_{|B},D(\uZ)_{|B})$.
Then we define
$$P\otimes_{\cE}\tilde G\in \Ext_{\cPic(\bS/B)}(\cE,D(\uZ)_{|B})$$ by
the diagram
\begin{equation}\label{qjhwdqwdqwdwqd}
\xymatrix{\cB\uT_{|B}\ar[d]&\cB\uT_{|B}\times_B\cB\uT_{|B}\ar[l]^{+}\ar[d]&&\\
P\otimes_\cE\tilde G&P\times_{\Z_{|B}}\tilde G\ar[d]\ar[l]\ar[r]&P\times_B \tilde G\ar[d]\\&\uZ_{|B}\ar[r]^{\diag}&\uZ_{|B}\times_B \uZ_{|B}}\ .
\end{equation}

The right lower square is cartesian,
and in the left upper square we take the fibre-wise quotient by the
anti-diagonal action of $\cB\uT_{|B}$.

\subsubsection{}

We have
$D(P)\in  \Ext_{\cPic(\bS/B)}(\uZ,D(\cE))$,
where
$$0\to \cB\uT_{|B}\to D(\cE)\to D(\cE)\to 0\ .$$
We define
$\overline{D(P)}$ to be the quotient of $D(P)$ by $\cB \uT_{|B}$ in the sense of \ref{dedwhedjewdn66} so that
$D(P)\to \overline{D(P)}$ is a gerbe with band $\uT$.

We define
$D(P)\otimes_{\overline{D(P)}} D(\tilde G)$ by
the diagram
\begin{equation}
\label{qjhwdqwdqwdwqd1}
\xymatrix{\cB\uT_{|B}\ar[d]&\cB\uT_{|B}\times_B\cB\uT_{|B}\ar[l]^{+}\ar[d]&&\\
D(P)\otimes_{\overline{D(P)}}D(\tilde G)&D(P)\times_{\Z_{|B}}D(\tilde G)\ar[d]\ar[l]\ar[r]&
D(P)\times_B D(\tilde G)\ar[d]\\ & \uZ_{|B} \ar[r]^{\diag}& \uZ_{|B} \times_B  \uZ_{|B} }\ .\end{equation}
Again, the right lower square is cartesian,
and in the left upper square we take the fibre-wise quotient by the
anti-diagonal action of $\cB\uT_{|B}$.

\subsubsection{}

\begin{prop}\label{whdwqjdhwdwqdqwuidiwqd}
We have an equivalence of Picard stacks
$$D(P\otimes_\cE \tilde G)\cong D(P)\otimes_{\overline{D(P)}} D(\tilde G)\ .$$
\end{prop}
\proof
The diagram \eqref{qjhwdqwdqwdwqd1} defines 
 $D(P)\otimes_{\overline{D(P)}} D(\tilde G)$ by forming the pull-back to the
 diagonal 
 and then taking the quotient of the anti-diagonal $\cB\uT_{|B}$-action in the
 fibre. One can obtain this diagram by dualizing  (\ref{qjhwdqwdqwdwqd}) and
 interchanging the order of pull-back and quotient.
\hB

\subsubsection{}

We can now finish the argument that
$\Psi$ is $H^3(B;\uZ)$-equivariant.
We let
$\Psi(P)=((E,H),(\hat E,\hat H),u)$ and
$\Psi(g+P)= ((E^\prime,H^\prime),(\hat E^\prime,\hat H^\prime),u^\prime)$.
Note that
$g+P=P\otimes_\cE \tilde G$.
An inspection of the construction of the first entry of $\Psi$ shows that
$E^\prime\cong E$ and $H^\prime\cong H\otimes_E \pr{E\to B}G$.
Proposition \ref{whdwqjdhwdwqdqwuidiwqd} shows that
$\overline{D(P\otimes_\cE \tilde G)}\cong \overline{D(P)}$.
This implies that $\hat E^\prime\cong \hat E$.
Furthermore, if we restrict
$D(P\otimes_\cE \tilde G)$ along  $\underline{\{1\}}_{|B}\to \uZ_{|B}$ we get
by Proposition \ref{whdwqjdhwdwqdqwuidiwqd} and the proof of Lemma
\ref{iudwqodhqwdqwdqwd} a diagram of stacks over $B$
$$\xymatrix{\hat H^{op}\otimes_{\hat E} G^{op}\ar[d]\ar[r]& \tilde{\hat H}\otimes_{R_\R^{\hat E}} G^{op}\ar[r]\ar[d]&D(P\otimes_\cE \tilde G)\ar[d]\ar[r]^\cong &D(P)\otimes_{\overline{D(P)}}D(\tilde G)\ar[d]\\
 \hat E\ar[r]^{\can}& R_\R^{\hat E}\ar[r]\ar[d]&\overline{D(P\otimes_\cE \tilde G)}\ar[d]\ar[r]^{\cong}&\overline{D(P)}\\&\underline{\{1\}}_{|B}\ar[r]&\uZ_{|B}&}$$
The map $u^\prime$ is induced by the evaluation
$(P\otimes_\cE  \tilde G)\times D(P\otimes_\cE \tilde G)\to \cB \uT$.
With the given identifications using the ``duality'' between (\ref{qjhwdqwdqwdwqd1}) 
and (\ref{qjhwdqwdqwdwqd}) we see that
 this evaluation the induced by the  product of the evaluations
$$(P\times  D(P))\times (\tilde G\times D(\tilde G))\to \cB\uT\times \cB \uT
\to \cB\uT\ .$$
After restriction to $\underline{\{1\}}_{|B}$ we see that
$u^\prime=u\otimes v$, where $v:G\otimes G^{op}\to \cB\uT$ is the canonical pairing.
This finishes the proof of the equivariance of $\Psi$.
\hB


\begin{theorem}\label{main}
The maps $\Psi$ and $\Phi$ are inverse to each other.
\end{theorem}
\proof 
We first show the assertion under the additional assumption that
$H^3(B;\uZ)=0$. In this case
an element $P\in Q_\cE$ is determined uniquely by the class $\widehat c(P)\in H^2(B;\uZ^n)$.
Similarly, a $T$-duality triple $t\in \Triple(B)$ with $c(t)=c(E)$ is uniquely determined by the class
$\widehat c(t)= c(\widehat E)$. Since for $P\in Q_\cE$ we have $\widehat c(\Phi\circ \Psi(P))\stackrel{(\ref{gghze673})}{=}\widehat c(\Psi(P))\stackrel{Lemma \:\ref{hewdjewhdf89}}{=}\widehat c(P)$ and
$\widehat c(\Psi\circ \Phi(t))\stackrel{Lemma \:\ref{hewdjewhdf89}}{=}\widehat c(\Phi(t))\stackrel{(\ref{gghze673})}{=}\widehat c(t)$  
this implies that
$\Phi\circ \Psi_{|Q_\cE}=\id_{Q_\cE}$ and
$\Psi\circ \Phi_{|s^{-1}(E)}=\id_{|s^{-1}(E)}$,
where $s:\Triple\to P$ is as in \ref{ewjdhdwed78283324}.
Note that $H^3(R_n;\uZ)=0$.
Therefore
$$\Psi(\Phi(t_{univ}))=t_{univ}\ ,\quad \Phi(\Psi(P_{univ}))\cong P_{univ}\ .$$

Now consider a general space $B\in \bS$. We first show that $\Psi\circ \Phi=\id$.
Let $t\in s^{-1}(E)$ be classified by the map $f_t:B\to R_n$, i.e. $f_t^*t_{univ}=t$.
Then we have
$\Psi(\Phi(t))=\Psi(f_t^* P_{univ})=f_t^*\Psi(P_{univ})=f_t^*t_{univ}=t$, i.e.
\begin{equation}\label{hfwjfwe728374}
\Psi\circ \Phi=\id\ .
\end{equation}

We consider the group $$\Gamma_\cE:=(\im(\alpha)+\im(s))/\im(\alpha)\subseteq H^3(B;\uZ)/\im(\alpha)\ .$$ This group is exactly the group $\ker(\pi^*)/C$ in the notation of  \cite[Theorem 2.24(3)]{math.GT/0501487}.
It follows from (\ref{gdhgqwdqw691990290}) that the action of $H^3(B;\uZ)$ on $Q_\cE$ induces an action of 
 $\Gamma_\cE$  on $Q_\cE$ which preserves the fibres of $f$.
In \cite[Theorem 2.24(3)]{math.GT/0501487} we have shown that
it also acts on $\Triple_E(B)$ and preserves the subsets $\Ext(E,h)$.

Let us fix a Picard stack $P\in Q_\cE$, and let $\widehat c:=\widehat c(P)$
and $h\in H^3(E;\uZ)$ be such that $f(P)\in \Ext(E,h)$. From (\ref{gdhgqwdqw691990290}) we see that
the group $\Gamma_\cE$
acts simply transitively on the  set
$$A_{\cE,\widehat c}:=\{Q\in f^{-1}(h)\mid\widehat c(Q)=\widehat c\}\ .$$
 By  \cite[Theorem 2.24(3)]{math.GT/0501487} it also acts simply transitively
on the set $$B_{\cE,\widehat c}:=\{t\in \Ext(E,h)\mid \widehat c(t)=\widehat c\}\ .$$
By Lemma \ref{hewdjewhdf89} we have $\Psi(A_{\cE,\widehat c})\subseteq B_{\cE,\widehat c}$.
By Lemma \ref{ghhwezfuwf78223} the map $\Psi$ is $\Gamma_\cE$-equivariant.
Hence it must induce a bijection between $A_{\cE,\widehat c}$ and $B_{\cE,\widehat c}$.
If we let $\widehat c$ run over all possible choices (solutions of $\widehat c\cup  c(E)=0$) we see that
$\Psi:Q_\cE\to \Triple_E(B)$ is a bijection. In view of (\ref{hfwjfwe728374}) we now also get $\Phi\circ \Psi=\id$. \hB \\[2cm]

{\em Acknowledgment:  We thank Tony Pantev for the crucial suggestion which initiated the work on this paper.}


\end{document}